\documentclass[review=flase,onefignum,onetabnum]{siamart220329}

\usepackage{amsmath}
\usepackage{amssymb}
\usepackage[bmargin=2.5cm,lmargin=2.5cm,rmargin=3.5cm,tmargin=2.5cm,asymmetric]{geometry}
\usepackage{mathrsfs}
\usepackage{mathtools}
\usepackage{amscd}
\usepackage[all]{xy}
\SelectTips{cm}{}
\usepackage{multirow}
\usepackage{stmaryrd}
\usepackage{enumitem}
\usepackage{textgreek}
\usepackage{hyperref}
\usepackage{cleveref}

\usepackage{tikz}
\tikzstyle{nodes1} = [circle, rounded corners, minimum width=1cm, minimum height=1cm,text centered, draw=black, fill=white!30]
\tikzstyle{arrow} = [thick,->,>=stealth]

\usepackage{xargs}                      
\usepackage{xcolor}  

\usepackage{algorithm}
\usepackage{algpseudocode}

\usepackage[colorinlistoftodos,prependcaption,textsize=tiny]{todonotes}
\newcommandx{\unsure}[2][1=]{\todo[linecolor=red,backgroundcolor=red!25,bordercolor=red,#1]{#2}}
\newcommandx{\change}[2][1=]{\todo[linecolor=blue,backgroundcolor=blue!25,bordercolor=blue,#1]{#2}}
\newcommandx{\info}[2][1=]{\todo[linecolor=OliveGreen,backgroundcolor=OliveGreen!25,bordercolor=OliveGreen,#1]{#2}}
\newcommandx{\improvement}[2][1=]{\todo[linecolor=Plum,backgroundcolor=Plum!25,bordercolor=Plum,#1]{#2}}
\newcommandx{\thiswillnotshow}[2][1=]{\todo[disable,#1]{#2}}

\newcommand{\reff}[1]{(\ref{#1})}

\newcommand{\tallopenbox}{\leavevmode
  \hbox to.6em{%
  \hfil\vrule
  \vbox to.65em{\hrule width.35em\vfil\hrule}%
  \vrule\hfil}}

\newcommand{\bdes}{\begin{description}}
	\newcommand{\edes}{\end{description}}
\newcommand{\bal}{\begin{align}}
	\newcommand{\eal}{\end{align}}
\newcommand{\bnum}{\begin{enumerate}}
	\newcommand{\enum}{\end{enumerate}}
\newcommand{\bit}{\begin{itemize}}
	\newcommand{\eit}{\end{itemize}}
\newcommand{\bea}{\begin{eqnarray}}
	\newcommand{\eea}{\end{eqnarray}}
\newcommand{\be}{\begin{equation}}
	\newcommand{\ee}{\end{equation}}
\newcommand{\baray}{\begin{array}}
	\newcommand{\earay}{\end{array}}
\newcommand{\bsry}{\begin{subarray}}
	\newcommand{\esry}{\end{subarray}}
\newcommand{\bca}{\begin{cases}}
	\newcommand{\eca}{\end{cases}}
\newcommand{\bcen}{\begin{center}}
	\newcommand{\ecen}{\end{center}}
\newcommand{\bbm}{\begin{bmatrix}}
	\newcommand{\ebm}{\end{bmatrix}}
\newcommand{\bmx}{\begin{matrix}}
	\newcommand{\emx}{\end{matrix}}
\newcommand{\bpm}{\begin{pmatrix}}
	\newcommand{\epm}{\end{pmatrix}}
\newcommand{\btab}{\begin{tabular}}
	\newcommand{\etab}{\end{tabular}}

\def\<#1,#2>{\langle #1,#2\rangle}
\def\<#1>{\langle #1\rangle}

\newtheorem{remark}[theorem]{Remark}

\makeatletter
\newcommand{\extp}{\@ifnextchar^\@extp{\@extp^{\,}}}

\def\@extp^#1{\mathop{\bigwedge\nolimits^{\!#1}}}
\makeatother

\renewcommand{\thesection}{\arabic{section}}
\renewcommand{\theequation}{\thesection.\arabic{equation}}
\numberwithin{equation}{section}

\begin{document}

\title{Asymptotic behavior of generic Riemann solutions to one-dimensional full compressible Navier-Stokes-Poisson equations
}

\author{
    Xin-Xiang Bian \thanks{School of Information and Mathematics, Yangtze University, (\email{bxx@yangtzeu.edu.cn})}
	\and Jiajin Shi \thanks{School of Mathematics and Statistics, Jiangsu Normal University, (\email{shijiajin@amss.ac.cn})}
    \and Lingling Xie\thanks{School of Mathematics and Statistics, Hefei Normal University, (\email{linglingxie@hfnu.edu.cn})}
}


\maketitle

\begin{abstract}
	This paper investigates the large-time asymptotic stability of the one-dimensional full compressible Navier-Stokes-Poisson equations when the corresponding quasi-neutral Euler equations admit a generic Riemann solution consisting of a shock wave, a contact discontinuity, and a rarefaction wave. We prove the existence and uniqueness of a global-in-time solution that converges in \(L^\infty\) norm, as time tends to infinity, to a composite wave profile consisting of a viscous shock wave with a time-dependent shift, a viscous contact wave, and an inviscid rarefaction wave. The main difficulty lies in the nonconservative lower-order source terms arising from the electric potential in the momentum and energy equations, which substantially complicate both the construction of individual wave profiles, particularly the viscous shock wave, and the derivation of relative entropy estimates. To overcome these difficulties, we exploit the structure of the Poisson equation governing the electric potential and combine the weighted relative entropy method with the \(a\)-contraction method with shifts to close the a priori estimates within an \(H^2\)-perturbation framework.
\end{abstract}




\section{Introduction}
\hspace*{\fill} \\

In dusty plasma or semiconductor devices \cite{Chen18,Markowich90}, the transport of charged particles under the influence of the self-consistent electric field of electrostatic potential force is governed by the one-dimensional unipolar non-isentropic compressible Navier-Stokes-Poisson system (CNSP), 
which in the Lagrangian coordinates reads  
\begin{subequations}\label{eq:CNSP}
	\begin{align}
		&v_t - u_x = 0, \label{eq:CNSPmass} \\
		&u_t + p_x = \mu \left(\frac{u_x}{v}\right)_x - \frac{\phi_x}{v}, \label{eq:CNSPmomentum} \\
		&E_t + (p u)_x = \kappa \left(\frac{\theta_x}{v}\right)_x + \mu \left(\frac{u u_x}{v}\right)_x - \frac{u\phi_x}{v},  \label{eq:CNSPenergy} \\
		&-\nu^2 \left(\frac{\phi_x}{v}\right)_x = 1 - v e^\phi, \label{eq:CNSPpoisson}
	\end{align}
\end{subequations}
for $t\ge 0$ and $x\in\mathbb{R}$. The unknown variables $v=1/n,u,\theta,\phi$ denote the specific volume,
the fluid velocity, the absolute temperature and the electric potential, respectively, where $n>0$ is the density of ions and $E = \epsilon + \frac{u^2}{2}$ is the specific total energy.
For the ideal polytropic gas, the internal energy $\epsilon$ and the pressure $p$ is given by 
\begin{align} \label{intro-p}
	\epsilon = \frac{R\theta}{\gamma - 1} + \text{const.},
	\quad
	p = p(v,\theta) = \frac{R\theta}{v}
\end{align}
with the gas constant $R>0$ and the adiabatic exponent $\gamma>1$. 
The constants $\mu>0, \kappa>0$ and $\nu>0$ correspond to the viscosity coefficient, the heat-conductivity coefficient and the Debye length.
In particular, the Poisson equaton \reff{eq:CNSPpoisson} is supplemented by the Boltzmann relation $n_e=e^\phi$ for the electron density, 
under the assumption justified by the physical fact that electrons reach equilibrium much faster than ions when the potential varies in a plasma.

In this paper, we are concerned with the large-time behavior of composite waves for the CNSP system \reff{eq:CNSP} with the associated initial data 
\begin{align} \label{CNSPID}
	(v,u,\theta,\phi)(t=0,x)=(v_0,u_0,\theta_0,\phi_0)(x) \rightarrow (v_\pm, u_\pm, \theta_\pm, \phi_\pm), 
	\quad \text{as}\quad x \rightarrow \pm\infty
\end{align}
for the given two different end states $(v_\pm, u_\pm, \theta_\pm, \phi_\pm)\in \mathbb{R}_+ \times \mathbb{R} \times \mathbb{R}_+ \times\mathbb{R}$
($\mathbb{R}_+:=(0,+\infty)$) and the quasi-neutral condition 
\begin{align} \label{quasi-neutral}
	\phi_\pm=-\ln v_\pm.
\end{align}

Formally speaking, the time-asymptotic behavior of the non-isentropic CNSP system \reff{eq:CNSP} with the initial data verifying \reff{CNSPID}
is expected to be determined by the associated quasi-neutral Euler equations
\begin{subequations}\label{eq:Euler}
	\begin{align}
		& v_t - u_x = 0, \label{eq:Eulermass} && \\
		& u_t + \left(p + \frac{1}{v}\right)_x = 0, \label{eq:Eulermomentum} && \\
		& E_t + (pu)_x + \left(\frac{1}{v}\right)_x u = 0, \label{eq:Eulerenergy} && \\[1mm]
		& \phi = - \ln v \label{eq:Eulerpoisson} &&
	\end{align}
\end{subequations}
with the Riemann initial data  
\begin{align} \label{EulerID}
	(v, u, \theta)(t = 0, x) &= 
	\begin{cases} 
	(v_-, u_-, \theta_-), & x < 0, \\ 
	(v_+, u_+, \theta_+), & x > 0 
	\end{cases} 
\end{align}
and the quasi-neutral initial data 
\begin{align*}
	\phi(t=0, x)=-\ln v(t=0, x)=
	\begin{cases}
		-\ln v_-, & x < 0, \\ 
		-\ln v_+, & x > 0.
	\end{cases}
\end{align*} 
The hyperbolic system \reff{eq:Eulermass}-\reff{eq:Eulerenergy} is strictly hyperbolic,  since it has three distinct real eigenvalues,
among which the first and third give rise to two genuinely nonlinear characteristic fields, whereas the second produces one linearly degenerate characteristic field. 
Therefore, the Riemann problem \reff{eq:Eulermass}-\reff{eq:Eulerenergy} with \reff{EulerID} has the self-similar solutions, called Riemann solution, 
consisting of a combination of at most three elementary solutions: shock, rarefaction and contact discontinuity (see also \cite{Dafermos79}).


We begin with the stability of viscous shock profiles, a topic that has been extensively studied in the literature. 
In 1960, Il'in and Ole\u{\i}nik \cite{Ilin60} first proved the stability of shock and rarefaction waves to the one-dimensional scalar Burgers equations.
Then Matsumura and Nishihara \cite{Matsumura85} extended the shock stability result to the one-dimensional compressible Navier-Stokes equations with physical viscosity, 
and independently Goodman \cite{Goodman86} obtained a similar result for a general system with artificial viscosity.
Both papers required the zero-mass condition in order to apply the anti-derivative method.
Duan et al.\cite{Duan20} successfully applied this method to establish the shock stability for the isentropic CNSP equations, with subsequent extensions to large initial density oscillations by Zhao \cite{Zhao23} and to space-periodic perturbations by Li et al.\cite{Li25}.
Meanwhile, the zero-mass condition was subsequently removed in \cite{Liu85,Liu15,Szepessy93} by means of the pointwise estimate method, and relaxed to a slightly weaker spectral condition in \cite{Mascia04}.
Notably, Li et al.\cite{Li26} and Duan and Li \cite{Duan26} have recently established the stability of shocks for the non-isentropic CNSP system, using the anti-derivative technique and the weighted relative entropy framework, respectively.
See \cite{Han23,Huang09,Kang-Kwon25,Kang21,Matsumura10,Wang25,Wang26} for more results on shocks.

Regarding the stability of rarefaction waves, one can employ an alternative approach based on direct energy methods.  
The stability result was first proved by Matsumura and Nishihara \cite{Matsumura86,Matsumura92} for the isentropic Navier-Stokes equations.
It was later generalized in \cite{Liu88,Nishihara04} to the full compressible Navier-Stokes-Fourier system and in \cite{Duan15,Duan13} to the isentropic CNSP system.  Further results can be found in \cite{Li18-1,Li18-2}.
For the linearly degenerate characteristic field, the viscous contact wave can be proved to be time-asymptotically stable for the one-dimensional non-physical \cite{Liu97,Xin96} and physical Navier-Stokes equations \cite{Huang08}.
In particular, the large-time behavior of composite waves consisting of rarefaction wave and viscous contact wave was established by Huang et al.\cite{Huang10} for the full compressible Navier-Stokes-Fourier system, and by Ruan et al.\cite{Ruan17} for the non-isentropic CNSP system, respectively.
Notably, Zhao and Zhao \cite{Zhao26} have considered the same case with large initial perturbations.

However, the stability of composite waves formed by different types of waves, especially shocks and rarefaction waves, has posed a major difficulty, since the anti-derivative technique for shocks and the elemental energy method for rarefaction and viscous contact waves cannot be applied simultaneously to such composite waves.
Fortunately, in recent years, this problem was first resolved by Kang et al.\cite{Kang23}. Their analysis employs the theory of $a$-contraction with shifts within a relative entropy framework, introduced by Dafermos \cite{Dafermos96} and Diperna \cite{DiPerna79}.
More recently, this approach has been further extended to establish the time-asymptotic stability of generic Riemann solutions for the full compressible Navier-Stokes-Fourier equations \cite{Kang25}, as well as that of composite waves for the relaxed compressible Navier-Stokes equations \cite{Guan26} and for the isentropic CNSP equations \cite{Shim26}.

This paper is devoted to the time-asymptotic stability of composite waves for the non-isentropic CNSP equations. 
More precisely, we extend the stability result of \cite{Kang25} to the wave configuration formed by the superposition of a viscous shock, a viscous contact wave, and an inviscid rarefaction wave for the one-dimensional non-isentropic CNSP system \reff{eq:CNSP}-\reff{CNSPID}.
Our proof also employs the $a$-contraction method with shifts and the relative entropy method to close the $H^2$-perturbation framework.
However, compared with \cite{Kang25}, the presence of the electric potential gives rise to several new difficulties:
\\[-1ex]

\noindent (i) $\,$ \hangindent=2em  The first difficulty stems from the existence of viscous shock waves for the non-isentropic CNSP system \reff{eq:viscous_shock}-\reff{far-field}. 
The essential problem lies in the lower-order source terms introduced by the electric potential into the momentum and energy equations. 
These terms are non-conservative and pose a great challenge to proving the existence of shock. 
However, we observe that the Poisson equation allows these terms to be transformed into a conservative form \reff{phi_relation} and \reff{phi_u_relation}. 
By further combining the center manifold theorem and the solvability of the associated ODEs, we establish the existence of shock waves and obtain the required estimates in Lemma \ref{VS-pro}, relying on the smallness of the shock strength. 
On the other hand, we note that \cite{Duan26,Li26} have recently proved the existence of shocks for the non-isentropic CNSP system, a result that was not available before our establishing the stability of composite waves. 
Since our approach differs slightly from theirs and the self-containedness of the paper is desirable, we include our proof in the Appendix for completeness. 
\\[-1ex]

\noindent (ii) $\,$ \hangindent=2em  Another difficulty lies in obtaining the sharp estimates of the dissipation term in Lemma \ref{Crucial}. 
First, the Taylor expansion reveals that these estimates are determined by the second derivative of the temperature with respect to the volume. 
A direct calculation inevitably produces singularities, as the singular terms appear to fail to obtain a leading-order term up to an error of the same order as the shock strength (see \eqref{phi_v}-\eqref{phi_v2}), since this would require the precise value of the exponents of the equivalent quantity.
We overcome this by combining Lemma \ref{VS-pro} with a shift of the far-field value to the origin, up to an error of the same order as the shock strength, thereby eliminating these singularities.
Note that the sharp estimates depend on the smallness of the shock strength, since a Taylor expansion is employed.
\\[-1ex]

\noindent (iii) $\,$ \hangindent=2em  Finally, compared with the weighted relative entropy method in \cite{Kang25}, the main difference arises from the last two terms of \reff{pf-vthetau}, which correspond to the source terms \reff{ele-add} induced by the electric potential. 
A direct relative entropy estimate without modification would produce error terms lacking smallness coefficient control, making the desired estimates unattainable. This is fatal to the present problem. 
Hence, these two terms must be separated into potentially problematic main terms and error terms controlled by a smallness coefficient. 
By exploiting the structure of the mass equation and the Poisson equation, we use the favorable terms inherent in the relative entropy to control the separated main terms, thereby successfully closing the estimates within the $H^2$-perturbation framework.
\\[-1ex]

Below, we state our main result on the nonlinear stability of the wave pattern composed of a combination of a viscous shock up to a time-dependent shift, a viscous contact wave and an inviscid rarefaction wave for the one-dimensional non-isentropic compressible CNSP system.
\\ 

\begin{theorem}  \label{Theorem-Asym}
	For a given constant state $(v_+, u_+, \theta_+)\in\mathbb{R}_+\times\mathbb{R}\times\mathbb{R}_+$, there exist constants $\delta_0, \varepsilon_0 > 0$ such that the following holds true.
	For any $(v^*, u^*, \theta^*) \in S_3(v_+, u_+, \theta_+)$, $(v_*, u_*, \theta_*)\in CD_2(v^*, u^*, \theta^*)$ and $(v_-, u_-, \theta_-)\in R_1(v_*, u_*, \theta_*)$ satisfying
	\[|v_+ - v^*| + |v^* - v_*| + |v_* - v_-| \le \delta_0,\]
	where $S_3(v_+, u_+, \theta_+)$ is given by the Rankine-Hugoniot condition \reff{RH} with $v_R$ replaced by $v_+$ and with $v<v_+$, and $CD_2(v^*, u^*, \theta^*)$ and $R_1(v_*, u_*, \theta_*)$ are the curves defined in \reff{CD_2_def} and \reff{R_1_def}, respectively.
	Denote $(v^r, u^r, \theta^r, \phi^r)(x/t)$ the $1$-rarefaction wave defined in \reff{RV} and \reff{RV-phi} with end states $(v_-, u_-, \theta_-)$ and $(v_*, u_*, \theta_*)$, 
	$(v^C, u^C, \theta^C, \phi^C)(x/\sqrt{1+t})$ the viscous $2$-contact wave to \reff{equ:VCD-ori} and \reff{eq:VCD} with end states $(v_*, u_*, \theta_*)$ and $(v^*, u^*, \theta^*)$,
	$(v^S, u^S, \theta^S, \phi^S)(x-\sigma t)$ the viscous $3$-shock solution to \reff{eq:CNSP} and \reff{eq:viscous_shock} with end states $(v^*, u^*, \theta^*)$ and $(v_+, u_+, \theta_+)$.
	Let $(v_0, u_0, \theta_0, \phi_0)$ be any initial data such that 
	\begin{align*}
		\begin{array}{l}
			\displaystyle \sum_{\pm} \| (v_0 - v_\pm, \, u_0 - u_\pm, \, \theta_0 - \theta_\pm, \, \phi_0-(-\ln v_\pm)) \|_{L^2(\mathbb{R}_\pm)} \\[4mm]
			\displaystyle  \qquad 
			+ \| (v_{0x},  u_{0x},  \theta_{0x}) \|_{H^1(\mathbb{R})} + \| \phi_{0x} \|_{H^2(\mathbb{R})} \le \varepsilon_0,
		\end{array}
	\end{align*}
	with $R_{-} : = -R_{+} = (-\infty, 0)$. 
	
	Then the one-dimensional non-isentropic CNSP system \reff{eq:CNSP} admits a unique global-in-time solution $(v, u, \theta, \phi)(t, x)$ for all $t\in\mathbb{R}_+$.
	Moreover, there exists an absolutely continuous shift $\mathbf{X}(t)$ (defined in \reff{Xt-def}) such that 
	\begin{align*}
		\begin{array}{l}
			\displaystyle v(t, x) - \Big( v^r(x/t) + v^C(x/\sqrt{1 + t}) + v^S(x - \sigma t - \mathbf{X}(t)) - v_* - v^* \Big) \in C(0, +\infty; H^2(\mathbb{R})), \\[4mm]
			\displaystyle u(t, x) - \Big( u^r(x/t) + u^C(x/\sqrt{1 + t}) + u^S(x - \sigma t - \mathbf{X}(t)) - u_* - u^* \Big) \in C(0, +\infty; H^2(\mathbb{R})), \\[4mm]
			\displaystyle \theta(t, x) - \Big( \theta^r(x/t) + \theta^C(x/\sqrt{1 + t}) + \theta^S(x - \sigma t - \mathbf{X}(t)) - \theta_* - \theta^* \Big) \in C(0, +\infty; H^2(\mathbb{R})), \\[4mm]
			\displaystyle \phi(t, x) - \Big( \phi^r(x/t) + \phi^C(x/\sqrt{1 + t}) + \phi^S(x - \sigma t - \mathbf{X}(t)) + \ln v_* + \ln v^* \Big) \in C(0, +\infty; H^3(\mathbb{R})),  \\[4mm]
			\displaystyle \big(u(t, x) - u^S( x - \sigma t - \mathbf{X}(t)) \big)_{xx} \in L^2(0, +\infty; H^1(\mathbb{R})), \\[4mm]
			\displaystyle \big(\theta(t, x) - \theta^S( x - \sigma t - \mathbf{X}(t)) \big)_{xx} \in L^2(0, +\infty; H^1(\mathbb{R})). 
		\end{array}
	\end{align*}
	In addition, as $t\to+\infty$, the solution $(v, u, \theta, \phi)(t,x)$ satisfies the following time-asymptotic behavior
	\begin{align} \label{Asym-beh}
		\begin{array}{l}
			\displaystyle \sup_{x \in \mathbb{R}} \Bigg| (v, u, \theta, \phi)(t, x) - \Big( v^r(x/t) + v^C(x/\sqrt{1+t}) + v^S(x - \sigma t - \mathbf{X}(t)) - v_* - v^*, \\[4mm]
			\displaystyle \qquad \qquad u^r(x/t) + u^C(x/\sqrt{1+t}) + u^S(x - \sigma t - \mathbf{X}(t)) - u_* - u^*, \\[4mm]
			\displaystyle \qquad \qquad \theta^r(x/t) + \theta^C(x/\sqrt{1+t}) + \theta^S(x - \sigma t - \mathbf{X}(t)) - \theta_* - \theta^* , \\[2mm]
			\displaystyle \qquad \qquad \phi^r(x/t) + \phi^C(x/\sqrt{1+t}) + \phi^S(x - \sigma t - \mathbf{X}(t)) + \ln v_* + \ln v^* \Big) \Bigg| \to 0,
		\end{array}
	\end{align}
	and it holds that
	\begin{align} \label{Xt-asym-theo}
		\lim_{t\to+\infty} |\dot{\mathbf{X}}(t)|\to 0.
	\end{align}
\end{theorem}
\hspace*{\fill}

\begin{remark}
	Theorem \ref{Theorem-Asym} implies that when the two far-field states in \reff{CNSPID} verifying the quasi-neutral condition \reff{quasi-neutral} are connected by the superposition of a shock, a contact discontinuity, and a rarefaction wave (i.e., \reff{Generic-Rie}), there exists a unique global-in-time solution to the non-isentropic CNSP system. 
	This solution converges in $L^\infty$ norm, as time goes to infinity, to the ansatz composed of a viscous shock with a time-dependent shift, a viscous contact wave, and an inviscid rarefaction wave. 
	The above composite wave also covers the cases of three single waves, namely, the shock, the rarefaction wave, and the contact discontinuity, since the single-wave structure corresponds to the degenerate cases of the composite wave.
	Moreover, compared with the results in \cite{Kang25}, we conclude that the presence of the electric potential satisfying the Poisson equation does not alter the asymptotic behavior of the composite wave, at least for the current wave configuration.
\end{remark}
\hspace*{\fill}

\begin{remark}
	From \reff{Xt-asym-theo}, the shift function $\mathbf{X}(t)$ grows at most sublinearly as $t\to+\infty$, i.e.,
	\begin{align*}
		\lim_{t\to+\infty} \dfrac{\mathbf{X}(t)}{t} = 0,
	\end{align*}
	which implies that the shifted viscous shock profile still keeps the original shape time-asymptotically.
\end{remark}
\hspace*{\fill}

The rest of the paper is organized as follows. In Section 2, we recall the well-known results on the generic Riemann problem for the quasi-neutral Euler equations \reff{eq:Euler}-\reff{EulerID}. 
In Section 3, we present some preliminaries on the approximate rarefaction wave together with the viscous versions of the contact wave and the shock wave, and state the existence of shock waves and the sharp estimates for the diffusion. 
The general set-up is laid out in Section 4, where we combine the local existence of solutions in Proposition \ref{Local Existence} with the a priori estimates in Proposition \ref{A prior estimate} to prove the global existence of solutions and derive their asymptotic behavior, thereby completing the proof of Theorem \ref{Theorem-Asym}. 
In Section 5, we combine the wave interaction estimates, the weighted relative entropy method, and the elliptic estimates to obtain the $L^2$ estimates for the perturbation. The final section is devoted to the $H^2$ estimates, while the Appendix is reserved for the existence of shock waves, the sharp diffusion estimates, and the basic estimates on the electric potential.
\\

\section{Riemann problem}
\hspace*{\fill} \\

We first recall the well-known results for the Riemann problem for the inviscid models \reff{eq:Euler}-\reff{EulerID}, first proposed and solved by Riemann \cite{Riemann60} in1860s and then generalized by Lax \cite{Lax57}.

This system of conservation laws \reff{eq:Eulermass}-\reff{eq:Eulerenergy} is strictly hyperbolic. Indeed, the derivative matrix of the flux function (though not given the explicit form) with respect to the conserved variables $(v, u, E)$ 
\begin{align} \label{Der-flux}
	\left(
	\begin{array}{ccc}
		0 & -1 & 0 \\[2mm]
		-\dfrac{p}{v} - \dfrac{1}{v^2} & -\dfrac{(\gamma-1)u}{v} & \dfrac{\gamma-1}{v} \\[3mm]
		-\dfrac{pu}{v} - \dfrac{u}{v^2} & p - \dfrac{(\gamma-1)u}{v} & \dfrac{(\gamma-1)u}{v}
	\end{array}
	\right)
\end{align}
is diagonalizable with distinct real eigenvalues.
The matrix characterizes the waves generated by the linearization of \reff{eq:Eulermass}-\reff{eq:Eulerenergy} about the fixed state $(v, u, \theta)\in\mathbb{R}^+\times\mathbb{R}\times\mathbb{R}^+$.
Its three eigenvalues are 
\begin{align} \label{Eigen}
	\sigma_1 = \sigma_1(v, \theta) = - \sqrt{\dfrac{\gamma R\theta + 1}{v^2}}, \quad \sigma_2 = 0, \quad \sigma_3 = \sigma_3(v, \theta) = \sqrt{\dfrac{\gamma R\theta + 1}{v^2}}
\end{align}
among which $\sigma_1$ and $\sigma_3$ give rise to two genuinely nonlinear characteristic fields, whereas $\sigma_2$ produces one linearly degenerate characteristic field. 
Therefore, the Riemann problem \reff{eq:Eulermass}-\reff{eq:Eulerenergy} with \reff{EulerID} has the self-similar solutions, called Riemann solutions, consisting of a combination
of at most three elementary solutions from the following five families: 1-rarefaction; 3-rarefaction; 1-shock; 3-shock and 2-contact discontinuity (see also \cite{Dafermos79}).
These families are determined by their associated curves in the states plane $\mathbb{R}^+\times\mathbb{R}\times\mathbb{R}^+$.
\\ 

\subsection{Rarefaction wave}
\hspace*{\fill} \\

For any state $(v_R, u_R, \theta_R)\in\mathbb{R}^+\times\mathbb{R}\times\mathbb{R}^+$,
the 1-rarefaction curve $R_1(v_R, u_R, \theta_R)$ generated by $\sigma_1$ is given by 
\begin{align} \label{R_1_def}
	R_1(v_R, u_R, \theta_R):=\Big\{
		(v,u,\theta)\Big| v<v_R, s(v,\theta)=s(v_R, \theta_R)=:s_R, u=u_R-\int_{v_R}^{v} \sigma_1(v^\prime, s_R) dv^\prime
	\Big\},
\end{align} 
where the entropy $s$ satisfies the Gibbs' relation $\theta ds=dE-udu+pdv$ and thus has the form $s=s(v,\theta)=\frac{R}{\gamma-1}\ln (\frac{R}{A}\theta v^{\gamma-1})$,
and $\sigma_1(v,s)=-\sqrt{\frac{\gamma p(v,s)}{v}+\frac{1}{v^2}}$ replacing $\theta$ in $p(v, \theta)$ by $s$.

For the Riemann problem \reff{EulerID} with $(v_-, u_-, \theta_-) = (v_L, u_L, \theta_L)$ and $(v_+, u_+, \theta_+) = (v_R, u_R, \theta_R)$, such that $(v_L, u_L, \theta_L)\in S_1(v_R, u_R, \theta_R)$.
Let $w_L=\sigma_1(v_L, \theta_L),\, w_R=\sigma_1(v_R, \theta_R)$, and then if $w_L<w_R$,
\begin{align} \label{RV-w}
	w^r(x/t)=\left\{ 
		\begin{array}{ll}
			w_L, & x/t < w_L, \\ [2mm] 
			x/t, & w_L\le x/t \le w_R, \\ [2mm]
			w_R, & x/t > w_R,
		\end{array}
	\right.
\end{align}
satisfies the inviscid Burgers equation 
\begin{align*}
    \left\{ 
       \begin{array}{l}
            w_t + w w_x = 0, \\ [2mm]
            w(0, x) = w_0(x): =\left\{
				\begin{array}{ll}
					w_L, & x<0, \\[1mm]
					w_R, & x>0.
				\end{array}
			\right. 
       \end{array}
    \right.
\end{align*}
The corresponding 1-rarefaction wave fan $(v^r, u^r, \theta^r)(t,x)=(v^r, u^r, \theta^r)(x/t)$ is defined as 
\begin{align} \label{RV}
	\begin{array}{l}
		\sigma_1(v^r(x/t), \theta^r(x/t)) = w^r(x/t), \\[3mm]
		z_1(v^r(x/t), u^r(x/t)) = z_1(v_L, u_L) = z_1(v_R, u_R), \\ [3mm] 
		s(v^r(x/t), \theta^r(x/t)) = s(v_L, \theta_L) = s(v_R, \theta_R), 
	\end{array}
\end{align}
and 
\begin{align}  \label{RV-phi}
	\phi^r(x/t) = - \ln v^r(x/t),
\end{align}
where $z_1(v,u)=u+\displaystyle{\int^v} \sigma_1(v^\prime, s) dv^\prime$ and $s$ are two 1-Riemann invariants to \reff{eq:Eulermass}-\reff{eq:Eulerenergy}.
It is easy to check that the approximate rarefaction wave $(v^r, u^r, \theta^r, \phi^r)(t,x)$ satisfies the quasi-neutral Euler equation 
\begin{align} \label{RV-equ}
	\left\{
	\begin{array}{l}
		v^r_t - u^r_x = 0, \\[2mm]
		u^r_t + \left(p(v^r,\theta^r)+\dfrac{1}{v^r} \right)_x = 0, \\[3mm]
		\dfrac{R}{\gamma - 1} \theta^r_t + p(v^r,\theta^r) u^r_x = 0, \\[4mm]
		\phi^r = -\ln v^r.
	\end{array}
	\right.
\end{align}
The same approach for third eigenvalue $\sigma_3$ yields 3-rarefaction curve and 3-rarefaction wave.
\\

\subsection{Contact discontinuity}
\hspace*{\fill} \\

The linearly degenerate eigenvalue $\sigma_2$ corresponds to the 2-contact discontinuity.
For some state  $(v_R, u_R, \theta_R)\in\mathbb{R}^+\times\mathbb{R}\times\mathbb{R}^+$, 
the 2-contact discontinuity curve is given by 
\begin{align} \label{CD_2_def}
	CD_2(v_R, u_R, \theta_R) := \Big\{ (v,u,\theta) \mid u = u_R, \quad p(v,\theta)+\frac{1}{v} = p(v_R,\theta_R)+\frac{1}{v_R} \Big\}.
\end{align}
The 2-contact discontinuity $(v^c, u^c, \theta^c)(t,x)$ connecting $(v_-, u_-, \theta_-)=(v_L,u_L,\theta_L)$ and $(v_+,u_+,\theta_+)=(v_R,u_R,\theta_R)$,
such that $(v_L, u_L, \theta_L)\in CD_2(v_R, u_R, \theta_R)$, is defined by 
\begin{align*}
	(v^c, u^c, \theta^c)(t,x) = \left\{ 
	\begin{array}{ll}
		(v_L, u_L, \theta_L), & x < 0, \\[3mm] 
		(v_R, u_R, \theta_R), & x > 0,
	\end{array} \right.
\end{align*}
which is a self-similar solution to the Riemann problem \reff{eq:Eulermass}-\reff{eq:Eulerenergy} and \reff{EulerID}, and 
\begin{align*}
	\phi^c(t,x)=-\ln v^c(t,x).
\end{align*} 
\\

\subsection{Inviscid shock}
\hspace*{\fill} \\

Next, we should give the pricise form of Rankine-Hugoniot condition for \reff{eq:Eulermass}-\reff{eq:Eulerenergy},
which is slightly different from the classical compressible Navier-Stokes equations.
First of all, by the mass equation \reff{eq:Eulermass} of quasi-neutral Euler system, it holds that 
\begin{align*}
	\left( \frac{1}{v} \right)_x u = \left( \frac{u}{v} \right)_x - \frac{u_x}{v} = \left( \frac{u}{v} \right)_x - \frac{v_t}{v}= \left( \frac{u}{v} \right)_x + (-\ln v)_t,
\end{align*}
which yields that \reff{eq:Eulerenergy} can be converted into  
\begin{align} \label{eq:Eulerenergy-new}
	\left( \frac{R\theta}{\gamma-1} + \frac{u^2}{2} - \ln v \right)_t + \left( \left( p(v,\theta) + \frac{1}{v} \right) u \right)_x  = 0.
\end{align}
For the Riemann problem \reff{EulerID} with $(v_+, u_+, \theta_+)=(v_R, u_R, \theta_R)$, combining with \reff{eq:Eulermass}, \reff{eq:Eulermomentum} and \reff{eq:Eulerenergy-new}, we can now obtain the Rankine-Hugoniot condition: 
\begin{align}\label{RH}
    \begin{array}{l}
        -\sigma (v_R - v) - (u_R - u) = 0, 
        \\[2mm]
        -\sigma (u_R - u) + \left[ \left(p(v_R, \theta_R) + \dfrac{1}{v_R}\right) - \left(p(v, \theta) +\dfrac{1}{v}\right) \right] = 0, 
        \\[4mm]
        -\sigma \left[ \dfrac{R(\theta_R - \theta)}{\gamma - 1} + \dfrac{u_R^2 - u^2}{2} - (\ln v_R - \ln v) \right] 
		\\[3mm]
		\qquad\qquad\qquad\qquad 
		+ \left[ \left(p(v_R, \theta_R) + \dfrac{1}{v_R}\right)u_R - \left(p(v, \theta) + \dfrac{1}{v}\right)u \right] = 0, 
    \end{array} 
\end{align}
which determines the shock speed of the form 
\begin{align*} 
	\sigma=:\sigma^\pm=\pm\sqrt{-\frac{\left( p(v_R, \theta_R) + \frac{1}{v_R} \right) - \left( p(v, {\theta}) + \frac{1}{v} \right)}{v_R - v}}.
\end{align*}

In this paper, we will later restrict ourselves to the case $\pmb{\sigma = \sigma^+}$.
The general theorey shows that this conditon determines two curves meeting at $(v_R, u_R, \theta_R)$.
The curve for $\sigma=\sigma^-$, called 1-shock curve $S_1(v_R, u_R, \theta_R)$, corresponds to admissible 1-shock for $v>v_R$, while the other for $\sigma=\sigma^+$ is 3-shock curve $S_3(v^R, u^R, \theta^R)$ with admissible 3-shock for $v<v_R$.
Whenever $(v_L, u_L, \theta_L)\in S_1(v_R, u_R, \theta_R)\cup S_3(v_R, u_R, \theta_R)$, 
the shock solution $(v^s, u^s, \theta^s)(t,x)$ to \reff{eq:Eulermass}-\reff{eq:Eulerenergy}, with $(v_-, u_-, \theta_-)=(v_L,u_L,\theta_L)$ and $(v_+,u_+,\theta_+)=(v_R,u_R,\theta_R)$,
is determined by the discontinuous travelling wave 
\begin{align*}
	(v^s, u^s, \theta^s)(t,x) = \left\{ 
	\begin{array}{ll}
		(v_L, u_L, \theta_L), & x < \sigma t, \\[3mm] 
		(v_R, u_R, \theta_R), & x > \sigma t,
	\end{array} \right.
\end{align*}
and 
\begin{align*}
	\phi^s(t,x)=-\ln v^s(t,x).
\end{align*}  

\subsection{Generic Riemann problem}
\hspace*{\fill} \\

For the generic Riemann problem with possibly nearby two states $(v_-,u_-,\theta_-), (v_+,u_+,\theta_+)\in\mathbb{R}^+\times\mathbb{R}\times\mathbb{R}^+$,  
it can be shown that there exist two unique middle states $(v_*, u_*, \theta_*), (v^*, u^*, \theta^*)\in\mathbb{R}^+\times\mathbb{R}\times\mathbb{R}^+$ such that: 
(i) $(v^*, u^*, \theta^*)$ belongs to a curve of the third family from $(v_+, u^+, \theta^+)$, either $R_3(v_+, u_+, \theta_+)$ or $S_3(v_+, u_+, \theta_+)$;
(ii) $(v_*, u_*, \theta_*)$ lies on a curve of the second family, i.e., $CD_2(v^*, u^*, \theta^*)$;
(iii) $(v_-, u_-, \theta_-)$ is located on a curve of the first family from $(v_*, u_*, \theta_*)$, either $R_1(v_*, u_*, \theta_*)$ or $S_1(v_*, u_*, \theta_*)$. 

Therefore, the Riemann solution $(v, u, \theta, \phi)$ to the quasi-neutral Euler equations \reff{eq:Euler}-\reff{EulerID} consists of a combination of the three related waves 
\begin{align} \label{Generic-sol}
	\begin{array}{l}
	(v, u, \theta, \phi)(t, x) = (v_1, u_1, \theta_1, -\ln v_1)(t, x) + (v_2, u_2, \theta_2, -\ln v_2)(t, x) + (v_3, u_3, \theta_3, -\ln v_3)(t, x) \\[2mm]
	\qquad\qquad\qquad\qquad - (v_*, u_*, \theta_*, -\ln v_*) - (v^*, u^*, \theta^*, -\ln v^*).
	\end{array}
\end{align}
The solution $(v_1, u_1, \theta_1)$ corresponding to $\sigma_1$ is 1-rarefaction wave if $(v_-, u_-, \theta_-)\in R_1(v_*, u_*, \theta_*)$, or 1-shock wave if $(v_-, u_-, \theta_-)\in S_1(v_*, u_*, \theta_*)$.
Then the solution $(v_2, u_2, \theta_2)$ associated with $\sigma_2$ is 2-contact discontinuity wave if $(v_*, u_*, \theta_*)\in CD_2(v^*, u^*, \theta^*)$.
Finally, the wave $(v_3, u_3, \theta_3)$ related to $\sigma_3$ is 3-rarefaction wave if $(v^*, u^*, \theta^*)\in R_3(v_+, u_+, \theta_+)$, or 3-shock wave if $(v^*, u^*, \theta^*)\in S_3(v_+, u_+, \theta_+)$.
This composite wave includes the pattern of three single waves, i.e., shock, rarefaction and contact discontinuity, since the single wave structure corresponds to its degenerate cases.
Without loss of generality, we will consider the generic Riemann problem with two states $(v_-,u_-,\theta_-)$ and $(v_+,u_+,\theta_+)$ such that 
\begin{align} \label{Generic-Rie}
	\begin{array}{l}
		(v_-, u_-, \theta_-)\in R_1(v_*, u_*, \theta_*), \\[2mm] 
		(v_*, u_*, \theta_*)\in CD_2(v^*, u^*, \theta^*), \\[2mm] 
		(v^*, u^*, \theta^*)\in S_3(v_+, u_+, \theta_+). \\[2mm]
	\end{array}
\end{align}

\section{Preliminaries} 
\hspace*{\fill} \\

As also shown above, the Riemann problem for the inviscid system \reff{eq:Euler}-\reff{EulerID} consists of three types of waves, i.e., shock, contact discontinuity and rarefaction wave, of which the first two are discontinuous and the last is only Lipschitz continuous. 
Since the rarefaction wave is not smooth enough for later use, we describe instead its approximate version, together with the viscous counterparts of contact discontinuity and shock for the CNSP system\reff{eq:CNSP}-\reff{CNSPID}.
In what follows, we gather the well-known results that will be used in the rest of the paper.
\\ 

\subsection{Approximate rarefaction wave}
\hspace*{\fill} \\

Consider two states $(v_-, u_-,\theta_-)$ and $(v_*, u_*, \theta_*)$ such that $(v_-, u_-, \theta_-)\in R_1(v_*, u_*, \theta_*)$, i.e., $\reff{Generic-Rie}_1$.
Set 
\begin{align*} 
	w_-=\sigma_1(v_-, \theta_-)=:\sigma_{1-}, \quad w_*=\sigma_1(v_*, \theta_*)=:\sigma_{1*}.
\end{align*}
When $w_-<w_*$, the smooth approximate 1-rarefaction wave $(v^R, u^R, \theta^R, \phi^R)(t, x)$ of the 1-rarefaction wave fan $(v^r, u^r, \theta^r, \phi^r)(x/t)$ defined in \reff{RV}-\reff{RV-w} is given by 
\begin{align} \label{ARV-sol}
	\begin{array}{l}
		\sigma_1(v^R, \theta^R)(t,x)=w(1+t,x), \\[2mm]
		z_1(v^R, u^R)(t,x)=z_1(v_-,u_-)=z_1(v_*,u_*), \\[2mm]
		s(v^R, \theta^R)(t,x)=s(v_-,\theta_-)=s(v_*,\theta_*), \\[2mm]
		\phi^R(t,x)=-\ln v^R(t,x),
	\end{array}
\end{align}
where $w(t,x)$ is the smooth solution to the Burgers equation
\begin{align*}
    \left\{ 
       \begin{array}{l}
            w_t + w w_x = 0, 
            \\ [3mm]
            w_0(x) = w(0, x) : = \dfrac{w_*+w_-}{2} + \dfrac{w_*-w_-}{2} \tanh x.
       \end{array}
    \right.
\end{align*}

On the other hand, it is easy to check that the approximate rarefaction wave $(v^R, u^R, \theta^R, \phi^R)(t,x)$ satisfies the quasi-neutral Euler equation 
\begin{align} \label{ARV}
	\left\{
	\begin{array}{l}
		v^R_t - u^R_x = 0, \\[2mm]
		u^R_t + \left(p(v^R,\theta^R)+\dfrac{1}{v^R} \right)_x = 0, \\[3mm]
		\dfrac{R}{\gamma - 1} \theta^R_t + p(v^R,\theta^R) u^R_x = 0, \\[4mm]
		\phi^R = -\ln v^R.
	\end{array}
	\right.
\end{align}
From \cite{Matsumura92}, the following lemma holds.
\\ 

\begin{lemma} \label{ARV-pro}
The smooth approximate 1-rarefaction wave $(v^R, u^R, \theta^R, \phi^R)(t, x)$ defined in \reff{ARV-sol} satisfies the following properties. 
Denote $\delta_R$ the rarefaction wave strength as $\delta_R := |v_* - v_-| \sim |u_* - u_-| \sim |\theta_* - \theta_-|$.
\\ 

	(1) $u^R_x = \dfrac{2 v^R (\gamma R\theta^R + 1)}{(\gamma+1)\gamma R\theta^R +2} w_x > 0$, 
		$v^R_x = \dfrac{v^R}{\sqrt{\gamma R\theta^R+1}} u^R_x > 0$,  
		$\theta^R_x = -\dfrac{(\gamma-1)\theta^R}{v^R} v^R_x < 0$, and\\
		
		$\qquad \phi^R = -\ln v^R$,
		$\quad \forall x \in \mathbb{R}, \ t \ge 0.$

	(2) The following estimates hold for all $t \ge 0$ and $p \in [1, +\infty]$:
	\begin{align*}
		\begin{array}{l}
			\|(v^R_x, u^R_x, \theta^R_x, \phi^R_x)\|_{L^p(\mathbb{R})} \le C \min\{\delta_R, \delta_R^{1/p} (1+t)^{-1+1/p}\}, \\[2mm]
			\|(v^R_{xx}, u^R_{xx}, \theta^R_{xx}, \phi^R_{xx})\|_{L^p(\mathbb{R})} \le C \min\{\delta_R, (1+t)^{-1}\}, \\[2mm]
			|(v^R_{xx}, u^R_{xx}, \theta^R_{xx}, \phi^R_{xx})| \le C |(v^R_x, u^R_x, \theta^R_x, \phi^R_x)|,  \quad \forall x \in \mathbb{R}.
		\end{array}
	\end{align*}

	(3) For $x \ge \sigma_{1*}(1+t), \ t \ge 0$, it holds that
	\begin{align*}
		\begin{array}{l}
			|(v^R, u^R, \theta^R, \phi^R)(t,x) - (v_*, u_*, \theta_*, -\ln v_*)| \le C \delta_R e^{-2|x - \sigma_{1*}(1+t)|} \\[2mm]
			|(v^R_x, u^R_x, \theta^R_x, \phi^R_x)(t,x)| \le C \delta_R e^{-2|x - \sigma_{1*}(1+t)|}.
		\end{array}
	\end{align*}

	(4) For $x \le \sigma_{1-}(1+t), \ t \ge 0$, it holds that
	\begin{align*}
		\begin{array}{l}
			|(v^R, u^R, \theta^R, \phi^R)(t,x) - (v_-, u_-, \theta_-, -\ln v_-)| \le C \delta_R e^{-2|x - \sigma_{1-}(1+t)|} \\[2mm]
			|(v^R_x, u^R_x, \theta^R_x, \phi^R_x)(t,x)| \le C \delta_R e^{-2|x - \sigma_{1-}(1+t)|}.
		\end{array}
	\end{align*}

	(5) $\underset{t \to +\infty}{\lim} \underset{x \in \mathbb{R}}{\sup} |(v^R, u^R, \theta^R, \phi^R)(t,x) - (v^r, u^r, \theta^r, \phi^r)(x/t)| = 0.$
\end{lemma}
\hspace*{\fill}

\subsection{Viscous contact wave}
\hspace*{\fill} \\

The viscous contact wave $(v^C, u^C, \theta^C, \phi^C)(t,x)$, i.e., the viscous version of the contact discontinuity connecting $(v_*, u_*, \theta_*)$ and $(v^*, u^*, \theta^*)$ such that $\reff{Generic-Rie}_2$ holds,
can be defined \cite{Huang08,Zhao26} by 
\begin{align} \label{equ:VCD-ori}
	\begin{array}{l}
		p(v^C, \theta^C) + \dfrac{1}{v^C} = p(v_*, \theta_*) + \dfrac{1}{v_*} = p(v^*, \theta^*) + \dfrac{1}{v^*}, \\[3mm]
		v^C_t - u^C_x = 0, \\[1mm]
		\dfrac{R}{\gamma - 1} \theta^C_t + p(v^C, \theta^C) u^C_x = \kappa \left( \dfrac{\theta^C_x}{v^C} \right)_x, \\[3mm]
		\phi^C = - \ln v^C,
	\end{array}
\end{align}
with the far-field condition 
\begin{align*}
	\theta^C(t, -\infty) = \theta_*, \quad \theta^C(t, +\infty) = \theta^*.
\end{align*}
Then we derive the absolute temperature $\theta^C$ satisfying
\begin{align} \label{theta_H}
	A(\theta^C) \theta_t^C = (B(\theta^C) \theta^C_x)_x,
\end{align}
where
\begin{align*}
	A(\theta^C) = \Big( \frac{1}{\gamma-1} + \frac{R \theta^C}{R \theta^C + 1} \Big) R, \quad B(\theta^C) = \frac{\kappa (p(v_*, \theta_*)+\dfrac{1}{v_*})}{R \theta^C + 1}.
\end{align*}
Define $H: \mathbb{R}^+ \to \mathbb{R}^+, \,\, x \mapsto \frac{\gamma R}{\gamma-1} x - \ln (R x + 1)$,
then $H^\prime(x) = ( \frac{1}{\gamma-1} + \frac{R x}{R x + 1} ) > 0$ implies the inverse function of $H(x)$.
Let 
\begin{align*}
\Lambda(t,x) = H(\theta^C(t,x)), \quad H_*: = H(\theta_*), \,\, H^*: = H(\theta^*),
\end{align*}
Then the equation \reff{theta_H} can be converted into
\begin{align*}
	\left\{
		\begin{array}{l}
			\Lambda_t = \left( \dfrac{B(H^{-1}(\Lambda))}{H^\prime (H^{-1}(\Lambda))} \Lambda_x \right)_x, \\[4mm]
			\Lambda(t, -\infty) = H_*, \Lambda(t, +\infty) = H^*.
		\end{array}
	\right.
\end{align*}
It follows from \cite{Hsiao93,van76} that the two-points boundary value problem admits a unique self-similar solution $\Lambda(t,x)=\Lambda(x/\sqrt{1+t})$ satisfying 
\begin{align*}
	(1+t)^{\frac{n}{2}}|\partial_x^n \Lambda| + |\Lambda-H_*| + |\Lambda-H^*| \le O(1) |H^*-H_*|\, e^{-\frac{C_1 x^2}{1+t}}, \quad x\rightarrow \pm\infty,
\end{align*}
with generic constant $C_1>0$. Therefore, we have 
\\ 

\begin{lemma} \label{VCD}
	Let $\delta_C$ denote the viscous contact wave strength as $\delta_C := |v^* - v_*| \sim |u^* - u_*| \sim |\theta^* - \theta_*|$.
	Then the viscous contact wave $(v^C,u^C,\theta^C,\phi^C)$ defined by the system \reff{equ:VCD-ori} satisfies
	\begin{align*}
		\begin{array}{l}
			|(v^C - v_*, u^C - u_*, \theta^C - \theta_*, \phi^C - (-\ln v_*))| =O(1) \delta_C \, e^{-\frac{C_1 x^2}{1+t}}, \quad \forall x < 0, \\[2mm]
			|(v^C - v^*, u^C - u^*, \theta^C - \theta^*, \phi^C - (-\ln v^*))| =O(1) \delta_C \, e^{-\frac{C_1 x^2}{1+t}}, \quad \forall x > 0, \\[3mm]
			|\partial_x^n (v^C, \theta^C, \phi^C)| =O(1) \delta_C \, (1+t)^{-\frac{n}{2}} \, e^{-\frac{C_1 x^2}{1+t}}, \quad \forall x \in \mathbb{R}, \, n \ge 1, \\[4mm]
			|\partial_x^n u^C| =O(1) \delta_C \, (1+t)^{-\frac{1+n}{2}} \, e^{-\frac{C_1 x^2}{1+t}}, \quad \forall x \in \mathbb{R}, \, n \ge 1, 
		\end{array}
	\end{align*}
	where $C_1>0$ is a generic constant.
\end{lemma}
\hspace*{\fill}

Then the viscous contact wave $(v^C, u^C, \theta^C, \phi^C)(t, x)$ defined in \reff{equ:VCD-ori} satisfies 
\begin{align} \label{eq:VCD}
	\left\{
	\begin{array}{l}
		v_t^C - u_x^C = 0, \\[2mm]
		u_t^C + \left( p(v^C, \theta^C) + \dfrac{1}{v^C} \right)_x = \mu \left( \dfrac{u_x^C}{v^C} \right)_x + \Phi(v^C, \phi^C)_x + Q_1^C, \\[3mm]
		\dfrac{R}{\gamma - 1} \theta_t^C + p(v^C, \theta^C) u_x^C = k \left( \dfrac{\theta_x^C}{v^C} \right)_x + \mu \dfrac{(u_x^C)^2}{v^C} + Q_2^C, \\[3mm]
		\phi^C = -\ln v^C,
	\end{array}
	\right.
\end{align}
where
\begin{align} \label{Q1C-Q2C}
	\begin{array}{l}
		Q_1^C = u_t^C - \mu \left( \dfrac{u_x^C}{v^C} \right)_x - \Phi(v^C, \phi^C)_x = O(1) \delta_C (1+t)^{-\frac{3}{2}} e^{-\frac{C_1 x^2}{1+t}}, \\[4mm]
		Q_2^C = -\mu \dfrac{(u_x^C)^2}{v^C} = O(1) \delta_C (1+t)^{-2} e^{-\frac{2 C_1 x^2}{1+t}},
	\end{array}
\end{align}
as $x \to \pm\infty$, due to Lemma \ref{VCD}.
\\ 

\subsection{Viscous shock}
\hspace*{\fill} \\

Next, we will consider the viscous 3-shock, the viscous counterpart of the inviscid 3-shock connecting two distinct states $(v^*, u^*, \theta^*)$ and $(v_+, u_+, \theta_+)$ such that $\reff{Generic-Rie}_3$ holds.
Before that, we need to determine the necessary relative size relation between two different constant states $(v, u, \theta)=(v^*, u^*, \theta^*)$ and $(v_R, u_R, \theta_R)=(v_+, u_+, \theta_+)$ in the Rankine-Hugoniot condition \reff{RH}.
By \reff{RH}, adding $(-u_+)$ times the second equality to the third equality, along with the first equality, gives 
\begin{align*}
	\frac{R(\theta_+ - {\theta^*})}{\gamma - 1} = (\ln v_+ - \ln v^*) - \left(p(v, \theta) + \frac{1}{v}\right)(v_+ - v^*) + \frac{\sigma^2}{2}(v_+ - v^*)^2,
\end{align*}
and then applying the Taylor expansion near the point $v^*$ yields
\begin{align*}
	\frac{R(\theta_+ - {\theta^*})}{\gamma - 1} = - p(v^*, {\theta^*})(v_+ - v^*) + O((v_+ - v^*)^2),
\end{align*}
which implies that $\theta_+-{\theta^*}$ and $v_+-v^*$ have opposite signs.
By the Lax entropy condition $\sigma_3(v_+, \theta_+) \le \sigma \le \sigma_3(v^*, {\theta^*})$ ($\sigma_3$ defined in \reff{Eigen}), it holds that 
\begin{align} \label{v_symbol}
	v^*< v_+,\,\, {\theta^*}>\theta_+,\,\, u^* > u_+,
\end{align}
the last of which is due to the first equality of \reff{RH}. 

Then we investigate the viscous shock profile, a special class of traveling wave solutions of the form $(v, u, \theta, \phi)(\xi)$ to CNSP \reff{eq:CNSP}, where $\xi=x-\sigma t$ and 
\begin{align} \label{sig-def}
	\sigma=\sqrt{-\frac{\left( p(v_+, \theta_+) + \frac{1}{v_+} \right) - \left( p(v^*, \theta^*) + \frac{1}{v^*} \right)}{v_+ - v^*}}
\end{align} 
is the 3-shock speed defined by replacing $(v, u, \theta)$ and $(v_R, u_R, \theta_R)$ in the Rankine-Hugoniot condition \reff{RH} by $(v^*, u^*, \theta^*)$ and $(v_+, u_+, \theta_+)$, respectively.
The viscous shock profile is governed by the ODEs
\begin{subequations}\label{eq:viscous_shock}
	\begin{align}
		&-\sigma v_\xi - u_\xi = 0, \label{eq:shock_a} \\[2mm]
		& -\sigma u_\xi + p_\xi = \mu \left(\frac{u_\xi}{v}\right)_\xi - \frac{\phi_\xi}{v}, \label{eq:shock_b} \\[2mm]
		& -\sigma \left(\frac{R\theta}{\gamma - 1} + \frac{u^2}{2}\right)_\xi + (pu)_\xi = \kappa \left(\frac{\theta_\xi}{v}\right)_\xi + \mu \left(\frac{u u_\xi}{v}\right)_\xi - \frac{\phi_\xi}{v} u, \label{eq:shock_c} \\[2mm]
		& -\nu^2 \left(\frac{\phi_\xi}{v}\right)_\xi = 1 - v e^\phi, \label{eq:shock_d}
	\end{align}
\end{subequations}
verifying the far-field condition
\begin{align} \label{far-field}
	(v, u, \theta, \phi)(-\infty)=(v^*, u^*, \theta^*, \phi^*), \qquad 
	(v, u, \theta, \phi)(+\infty)=(v_+, u_+, \theta_+, \phi_+)
\end{align}
with $\phi^*=-\ln v^*, \,\, \phi_+=-\ln v_+.$

Note that, after a little calculation, we can derive a simpler form of the equations satisfied by the viscous shock wave.
First, integrating \reff{eq:shock_a} over $[-\infty, \xi]$ and $[\xi, +\infty]$ respectively yields
\begin{align}\label{u_relation}
	u = (u^* + \sigma v^*) - \sigma v = (u_+ + \sigma v_+) - \sigma v.
\end{align}
From the Poisson equation \reff{eq:shock_d}, we derive the relation
\begin{align}\label{phi_relation}
	\frac{\phi_\xi}{v} = - \left[ \frac{\nu^2}{2} \left(\frac{\phi_\xi}{v}\right)^2 - e^\phi \right]_\xi.
\end{align}
Plugging \reff{phi_relation} back into \reff{eq:shock_b}, integrating over $[-\infty, \xi]$, and applying \reff{u_relation} gives 
\begin{align}\label{eq:v_xi}
	v_\xi = - \frac{v}{\mu \sigma} \left[ \sigma^2 (v - v^*) + (p(v, \theta) - p(v^*, \theta^*)) - \frac{\nu^2}{2} \left(\frac{\phi_\xi}{v}\right)^2 + (e^\phi - e^{\phi^*}) \right].
\end{align}

With \reff{u_relation} and \reff{phi_relation}, it holds that 
\begin{align} \label{phi_u_relation}
	\begin{array}{ll} 
		\dfrac{\phi_\xi}{v} u 
	 	&= \dfrac{\phi_\xi}{v} \left( (u^* + \sigma v^*) - \sigma v \right) = \dfrac{\phi_\xi}{v} (u^* + \sigma v^*) - \sigma \phi_\xi \\[3mm]
	 	&= - (u^* + \sigma v^*) \left[ \dfrac{\nu^2}{2} \left(\dfrac{\phi_\xi}{v}\right)^2 - e^\phi \right]_\xi - \sigma \phi_\xi.
	\end{array}
\end{align}
Substituting the above relation into \reff{eq:shock_c} and integrating over $[-\infty, \xi]$ yields
\begin{align*}
	\begin{array}{l}
		-\sigma \left( \dfrac{R(\theta - {\theta^*})}{\gamma - 1} + \dfrac{u^2 - (u^*)^2}{2} \right) + ( p(v,\theta) u - p(v^*, {\theta^*})u^* ) = \kappa \dfrac{\theta_\xi}{v} + \mu \dfrac{u u_\xi}{v} \\[3mm]
		\qquad\qquad + (u^* + \sigma v^*) \left[ \dfrac{\nu^2}{2} \left(\dfrac{\phi_\xi}{v}\right)^2 - e^\phi + e^{\phi^*} \right] + \sigma (\phi - \phi^*).
	\end{array}
\end{align*}
Then replacing $u_\xi$ by \reff{eq:shock_a} and \reff{eq:v_xi}, it follows that 
\begin{align}\label{eq:theta_xi}	
	\begin{array}{l}
		\theta_\xi =
		- \dfrac{\sigma v}{\kappa} \left[ \dfrac{R(\theta - {\theta^*})}{\gamma - 1} - \dfrac{\sigma^2 }{2} (v - v^*)^2 + p(v^*, {\theta^*})(v - v^*) + \dfrac{\nu^2}{2} \dfrac{(\phi_\xi)^2}{v}  \right. \\[3mm]
		\quad \left. \vphantom{\frac{R(\theta - {\theta^*})}{\gamma - 1}}\qquad\qquad 
		- v (e^\phi - e^{\phi^*}) + (\phi - \phi^*) \right].
	\end{array}
\end{align}
Finally, Substituting \reff{eq:v_xi} into the Poisson equation \reff{eq:shock_d} leads to 
\begin{align}\label{eq:phi_xi_xi}
	\phi_{\xi\xi} 
	= - \frac{v}{\nu^2} (1 - v e^\phi) - \frac{\phi_\xi}{\mu \sigma} \left[ \sigma^2 (v - v^*) + (p - p(v^*, {\theta^*})) - \frac{\nu^2}{2} \left(\frac{\phi_\xi}{v}\right)^2 + (e^\phi - e^{\phi^*}) \right].
\end{align}
Then the existence of viscous shock profile to the ODEs \reff{eq:viscous_shock} is equivalent to the existence of solution to \reff{eq:v_xi}-\reff{eq:phi_xi_xi} with \reff{u_relation}, which will be proved in Appendix A.
The properties of the 3-viscous shock wave can be listed as follows.
\\

\begin{lemma} \label{VS-pro} 
	For any state $(v^*, u^*, \theta^*)$ such that $(v^*, u^*, \theta^*)\in S_3(v_+, u_+, \theta_+)$.
	Let $\delta_S$ denote the viscous shock wave strength as $\delta_S := |v_+ - v^*| \sim |u_+ - u^*| \sim |\theta_+ - \theta^*|$. 
	Then there exists a unique (up to a shift) shock profile $(v^S, u^S, \theta^S, \phi^S)(\xi)$ to \reff{eq:viscous_shock}-\reff{far-field},
	and it holds that 
	\begin{align}\label{exist_all}
		\begin{array}{lll}
			&v^S_\xi > 0,\ u^S_\xi < 0,\ \theta^S_\xi < 0,\ \phi^S_\xi < 0, &\quad \forall \xi \in \mathbb{R},  \\[3mm]
			&|v^S_\xi| \sim |u^S_\xi| \sim |\theta^S_\xi| \sim |\phi^S_\xi|, &\quad \forall \xi \in \mathbb{R},  \\[3mm]
			&\left| ({v^S} - v^*, u^S - u^*, {\theta^S} - \theta^*, \phi^S - \phi^* )\right| \le C \delta_S e^{-C \delta_S |\xi|}, &\quad \forall \xi < 0,\\[3mm]
			&\left| ({v^S} - v_+, u^S - u_+, {\theta^S} - \theta_+, \phi^S - \phi_+ )\right| \le C \delta_S e^{-C \delta_S |\xi|}, &\quad \forall \xi > 0,\\[3mm]
			&\left| (v^S_\xi, u^S_\xi, \theta^S_\xi, \phi^S_\xi) \right| \le C \delta_S^2 e^{-C \delta_S |\xi|}, &\quad \forall \xi \in \mathbb{R}, \\[3mm]
			&\left| ({v^S}_{\xi\xi}, u^S_{\xi\xi}, {\theta^S}_{\xi\xi}, \phi^S_{\xi\xi}) \right| \le C \delta_S \left| (v^S_\xi, u^S_\xi, \theta^S_\xi, \phi^S_\xi) \right|, &\quad \forall \xi \in \mathbb{R}. \\[2mm]
		\end{array}
	\end{align}
	Moreover, the solution satisfies $v^S_\xi(\xi=0)=C\delta_S^2$ and $\partial_{\xi}^n\phi^S(\xi=0)=C\delta_S^{1+n}$. In particular, it yields
	\begin{align} \label{exist_uv}
		|u^S_\xi + \sigma^* v^S_\xi| \le C\delta_S |v^S_\xi|, \qquad\qquad  \forall \xi\in\mathbb{R},
	\end{align}
	and 
	\begin{align} \label{exist_theta_v}
		\left|\theta^S_\xi + \dfrac{(\gamma-1)p(v^*, {\theta^*})}{R} v^S_\xi \right| \le C\delta_S |v^S_\xi|, \qquad\qquad \forall \xi \in \mathbb{R},
	\end{align}
	where 
	\begin{align*} 
		\sigma^*:=\sigma_3(v^*, {\theta^*})= \sqrt{\dfrac{\gamma R{\theta^*} + 1}{(v^*)^2}}
	\end{align*}
	which satisfies
	\begin{align} \label{sig_sig-}
		|\sigma-\sigma^*|\le C\delta_S.
	\end{align}
\end{lemma}

\begin{remark}
	The presence of the electric potential introduces an additional non-conservative source term into the momentum equation \reff{eq:shock_b}. 
	Note that, however, the electric potential satisfies the Poisson equation \reff{eq:shock_d}, this term can be equivalently transformed into a conservative form \reff{phi_relation}. 
	By combining this with the center manifold theorem and the solvability of the associated simple ODEs \reff{eta_xi} (see \cite{Duan20,Majda85}), one obtains the solvability of ODEs \reff{eq:v_xi}-\reff{eq:phi_xi_xi} with \reff{u_relation}. 
	The above proof relies primarily on the smallness of the shock strength and the Taylor expansions near the critical points.
\end{remark}
\hspace*{\fill}

We proceed to obtain sharp estimates for the diffusion, which are central to the relative entropy estimates that follow and exploit the key properties of the viscous shock profiles.
The proof is later given in Appendix B.
\\  

\begin{lemma} \label{Crucial}
	For the viscous shock profile $({v^S}, u^S, {\theta^S}, \phi^S)$ (defined in Lemma \ref{VS-pro}) connecting 
	two end states $(v^*, u^*, \theta^*, \phi^*)$ and $(v_+, u_+, \theta_+, \phi_+)$ with the quasi-neutral condition \reff{quasi-neutral}, define 
	\begin{align*}
		\hat{p}:=\hat{p}({v^S}, {\theta^S})=p({v^S}, {\theta^S})+\frac{1}{{v^S}}=\frac{R{\theta^S} + 1}{{v^S}}, \quad 
		\hat{p}_+:=\hat{p}(v_+, \theta_+),\quad \hat{p}^*:=\hat{p}(v^*, {\theta^*}),
	\end{align*} 
	and let $\delta_S:=v_+-v^*>0$ (by \reff{v_symbol}) denote the viscous shock wave amplitude.
	Then it holds that 
	\begin{align} \label{Crucial_est}
		\left| \frac{\hat{p}-\hat{p}_+}{{v^S} - v_+} - \frac{\hat{p}-\hat{p}^*}{{v^S} - v^*} - \sigma^* \alpha^* \beta^* \delta_S \right| \leq C \delta_S^2,
	\end{align}
	where 
	\begin{align*}
		\sigma^* := \sigma_3(v^*, {\theta^*}) = \sqrt{\frac{\gamma R {\theta^*} + 1}{(v^*)^2}}, 
		\quad\quad \alpha^* := \left( \frac{\gamma(\gamma+1)p(v^*, {\theta^*})}{2} + \frac{1}{v^*} \right) \frac{1}{(v^*)^2 \sigma^*},
	\end{align*}
	and 
	\begin{align*}
		\beta^*  := \dfrac{R \mu \left(\gamma p(v^*, {\theta^*}) + \frac{1}{v^*}\right)}{\kappa(\gamma-1)^2 p(v^*, {\theta^*}) + R \mu \left(\gamma p(v^*, {\theta^*}) + \frac{1}{v^*}\right)}.
	\end{align*}
\end{lemma}
\hspace*{\fill}

\begin{remark}
	The optimal estimate of the diffusion term is determined by the second derivative of the temperature with respect to the volume \reff{equi_form_p}, as revealed by a Taylor expansion. 
	This ultimately depends on the derivatives of the electric potential. We emphasize that the presence of the electric potential complicates this computation. 
	A direct calculation may exhibit singularities at the far-field state, since the singular terms appear to admit no equivalent quantity in terms of $\delta_S$ (see \eqref{phi_v}-\eqref{phi_v2}), which would require the precise values of the exponents of the equivalent forms. 
	We overcome this by using the estimates $v^S_\xi(\xi=0)=C\delta_S^2$ and $\partial_{\xi}^n\phi^S(\xi=0)=C\delta_S^{1+n}$ in Lemma \ref{VS-pro}: a Taylor expansion shifts $v^*$ to $v^S(\xi=0)$ up to an error of order $\delta_S$, thereby eliminating the singularity.
\end{remark}
\hspace*{\fill}

\section{Proof of Theorem \ref{Theorem-Asym}}
\hspace*{\fill} \\ 

\subsection{Local in time Estimates}
\hspace*{\fill} \\

We first rewrite the CNSP system \reff{eq:CNSP} in a divergence form, which is curcial in our later analysis.
From the Poisson equation \reff{eq:CNSPpoisson}, we derive the following relation 
\begin{align} \label{ele-add}
	\dfrac{\phi_x}{v}
	=\left(\dfrac{1}{v} + \dfrac{\nu^2}{v}\left(\dfrac{\phi_x}{v}\right)_x -\dfrac{\nu^2}{2} \left(\dfrac{\phi_x}{v}\right)^2\right)_x
	=:\left(\dfrac{1}{v}-\Phi(v,\phi)\right)_x. 
\end{align}
Plugging it into the CNSP system \reff{eq:Euler} yields 
\begin{align*}
	\begin{array}{l}
	v_t - u_x = 0, \\[2mm]
	u_t + \hat{p}(v, \theta)_x = \mu \left( \dfrac{u_x}{v} \right)_x + \Phi(v, \phi)_x, \\[2mm]
	\dfrac{R}{\gamma - 1} \theta_t + p(v, \theta) u_x = \kappa \left( \dfrac{\theta_x}{v} \right)_x + \mu \dfrac{u^2_x}{v}, \\[2mm]
	- \nu^2 \left( \dfrac{\phi_x}{v} \right)_x = 1 - v e^\phi.
	\end{array}
\end{align*}
where 
\begin{align} \label{hat-p-Phi}
	\hat{p}(v,\theta) : = p(v,\theta)+\dfrac{1}{v}, 
	\qquad 
	\Phi(v, \phi) : = \dfrac{\nu^2}{2} \left(\dfrac{\phi_x}{v}\right)^2-\dfrac{\nu^2}{v}\left(\dfrac{\phi_x}{v}\right)_x. 
\end{align}
Based on the change of variable associated to the speed of propagation of the shock $(t, x) \mapsto (t, \xi=x-\sigma t)$ ($\sigma$ defined in \reff{sig-def}), 
the above CNSP system is rewritten as 
\begin{align} \label{eq:CNSP-change}
	\left\{
	\begin{array}{l}
		v_t - \sigma v_\xi - u_\xi = 0, \\[2mm]
		u_t - \sigma u_\xi + \hat{p}(v, \theta)_\xi = \mu \left( \dfrac{u_\xi}{v} \right)_\xi + \Phi(v, \phi)_\xi, \\[2mm]
		\dfrac{R}{\gamma-1} \theta_t - \dfrac{\sigma R}{\gamma-1} \theta_\xi + p(v, \theta) u_\xi = \kappa \left( \dfrac{\theta_\xi}{v} \right)_\xi + \mu \dfrac{u_\xi^2}{v}, \\[2mm]
		- \nu^2 \left( \dfrac{\phi_\xi}{v} \right)_\xi = 1 - v e^\phi.
	\end{array}
	\right.
\end{align}
It follows from \reff{ARV} that the approximate rarefaction wave $(v, u, \theta, \phi)(t, \xi) = (v^R, u^R, \theta^R, \phi^R) (t, \xi+\sigma t)$ satisfies
\begin{align} \label{eq:ARV-change}
	\left\{
	\begin{array}{l}
		v_t - \sigma v_\xi - u_\xi = 0, \\[2mm]
		u_t - \sigma u_\xi + \hat{p}(v,\theta)_\xi = 0, \\[3mm]
		\dfrac{R}{\gamma - 1} \theta_t - \dfrac{\sigma R}{\gamma-1} \theta_\xi + p(v,\theta) u_\xi = 0, \\[4mm]
		\phi = -\ln v.
	\end{array}
	\right.
\end{align}
and from \reff{eq:VCD} that the viscous contact wave $(v, u, \theta, \phi)(t, \xi)=(v^C, u^C, \theta^C, \phi^C)(t, \xi+\sigma t)$ verifies
\begin{align} \label{eq:VCD-change}
	\left\{
	\begin{array}{l}
	v_t - \sigma v_\xi - u_\xi = 0, \\[2mm]
	u_t - \sigma u_\xi + \hat{p}(v, \theta)_\xi = \mu \left( \dfrac{u_\xi}{v} \right)_\xi + \Phi(v, \phi)_\xi + Q_1^C, \\[3mm]
	\dfrac{R}{\gamma-1} \theta_t - \dfrac{\sigma R}{\gamma-1} \theta_\xi + p(v, \theta) u_\xi = \kappa \left( \dfrac{\theta_\xi}{v} \right)_\xi + \mu \dfrac{u_\xi^2}{v} + Q_2^C, \\[4mm]
	\phi = - \ln v,
	\end{array}
	\right.
\end{align}
with the error terms $Q^C_1, Q^C_2$ defined in \reff{Q1C-Q2C}.
Using \reff{eq:viscous_shock} and \reff{ele-add}, the viscous shock shifted by $\mathbf{X}(t)$ (defined in \reff{Xt-def}) $(v^S, u^S, \theta^S, \phi^S)(\xi- \mathbf{X}(t))$ satisfies the equations
\begin{align}  \label{eq:VShock-change}
	\left\{
	\begin{array}{l}
		v_t^S - \sigma v_\xi^S + \mathbf{\dot{X}}(t) v_\xi^S - u_\xi^S = 0, \\[2mm]
		u_t^S - \sigma u_\xi^S + \mathbf{\dot{X}}(t) u_\xi^S + \hat{p}(v^S, \theta^S)_\xi = \mu \left( \dfrac{u_\xi^S}{v^S} \right)_\xi + \Phi(v^S, \phi^S)_\xi, \\[2mm]
		\dfrac{R}{\gamma-1} \theta_t^S - \dfrac{\sigma R}{\gamma-1} \theta_\xi^S + \mathbf{\dot{X}}(t) \dfrac{R}{\gamma-1} \theta_\xi^S + p(v^S, \theta^S) u_\xi^S = \kappa \left( \dfrac{\theta_\xi^S}{v^S} \right)_\xi + \mu \dfrac{(u_\xi^S)^2}{v^S}, \\[2mm]
		- \nu^2 \left( \dfrac{\phi_\xi^S}{v^S} \right)_\xi = 1 - v^S e^{\phi^S}.
	\end{array}
	\right.
\end{align}
Later, \textbf{denote $\mathbf{(v^S, u^S, \theta^S, \phi^S)}$ be the viscous shock shifted by $\mathbf{X}(t)$}.

We will consider the stability of the solution to \reff{eq:CNSP-change} around the composite wave consisting of the approximate 1-rarefaction wave, the viscous 2-contact wave and the viscous 3-shock wave shifted by $\mathbf{X}(t)$:
\begin{align} \label{CompositeWave}
	\begin{array}{ll}
		(\bar{v}, \bar{u}, \bar{\theta}, \bar{\phi})(t, \xi) : = \Big( 
		&v^R(t, \xi + \sigma t) + v^C(t, \xi + \sigma t) + v^S(\xi - \mathbf{X}(t)) - v_* - v^*, \\[2mm]
		&u^R(t, \xi + \sigma t) + u^C(t, \xi + \sigma t) + u^S(\xi - \mathbf{X}(t)) - u_* - u^*, \\[2mm]
		&\theta^R(t, \xi + \sigma t) + \theta^C(t, \xi + \sigma t) + \theta^S(\xi - \mathbf{X}(t)) - \theta_* - \theta^*, \\[2mm]
		&\phi^R(t, \xi + \sigma t) + \phi^C(t, \xi + \sigma t) + \phi^S(\xi - \mathbf{X}(t)) - ( - \ln v_* ) - (- \ln v^*) 
		\Big).
	\end{array}
\end{align}
It follows from \reff{eq:ARV-change}-\reff{eq:VShock-change} that $(\bar{v}, \bar{u}, \bar{\theta}, \bar{\phi})$ satisfies the following equations 
\begin{align} \label{RCDS}
	\left\{
		\begin{array}{l}
		\bar{v}_t - \sigma \bar{v}_\xi + \mathbf{\dot{X}}(t) v^S_\xi - \bar{u}_\xi = 0, \\[2mm]
		\bar{u}_t - \sigma \bar{u}_\xi + \mathbf{\dot{X}}(t) u^S_\xi + \hat{p}(\bar{v}, \bar{\theta})_\xi = \mu \left( \dfrac{\bar{u}_\xi}{\bar{v}} \right)_\xi + \Phi(\bar{v}, \bar{\phi})_\xi + Q_1, \\[3mm]
		\dfrac{R}{\gamma-1} \bar{\theta}_t - \dfrac{\sigma R}{\gamma-1} \bar{\theta}_\xi + \mathbf{\dot{X}}(t) \dfrac{R}{\gamma-1} \theta^S_\xi + p(\bar{v}, \bar{\theta}) \bar{u}_\xi = \kappa \left( \dfrac{\bar{\theta}_\xi}{\bar{v}} \right)_\xi + \mu \dfrac{\bar{u}_\xi^2}{\bar{v}} + Q_2, \\[3mm]
		- \nu^2 \left( \dfrac{\bar{\phi}_\xi}{\bar{v}} \right)_\xi = 1 - \bar{v} e^{\bar{\phi}} + Q_3,
	\end{array}
	\right.
\end{align}
where the error terms are 
\begin{align} \label{Q123}
	Q_i := Q_i^I + Q_i^R + Q_i^C, \qquad i=1,2,3,
\end{align}
with wave interaction terms 
\begin{align} \label{Q123_I}
	\begin{array}{l}
		Q_1^I := \Big( \hat{p}(\bar{v}, \bar{\theta}) - \hat{p}(v^R, \theta^R) - \hat{p}(v^C, \theta^C) - \hat{p}(v^S, \theta^S) \Big)_\xi \\[2mm]
		\qquad\quad - \mu \left( \dfrac{\bar{u}_\xi}{\bar{v}} - \dfrac{u_\xi^C}{v^C} - \dfrac{u_\xi^S}{v^S} - \dfrac{u_\xi^R}{v^R} \right)_\xi \\[4mm]
		\qquad\quad - \Big( \Phi(\bar{v}, \bar{\phi}) - \Phi(v^C, \phi^C) - \Phi(v^S, \phi^S) - \Phi(v^R, \phi^R) \Big)_\xi, \\[4mm]
		Q_2^I := \Big( p(\bar{v}, \bar{\theta}) \bar{u}_\xi - p(v^R, \theta^R) u_\xi^R - p(v^C, \theta^C) u_\xi^C - p(v^S, \theta^S) u_\xi^S \Big) \\[2mm]
		\qquad - \kappa \left( \dfrac{\bar{\theta}_\xi}{\bar{v}} - \dfrac{\theta_\xi^C}{v^C} - \dfrac{\theta_\xi^S}{v^S} - \dfrac{\theta_\xi^R}{v^R} \right)_\xi \\[4mm]
		\qquad - \mu \left( \dfrac{\bar{u}_\xi^2}{\bar{v}} - \dfrac{(u_\xi^C)^2}{v^C} - \dfrac{(u_\xi^S)^2}{v^S} - \dfrac{(u_\xi^R)^2}{v^R} \right), \\[4mm]
		Q_3^I := (\bar{v} e^{\bar{\phi}} - v^S e^{\phi^S}) - \nu^2 \left( \dfrac{\bar{\phi}_\xi}{\bar{v}} - \dfrac{\phi_\xi^C}{v^C} - \dfrac{\phi_\xi^S}{v^S} - \dfrac{\phi_\xi^R}{v^R} \right)_\xi.
	\end{array}
\end{align}
and the error terms due to the inviscid approximate rarefaction wave 
\begin{align} \label{Q123_R}
	\begin{array}{l}
		Q_1^R = - \mu \left( \dfrac{u_\xi^R}{v^R} \right)_\xi - \Phi(v^R, \phi^R)_\xi, \\[2mm]
		Q_2^R = - \kappa \left( \dfrac{\theta_\xi^R}{v^R} \right)_\xi - \mu \dfrac{(u_\xi^R)^2}{v^R}, \\[2mm]
		Q_3^R = - \nu^2 \left( \dfrac{\phi_\xi^R}{v^R} \right)_\xi = \nu^2 \left( \dfrac{v_\xi^R}{(v^R)^2} \right)_\xi, \quad (\text{by} \,\, \reff{eq:ARV-change}_4).
	\end{array}
\end{align}
and the error terms associated with the viscous contact wave from \reff{Q1C-Q2C}  
\begin{align} \label{Q123_C}
	\begin{array}{l}
		Q_1^C = u_t^C - \mu \left( \dfrac{u_x^C}{v^C} \right)_x - \Phi(v^C, \phi^C)_x = O(1) \delta_C (1+t)^{-\frac{3}{2}} e^{-\frac{C_1 x^2}{1+t}}, \\[3mm]
		Q_2^C = -\mu \dfrac{(u_x^C)^2}{v^C} = O(1) \delta_C (1+t)^{-2} e^{-\frac{2 C_1 x^2}{1+t}}, \\[3mm]
		Q_3^C = - \nu^2 \left( \dfrac{\phi_\xi^C}{v^C} \right)_\xi = \nu^2 \left( \dfrac{v_\xi^C}{(v^C)^2} \right)_\xi = O(1) \delta_C (1+t)^{-1} e^{-\frac{C_1 x^2}{1+t}}.
	\end{array}
\end{align}
The last equality is due to $\reff{eq:VCD-change}_4$ and Lemma \ref{VCD}. 
From \reff{eq:CNSP-change} and \reff{RCDS}, we have the perturbation system
\begin{align} \label{PS}
	\left\{
	\begin{array}{l}
	(v - \bar{v})_t - \sigma (v - \bar{v})_\xi - \mathbf{\dot{X}}(t) v^S_\xi - (u - \bar{u})_\xi = 0, \\[2mm]
	(u - \bar{u})_t - \sigma (u - \bar{u})_\xi - \mathbf{\dot{X}}(t) u^S_\xi + (\hat{p}(v, \theta) - \hat{p}(\bar{v}, \bar{\theta}))_\xi  \\[3mm]
	\qquad\qquad\qquad\qquad\qquad\qquad = \mu \left( \dfrac{u_\xi}{v} - \dfrac{\bar{u}_\xi}{\bar{v}} \right)_\xi + \left(\Phi(v, \phi) - \Phi(\bar{v}, \bar{\phi})\right)_\xi - Q_1, \\[3mm]
	\dfrac{R}{\gamma-1} (\theta - \bar{\theta})_t - \dfrac{\sigma R}{\gamma-1} (\theta - \bar{\theta})_\xi - \mathbf{\dot{X}}(t) \dfrac{R}{\gamma-1} \theta^S_\xi + (p(v,\theta) u_\xi - p(\bar{v}, \bar{\theta}) \bar{u}_\xi) \\[3mm]
	\qquad\qquad\qquad\qquad\qquad\qquad = \kappa \left( \dfrac{\theta_\xi}{v} - \dfrac{\bar{\theta}_\xi}{\bar{v}} \right)_\xi + \mu \left( \dfrac{u_\xi^2}{v} - \dfrac{\bar{u}_\xi^2}{\bar{v}} \right) - Q_2, \\[3mm]
	- \nu^2 \left( \dfrac{\phi_\xi}{v} - \dfrac{\bar{\phi}_\xi}{\bar{v}} \right)_\xi = - e^{\bar{\phi}} \Big( \bar{v} \tilde{\phi} + (v-\bar{v})\Big) + e^{\bar{\phi}} N_p - Q_3. 
	\end{array}
	\right.
\end{align}
where we set $\tilde{\phi} :=  \phi - \bar{\phi}$ and 
\begin{align} \label{Np-def}
	\displaystyle N_p := \bar{v} (1 + \tilde{\phi} - e^{\tilde{\phi}}) + (v-\bar{v}) (1 - e^{\tilde{\phi}}).
\end{align}

We first recall local-in-time estimates for strong solutions to CNSP \reff{eq:CNSP} and so also for \reff{eq:CNSP-change}.
The proof follows from the standard iteration method (see \cite{Nash62,Matsumura80,Solonnikov76}) and is omitted here for brevity.
\\ 

\begin{proposition} \label{Local Existence}
	Let $\underline{v}, \underline{u}, \underline{\theta}, \underline{\phi}$ be smooth monotone functions such that
	\begin{align}  \label{ini_und_vu}
		\begin{array}{l}
			\underline{v}(x) = v_\pm, \quad \underline{u}(x) = u_\pm, \quad \underline{\theta}(x) = \theta_\pm, \quad \underline{\phi}(x) = \phi_\pm, \quad \text{for} \,\, |x| \ge 1.
		\end{array}
	\end{align}
	Given any constants $M_0, M_1, \underline{k}_0, \overline{k}_0, \underline{k}_1, \overline{k}_1$, with $M_1 > M_0 > 0$ and $\overline{k}_1 > \overline{k}_0 > \underline{k}_0 > \underline{k}_1 > 0$, there exists a finite time $T_0 > 0$ such that if
	\begin{align*}
		\begin{array}{l}
			\| (v_0 - \underline{v}, u_0 - \underline{u}, \theta_0 - \underline{\theta}) \|_{H^2(\mathbb{R})} + \| \phi_0 - \underline{\phi} \|_{H^3(\mathbb{R})} \le M_0, \\[3mm]
			0 < \underline{k}_0 \le v_0(x), \ \theta_0(x) \le \overline{k}_0, \quad \forall x \in \mathbb{R},
		\end{array}
	\end{align*}
	then the Cauchy problem to the CNSP system \reff{eq:CNSP} has a unique solution $(v, u, \theta, \phi)$ on $[0, T_0]$ satisfying
	\begin{align*}
		\begin{array}{l}
			v - \underline{v} \in C([0, T_0]; H^2(\mathbb{R})), \\[3mm]
			(u - \underline{u}, \theta - \underline{\theta}) \in C([0, T_0]; H^2(\mathbb{R})) \cap L^2(0, T_0; H^3(\mathbb{R})), \\[3mm]
			\phi - \underline{\phi} \in C([0, T_0]; H^3(\mathbb{R})),
		\end{array}
	\end{align*}
	and
	\begin{align*}
		\begin{array}{l}
			\| (v - \underline{v}, u - \underline{u}, \theta - \underline{\theta}) \|_{L^\infty(0, T_0; H^2(\mathbb{R}))} + \| \phi - \underline{\phi} \|_{L^\infty(0, T_0; H^3(\mathbb{R}))} \le M_1, \\[3mm]
			0<\underline{k}_1 \le v(t, x), \ \theta(t, x) \le \overline{k}_1, \quad \forall (t, x) \in [0, T_0] \times \mathbb{R}.
		\end{array}
	\end{align*}
\end{proposition}
\hspace*{\fill}  

\subsection{Construction of shift}
\hspace*{\fill} \\ 

To prove the Theorme \ref{Theorem-Asym}, we use the continuity argument based on the local-in-time existence of Proposition \ref{Local Existence} and the a priori estimates in Proposition \ref{A prior estimate}.
Since these estimates depend on the shift function $\mathbf{X}(t)$, then we give its definition. 
The definition depends on the weight function $a:\mathbb{R} \mapsto \mathbb{R}$ defined in \reff{a-def} and $||a||_{C^1(\mathbb{R})}\le 2$.
Define the shift $\mathbf{X}(t)$ as a solution to the following ODE
\begin{align} \label{Xt-def}
	\left\{
	\begin{array}{l}
		\mathbf{\dot{X}}(t) = -\dfrac{M}{\delta_S}  \displaystyle \int_{\mathbb{R}} a(\xi-\mathbf{X}(t)) \bigg[ \dfrac{\hat{p}(\bar{v}, \bar{\theta}) \, v_\xi^S(\xi-\mathbf{X}(t))}{\bar{v}} (v - \bar{v}) + u_\xi^S(\xi-\mathbf{X}(t)) (u - \bar{u}) \\[2mm]
		\qquad\qquad\qquad\qquad\qquad\qquad\qquad
		+ \dfrac{R \, \theta_\xi^S(\xi-\mathbf{X}(t))}{(\gamma - 1) \bar{\theta}} (\theta - \bar{\theta}) \bigg] d\xi, \\[2mm]
		\mathbf{X}(0) = 0,
	\end{array}
	\right.
\end{align}
where $M$ is a specific constant chosen as 
$M := \dfrac{3 \alpha^* \beta^*}{2 (\sigma^*)^2} \left( 1 + \dfrac{2 \kappa (\gamma - 1)^2 p(v^*, \theta^*)}{R \mu (\gamma p(v^*, \theta^*) + \frac{1}{v^*})} \right)$ in \reff{M-def},
with $\alpha^*, \beta^*, \sigma^*$ defined in Lemma \ref{Crucial}.

The following lemma proves that the ODE \reff{Xt-def} has a unique absolutely continuous solution $\mathbf{X}(t)$ defined on any interval in time $[0,T]$.
The proof of this lemma relies on the well-known Cauchy-Lipschitz theorem, and we refer to \cite{Kang25} for a detailed argument.
\\ 

\begin{lemma} \label{Cauchy-Lip}
    For any $c_1,c_2>0$, there exists a constant $C>0$ such that the following holds.
    For any $T>0$ and any functions $v,\theta \in L^\infty((0,T)\times\mathbb{R})$ verifying
    \begin{align*}
        c_1\le v(t, x), \theta(t, x)\le c_2, \qquad \forall (t,x)\in[0,T]\times\mathbb{R},
    \end{align*}
    the ODE \eqref{Xt-def} has a unique absolutely continuous solution $X(t)$ on $[0,T]$.
    Moreover,
    \begin{align*}
        |\mathbf{X}(t)|\le Ct, \qquad \forall t\le T.
    \end{align*} 
\end{lemma}

\begin{proof}
	By the well-known Cauchy-Lipschitz theorem, it only suffices to prove the following two estimates.
	Let $F(t, \mathbf{X}(t))$ denote the right-hand side of the ODE \reff{Xt-def}. 
	Thanks to the facts that $\| a \|_{C^1(\mathbb{R})}\le 2$ and by Lemma \ref{VS-pro}, we have 
	\begin{align*}
		\underset{\mathbf{X}\in\mathbb{R}}{\sup}|F(t, \mathbf{X}(t))| \le \dfrac{C}{\delta_S} \| (v-\bar{v}, u-\bar{u}, \theta-\bar{\theta}) \|_{L^\infty(\mathbb{R})} \displaystyle{\int_\mathbb{R}} |(v^S_\xi, u^S_\xi, \theta^S_\xi)| d\xi \le C,
	\end{align*}
	and by \reff{par-a}, 
	\begin{align*}
		\underset{\mathbf{X}\in\mathbb{R}}{\sup}|\partial_X F(t,\mathbf{X}(t))| \le \dfrac{C}{\delta_S} \left\{ \displaystyle \int_{\mathbb{R}} \left| (v_\xi^S, u_\xi^S, \theta_\xi^S) \right|^2 d\xi + \| (v - \bar{v}, u - \bar{u}, \theta - \bar{\theta}) \|_{L^\infty(\mathbb{R})} \right. \\[2mm]
		\left. \times\displaystyle \int_{\mathbb{R}} \left( \left| (v_\xi^S, u_\xi^S, \theta_\xi^S) \right| + \left| (v_{\xi\xi}^S, u_{\xi\xi}^S, \theta_{\xi\xi}^S) \right| \right) d\xi \right\} \le C. 
	\end{align*}
	Then it yields $|\mathbf{X}(t)|\le Ct$ from the above first estimate.
\end{proof}
\hspace*{\fill} 

\subsection{A prior estimates}
\hspace*{\fill}\\ 

Below, we state the a prior estimates in Proposition \ref{A prior estimate}, which is a key point of the subsequent estimates.
\\ 

\begin{proposition} \label{A prior estimate}
	Suppose that $(v, u, \theta, \phi)$ is the solution to \reff{eq:CNSP-change} on $[0,T]$ for some $T>0$, 
	and $(\bar{v}, \bar{u}, \bar{\theta}, \bar{\phi})$ is the composite wave defined in \reff{CompositeWave} with the shift $\mathbf{X}(t)$ being the absolutely continuous solution to \reff{Xt-def} with the weight function $a$ defined in \reff{a-def}.   
	Then there exists positive constants $C_0,\delta_0,\varepsilon_1$ such that if all the wave strength satisfies $\delta_R, \delta_C, \delta_S<\delta_0$ and 
	\begin{align*}
		\begin{array}{l}
			v - \bar{v} \in C([0, T]; H^2(\mathbb{R})), \\[2mm]
			(u - \bar{u}, \theta - \bar{\theta}) \in C([0, T]; H^2(\mathbb{R})) \cap L^2(0, T; H^3(\mathbb{R})), \\[2mm]
			\phi - \bar{\phi} \in C([0, T]; H^3(\mathbb{R})),
		\end{array}
	\end{align*}
	and
	\begin{align} \label{ProAss}
		\begin{array}{l}
			\| (v - \bar{v}, u - \bar{u}, \theta - \bar{\theta}) \|_{L^\infty(0, T; H^2(\mathbb{R}))} + \|\phi - \bar{\phi}\|_{L^\infty(0, T; H^3(\mathbb{R}))} \le \varepsilon_1,
		\end{array}
	\end{align}
	then it holds that for all $t\le T$,
	\begin{align} \label{ProFin}
		\begin{array}{l}
			\displaystyle \sup_{t \in [0, T]} \left( \| (v - \bar{v}, \, u - \bar{u}, \, \theta - \bar{\theta}) \|_{H^2(\mathbb{R})}^2 + \| \phi-\bar{\phi} \|_{H^3(\mathbb{R})}^2 \right) \\[4mm]
			\displaystyle + \int_0^t \left( \delta_S |\mathbf{\dot{X}}(\tau)|^2  + \mathcal{G}^R(U) + \mathcal{G}^S(U) + \mathcal{G}^C(U) + \| (v - \bar{v})_\xi \|_{H^1(\mathbb{R})}^2 \right. \\[4mm]
			\displaystyle \left. \qquad\qquad + \| (u - \bar{u}, \theta - \bar{\theta})_\xi \|_{H^2(\mathbb{R})}^2 \right) d\tau \\[4mm]
			\displaystyle \le C_0 \| (v_0 - \bar{v}(0, \cdot), \, u_0 - \bar{u}(0, \cdot), \, \theta_0 - \bar{\theta}(0, \cdot)) \|_{H^2(\mathbb{R})}^2 + C_0 \delta_0^{1/3},
		\end{array}
	\end{align}
	where 
	\begin{align} \label{ProG-RSC}
		\begin{array}{l}
			\displaystyle \mathcal{G}^R = \int_{\mathbb{R}} |v_\xi^R| \, |(v - \bar{v}, \theta - \bar{\theta})|^2 \, d\xi, \\[3mm]
			\displaystyle \mathcal{G}^S = \int_{\mathbb{R}} |v_\xi^S| \, |(v - \bar{v}, u - \bar{u}, \theta - \bar{\theta})|^2 \, d\xi, \\[3mm]
			\displaystyle \mathcal{G}^C = (1 + \tau)^{-1} \, \int_{\mathbb{R}} e^{ - \frac{C_1 |\xi + \sigma \tau|^2}{ 1 + \tau } } \,|(v - \bar{v}, u - \bar{u}, \theta - \bar{\theta})|^2 \, d\xi.
		\end{array}
	\end{align}
	Moreover, by \reff{Xt-def},
	\begin{align} \label{Xt-up}
		\begin{array}{l}
			|\dot{\mathbf{X}}(t)| \le C_0 \| (v - \bar{v}, u - \bar{u}, \theta - \bar{\theta})(t, \cdot) \|_{L^\infty(\mathbb{R})}, \qquad t \le T.
		\end{array}
	\end{align}
\end{proposition}
\hspace*{\fill}

The proof of this central proposition is postponed to the next sections. 
To this end, we first prove how Proposition \ref{Local Existence} and Proposition \ref{A prior estimate} together lead to Theorem \ref{Theorem-Asym}.
\\ 

\subsection{Global in time Existence}
\hspace*{\fill}\\ 

We can use the local-in-time existence and the standard continuity argument to derive the global in time existence of the solution perturbations.
We shall abuse the notation of $x$ and $\xi$ in what follows, but the meaning will be unambiguous from the context.
Let us choose the same positive constants $\delta_0, \varepsilon_1, C_0$ defined in Proposition \ref{A prior estimate}. 
Note that $\delta_0$ in Theorem \ref{Theorem-Asym} can be chosen even smaller than this proposition.
First of all, we take the smooth monotone functions $(\underline{v}, \underline{u}, \underline{\theta}, \underline{\phi})(x)$ defined in \reff{ini_und_vu} satisfying that for some $C_*>0$,
\begin{align} \label{vu1}
	\begin{array}{l}
		\displaystyle \sum_{\pm} \| (\underline{v} - v_\pm, \, \underline{u} - u_\pm, \, \underline{\theta} - \theta_\pm, \, \underline{\phi} - \phi_\pm) \|_{L^2(\mathbb{R_{\pm}})} \\[4mm]
		\displaystyle \qquad + \| ( \partial_x\underline{v}, \, \partial_x\underline{u}, \, \partial_x\underline{\theta} ) \|_{H^1(\mathbb{R})} + \| \partial_x\underline{\phi} \|_{H^2(\mathbb{R})} \le C_* \delta_0.
	\end{array}
\end{align}
This together with Lemma \ref{ARV-pro} and \ref{VS-pro} implies that for some $C_1>0$,
\begin{align} \label{vu2}
	\begin{array}{l}
		\displaystyle \| (\underline{v} - \bar{v}(0, \cdot), \, \underline{u} - \bar{u}(0, \cdot), \, \underline{\theta} - \bar{\theta}(0, \cdot)) \|_{H^2(\mathbb{R})} + \| \underline{\phi} - \bar{\phi}(0, \cdot) \|_{H^3(\mathbb{R})} \le C_1 \sqrt{\delta_0}.
	\end{array}
\end{align}
Let us choose the positive constants $\varepsilon_0, \varepsilon_*$ by the smallness of $\delta_0$,
\begin{align*}
	\displaystyle \varepsilon_* := \dfrac{\dfrac{\varepsilon_1}{2} - C_0^{1/2} \delta_0^{1/6}}{C_0 + 1} - C_1 \sqrt{\delta_0}, 
	\qquad 
	\varepsilon_0 := \varepsilon_* - C_* \delta_0.
\end{align*}

Consider the initial data $(v_0, u_0, \theta_0, \phi_0)$ verifying the hypothesis of Theorem \ref{Theorem-Asym}, i.e.,
\begin{align*}
	\begin{array}{l}
		\displaystyle \sum_{\pm} \| (v_0 - v_\pm, \, u_0 - u_\pm, \, \theta_0 - \theta_\pm, \, \phi_0-(-\ln v_\pm)) \|_{L^2(\mathbb{R}_\pm)} \\[4mm]
		\displaystyle  \qquad 
		+ \| (v_{0x},  u_{0x},  \theta_{0x}) \|_{H^1(\mathbb{R})} + \| \phi_{0x} \|_{H^2(\mathbb{R})} \le \varepsilon_0,
	\end{array}
\end{align*}
which together with \reff{vu1} implies 
\begin{align} \label{vu3}
	\begin{array}{l}
		\displaystyle \| (v_0 - \underline{v}, \, u_0 - \underline{u}, \, \theta_0 - \underline{\theta}) \|_{H^2(\mathbb{R})} 
		+ \| \phi_0 - \underline{\phi} \|_{H^3(\mathbb{R})} \le \varepsilon_0 + C_* \delta_0 = \varepsilon_*.
	\end{array}
\end{align}
Then it follows from the classical one-dimensional Sobolev inequality that 
\[\| (v_0 - \underline{v}, \, \theta_0 - \underline{\theta}) \|_{L^\infty(\mathbb{R})} \le C \varepsilon_*,\]
which together with the smallness of $\varepsilon_*$ means that
\begin{align} \label{v-bound}
	0 < C^{-1} < v_0(\xi), \ \theta_0(\xi) < C, \quad \forall \xi \in \mathbb{R}.
\end{align}
Since $0 < \varepsilon_* < \dfrac{\varepsilon_1}{2}$, Proposition \ref{Local Existence} with \reff{vu3} and \reff{v-bound} yields that 
there exists $T_0 > 0$ such that the CNSP system \reff{eq:CNSP-change} has a unique local solution $(v, u, \theta, \phi)$ on $[0, T_0]$ satisfying
\begin{align} \label{vu4-sol}
	\| (v - \underline{v}, \, u - \underline{u}, \, \theta - \underline{\theta}) \|_{L^\infty(0, T_0; H^2(\mathbb{R}))} + \| \phi - \underline{\phi} \|_{L^\infty(0, T_0; H^3(\mathbb{R}))} \le \dfrac{\varepsilon_1}{2},
\end{align}
and 
\begin{align*}
	0 < \dfrac{C^{-1}}{2} < v(t, \xi), \ \theta(t, \xi) < 2C, \quad \forall (t, \xi) \in [0, T_0] \times \mathbb{R}.
\end{align*}
Then using the same argument as in \reff{vu2}, and applying Lemma \ref{ARV-pro} and \ref{Cauchy-Lip},
we have 
\begin{align*}
	\displaystyle \| (\underline{v} - \bar{v}(t, \cdot), \, \underline{u} - \bar{u}(t, \cdot), \, \underline{\theta} - \bar{\theta}(t, \cdot)) \|_{H^2(\mathbb{R})} + \| \underline{\phi} - \bar{\phi}(t, \cdot) \|_{H^3(\mathbb{R})} \le C \sqrt{\delta_0} \, (1 + \sqrt{t}).
\end{align*}
Using the smallness of $\delta_0$, and taking $T_1 \in (0,T_0)$ small enough such that $C \sqrt{\delta_0} \, (1 + \sqrt{T_1}) \le \dfrac{\varepsilon_1}{2},$
we obtain 
\begin{align*}
	\displaystyle \| (v - \bar{v}, \, u - \bar{u}, \, \theta - \bar{\theta}) \|_{L^\infty(0, T_1; H^2(\mathbb{R}))} + \| \phi - \bar{\phi} \|_{L^\infty(0, T_1; H^3(\mathbb{R}))} \le \varepsilon_1.
\end{align*}
By Proposition \ref{Local Existence} and the absolutely continuous $\mathbf{\dot{X}}(t)$ in Lemma \ref{Cauchy-Lip}, 
we have 
\begin{align*}
	(v - \bar{v}, \, u - \bar{u}, \, \theta - \bar{\theta}) \in C([0, T_1]; H^2(\mathbb{R})),
\end{align*}
and
\begin{align*}
	\phi - \bar{\phi} \in C([0, T_1]; H^3(\mathbb{R})).
\end{align*}

Define the maximal existence time
\begin{align*}
	\begin{array}{l}
		\displaystyle T_M := \sup \left\{ t > 0 \ \middle| \ \sup_{\tau \in [0, t]} \left( \| (v - \bar{v}, \, u - \bar{u}, \, \theta - \bar{\theta}) \|_{H^2(\mathbb{R})} + \| \phi - \bar{\phi} \|_{H^3(\mathbb{R})} \right) \le \varepsilon_1 \right\}.
	\end{array}
\end{align*}
Obviously, $T_M \ge T_1$. We can prove $T_M = +\infty$. Otherwise, if $T_M < +\infty$, then
\begin{align} \label{vu5}
	\begin{array}{l}
		\displaystyle \sup_{\tau \in [0, T_M]} \left( \| (v - \bar{v}, \, u - \bar{u}, \, \theta - \bar{\theta}) \|_{H^2(\mathbb{R})} + \| \phi - \bar{\phi} \|_{H^3(\mathbb{R})} \right) = \varepsilon_1.
	\end{array}
\end{align}
However, it follows from \reff{vu2} and \reff{vu3} that
\begin{align*}
	\begin{array}{l}
		\displaystyle \| (v_0 - \bar{v}(0, \cdot), \, u_0 - \bar{u}(0, \cdot), \, \theta_0 - \bar{\theta}(0, \cdot)) \|_{H^2(\mathbb{R})} + \| \phi_0 - \bar{\phi}(0, \cdot) \|_{H^3(\mathbb{R})} \\[4mm]
		\displaystyle \qquad \le \varepsilon_* + C_1 \sqrt{\delta_0} = \dfrac{\dfrac{\varepsilon_1}{2} - C_0^{1/2} \delta_0^{1/6}}{C_0 + 1},
	\end{array}
\end{align*}
which together with Proposition \ref{A prior estimate} implies 
\begin{align*}
	\begin{array}{l}
		\displaystyle \sup_{\tau \in [0, T_M]} \left( \| (v - \bar{v}, \, u - \bar{u}, \, \theta - \bar{\theta}) \|_{H^2(\mathbb{R})} + \| \phi - \bar{\phi} \|_{H^3(\mathbb{R})} \right) \\[4mm]
		\displaystyle \qquad \le C_0 \dfrac{\dfrac{\varepsilon_1}{2} - C_0^{1/2} \delta_0^{1/6}}{C_0 + 1} + C_0^{1/2} \delta_0^{1/6} < \dfrac{\varepsilon_1}{2},
	\end{array}
\end{align*}
which contradicts with \reff{vu5}. 

Therefore, the estimates \reff{ProFin} and \reff{Xt-up} hold for $T_M = +\infty$, and 
\[(v - \bar{v}, \, u - \bar{u}, \, \theta - \bar{\theta}) \in C(0, +\infty; H^2(\mathbb{R})),\]
and 
\[\phi - \bar{\phi} \in C(0, +\infty; H^3(\mathbb{R})).\]
Moreover, since the rarefaction wave $(v^r, u^r, \theta^r, \phi^r)(x/t)$ is Lipschitz continuous in $x$ for all $t>0$ and from Lemma \ref{ARV-pro}, we have 
\begin{align*}
	\begin{array}{l}
		\displaystyle v(t, x) - \left( v^r(x/t) + v^C(x/\sqrt{1 + t}) + v^S(x - \sigma t - \mathbf{X}(t)) - v_* - v^* \right) \in C(0, +\infty; H^2(\mathbb{R})), \\[4mm]
		\displaystyle u(t, x) - \left( u^r(x/t) + u^C(x/\sqrt{1 + t}) + u^S(x - \sigma t - \mathbf{X}(t)) - u_* - u^* \right) \in C(0, +\infty; H^2(\mathbb{R})), \\[4mm]
		\displaystyle \theta(t, x) - \left( \theta^r(x/t) + \theta^C(x/\sqrt{1 + t}) + \theta^S(x - \sigma t - \mathbf{X}(t)) - \theta_* - \theta^* \right) \in C(0, +\infty; H^2(\mathbb{R})), \\[4mm]
		\displaystyle \phi(t, x) - \left( \phi^r(x/t) + \phi^C(x/\sqrt{1 + t}) + \phi^S(x - \sigma t - \mathbf{X}(t)) + \ln v_* + \ln v^* \right) \in C(0, +\infty; H^3(\mathbb{R})).
	\end{array}
\end{align*}
It follows from the estimates \reff{ProFin} with $T_M = +\infty$ that 
$(u - \bar{u}, \, \theta - \bar{\theta})_{\xi\xi} \in L^2(0, +\infty; H^1(\mathbb{R})).$
Then combining with $(u^R, \theta^R)_{\xi\xi} \in L^2(0, +\infty; H^1(\mathbb{R}))$ by Lemma \ref{ARV-pro}
and $(u^C, \theta^C)_{\xi\xi} \in L^2(0, +\infty; H^1(\mathbb{R}))$ by Lemma \ref{VCD} yields that 
\begin{align*}
	\begin{array}{l}
		\big(u(t, x) - u^S( x - \sigma t - \mathbf{X}(t)) \big)_{xx} \in L^2(0, +\infty; H^1(\mathbb{R})), \\[4mm]
		\big(\theta(t, x) - \theta^S( x - \sigma t - \mathbf{X}(t)) \big)_{xx} \in L^2(0, +\infty; H^1(\mathbb{R})). 
	\end{array}
\end{align*}

Especially, since the right-hand side of \reff{ProFin} with $T_M=+\infty$ can be small enough, if follows from the Sobolev inequality that 
\[0 < \dfrac{C^{-1}}{3} < v(t, \xi), \theta(t, \xi)< 3C.\]
Using the above computation, we can obtain the long-time asymptotic behavior as follows.
\\

\subsection{Asymptotic behaviors}
\hspace*{\fill}\\ 

Next, we turn to the time-asymptotic behavior of the perturbation $(v - \bar{v}, \, u - \bar{u}, \, \theta - \bar{\theta}, \, \phi - \bar{\phi})$. We claim that
\begin{align} \label{BV}
	\displaystyle \int_0^{+\infty} \left[ |g(t)| + |g'(t)| \right] dt < +\infty, 
\end{align}
in which 
\begin{align*}
	\begin{array}{l}
		\displaystyle g(t) := \| (v - \bar{v}, \, u - \bar{u}, \, \theta - \bar{\theta}, \, \phi - \bar{\phi}) \|_{L^2(\mathbb{R})}^2.
	\end{array}
\end{align*}
Since the estimates in Proposition \ref{A prior estimate} hold for $T=+\infty$, it only suffices to check that
\begin{align*}
	\begin{array}{l}
		\displaystyle \int_0^{+\infty} |g'(t)| \, dt < +\infty.
	\end{array}
\end{align*}
Then by $\reff{PS}_1$, \reff{PS-u-form}, and $\reff{PS}_3$, and integration by parts, if follows 
\begin{align*}
	\begin{array}{l}
		\displaystyle \int_0^{+\infty} |g'(t)| \, dt = \int_0^{+\infty} \left| \dfrac{d}{dt} \Big| (v - \bar{v}, \, u - \bar{u}, \, \theta - \bar{\theta}, \, \phi - \bar{\phi}) \Big|_{L^2(\mathbb{R})}^2 \right| dt \\[3mm]
		\displaystyle  = 2 \int_0^{+\infty} \Big| \int_{\mathbb{R}} \left[ (v - \bar{v})_\xi \cdot (v - \bar{v})_{\xi t} + (u - \bar{u})_\xi \cdot (u - \bar{u})_{\xi t} + (\theta - \bar{\theta})_\xi \cdot (\theta - \bar{\theta})_{\xi t} \right.  \\[3mm]
		\displaystyle \qquad \qquad \qquad \quad  \left. + (\phi - \bar{\phi})_\xi \cdot (\phi - \bar{\phi})_{\xi t} \right] d\xi \Big| dt  \\[3mm]
		\displaystyle = 2 \int_0^{+\infty} \Big| \int_{\mathbb{R}} \left[ (v - \bar{v})_{\xi\xi} \, (v - \bar{v})_t + (u - \bar{u})_{\xi\xi} \, (u - \bar{u})_t + (\theta - \bar{\theta})_{\xi\xi} \, (\theta - \bar{\theta})_t \right. \\[3mm]
		\displaystyle \qquad \qquad \qquad \quad \left. + (\phi - \bar{\phi})_{\xi\xi} \, (\phi - \bar{\phi})_t \right] d\xi \Big| dt \\[3mm]
		\displaystyle = 2 \int_0^{+\infty} \Big| \int_{\mathbb{R}} \Big\{ (v - \bar{v})_{\xi\xi} \left[ \sigma (v - \bar{v})_\xi + \mathbf{\dot{X}}(t) v_\xi^S + (u - \bar{u})_\xi \right]  \\[3mm]
		\displaystyle \qquad \qquad  + (u - \bar{u})_{\xi\xi} \left[ \sigma (u - \bar{u})_\xi + \mathbf{\dot{X}}(t) u_\xi^S - (p(v, \theta) - p(\bar{v}, \bar{\theta}))_\xi \right. \\[3mm]
		\displaystyle \qquad \qquad \qquad \qquad \qquad \qquad \left. + \mu \left( \dfrac{u_\xi}{v} - \dfrac{\bar{u}_\xi}{\bar{v}} \right)_\xi - \left( \dfrac{\phi_\xi}{v} - \dfrac{\bar{\phi}_\xi}{\bar{v}} \right) - Q_1 \right]  \\[3mm]
		\displaystyle \qquad \qquad  + \dfrac{\gamma - 1}{R} (\theta - \bar{\theta})_{\xi\xi}  \left[ \dfrac{\sigma R}{\gamma - 1} (\theta - \bar{\theta})_\xi + \mathbf{\dot{X}}(t) \dfrac{R}{\gamma - 1} \theta_\xi^S - (p(v, \theta) u_\xi - p(\bar{v}, \bar{\theta}) \bar{u}_\xi) \right.  \\[3mm]
		\displaystyle \qquad \qquad \qquad \qquad \qquad \qquad \left. + \kappa \left( \dfrac{\theta_\xi}{v} - \dfrac{\bar{\theta}_\xi}{\bar{v}} \right)_\xi + \mu \left( \dfrac{u_\xi^2}{v} - \dfrac{\bar{u}_\xi^2}{\bar{v}} \right) - Q_2 \right] \\[3mm]
		\displaystyle \qquad \qquad  + (\phi - \bar{\phi})_{\xi\xi} \, (\phi - \bar{\phi})_t \Big\} d\xi \Big| dt \\[6mm]
		\displaystyle \le C \int_0^{+\infty} \left( \delta_S |\mathbf{\dot{X}}(t)|^2 + \mathcal{G}^S + \mathcal{G}^R + \mathcal{G}^C + \| (v - \bar{v}, \, u - \bar{u}, \, \theta - \bar{\theta})_\xi \|_{H^1(\mathbb{R})}^2 \right. \\[3mm]
		\displaystyle \qquad \left. + \| ((\phi-\bar{\phi})_t, (\phi-\bar{\phi})_\xi, (\phi-\bar{\phi})_{\xi\xi}) \|_{L^2(\mathbb{R})}^2 + \| Q_1 \|_{L^2(\mathbb{R})}^2 + \| Q_2 \|_{L^2(\mathbb{R})}^2 \right) dt.
	\end{array}
\end{align*}
This together with Lemma \ref{phi-hig}, \reff{Q12-Est-all}, and Proposition \ref{A prior estimate} together with $T=+\infty$ yields 
\begin{align*}
	\begin{array}{l}
		\displaystyle \int_0^{+\infty} |g'(t)| \, dt \\[3mm]
		\displaystyle \le C \int_0^{+\infty} \left[ \delta_S |\dot{\mathbf{X}}(t)|^2 + \mathcal{G}^S + \mathcal{G}^R + \mathcal{G}^C + \| (v - \bar{v}, \, u - \bar{u}, \, \theta - \bar{\theta})_\xi \|_{H^1(\mathbb{R})}^2 \right] dt + C \delta_0 \\[3mm]
		\displaystyle < +\infty.
	\end{array}
\end{align*}
Therefore, the classical estimates \reff{BV} yields 
\[\lim_{t \to +\infty} \| (v - \bar{v}, \, u - \bar{u}, \, \theta - \bar{\theta}, \, \phi - \bar{\phi})_\xi \|_{L^2(\mathbb{R})} = 0,\]
which follows from the Sobolev inequality that
\[\lim_{t \to +\infty} \sup_{\xi \in \mathbb{R}} | (v - \bar{v}, \, u - \bar{u}, \, \theta - \bar{\theta}, \, \phi - \bar{\phi}) | = 0.\]
Thus, combining with the fifth estimate in Lemma \ref{ARV-pro} and returning to the original coordinates via the shock-adapted transformation $(t, \xi) \mapsto (t, x = \xi + \sigma t)$ gives the time-asymptotic behavior in Theorem \ref{Theorem-Asym}, and it holds from \reff{Xt-up} 
\[|\dot{\mathbf{X}}(t)| \le C \| (v - \bar{v}, \, u - \bar{u}, \, \theta - \bar{\theta}) \|_{L^\infty(\mathbb{R})} \to 0, \quad \text{as} \quad t \to +\infty.\]
Hence, we complete the proof of Theorem \ref{Theorem-Asym}.
Later, \textbf{denote $\mathbf{(v^S, u^S, \theta^S, \phi^S)}$ be the viscous shock shifted by $\mathbf{X}(t)$}.
\\

\section{Relative Entropy Estimates}
\hspace*{\fill}\\ 

This section is concentrated on the proof of the following key lemma for $L^2$ estimate around the viscous superposition ansatz \reff{CompositeWave}, which is based on the modulated relative entropy introduced by Dafermos \cite{Dafermos96} and DiPerna \cite{DiPerna79},
and additionally treats the structure of the approximate rarefaction wave and the viscous contact wave with the heat kernel estimates \cite{Huang10} and wave interaction estimates between three waves.
\\ 

\begin{lemma} \label{Low-Est-lem}
	Under the hypotheses of Proposition \ref{A prior estimate}, there exists a generic constant $C>0$ independent of $\delta_0, \varepsilon_1, T$ such that for all $t\in[0,T]$, 
	\begin{align*}
		\begin{array}{l}
			\displaystyle \sup_{t\in[0,T]} \| \left(v-\bar{v}, \, u-\bar{u}, \, \theta-\bar{\theta}, \, (\phi-\bar{\phi})_\xi, \, (\phi-\bar{\phi})_{\xi\xi} \right) \|_{L^2(\mathbb{R})}^2 \\[3mm]
			\displaystyle \qquad \qquad  + \int_0^t \left( \delta_S |\dot{\mathbf{X}}(\tau)|^2 + \mathcal{G}^R + \mathcal{G}^S + \mathcal{G}^C + D \right) d\tau \\[4mm]
			\displaystyle \le \| \left(v-\bar{v}, \, u-\bar{u}, \, \theta-\bar{\theta}, \, (\phi-\bar{\phi})_\xi, \, (\phi-\bar{\phi})_{\xi\xi} \right)(t=0) \|_{L^2(\mathbb{R})}^2 + C \delta_0^{1/3} \\[4mm]
			\displaystyle \qquad \qquad + C (\lambda + \delta_0^{1/3} + \varepsilon_1) \int_0^t \left( \| (v - \bar{v})_\xi \|_{H^1(\mathbb{R})}^2 + \| (u - \bar{u})_{\xi\xi} \|_{L^2(\mathbb{R})}^2 \right) d\tau,
		\end{array}
	\end{align*}
	where $\mathcal{G}^R, \mathcal{G}^S, \mathcal{G}^C$ is defined in \reff{ProG-RSC}.
\end{lemma}

\hspace*{\fill}

\subsection{Wave interaction estimates}
\hspace*{\fill}\\

We will present some wave interaction estimates between the three wave profiles, including not only $Q_i^I$ but also some useful lower-order terms for later use.
Note that a prior assumption \reff{ProAss} with the help of the Sobolev inequality implies 
\begin{align} \label{Ass-Inf}
	\| (v - \bar{v}, \, u - \bar{u}, \, \theta - \bar{\theta}) \|_{L^\infty((0, T) \times \mathbb{R})} \le C \varepsilon_1, \quad \text{and so} \quad v, \, u, \, \theta \in L^\infty((0, T) \times \mathbb{R}).
\end{align}
which togther with the ODE \reff{Xt-def} means that 
\begin{align} \label{WI-Xt-up}
	\begin{array}{l}
		\displaystyle |\dot{\mathbf{X}}(t)| \le \dfrac{C}{\delta_S} \| (v - \bar{v}, \, u - \bar{u}, \, \theta - \bar{\theta}) \|_{L^\infty(\mathbb{R})} \int_{\mathbb{R}} |(v^S_\xi, u^S_\xi, \theta^S_\xi)| \, d\xi \\[3mm]
		\qquad \quad \le C \| (v - \bar{v}, \, u - \bar{u}, \, \theta - \bar{\theta}) \|_{L^\infty(\mathbb{R})} \\[3mm]
		\qquad \quad \le C\varepsilon_1.
	\end{array}
\end{align}
\hspace*{\fill}

\begin{lemma} \label{WI-est}
	Let $\mathbf{X}(t)$ be the shifted defined by \reff{Xt-def} with \reff{Ass-Inf}-\reff{WI-Xt-up} and $Q_i^I (i=1,2,3)$ is defined in \reff{Q123_I}.
	Then there exists a generic constant $C > 0$, independent of $T, \delta_R, \delta_S, \delta_C$ such that for all $t\le T$,
	\begin{align} \label{WI-1}
		\begin{array}{l}
			\| |(v_\xi^R, \, \theta_\xi^R) | \cdot | (v^C - v_*, \, \theta^C - \theta_*, \, v^S - v^*, \, \theta^S - \theta^*) | \|_{L^2(\mathbb{R})} \\[3mm]
			\qquad + \| |(v_\xi^C, \, \theta_\xi^C) | \cdot | (v^R - v_*, \, \theta^R - \theta_*, \, v^S - v^*, \, \theta^S - \theta^*) | \|_{L^2(\mathbb{R})} \\[3mm]
			\le C \delta_S (\delta_R + \delta_C) e^{-C \delta_S t} + C \delta_R \delta_C e^{-Ct}, \\[3mm]
			\| (v_\xi^S, \, \theta_\xi^S) | \cdot | (v^R - v_*, \, \theta^R - \theta_*, \, v^C - v^*, \, \theta^C - \theta^*) | \|_{L^2(\mathbb{R})} \\[3mm]
			\le C \delta_S^{3/2} (\delta_R + \delta_C) e^{-C \delta_S t},
		\end{array}
	\end{align}
	and 
	\begin{align} \label{WI-2}
		\begin{array}{l}
			\| (v^R - v_*) (v^S - v^*) \|_{L^2(\mathbb{R})} \le C \delta_S^{1/2} \delta_R e^{-C \delta_S t} + C \delta_S \delta_R e^{-Ct}, \\[3mm]
			\| (v^C - v^*) (v^S - v^*) \|_{L^2(\mathbb{R})} \le C \delta_S^{1/2} \delta_C e^{-C \delta_S t} + C \delta_S \delta_C e^{-Ct}, \\[3mm]
			\| (v^R - v_*) (v^C - v_*) \|_{L^2(\mathbb{R})} \le C \delta_R \delta_C e^{-Ct}.
		\end{array}
	\end{align}
	In particular,
	\begin{align} \label{WI-3}
		\begin{array}{l}
			\| (Q_1^I, \, Q_2^I) \|_{L^2(\mathbb{R})} \le C \delta_S (\delta_R + \delta_C) e^{-C \delta_S t} + C \delta_R \delta_C e^{-Ct}, \\[3mm]
			\| Q_3^I \|_{L^2(\mathbb{R})} \le C\delta_S^{1/2} (\delta_R + \delta_C) e^{-C \delta_S t} + C ( \delta_S \delta_C + \delta_S \delta_R + \delta_R \delta_C ) e^{-C t}, \\[3mm]
			\| ((Q_3^I)_\xi, \bar{v}_\xi Q_3^I) \|_{L^2(\mathbb{R})} \le C\delta_S (\delta_R + \delta_C) e^{-C \delta_S t} + C \delta_R \delta_C  e^{-C t}. \\[3mm]
		\end{array}
	\end{align}
\end{lemma}

\begin{remark}
	Although the estimates for the first two in \reff{WI-2} and the middle one in\reff{WI-3} appear to be singular in $\delta_S$ when integrated over time in $ [0, +\infty)$,
	no such issues arises, since the subsequent analysis only requires the $L^2$ integral in time, namely,
	\[\int_{0}^{+\infty} (\| (v^C-v^*, v^R-v_*)(v^S-v^*) \|^2_{L^2(\mathbb{R})} + \|Q_3^I\|^2_{L^2(\mathbb{R})} ) dt,\]
	which is not singular in $\delta_S$.
\end{remark}
\hspace*{\fill}

\begin{proof}[Proof]
	For brevity, we only prove the estimates in \reff{WI-2} and the last two in \reff{WI-3}, since the remaining cases are covered in Lemma 4.2 of \cite{Kang25}.
	From \reff{WI-Xt-up}, it follows that 
	\begin{align*}
		|\mathbf{X}(t)| \le C\varepsilon_1 t \le \frac{\sigma t}{4}, \quad \forall t \in [0,T],
	\end{align*}
	where $\varepsilon_1>0$ is a small constant and $\sigma>0$ is defined in \reff{sig-def}.
	\\ 

	\noindent$\bullet$ 
	Since 
	\begin{align*}
		\forall \xi \le -\dfrac{\sigma t}{2} \implies \xi - \mathbf{X}(t) \le -\dfrac{\sigma t}{4} < 0 \implies |\xi - \mathbf{X}(t)| \ge \dfrac{\sigma t}{4},
	\end{align*}
	by Lemma \ref{VS-pro}, it holds that $\forall \xi \le -\dfrac{\sigma t}{2}$,
	\begin{align} \label{v1-1}
		|v^S - v^*| \le C \delta_S \, e^{-C \delta_S |\xi - \mathbf{X}(t)|} \le C \delta_S \, e^{-\frac{C \delta_S |\xi - \mathbf{X}(t)|}{2}} e^{-\frac{C \delta_S \sigma t}{8}}.
	\end{align}
	On the other hand, due to $\sigma_{1*} := \sigma_1(v_*, \theta_*) < 0$ by \reff{Eigen},
	\begin{align*}
		\forall \xi > -\dfrac{\sigma t}{2} \implies \xi + \sigma t > \dfrac{\sigma t}{2} > 0 > \dfrac{\sigma_{1*}}{2} (1+t) > \sigma_{1*} (1+t),
	\end{align*}
	it follows from Lemma \ref{ARV-pro} that $\forall \xi > -\dfrac{\sigma t}{2}$,
	\begin{align} \label{v1-2}
		|v^R - v_*| \le C \delta_R \, e^{-2|\xi + \sigma t - \sigma_{1*} (1+t)|} = C \delta_R \, e^{-2 |\xi + \sigma t - \frac{\sigma_{1*}}{2}(1+t)|} e^{-|\sigma_{1*}|(1+t)}.
	\end{align}
	Using \reff{v1-1} and \reff{v1-2} together with
	\begin{align*}
		\| v^S - v^* \|_{L^\infty(\mathbb{R})} \le C \delta_S, \quad \| v^R - v_* \|_{L^\infty(\mathbb{R})} \le C \delta_R,
	\end{align*}
	from Lemma \ref{VS-pro} and \ref{ARV-pro}, we obtain 
	\begin{align*}
		\begin{array}{l}
			\displaystyle \| (v^R - v_*) (v^S - v^*) \|_{L^2(\mathbb{R})}^2 = \left( \int_{\left\{ \xi \le -\frac{\sigma t}{2} \right\}} + \int_{\left\{ \xi > -\frac{\sigma t}{2} \right\}} \right) (v^R - v_*)^2 (v^S - v^*)^2 \, d\xi \\[4mm]
			\displaystyle \le \int_{\left\{ \xi \le -\frac{\sigma t}{2} \right\}} (v^S - v^*)^2 \, d\xi \cdot \| v^R - v_* \|_{L^\infty(\mathbb{R})}^2 + \int_{\left\{ \xi > -\frac{\sigma t}{2} \right\}} (v^R - v_*)^2 \, d\xi \cdot \| v^S - v^* \|_{L^\infty(\mathbb{R})}^2 \\[5mm]
			\displaystyle \le C \delta_S \delta_R^2 \, e^{-C \delta_S t} + \delta_S^2 \delta_R^2 \, e^{-Ct}.
		\end{array}
	\end{align*}
	Thus,
	\begin{align} \label{WI-first}
		\displaystyle \| (v^R - v_*) (v^S - v^*) \|_{L^2(\mathbb{R})} \le C \delta_S^{1/2} \delta_R \, e^{-C \delta_S t} + C \delta_S \delta_R \, e^{-Ct}. 
	\end{align}
	\hspace*{\fill}

	\noindent$\bullet$
	Since
	\[\xi > -\dfrac{\sigma t}{2} \implies \xi + \sigma t > \dfrac{\sigma t}{2} > 0,\]
	by Lemma \ref{VCD},
	\begin{align} \label{v2-1}
		|v^C - v^*| \le C \delta_C \, e^{-\frac{C_1 |\xi + \sigma t|^2}{1+t}} \le C \delta_C \, e^{-\frac{C_1 |\xi + \sigma t|^2}{2(1+t)}} e^{-Ct}. 
	\end{align}
	Using \reff{v1-1} and \reff{v2-1} together with 
	\begin{align*}
		\|v^C-v^*\|_{L^\infty(\mathbb{R})} \le C\delta_C, \qquad \|v^S-v^*\|_{L^\infty(\mathbb{R})} \le C\delta_S,
	\end{align*}
	we get
	\begin{align*}
		\begin{array}{l}
			\displaystyle \| (v^C - v^*) (v^S - v^*) \|_{L^2(\mathbb{R})}^2 = \left( \int_{\left\{ \xi \le -\frac{\sigma t}{2} \right\}} + \int_{\left\{ \xi > -\frac{\sigma t}{2} \right\}} \right) (v^C - v^*)^2 (v^S - v^*)^2 \, d\xi \\[4mm]
			\displaystyle \le \int_{\left\{ \xi \le -\frac{\sigma t}{2} \right\}} (v^S - v^*)^2 \, d\xi \cdot \| v^C - v^* \|_{L^\infty(\mathbb{R})}^2 + \int_{\left\{ \xi > -\frac{\sigma t}{2} \right\}} (v^C - v^*)^2 \, d\xi \cdot \| v^S - v^* \|_{L^\infty(\mathbb{R})}^2 \\[5mm]
			\displaystyle \le C \delta_S \delta_C^2 \, e^{-C \delta_S t} + \delta_C^2 \delta_S^2 \, e^{-Ct}.
		\end{array}
	\end{align*}
	Thus,
	\begin{align} \label{WI-second}
		\| (v^C - v^*) (v^S - v^*) \|_{L^2(\mathbb{R})} \le C \delta_S^{1/2} \delta_C \, e^{-C \delta_S t} + C \delta_S \delta_C \, e^{-Ct}. 
	\end{align}
	\hspace*{\fill}

	\noindent$\bullet$
	Since
	\begin{align*}
		\xi + \sigma t > \dfrac{\sigma_{1*}}{2}(1+t) \implies \xi + \sigma t - \sigma_{1*}(1+t) > -\dfrac{\sigma_{1*}}{2}(1+t) = \dfrac{|\sigma_{1*}|}{2}(1+t) > 0,
	\end{align*}
	it holds from Lemma \ref{ARV-pro} that
	\begin{align} \label{v3-1}
		|v^R - v_*| \le C \delta_R \, e^{-2|\xi + \sigma t - \sigma_{1*}(1+t)|} \le C \delta_R \, e^{-2\left|\xi + \sigma t - \frac{\sigma_{1*}}{2}(1+t)\right|} e^{-|\sigma_{1*}|(1+t)}. 
	\end{align}
	On the other hand, when
	\begin{align*}
		\xi + \sigma t \le \dfrac{\sigma_{1*}}{2} (1+t) < 0 \implies |\xi + \sigma t| \ge \dfrac{|\sigma_{1*}|}{2}(1+t),
	\end{align*}
	immediately we get
	\begin{align} \label{v3-2}
		|v^C - v_*| \le C \delta_C \, e^{-\frac{C_1 |\xi + \sigma t|^2}{1+t}} \le C \delta_C \, e^{-\frac{C_1 |\xi + \sigma t|^2}{2(1+t)}} e^{-\frac{|\sigma_{1*}|^2}{8}(1+t)}. 
	\end{align}
	Using \reff{v3-1} an \reff{v3-2} together with 
	\begin{align*}
		\| v^R - v_* \|_{L^\infty(\mathbb{R})} \le C \delta_R, \quad \| v^C - v_* \|_{L^\infty(\mathbb{R})} \le C \delta_C,
	\end{align*}
	it follows that 
	\begin{align*}
		\begin{array}{l}
			\displaystyle \| (v^R - v_*) (v^C - v_*) \|_{L^2(\mathbb{R})}^2 \\[3mm]
			\displaystyle = \left( \int_{\left\{ \xi + \sigma t \le \frac{\sigma_{1*}}{2}(1+t) \right\}} + \int_{\left\{ \xi + \sigma t > \frac{\sigma_{1*}}{2}(1+t) \right\}} \right) (v^R - v_*)^2 (v^C - v_*)^2 \, d\xi \\[4mm]
			\displaystyle \le \int_{\left\{ \xi + \sigma t \le \frac{\sigma_{1*}(1+t)}{2} \right\}} (v^C - v_*)^2 \, d\xi \cdot \| v^R - v_* \|_{L^\infty(\mathbb{R})}^2 \\[4mm]
			\displaystyle \qquad\qquad\qquad + \int_{\left\{ \xi + \sigma t > \frac{\sigma_{1*}(1+t)}{2} \right\}} (v^R - v_*)^2 \, d\xi \cdot \| v^C - v_* \|_{L^\infty(\mathbb{R})}^2 \\[5mm]
			\le C \delta_R^2 \delta_C^2 \, e^{-Ct},
		\end{array}
	\end{align*}
	which implies 
	\begin{align}  \label{WI-third}
		\| (v^R - v_*) (v^C - v_*) \|_{L^2(\mathbb{R})} \le C \delta_R \delta_C \, e^{-Ct}. 
	\end{align}
	\hspace*{\fill}

	\noindent$\bullet$
	Finally, from \reff{Q123_I},
	\begin{align*}
		\displaystyle Q_3^I = (\bar{v} e^{\bar{\phi}} - v^S e^{\phi^S}) - \nu^2 \left[ \dfrac{\bar{\phi}_\xi}{\bar{v}} - \dfrac{\phi_\xi^S}{v^S} - \dfrac{\phi_\xi^C}{v^C} - \dfrac{\phi_\xi^R}{v^R} \right]_\xi,
	\end{align*}
	we only need to focus on the lower-order term $\bar{v} e^{\bar{\phi}} - v^S e^{\phi^S}$, since the higher-order term can be easily controlled by Lemma \ref{ARV-pro}-\ref{VS-pro} and \reff{WI-1}.
	From the composite wave $\bar{\phi}$ defined in \reff{CompositeWave}, and 
	$\phi^R = -\ln v^R, \, \phi^C = -\ln v^C$  by $\reff{eq:ARV-change}_4$ and $\reff{eq:VCD-change}_4$,
	we have
	\begin{align*}
		\begin{array}{l}
			\displaystyle \bar{v} e^{\bar{\phi}} - v^S e^{\phi^S} = e^{\bar{\phi}} (\bar{v} - v^S e^{\phi^S - \bar{\phi}}) = e^{\bar{\phi}} \left[ (v^R - v_*) + (v^C - v^*) + v^S \left( 1 - \dfrac{v^R v^C}{v_* v^*} \right) \right] \\[3mm]
			\displaystyle  = e^{\bar{\phi}} \left\{ (v^R - v_*) + (v^C - v^*) + v^S \left[ \left( 1 - \dfrac{v^C}{v^*} \right) + \dfrac{v^C}{v^*} \left( 1 - \dfrac{v^R}{v_*} \right) \right] \right\} \\[3mm]
			\displaystyle  = e^{\bar{\phi}} \left[ \dfrac{1}{v^*} (v^C - v^*) (v^S - v^*) + \dfrac{1}{v^*} (v^R - v_*) (v^S - v^*) + \dfrac{v^S}{v_* v^*} (v^R - v_*) (v^C - v_*) \right]. 
		\end{array}
	\end{align*}
	By \reff{WI-first}, \reff{WI-second} and \reff{WI-third}, we have
	\begin{align*}
		\| Q_3^I \|_{L^2(\mathbb{R})} \le C \delta_S^{1/2} (\delta_R + \delta_C) e^{-C \delta_S t} +C ( \delta_S \delta_C + \delta_S \delta_R + \delta_R \delta_C ) e^{-Ct}.
	\end{align*}
	Similarly, from \reff{WI-1},
	\begin{align*}
		\| \left((Q_3)_\xi, \bar{v}_\xi Q_3\right) \|_{L^2(\mathbb{R})} \le C \delta_S (\delta_R + \delta_C) e^{-C \delta_S t} + C \delta_R \delta_C e^{-Ct}.
	\end{align*}
\end{proof}
\hspace*{\fill}

\subsection{Construction of weight function}
\hspace*{\fill}\\ 

Deine the weight function 
\begin{align} \label{a-def}
	a^{-\mathbf{X}}(\xi)=a(\xi-\mathbf{X}(t)):=1+\dfrac{\lambda}{\delta_S} \big(v^S(\xi-\mathbf{X}(t))-v^* \big),
\end{align}
where the parameter $\lambda$ is chosen to be sufficiently small but far bigger than $\delta_S:=|v_+-v^*|$ such that 
\begin{align} \label{lam-range}
	\delta_S << \lambda \le C\sqrt{\delta_S},
\end{align}
It follows that 
\begin{align} \label{a-range}
	1 < a^{-\mathbf{X}}(\xi) < 1 + \lambda,
\end{align}
and by Lemma \reff{VS-pro},
\begin{align} \label{par-a}
	0 < a^{-\mathbf{X}}_\xi = \dfrac{\lambda}{\delta_S} v^S_\xi \sim \dfrac{\lambda}{\delta_S} |u^S_\xi| \sim \dfrac{\lambda}{\delta_S} |\theta^S_\xi| \sim \dfrac{\lambda}{\delta_S} |\phi^S_\xi|.
\end{align}
\hspace*{\fill}

\subsection{Relative entropy method}
\hspace*{\fill}\\ 

Let 
\begin{align}  \label{U-nota}
	U=(v, \, u, \, \theta)^t, \qquad \bar{U}=(\bar{v}, \bar{u}, \bar{\theta}),
\end{align}
where the former is the solution to the reformulated CNSP system \reff{eq:CNSP-change} and the latter is the composite wave defined in \reff{CompositeWave}, consisting of the approximate 1-rarefaction wave, the viscous 2-contact wave and the viscous 3-shock wave.

For the one-dimensional compressible Navier-Stokes system, it has the relative entropy of $U$ and $\bar{U}$ (See Appendix A in \cite{Kang25}) given by 
\begin{align} \label{eta-class}
	\eta(U | \bar{U}) = R \Psi\left( \dfrac{v}{\bar{v}} \right) + \dfrac{R}{\gamma - 1} \Psi\left( \dfrac{\theta}{\bar{\theta}} \right) + \dfrac{(u - \bar{u})^2}{2\bar{\theta}},
\end{align}  
with the convex function $\Psi(z):=z-1-\ln z$. Next, we first compute the relative entropy fucntional weighted by $a^{-\mathbf{X}}(\xi) \bar{\theta}(t,\xi)$ defined by 
\begin{align*}
	\displaystyle \int_{\mathbb{R}} a^{-\mathbf{X}}(\xi) \, \bar{\theta}(t, \xi) \, \eta(U(t, \xi) | \bar{U}(t, \xi)) \, d\xi.
\end{align*}

\begin{lemma} \label{REM-lem}
	Let $a^{-\mathbf{X}}$ be the weight function defined in \reff{a-def}. Let $(U, \phi)$ and $(\bar{U}, \bar{\phi})$ be the solution to the reformulated CNSP system \reff{eq:CNSP-change} and the superposition wave defined in \reff{CompositeWave} respectively.
	Set the modulated relative entropy
	\begin{align} \label{eta-new}
		\begin{array}{l}
			\tilde{\eta}\big((U, \phi)|(\bar{U}, \bar{\phi})\big) := \bar{\theta} \eta( U | \bar{U} ) + \dfrac{(v - \bar{v})^2}{2 \bar{v}^2} - \dfrac{(v - \bar{v})^3}{3 \bar{v}^3} - \dfrac{\nu^2}{\bar{v}^2} (v - \bar{v}) \tilde{\phi}_{\xi\xi} + \dfrac{\nu^4}{2 \bar{v}^3 e^{\bar{\phi}}} \tilde{\phi}_{\xi\xi}^2 + \dfrac{\nu^2}{2 \bar{v}} \tilde{\phi}_\xi^2.
		\end{array}
	\end{align}
	Then 
	\begin{align} \label{RE-ori}
		\begin{array}{l}
			\displaystyle \dfrac{d}{dt} \int_{R} a^{-\mathbf{X}}(\xi) \, \tilde{\eta}\big( (U, \phi) \, | \, (\bar{U}, \bar{\phi}) \big) \,d\xi =\mathbf{\dot{X}}(t) \sum_{i=1}^{4} \mathbf{Y}_i(U) - \mathbf{G}(U)- \mathbf{D}(U) \\[2mm]
			\displaystyle \qquad\qquad\qquad\qquad\,\, + \mathbf{B}_1(U) + \mathbf{B}_2(U) + \mathbf{E}_1(U) + \mathbf{E}_2(U),
		\end{array}
	\end{align}
	where 
	\begin{align*} 
		\begin{array}{l}
			\displaystyle \mathbf{Y}_1(U):=\int_{\mathbb{R}} a^{-\mathbf{X}} \left[ \dfrac{\hat{p}(\bar{v}, \bar{\theta}) v_\xi^S}{\bar{v}} (v - \bar{v}) + u_\xi^S (u - \bar{u}) + \dfrac{R \theta_\xi^S}{(\gamma - 1) \bar{\theta}} (\theta - \bar{\theta}) \right] d\xi,  \\[4mm]
			\displaystyle \mathbf{Y}_2(U):=-\int_{\mathbb{R}} a^{-\mathbf{X}} \left[  R \theta_\xi^S \, \Psi\left( \dfrac{v}{\bar{v}} \right) + \dfrac{R \theta_\xi^S}{\gamma - 1} \, \Psi\left( \dfrac{\theta}{\bar{\theta}} \right) + \dfrac{v_\xi^S}{\bar{v}^4} (v - \bar{v})^3 \right] d\xi, \\[4mm]
			\displaystyle \mathbf{Y}_3(U):=-\int_{\mathbb{R}} a^{-\mathbf{X}}_\xi \left[ \bar{\theta}\eta(U|\bar{U}) + \dfrac{(v - \bar{v})^2}{2 \bar{v}^2} - \dfrac{(v - \bar{v})^3}{3 \bar{v}^3} - \dfrac{\nu^2}{\bar{v}^2} (v - \bar{v}) \tilde{\phi}_{\xi\xi} + \dfrac{\nu^4}{2 \bar{v}^3 e^{\bar{\phi}}} \tilde{\phi}_{\xi\xi}^2 + \dfrac{\nu^2}{2 \bar{v}} \tilde{\phi}_\xi^2 \right] d\xi, \\[4mm]
			\displaystyle \mathbf{Y}_4(U):=\int_{\mathbb{R}} a^{-\mathbf{X}} \left[ - 2 \nu^2 \dfrac{v_\xi^S}{\bar{v}^3} (v - \bar{v}) \tilde{\phi}_{\xi\xi} - \nu^2 \dfrac{v_\xi^S}{\bar{v}^2} \tilde{\phi}_{\xi\xi} + \dfrac{\nu^4}{2} \left( \dfrac{3 v_\xi^S}{\bar{v}^4 e^{\bar{\phi}}} + \dfrac{\phi_\xi^S}{\bar{v}^3 e^{\bar{\phi}}} \right) \tilde{\phi}_{\xi\xi}^2 + \dfrac{\nu^2 v_\xi^S}{\bar{v}^3 e^{\bar{\phi}}} \tilde{\phi}_\xi^2 \right] d\xi,
		\end{array}
	\end{align*} 
	\begin{align*} 
		\begin{array}{l}
			\displaystyle \mathbf{G}(U):=\int_{\mathbb{R}} \sigma a^{-\mathbf{X}}_\xi \left[ \bar{\theta}\eta(U|\bar{U}) + \dfrac{(v-\bar{v})^2}{2\bar{v}^2} - \dfrac{(v-\bar{v})^3}{3\bar{v}^3}\right] \, d\xi, \\[4mm]
			\displaystyle \mathbf{D}(U):=\int_{\mathbb{R}} a^{-\mathbf{X}} \left[ \dfrac{\mu}{v} |(u - \bar{u})_\xi|^2 + \dfrac{k}{v \theta} |(\theta - \bar{\theta})_\xi|^2 \right] \, d\xi,
		\end{array}
	\end{align*}
	\begin{align*} 
		\begin{array}{l}
			\displaystyle \mathbf{B}_1(U):= \int_{\mathbb{R}} a^{-\mathbf{X}} \Bigg\{ \Big[ R \left( (\theta_t^R - \sigma \theta_\xi^R) + (\theta_t^C - \sigma \theta_\xi^C) - \sigma \theta_\xi^S \right) \Psi\left( \dfrac{v}{\bar{v}} \right) - \dfrac{p(\bar{v}, \bar{\theta}) \bar{u}_\xi}{v \bar{v}} (v - \bar{v})^2\\[3mm] 
			\displaystyle \qquad\qquad\qquad\qquad -\left( (v^R_t - \sigma v^R_\xi ) + ( v^C_t - \sigma v^C_\xi ) - \sigma v^S_\xi \right)  \dfrac{(v - \bar{v})^2}{\bar{v}^3} \Big]  \\[4mm]
			\displaystyle \qquad\qquad\qquad\,\,  + \Big[ \dfrac{R}{\gamma - 1} \left( (\theta_t^R - \sigma \theta_\xi^R) + (\theta_t^C - \sigma \theta_\xi^C) - \sigma \theta_\xi^S \right) \Psi\left( \dfrac{\theta}{\bar{\theta}} \right) \\[4mm]
			\displaystyle \qquad\qquad\qquad\qquad - \dfrac{\bar{u}_\xi}{\theta} (\theta - \bar{\theta}) (p(v, \theta) - p(\bar{v}, \bar{\theta})) + \dfrac{p(\bar{v}, \bar{\theta}) \bar{u}_\xi}{\theta \bar{\theta}} (\theta - \bar{\theta})^2 \Big] \Bigg\} \, d\xi, \\[4mm]
			\displaystyle \mathbf{B}_2(U):= \int_{\mathbb{R}} a^{-\mathbf{X}}_\xi (u - \bar{u}) (\hat{p}(v, \theta) - \hat{p}(\bar{v}, \bar{\theta})) \, d\xi, 
		\end{array}
	\end{align*}
	and 
	\begin{align*}
		\begin{array}{l}
			\displaystyle \mathbf{E}_1(U):= \int_{\mathbb{R}} a^{-\mathbf{X}} \left\{ \left[ \left( v_t^R + v_t^C - \sigma \bar{v}_\xi \right) \dfrac{(v - \bar{v})^3}{\bar{v}^4} - \dfrac{(v - \bar{v})^3}{\bar{v}^3 v} (u - \bar{u})_\xi \right] \right. \\[3mm]
			\displaystyle \qquad + \left[ \kappa \dfrac{\theta - \bar{\theta}}{\theta} \Big( \dfrac{\theta_\xi}{v} - \dfrac{\bar{\theta}_\xi}{\bar{v}} \Big) + \mu (u - \bar{u}) \left( \dfrac{u_\xi}{v} - \dfrac{\bar{u}_\xi}{\bar{v}} \right) \right]_\xi + \Bigg[ -\kappa \dfrac{(\theta - \bar{\theta})_\xi}{\theta}  \Big( \dfrac{1}{v} - \dfrac{1}{\bar{v}} \Big) \bar{\theta}_\xi \\[3mm]
			\displaystyle \qquad \qquad + \kappa \dfrac{(\theta - \bar{\theta}) \theta_\xi}{\theta^2} \Big( \dfrac{\theta_\xi}{v} - \dfrac{\bar{\theta}_\xi}{\bar{v}} \Big) + \mu \dfrac{\theta - \bar{\theta}}{\theta} \Big( \dfrac{u_\xi^2}{v} - \dfrac{\bar{u}_\xi^2}{\bar{v}} \Big) - \dfrac{(\theta - \bar{\theta})^2}{\theta \bar{\theta}} \Big( \kappa \big( \dfrac{\bar{\theta}_\xi}{\bar{v}} \big)_\xi + \mu \dfrac{\bar{u}_\xi^2}{\bar{v}} \Big) \Bigg] \\[3mm]
			\displaystyle \qquad \left. - \left[ (u - \bar{u}) Q_1 + \dfrac{\theta - \bar{\theta}}{\theta} Q_2 \right] \right\} d\xi, 
		\end{array}
	\end{align*}
	\begin{align*}
		\begin{array}{l}
			\displaystyle \mathbf{E}_2(U) := \int_{\mathbb{R}} \Bigg\{ - (a^{-\mathbf{X}} (u - \bar{u}))_\xi \left[ \dfrac{\nu^2}{2} \Big( \big( \dfrac{\phi_\xi}{v} \big)^2 - \big( \dfrac{\bar{\phi}_\xi}{\bar{v}} \big)^2 \Big) + \nu^2 \Big( \dfrac{\phi_\xi v_\xi}{v^3} - \dfrac{\bar{\phi}_\xi \bar{v}_\xi}{\bar{v}^3} \Big) - \nu^2 \Big( \dfrac{1}{v^2} - \dfrac{1}{\bar{v}^2} \Big) \phi_{\xi\xi} \right] \\[3mm]
			\displaystyle \qquad + \left[\nu^2 a^{-\mathbf{X}}_\xi (u - \bar{u}) \dfrac{\tilde{\phi}_{\xi\xi}}{\bar{v}^2} + 2 \nu^2 a^{-\mathbf{X}} \left( v_t^R + v_t^C \right) (v - \bar{v}) \dfrac{\tilde{\phi}_{\xi\xi}}{\bar{v}^3} + \nu^2 \left( \dfrac{a^{-\mathbf{X}}}{\bar{v}^2} \right)_\xi (v - \bar{v}) \tilde{\phi}_{\xi t} \right]\\[4mm]
			\displaystyle \qquad + \Bigg[ a^{-\mathbf{X}} \dfrac{\nu^4}{2} \bigg( - \dfrac{3 (v_t^R + v_t^C)}{\bar{v}^4 e^{\bar{\phi}}} - \dfrac{\phi_t^R + \phi_t^C}{\bar{v}^3 e^{\bar{\phi}}} \bigg) \tilde{\phi}_{\xi\xi}^2 - \nu^2 a^{-\mathbf{X}} \dfrac{v_t^R + v_t^C}{2 \bar{v}^2} \tilde{\phi}_\xi^2  \\[3mm]
			\displaystyle \qquad \qquad - \dfrac{\nu^4}{\bar{v} e^{\bar{\phi}}} \Big( \dfrac{a^{-\mathbf{X}}}{\bar{v}^2} \Big)_\xi \tilde{\phi}_{\xi t} \cdot \tilde{\phi}_{\xi\xi} + \Big( \dfrac{\nu^4 a^{-\mathbf{X}}}{\bar{v}^3 e^{\bar{\phi}}} \Big)_\xi \dfrac{\sigma \tilde{\phi}_{\xi\xi}^2}{2} - \dfrac{\sigma \nu^4 a^{-\mathbf{X}}}{\bar{v}^2} \Big( \dfrac{1}{\bar{v} e^{\bar{\phi}}} \Big)_\xi \tilde{\phi}_{\xi\xi}^2 - \Big(\dfrac{a^{-\mathbf{X}}}{\bar{v}} \Big)_\xi \dfrac{\sigma \nu^2 \tilde{\phi}_\xi^2}{2} \Bigg] \\[3mm]
			\displaystyle \qquad + \nu^2 \dfrac{a^{-\mathbf{X}}}{\bar{v}^2} \left( \tilde{\phi}_{\xi t} - \sigma \tilde{\phi}_{\xi\xi} \right) \left[ \Big( - \nu^2 e^{\bar{\phi}} \tilde{\phi}_\xi \dfrac{\bar{v}_\xi}{\bar{v}^2} - \nu^2 e^{-\bar{\phi}} \Big( \dfrac{v - \bar{v}}{\bar{v} v} \phi_\xi \Big)_\xi + N_p \Big)_\xi -\bar{v}_\xi \tilde{\phi} \right] \\[3mm]
			\displaystyle \qquad - \nu^2 \dfrac{a^{-\mathbf{X}}}{\bar{v}^2} \left( \tilde{\phi}_{\xi t} - \sigma \tilde{\phi}_{\xi\xi} \right) (e^{-\bar{\phi}} Q_3)_\xi \Bigg\} d\xi,    
		\end{array}
	\end{align*}
	with $Q_i \, (i=1,2,3)$ defined in \reff{Q123} and $N_p$ defined in \reff{Np-def}.
\end{lemma}
\hspace*{\fill}


\begin{remark}
	For the sake of more straightforward estimates in what follows, we directly distinguish five groups of terms in the current lemma: the terms governing the definition of $\mathbf{\dot{X}}(t)$, the good terms $\mathbf{G}$ and $\mathbf{D}$, the bad terms $\mathbf{B}_1$ and $\mathbf{B}_2$, and the error terms $\mathbf{E}_1$ involving $U$ and $\mathbf{E}_2$ relying on $\phi$, respectively.
\end{remark}
\hspace*{\fill}

\begin{proof} [Proof]
	We compute the evolution of the relative entropy functional
	\begin{align} \label{pf-RE}
		\begin{array}{l}
			\displaystyle \dfrac{d}{dt} \int_{\mathbb{R}} a^{-\mathbf{X}}(\xi) \bar{\theta}(t, \xi) \, \eta(U | \bar{U}) \, d\xi \\[3mm]
			\displaystyle = - \dot{\mathbf{X}}(t) \int_{\mathbb{R}} a_{\xi}^{-\mathbf{X}} \bar{\theta} \, \eta(U | \bar{U}) \, d\xi + \int_{\mathbb{R}} a^{-\mathbf{X}} \partial_t \left[ R \bar{\theta} \, \Phi\left( \dfrac{v}{\bar{v}} \right) + \dfrac{R \bar{\theta}}{\gamma - 1} \Phi\left( \dfrac{\theta}{\bar{\theta}} \right) + \dfrac{(u - \bar{u})^2}{2} \right] d\xi. 
		\end{array}
	\end{align}
	We next focus on calculating the second term on the right-hand side of \reff{pf-RE}.
	Since $\Psi'(z) = 1 - \dfrac{1}{z}$, $\reff{PS}_1$ and $\reff{RCDS}_1$, it follows that 
	\begin{align*}
		\begin{array}{l}
			\displaystyle \partial_t \left[ \bar{\theta} \, \Psi\left( \dfrac{v}{\bar{v}} \right) \right] = \bar{\theta}_t \, \Psi\left( \dfrac{v}{\bar{v}} \right) + \bar{\theta} \, \Psi'\left( \dfrac{v}{\bar{v}} \right) \left( \dfrac{v}{\bar{v}} \right)_t \\[3mm]
			\displaystyle = \bar{\theta}_t \, \Psi\left( \dfrac{v}{\bar{v}} \right) + \bar{\theta} \left( 1 - \dfrac{\bar{v}}{v} \right) \left( \dfrac{v_t}{\bar{v}} - \dfrac{v \bar{v}_t}{\bar{v}^2} \right) \\[3mm]
			\displaystyle = \bar{\theta}_t \, \Psi\left( \dfrac{v}{\bar{v}} \right) + \bar{\theta} \left( \dfrac{1}{\bar{v}} - \dfrac{1}{v} \right) \left[ (v - \bar{v})_t + \bar{v}_t \left( 1 - \dfrac{v}{\bar{v}} \right) \right] \\[3mm]
			\displaystyle = \bar{\theta}_t \, \Psi\left( \dfrac{v}{\bar{v}} \right) + \bar{\theta} \left( \dfrac{1}{\bar{v}} - \dfrac{1}{v} \right) \left[ \sigma (v - \bar{v})_\xi + \dot{\mathbf{X}}(t) v_\xi^S + (u - \bar{u})_\xi \right] \\[3mm]
			\displaystyle \qquad + \bar{\theta} \left( \dfrac{1}{\bar{v}} - \dfrac{1}{v} \right) \left( 1 - \dfrac{v}{\bar{v}} \right) \left[ \sigma \bar{v}_\xi - \dot{\mathbf{X}}(t) v_\xi^S + \bar{u}_\xi \right].
		\end{array}
	\end{align*}
	Since
	\begin{align*}
		\begin{array}{l}
			\displaystyle \partial_\xi \left[ \sigma \bar{\theta} \, \Psi\left( \dfrac{v}{\bar{v}} \right) \right] = \sigma \bar{\theta}_\xi \, \Psi\left( \dfrac{v}{\bar{v}} \right) + \sigma \bar{\theta} \, \Psi'\left( \dfrac{v}{\bar{v}} \right) \left( \dfrac{v}{\bar{v}} \right)_\xi \\[3mm]
			\displaystyle = \sigma \bar{\theta}_\xi \, \Psi\left( \dfrac{v}{\bar{v}} \right) + \sigma \bar{\theta} \left( \dfrac{1}{\bar{v}} - \dfrac{1}{v} \right) \left[ (v - \bar{v})_\xi + \bar{v}_\xi \left( 1 - \dfrac{v}{\bar{v}} \right) \right].
		\end{array}
	\end{align*}
	Then
	\begin{align} \label{pf-v}
		\begin{array}{l}
			\displaystyle \partial_t \left[ R \bar{\theta} \, \Psi\left( \dfrac{v}{\bar{v}} \right) \right] = \partial_\xi \left[ \sigma R \bar{\theta} \, \Psi\left( \dfrac{v}{\bar{v}} \right) \right] + R (\bar{\theta}_t - \sigma \bar{\theta}_\xi) \Psi\left( \dfrac{v}{\bar{v}} \right) \\[3mm]
			\displaystyle \qquad + R \bar{\theta} \left( \dfrac{1}{\bar{v}} - \dfrac{1}{v} \right) \dfrac{v}{\bar{v}} v_\xi^S \dot{\mathbf{X}}(t) + R \bar{\theta} \left( \dfrac{1}{\bar{v}} - \dfrac{1}{v} \right) (u - \bar{u})_\xi \\[3mm]
			\displaystyle \qquad + R \bar{\theta} \left( \dfrac{1}{\bar{v}} - \dfrac{1}{v} \right) \left( 1 - \dfrac{v}{\bar{v}} \right) \bar{u}_\xi \\[3mm]
			\displaystyle = \partial_\xi \left[ \sigma R \bar{\theta} \, \Psi\left( \dfrac{v}{\bar{v}} \right) \right] + R (\bar{\theta}_t - \sigma \bar{\theta}_\xi) \Psi\left( \dfrac{v}{\bar{v}} \right) + \dot{\mathbf{X}}(t) \dfrac{p(\bar{v}, \bar{\theta})}{\bar{v}} (v - \bar{v}) v_\xi^S \\[3mm]
			\displaystyle \qquad - \dfrac{p(\bar{v}, \bar{\theta})}{\bar{v} v} (v - \bar{v})^2 \bar{u}_\xi + R \bar{\theta} \left( \dfrac{1}{\bar{v}} - \dfrac{1}{v} \right) (u - \bar{u})_\xi.
		\end{array}
	\end{align}
	Similarly, using $\reff{PS}_3$ and $\reff{RCDS}_3$, we obtain
	\begin{align*}
		\begin{array}{l}
			\displaystyle \partial_t \left[ \dfrac{R}{\gamma - 1} \bar{\theta} \, \Psi\left( \dfrac{\theta}{\bar{\theta}} \right) \right] = \dfrac{R}{\gamma - 1} \bar{\theta}_t \, \Psi\left( \dfrac{\theta}{\bar{\theta}} \right) + \dfrac{R \bar{\theta}}{\gamma - 1} \Psi'\left( \dfrac{\theta}{\bar{\theta}} \right) \left( \dfrac{\theta}{\bar{\theta}} \right)_t \\[3mm]
			\displaystyle = \dfrac{R}{\gamma - 1} \bar{\theta}_t \, \Psi\left( \dfrac{\theta}{\bar{\theta}} \right) + \bar{\theta} \left( \dfrac{1}{\bar{\theta}} - \dfrac{1}{\theta} \right) \left[ \dfrac{R}{\gamma - 1} (\theta - \bar{\theta})_t + \dfrac{R}{\gamma - 1} \bar{\theta}_t \left( 1 - \dfrac{\theta}{\bar{\theta}} \right) \right] \\[3mm]
			\displaystyle = \dfrac{R}{\gamma - 1} \bar{\theta}_t \, \Psi\left( \dfrac{\theta}{\bar{\theta}} \right) + \bar{\theta} \left( \dfrac{1}{\bar{\theta}} - \dfrac{1}{\theta} \right) \left[ \dfrac{\sigma R}{\gamma - 1} (\theta - \bar{\theta})_\xi + \dot{\mathbf{X}}(t) \dfrac{R}{\gamma - 1} \theta_\xi^S - (p(v, \theta) u_\xi - p(\bar{v}, \bar{\theta}) \bar{u}_\xi) \right. \\[3mm]
			\displaystyle \qquad \left. + \kappa \left( \dfrac{\theta_\xi}{v} - \dfrac{\bar{\theta}_\xi}{\bar{v}} \right)_\xi + \mu \left( \dfrac{u_\xi^2}{v} - \dfrac{\bar{u}_\xi^2}{\bar{v}} \right) - Q_2 \right] \\[3mm]
			\displaystyle \qquad + \bar{\theta} \left( \dfrac{1}{\bar{\theta}} - \dfrac{1}{\theta} \right) \left( 1 - \dfrac{\theta}{\bar{\theta}} \right) \left[ \dfrac{\sigma R}{\gamma - 1} \bar{\theta}_\xi - \dot{\mathbf{X}}(t) \dfrac{R}{\gamma - 1} \theta_\xi^S - p(\bar{v}, \bar{\theta}) \bar{u}_\xi + \kappa \left( \dfrac{\bar{\theta}_\xi}{\bar{v}} \right)_\xi + \mu \dfrac{\bar{u}_\xi^2}{\bar{v}} + Q_2 \right].
		\end{array}
	\end{align*}
	and
	\begin{align*}
		\begin{array}{l}
			\displaystyle \partial_\xi \left[ \dfrac{\sigma R}{\gamma - 1} \bar{\theta} \, \Psi\left( \dfrac{\theta}{\bar{\theta}} \right) \right] = \dfrac{\sigma R}{\gamma - 1} \bar{\theta}_\xi \, \Psi\left( \dfrac{\theta}{\bar{\theta}} \right) + \dfrac{\sigma R}{\gamma - 1} \bar{\theta} \, \Psi'\left( \dfrac{\theta}{\bar{\theta}} \right) \left( \dfrac{\theta}{\bar{\theta}} \right)_\xi \\[3mm]
			\displaystyle = \dfrac{\sigma R}{\gamma - 1} \bar{\theta}_\xi \, \Psi\left( \dfrac{\theta}{\bar{\theta}} \right) + \dfrac{\sigma R}{\gamma - 1} \bar{\theta} \left( \dfrac{1}{\bar{\theta}} - \dfrac{1}{\theta} \right) \left[ (\theta - \bar{\theta})_\xi + \bar{\theta}_\xi \left( 1 - \dfrac{\theta}{\bar{\theta}} \right) \right].
		\end{array}
	\end{align*}
	Then it holds
	\begin{align} \label{pf-theta}
		\begin{array}{l}
			\displaystyle \partial_t \left[ \dfrac{R}{\gamma - 1} \bar{\theta} \, \Psi\left( \dfrac{\theta}{\bar{\theta}} \right) \right] = \partial_\xi \left[ \dfrac{\sigma R}{\gamma - 1} \bar{\theta} \, \Psi\left( \dfrac{\theta}{\bar{\theta}} \right) \right] + \dfrac{R}{\gamma - 1} (\bar{\theta}_t - \sigma \bar{\theta}_\xi) \Psi\left( \dfrac{\theta}{\bar{\theta}} \right) \\[3mm]
			\displaystyle \qquad + \dot{\mathbf{X}}(t) \dfrac{R}{(\gamma - 1) \bar{\theta}} (\theta - \bar{\theta}) \theta_\xi^S - \dfrac{R}{v} (\theta - \bar{\theta}) (u - \bar{u})_\xi - \dfrac{\bar{u}_\xi}{\theta} (\theta - \bar{\theta}) (p(v, \theta) - p(\bar{v}, \bar{\theta})) \\[3mm]
			\displaystyle \qquad + \dfrac{p(\bar{v}, \bar{\theta}) (\theta - \bar{\theta})^2}{\bar{\theta} \theta} \bar{u}_\xi + \kappa \left[ \dfrac{\theta - \bar{\theta}}{\theta} \left( \dfrac{\theta_\xi}{v} - \dfrac{\bar{\theta}_\xi}{\bar{v}} \right) \right]_\xi - \dfrac{\kappa}{\theta v} |(\theta - \bar{\theta})_\xi|^2 \\[3mm]
			\displaystyle \qquad - \kappa \dfrac{(\theta - \bar{\theta})_\xi}{\theta} \left( \dfrac{1}{v} - \dfrac{1}{\bar{v}} \right) \bar{\theta}_\xi + \kappa \dfrac{(\theta - \bar{\theta}) \theta_\xi}{\theta^2} \left( \dfrac{\theta_\xi}{v} - \dfrac{\bar{\theta}_\xi}{\bar{v}} \right) + \mu \dfrac{\theta - \bar{\theta}}{\theta} \left( \dfrac{u_\xi^2}{v} - \dfrac{\bar{u}_\xi^2}{\bar{v}} \right) \\[3mm]
			\displaystyle \qquad - \dfrac{\theta - \bar{\theta}}{\theta} Q_2 - \dfrac{(\theta - \bar{\theta})^2}{\theta \bar{\theta}} \left[ \kappa \left( \dfrac{\bar{\theta}_\xi}{\bar{v}} \right)_\xi + \mu \dfrac{\bar{u}_\xi^2}{\bar{v}} + Q_2 \right].
		\end{array}
	\end{align}
	Similarly, multiplying $\reff{PS}_2$ by $(u - \bar{u})$ yields
	\begin{align} \label{pf-u}
		\begin{array}{l}
			\displaystyle \partial_t \left[ \dfrac{(u - \bar{u})^2}{2} \right] = \partial_\xi \left[ \dfrac{\sigma (u - \bar{u})^2}{2} \right] + \dot{\mathbf{X}}(t) (u - \bar{u}) u_\xi^S - (u - \bar{u}) (\hat{p}(v, \theta) - \hat{p}(\bar{v}, \bar{\theta}))_\xi \\[3mm]
			\displaystyle \qquad + \mu (u - \bar{u}) \left( \dfrac{u_\xi}{v} - \dfrac{\bar{u}_\xi}{\bar{v}} \right)_\xi + (u - \bar{u}) \left( \Phi(v, \phi) - \Phi(\bar{v}, \bar{\phi}) \right)_\xi - (u - \bar{u}) Q_1 \\[3mm]
			\displaystyle = \partial_\xi \left[ \dfrac{\sigma (u - \bar{u})^2}{2} \right] + \dot{\mathbf{X}}(t) (u - \bar{u}) u_\xi^S - \partial_\xi \left[ (u - \bar{u}) (\hat{p}(v, \theta) - \hat{p}(\bar{v}, \bar{\theta})) \right] \\[3mm]
			\displaystyle \qquad + (\hat{p}(v, \theta) - \hat{p}(\bar{v}, \bar{\theta})) (u - \bar{u})_\xi + \mu \left[ (u - \bar{u}) \left( \dfrac{u_\xi}{v} - \dfrac{\bar{u}_\xi}{\bar{v}} \right) \right]_\xi - \dfrac{\mu}{v} |(u - \bar{u})_\xi|^2\\[3mm]
			\displaystyle \qquad - \mu (u - \bar{u})_\xi \bar{u}_\xi \left( \dfrac{1}{v} - \dfrac{1}{\bar{v}} \right) + \left( \Phi(v, \phi) - \Phi(\bar{v}, \bar{\phi}) \right)_\xi (u - \bar{u}) - (u - \bar{u}) Q_1. 
		\end{array}
	\end{align}
	Since
	\begin{align*}
		\begin{array}{l}
			\displaystyle R \bar{\theta} \left( \dfrac{1}{v} - \dfrac{1}{\bar{v}} \right) (u - \bar{u})_\xi - \dfrac{R}{v} (\theta - \bar{\theta}) (u - \bar{u})_\xi + (\hat{p}(v, \theta) - \hat{p}(\bar{v}, \bar{\theta})) (u - \bar{u})_\xi = \left( \dfrac{1}{v} - \dfrac{1}{\bar{v}} \right) (u - \bar{u})_\xi,
		\end{array}
	\end{align*}
	and
	\begin{align*}
		\bar{\theta}_t - \sigma \bar{\theta}_\xi = (\theta_t^R - \sigma \theta_\xi^R) + (\theta_t^C - \sigma \theta_\xi^C) - \sigma \theta_\xi^S - \dot{\mathbf{X}}(t) \theta_\xi^S,
	\end{align*}
	it follows from \reff{pf-v}-\reff{pf-u} that
	\begin{align} \label{pf-vthetau}
		\begin{array}{l}
			\displaystyle \partial_t \left[ R \bar{\theta} \, \Psi\left( \dfrac{v}{\bar{v}} \right) + \dfrac{R \bar{\theta}}{\gamma - 1} \, \Psi\left( \dfrac{\theta}{\bar{\theta}} \right) + \dfrac{(u - \bar{u})^2}{2} \right] \\[4mm]
			\displaystyle = \partial_\xi \left[ \sigma a^{-\mathbf{X}} \left( R \bar{\theta} \, \Psi\left( \dfrac{v}{\bar{v}} \right) + \dfrac{R \bar{\theta}}{\gamma - 1} \, \Psi\left( \dfrac{\theta}{\bar{\theta}} \right) + \dfrac{(u - \bar{u})^2}{2} \right) - a^{-\mathbf{X}} (u - \bar{u}) (\hat{p}(v, \theta) - \hat{p}(\bar{v}, \bar{\theta})) \right] \\[4mm]
			\displaystyle \quad - \sigma a_\xi^{-\mathbf{X}} \left( R \bar{\theta} \, \Psi\left( \dfrac{v}{\bar{v}} \right) + \dfrac{R \bar{\theta}}{\gamma - 1}  \, \Psi\left( \dfrac{\theta}{\bar{\theta}} \right) + \dfrac{(u - \bar{u})^2}{2} \right) + a_\xi^{-\mathbf{X}} (u - \bar{u}) (\hat{p}(v, \theta) - \hat{p}(\bar{v}, \bar{\theta})) \\[4mm]
			\displaystyle \quad - a^{-\mathbf{X}} \left[ \dfrac{\kappa}{v \theta} |(\theta - \bar{\theta})_\xi|^2 + \dfrac{\mu}{v} |(u - \bar{u})_\xi|^2 \right] \\[4mm]
			\displaystyle \quad + \dot{\mathbf{X}}(t) a^{-\mathbf{X}} \left\{ \left[ \dfrac{p(\bar{v}, \bar{\theta}) v^S_\xi}{\bar{v}} (v - \bar{v}) + \dfrac{R \theta_\xi^S}{(\gamma - 1) \bar{\theta}} (\theta - \bar{\theta}) + u_\xi^S (u - \bar{u}) \right] \right. \\[4mm]
			\displaystyle \qquad\qquad\qquad\quad  \left.+ \left[ - R \theta_\xi^S \, \Psi\left( \dfrac{v}{\bar{v}} \right) - \dfrac{R}{\gamma - 1} \theta_\xi^S \, \Psi\left( \dfrac{\theta}{\bar{\theta}} \right) \right] \right\} \\[4mm]
			\displaystyle \quad + a^{-\mathbf{X}} \left\{ \left[ R \left( (\theta_t^R - \sigma \theta_\xi^R) + (\theta_t^C - \sigma \theta_\xi^C) - \sigma \theta_\xi^S \right) \Psi\left( \dfrac{v}{\bar{v}} \right) - \dfrac{p(\bar{v}, \bar{\theta}) \bar{u}_\xi}{\bar{v} v} (v - \bar{v})^2 \right] \right. \\[4mm]
			\displaystyle \qquad\qquad\quad + \left[ \dfrac{R}{\gamma - 1} \left( (\theta_t^R - \sigma \theta_\xi^R) + (\theta_t^C - \sigma \theta_\xi^C) - \sigma \theta_\xi^S \right) \Psi\left( \dfrac{\theta}{\bar{\theta}} \right) \right. \\[4mm]
			\displaystyle \qquad\qquad\qquad  \left.\left. - \dfrac{\bar{u}_\xi}{\theta} (\theta - \bar{\theta}) (p(v, \theta) - p(\bar{v}, \bar{\theta})) + \dfrac{p(\bar{v}, \bar{\theta}) \bar{u}_\xi}{\bar{\theta} \theta} (\theta - \bar{\theta})^2 \right] \right\} \\[4mm]
			\displaystyle \quad + a^{-\mathbf{X}} \left\{ \left[ \kappa \dfrac{\theta - \bar{\theta}}{\theta} \left( \dfrac{\theta_\xi}{v} - \dfrac{\bar{\theta}_\xi}{\bar{v}} \right) + \mu (u - \bar{u}) \left( \dfrac{u_\xi}{v} - \dfrac{\bar{u}_\xi}{\bar{v}} \right) \right]_\xi + \left[ - \kappa \dfrac{(\theta - \bar{\theta})_\xi}{\theta} \left( \dfrac{1}{v} - \dfrac{1}{\bar{v}} \right) \bar{\theta}_\xi \right. \right. \\[4mm]
			\displaystyle \quad \qquad \left. + \kappa \dfrac{(\theta - \bar{\theta}) \theta_\xi}{\theta^2} \left( \dfrac{\theta_\xi}{v} - \dfrac{\bar{\theta}_\xi}{\bar{v}} \right) + \mu \dfrac{\theta - \bar{\theta}}{\theta} \left( \dfrac{u_\xi^2}{v} - \dfrac{\bar{u}_\xi^2}{\bar{v}} \right) - \dfrac{(\theta - \bar{\theta})^2}{\theta \bar{\theta}} \left( \kappa \left( \dfrac{\bar{\theta}_\xi}{\bar{v}} \right)_\xi + \mu \dfrac{\bar{u}_\xi^2}{\bar{v}} + Q_2 \right) \right] \\[4mm]
			\displaystyle \quad \qquad \left. - \left[ (u - \bar{u}) Q_1 + \dfrac{\theta - \bar{\theta}}{\theta} Q_2 \right] \right\} \\[4mm]
			\displaystyle \quad + a^{-\mathbf{X}} \left( \dfrac{1}{v} - \dfrac{1}{\bar{v}} \right) (u - \bar{u})_\xi + a^{-\mathbf{X}} \left( \Phi(v, \phi) - \Phi(\bar{v}, \bar{\phi}) \right)_\xi (u - \bar{u}). 
		\end{array}
	\end{align}
	We will rewrite the last two terms on the right-hand side of \reff{pf-vthetau}.
	By $\reff{PS}_1$, $\bar{v}_t = v_t^R + v_t^C - \dot{\mathbf{X}}(t) v_\xi^S$ and $ \bar{v}_\xi = v^R_\xi + v^C_\xi + v^S_\xi $, it yields 
	\begin{align} \label{pf-hat-p}
		\begin{array}{l}
			\displaystyle a^{-\mathbf{X}} \left( \dfrac{1}{v} - \dfrac{1}{\bar{v}} \right) (u - \bar{u})_\xi = - a^{-\mathbf{X}} \dfrac{v - \bar{v}}{v \bar{v}} (u - \bar{u})_\xi \\[4mm]
			\displaystyle = - a^{-\mathbf{X}} \left[ \dfrac{v - \bar{v}}{\bar{v}^2} - \dfrac{(v - \bar{v})^2}{\bar{v}^3} \right] (u - \bar{u})_\xi - a^{-\mathbf{X}} \dfrac{(v - \bar{v})^3}{\bar{v}^3 v} (u - \bar{u})_\xi \\[4mm]
			\displaystyle = - a^{-\mathbf{X}} \left[ \dfrac{v - \bar{v}}{\bar{v}^2} - \dfrac{(v - \bar{v})^2}{\bar{v}^3} \right] \left( (v - \bar{v})_t - \sigma (v - \bar{v})_\xi - \dot{\mathbf{X}}(t) v_\xi^S \right) - a^{-\mathbf{X}} \dfrac{(v - \bar{v})^3}{\bar{v}^3 v} (u - \bar{u})_\xi \\[4mm]
			\displaystyle = - \left[ \dfrac{a^{-\mathbf{X}}}{\bar{v}^2} \left( \dfrac{(v - \bar{v})^2}{2} \right)_t  - \dfrac{a^{-\mathbf{X}}}{\bar{v}^3} \left( \dfrac{(v - \bar{v})^3}{3} \right)_t \right] + \sigma \left[ \dfrac{a^{-\mathbf{X}}}{\bar{v}^2} \left( \dfrac{(v - \bar{v})^2 }{2} \right)_\xi - \dfrac{a^{-\mathbf{X}}}{\bar{v}^3} \left( \dfrac{(v - \bar{v})^3 }{3} \right)_\xi \right] \\[4mm]
			\displaystyle \qquad + \dot{\mathbf{X}}(t) a^{-\mathbf{X}} v_\xi^S \left( \dfrac{v - \bar{v}}{\bar{v}^2} - \dfrac{(v - \bar{v})^2}{\bar{v}^3} \right) - a^{-\mathbf{X}} \dfrac{(v - \bar{v})^3}{\bar{v}^3 v} (u - \bar{u})_\xi \\[4mm]
			\displaystyle = - \partial_t \left[ a^{-\mathbf{X}} \left( \dfrac{(v - \bar{v})^2}{2 \bar{v}^2} - \dfrac{(v - \bar{v})^3}{3 \bar{v}^3} \right) \right] + \partial_\xi \left[ \sigma a^{-\mathbf{X}} \left( \dfrac{(v - \bar{v})^2}{2 \bar{v}^2} - \dfrac{(v - \bar{v})^3}{3 \bar{v}^3} \right) \right] \\[4mm]
			\displaystyle \qquad + \dot{\mathbf{X}}(t) \left\{ a_\xi^{-\mathbf{X}} \left( - \dfrac{(v - \bar{v})^2}{2 \bar{v}^2} + \dfrac{(v - \bar{v})^3}{3 \bar{v}^3} \right) + a^{-\mathbf{X}} \left( \dfrac{v_\xi^S}{\bar{v}^2} (v - \bar{v}) - \dfrac{v_\xi^S}{\bar{v}^4} (v - \bar{v})^3 \right) \right\} \\[4mm]
			\displaystyle \qquad - a^{-\mathbf{X}} ( ( v_t^R - \sigma v_\xi^R ) + ( v_t^C - \sigma v_\xi^C ) - \sigma v^S_\xi ) \dfrac{(v - \bar{v})^2}{\bar{v}^3} - \sigma a_\xi^{-\mathbf{X}} \left( \dfrac{(v - \bar{v})^2}{2 \bar{v}^2} - \dfrac{(v - \bar{v})^3}{3 \bar{v}^3} \right) \\[4mm]
			\displaystyle \qquad + a^{-\mathbf{X}} \left[ \left( v_t^R + v_t^C - \sigma \bar{v}_\xi \right) \dfrac{(v - \bar{v})^3}{\bar{v}^4} - \dfrac{(v - \bar{v})^3}{\bar{v}^3 v} (u - \bar{u})_\xi \right]. 
		\end{array}
	\end{align}
	By $\reff{hat-p-Phi}_2$,
	\begin{align} \label{pf-Phi}
		\begin{array}{l}
			\displaystyle a^{-\mathbf{X}} \left( \Psi(v, \phi) - \Psi(\bar{v}, \bar{\phi}) \right)_\xi (u - \bar{u}) \\[4mm]
			\displaystyle = \partial_\xi \left[ a^{-\mathbf{X}} (u - \bar{u}) (\Psi(v, \phi) - \Psi(\bar{v}, \bar{\phi})) \right] - (a^{-\mathbf{X}} (u - \bar{u}))_\xi (\Psi(v, \phi) - \Psi(\bar{v}, \bar{\phi})) \\[4mm]
			\displaystyle = \partial_\xi \left[ a^{-\mathbf{X}} (u - \bar{u}) (\Psi(v, \phi) - \Psi(\bar{v}, \bar{\phi})) \right] - (a^{-\mathbf{X}} (u - \bar{u}))_\xi \bigg[ \dfrac{\nu^2}{2} \Big( \big( \dfrac{\phi_\xi}{v} \big)^2 - \big( \dfrac{\bar{\phi}_\xi}{\bar{v}} \big)^2 \Big) \\[4mm]
			\displaystyle \quad + \nu^2 \left( \dfrac{\phi_\xi v_\xi}{v^3} - \dfrac{\bar{\phi}_\xi \bar{v}_\xi}{\bar{v}^3} \right) - \nu^2 \left( \dfrac{1}{v^2} - \dfrac{1}{\bar{v}^2} \right) \phi_{\xi\xi} \bigg]  + \nu^2 a_\xi^{-\mathbf{X}} (u - \bar{u}) \dfrac{\tilde{\phi}_{\xi\xi}}{\bar{v}^2} + \nu^2 a^{-\mathbf{X}} (u - \bar{u})_\xi \dfrac{\tilde{\phi}_{\xi\xi}}{\bar{v}^2}. 
		\end{array}
	\end{align}
	From $\reff{PS}_1$, the last term of the above equality is the main difficult term and can be represented as
	\begin{align}  \label{pf-Phi1}
		\begin{array}{l}
			\displaystyle \nu^2 a^{-\mathbf{X}} (u - \bar{u})_\xi \dfrac{\tilde{\phi}_{\xi\xi}}{\bar{v}^2} = \nu^2 a^{-\mathbf{X}} \dfrac{\tilde{\phi}_{\xi\xi}}{\bar{v}^2} \left( (v - \bar{v})_t - \sigma (v - \bar{v})_\xi - \dot{\mathbf{X}}(t) v_\xi^S \right) \\[4mm]
			\displaystyle = \partial_t \left[ \nu^2 a^{-\mathbf{X}} (v - \bar{v}) \dfrac{\tilde{\phi}_{\xi\xi}}{\bar{v}^2} \right] - \nu^2 \dfrac{a^{-\mathbf{X}}}{\bar{v}^2} (v - \bar{v}) \tilde{\phi}_{\xi\xi t} - \nu^2 \left( \dfrac{a^{-\mathbf{X}}}{\bar{v}^2} \right)_t (v - \bar{v}) \tilde{\phi}_{\xi\xi} \\[4mm]
			\displaystyle \quad - \sigma \nu^2 \dfrac{a^{-\mathbf{X}}}{\bar{v}^2} (v - \bar{v})_\xi \tilde{\phi}_{\xi\xi} - \dot{\mathbf{X}}(t) \nu^2 a^{-\mathbf{X}} v_\xi^S \dfrac{\tilde{\phi}_{\xi\xi}}{\bar{v}^2} \\[4mm]
			\displaystyle = \partial_t \left[ \nu^2 a^{-\mathbf{X}} (v - \bar{v}) \dfrac{\tilde{\phi}_{\xi\xi}}{\bar{v}^2} \right] - \partial_\xi \left[ \nu^2 a^{-\mathbf{X}} (v - \bar{v}) \dfrac{\tilde{\phi}_{\xi\xi}}{\bar{v}^2} \right] \\[4mm]
			\displaystyle \quad + \dot{\mathbf{X}}(t) \left[ \nu^2 a_\xi^{-\mathbf{X}} (v - \bar{v}) \dfrac{\tilde{\phi}_{\xi\xi}}{\bar{v}^2} - a^{-\mathbf{X}} \left( 2 \nu^2 v_\xi^S (v - \bar{v}) \dfrac{\tilde{\phi}_{\xi\xi}}{\bar{v}^3} + \nu^2 v_\xi^S \dfrac{\tilde{\phi}_{\xi\xi}}{\bar{v}^2} \right) \right] \\[4mm]
			\displaystyle \quad + 2 \nu^2 a^{-\mathbf{X}} (v_t^R + v_t^C) (v - \bar{v}) \dfrac{\tilde{\phi}_{\xi\xi}}{\bar{v}^3} + \nu^2 \left( \dfrac{a^{-\mathbf{X}}}{\bar{v}^2} \right)_\xi (v - \bar{v}) \tilde{\phi}_{\xi t} + \nu^2 \dfrac{a^{-\mathbf{X}}}{\bar{v}^2} (v - \bar{v})_\xi (\tilde{\phi}_{\xi t} - \sigma \tilde{\phi}_{\xi\xi}). 
		\end{array}
	\end{align}
	By $\reff{PS}_4$, we have  
	\begin{align*}
		\displaystyle v - \bar{v} = \nu^2 e^{-\bar{\phi}} \dfrac{\tilde{\phi}_{\xi\xi}}{\bar{v}} - \bar{v} \tilde{\phi} - \nu^2 e^{-\bar{\phi}} \dfrac{\tilde{\phi}_\xi \bar{v}_\xi}{\bar{v}^2} - \nu^2 e^{-\bar{\phi}} \left( \dfrac{v - \bar{v}}{\bar{v} v} \phi_\xi \right)_\xi + N_p - e^{-\bar{\phi}} Q_3,
	\end{align*}
	which yields 
	\begin{align*}
		\displaystyle (v - \bar{v})_\xi = \nu^2 \left(e^{-\bar{\phi}} \dfrac{\tilde{\phi}_{\xi\xi}}{\bar{v}}\right)_\xi - \bar{v} \tilde{\phi}_\xi + \left[ - \nu^2 e^{-\bar{\phi}} \tilde{\phi}_\xi \cdot \dfrac{\bar{v}_\xi}{\bar{v}^2} - \nu^2 e^{-\bar{\phi}} \left( \dfrac{v - \bar{v}}{\bar{v} v} \phi_\xi \right)_\xi + N_p - e^{-\bar{\phi}} Q_3 \right]_\xi - \bar{v}_\xi \tilde{\phi}.
	\end{align*}
	Then it follows  
	\begin{align}  \label{pf-Phi2}
		\begin{array}{l}
			\displaystyle \dfrac{a^{-\mathbf{X}}}{\bar{v}^2} (v - \bar{v})_\xi (\tilde{\phi}_{\xi t} - \sigma \tilde{\phi}_{\xi\xi}) \\[4mm]
			\displaystyle = \dfrac{a^{-\mathbf{X}}}{\bar{v}^2} (\tilde{\phi}_{\xi t} - \sigma \tilde{\phi}_{\xi\xi}) \, \nu^2 \left( e^{-\bar{\phi}} \dfrac{\tilde{\phi}_{\xi\xi}}{\bar{v}} \right)_\xi - \dfrac{a^{-\mathbf{X}}}{\bar{v}} (\tilde{\phi}_{\xi t} - \sigma \tilde{\phi}_{\xi\xi}) \tilde{\phi}_\xi \\[4mm]
			\displaystyle \qquad + \dfrac{a^{-\mathbf{X}}}{\bar{v}^2} (\tilde{\phi}_{\xi t} - \sigma \tilde{\phi}_{\xi\xi}) \left\{ \left[ -\nu^2 e^{-\bar{\phi}} \tilde{\phi}_\xi \cdot \dfrac{\bar{v}_\xi}{\bar{v}^2} - \nu^2 e^{-\bar{\phi}} \left( \dfrac{v - \bar{v}}{\bar{v} v} \phi_\xi \right)_\xi + N_p - e^{-\bar{\phi}} Q_3 \right]_\xi - \bar{v}_\xi \tilde{\phi} \right\} \\[4mm]
			\displaystyle = \partial_t \left[ - \dfrac{\nu^2 a^{-\mathbf{X}}}{2 \bar{v}^3 e^{\bar{\phi}}} \tilde{\phi}_{\xi\xi}^2 - \dfrac{a^{-\mathbf{X}}}{2 \bar{v}} \tilde{\phi}_\xi^2 \right] + \partial_\xi \left[ \dfrac{\nu^2 a^{-\mathbf{X}}}{\bar{v}^3 e^{\bar{\phi}}} \tilde{\phi}_{\xi t} \cdot \tilde{\phi}_{\xi\xi} - \dfrac{\sigma \nu^2 a^{-\mathbf{X}}}{2 \bar{v}^3 e^{\bar{\phi}}} \tilde{\phi}_{\xi\xi}^2 + \dfrac{\sigma a^{-\mathbf{X}}}{2 \bar{v}} \tilde{\phi}_\xi^2 \right] \\[4mm]
			\displaystyle \qquad + \dot{\mathbf{X}}(t) \left\{ a_\xi^{-\mathbf{X}} \left( - \dfrac{\nu^2}{2 \bar{v}^3 e^{\bar{\phi}}} \tilde{\phi}_{\xi\xi}^2 - \dfrac{1}{2 \bar{v}} \tilde{\phi}_\xi^2 \right) + a^{-\mathbf{X}} \left[ \dfrac{\nu^2}{2} \left(\dfrac{3 v_\xi^S}{\bar{v}^4 e^{\bar{\phi}}} + \dfrac{\phi_\xi^S}{\bar{v}^3 e^{\bar{\phi}}} \right) \tilde{\phi}_{\xi\xi}^2 + \dfrac{v_\xi^S}{2 \bar{v}^2} \tilde{\phi}_\xi^2 \right] \right\} \\[4mm]
			\displaystyle \qquad + \left\{ a^{-\mathbf{X}} \dfrac{\nu^2}{2} \left( - \dfrac{3 (v_t^R + v_t^C)}{\bar{v}^4 e^{\bar{\phi}}} - \dfrac{\phi_t^R + \phi_t^C}{\bar{v}^3 e^{\bar{\phi}}} \right) \tilde{\phi}_{\xi\xi}^2 - a^{-\mathbf{X}} \dfrac{v_t^R + v_t^C}{2 \bar{v}^2} \tilde{\phi}_\xi^2 \right. \\[4mm]
			\displaystyle \qquad \qquad \left. - \dfrac{\nu^2}{\bar{v} e^{\bar{\phi}}} \left( \dfrac{a^{-\mathbf{X}}}{\bar{v}^2} \right)_\xi \tilde{\phi}_{\xi t} \cdot \tilde{\phi}_{\xi\xi} + \left( \dfrac{\nu^2 a^{-\mathbf{X}}}{\bar{v}^3 e^{\bar{\phi}}} \right)_\xi \dfrac{\sigma \tilde{\phi}_{\xi\xi}^2}{2} - \dfrac{\sigma \nu^2 a^{-\mathbf{X}}}{\bar{v}^2} \left( \dfrac{1}{\bar{v} e^{\bar{\phi}}} \right)_\xi \tilde{\phi}_{\xi\xi}^2 - \left( \dfrac{a^{-\mathbf{X}}}{\bar{v}} \right)_\xi \dfrac{\sigma \tilde{\phi}_\xi^2}{2} \right\} \\[4mm]
			\displaystyle \qquad + \dfrac{a^{-\mathbf{X}}}{\bar{v}^2} (\tilde{\phi}_{\xi t} - \sigma \tilde{\phi}_{\xi\xi}) \left\{ \left[ -\nu^2 e^{-\bar{\phi}} \tilde{\phi}_\xi \cdot \dfrac{\bar{v}_\xi}{\bar{v}^2} - \nu^2 e^{-\bar{\phi}} \left( \dfrac{v - \bar{v}}{\bar{v} v} \phi_\xi \right)_\xi + N_p - e^{-\bar{\phi}} Q_3 \right]_\xi - \bar{v}_\xi \tilde{\phi} \right\}.
		\end{array}
	\end{align}
	combining \reff{pf-hat-p}-\reff{pf-Phi2} with \reff{pf-vthetau} and classifying the terms into five categories, we get the desired representation of this lemma.
\end{proof}
\hspace*{\fill}

From the construction of shift \reff{Xt-def}, it implies 
\begin{align*}
	\dot{\mathbf{X}}(t) = - \dfrac{M}{\delta_S} \mathbf{Y}_1(U),
\end{align*}
which yields
\begin{align*}
	\displaystyle \dot{\mathbf{X}}(t) \sum_{i=1}^{4}\mathbf{Y}_i(U) = - \dfrac{\delta_S}{M} |\dot{\mathbf{X}}(t)|^2 + \dot{\mathbf{X}}(t) \sum_{i=2}^4 \mathbf{Y}_i(U).
\end{align*}
Then it follows from \reff{RE-ori} in \ref{REM-lem} that 
\begin{align}  \label{RE-mid1}
	\begin{array}{l}
		\displaystyle \dfrac{d}{dt} \int_{\mathbb{R}} a^{-\mathbf{X}}(\xi) \, \tilde{\eta}\big( (U, \phi) \, | \, (\bar{U}, \bar{\phi}) \big) \, d\xi = \underset{\mathbf{MT}}{\underbrace{- \dfrac{\delta_S}{2M} |\dot{\mathbf{X}}(t)|^2 + \mathbf{B}_1 + \mathbf{B}_2 - \mathbf{G} - \dfrac{3}{4} \mathbf{D}}} \\[4mm]
		\displaystyle \qquad - \dfrac{\delta_S}{2M} |\dot{\mathbf{X}}(t)|^2 + \dot{\mathbf{X}}(t) \sum_{i=2}^4 \mathbf{Y}_i - \dfrac{1}{4} \mathbf{D} + \mathbf{E}_1 + \mathbf{E}_2.
	\end{array}
\end{align}
where the $\mathbf{MT}$ term is the leading order and the remaining terms are errors that need to be controlled.
\\ 

\subsection{Leading order estimates}
\hspace*{\fill}\\ 

We will give the estimates for the leading order term $\mathbf{MT}$ throught the below Poincare$\acute{}$-type inequality.
\\ 

\begin{lemma}[\cite{Kang23}  Appendix B] \label{Poincare} 
    For any $f:[0,1]\rightarrow\mathbb{R}$ satisfying $\displaystyle \int_0^1 y(1-y)|f^\prime|^2 dy<\infty$,
    \begin{align*}
        \int_0^1 \left|f-\int_0^1 f dy \right|^2 dy \le \dfrac{1}{2}\int_0^1 y(1-y)|f^\prime|^2 dy.
    \end{align*}
\end{lemma}

The leading order estimates are as follows.
\\ 

\begin{lemma} \label{LO-lem}
	There exists a generic constant $C > 0$ such that
	\begin{align*}
		\begin{array}{l}
			\displaystyle - \dfrac{\delta_S}{2M} |\dot{\mathbf{X}}(t)|^2 + \mathbf{B}_1 + \mathbf{B}_2 - \mathbf{G} - \dfrac{3}{4} \mathbf{D} \\[4mm]
			\displaystyle \le - C \mathcal{G}^R - C \mathcal{G}^S + C \delta_S^{4/3} (\delta_R + \delta_C)^{4/3} e^{-C \delta_S t} \\[4mm]
			\displaystyle \qquad + C \delta_C (1+t)^{-1} \int_\mathbb{R} e^{-\frac{C_1 |\xi + \sigma t|^2}{ ( 1 + t) }} |(v-\bar{v}, \theta-\bar{\theta})|^2 \, d\xi,
		\end{array}
	\end{align*}
	where
	\begin{align}  \label{LO-G_SR}
		\begin{array}{l}
			\displaystyle \mathcal{G}^R(U) := \int_{\mathbb{R}} |v_\xi^R| \cdot |(v - \bar{v}, \theta - \bar{\theta})|^2 \, d\xi, \\[3mm]
			\displaystyle \mathcal{G}^S(U) := \int_{\mathbb{R}} |v_\xi^S| \cdot |(v - \bar{v}, u - \bar{u}, \theta - \bar{\theta})|^2 \, d\xi. \\[3mm]
		\end{array}
	\end{align}
\end{lemma}

\begin{proof}[Proof]
	Before proving this lemma, we first define the two new variables
	\begin{align} \label{wy-def}
		w:=u-\bar{u} \qquad  \text{and}  \qquad y:=\dfrac{v^S(\xi-\mathbf{X}(t))-v^*}{\delta_S} \in (0,1),
	\end{align}
	which together with Lemma \ref{VS-pro} yields 
	\begin{align} \label{par-y}
		\dfrac{dy}{d\xi} = \dfrac{v^S_\xi}{\delta_S} > 0.
	\end{align}
	It follows from \reff{a-def} and \reff{par-a} that $ a^{-\mathbf{X}}(\xi) = 1 + \lambda y$ and $ a^{-\mathbf{X}}_\xi = \lambda \dfrac{dy}{d\xi}$.
	Next, we will use the $O(1)$ constants 
	\begin{align} \label{left-p-sig}
		p^*:=p(v^*, \theta^*), \quad \hat{p}^*:=\hat{p}(v^*, \theta^*), \quad \sigma^*:=\sigma_3(v^*, \theta^*)
	\end{align}
	with \reff{intro-p}, \reff{hat-p-Phi} and \reff{Eigen}.
	Note that these constants are independent of the small constants $\delta_0, \varepsilon_1$, since $\frac{v_+}{2} \le v^* \le v_+$ and $\frac{\theta}{2} \le \theta^* \le \theta_+$.
	the proof of this lemma is separated into five steps. 
	\\ 

	\noindent $\bullet$ {\bf Step 1: Estimates on $\mathbf{B}_2 - \mathbf{G}$.}
	
	First, it follows from \reff{hat-p-Phi} that $\hat{p}(v, \theta) - \hat{p}(\bar{v}, \bar{\theta}) = \dfrac{R}{v} (\theta - \bar{\theta}) - \dfrac{\hat{p}(\bar{v}, \bar{\theta})}{v} (v - \bar{v})$, then
	\begin{align*}
		\begin{array}{l}
			\displaystyle (u - \bar{u}) (\hat{p}(v, \theta) - \hat{p}(\bar{v}, \bar{\theta})) - \sigma \left[ \bar{\theta}\eta( U | \bar{U} )  + \dfrac{(v - \bar{v})^2}{2 \bar{v}^2} - \dfrac{(v - \bar{v})^3}{3 \bar{v}^3} \right] \\[4mm]
			\displaystyle = (u - \bar{u}) \left( \dfrac{R}{v} (\theta - \bar{\theta}) - \dfrac{\hat{p}(\bar{v}, \bar{\theta})}{v} (v - \bar{v}) \right) \\[4mm] 
			\displaystyle \qquad - \sigma \left[ R \bar{\theta} \, \Psi\left( \dfrac{v}{\bar{v}} \right) + \dfrac{R}{\gamma - 1} \bar{\theta} \, \Psi\left( \dfrac{\theta}{\bar{\theta}} \right) + \dfrac{(u - \bar{u})^2}{2} + \dfrac{(v - \bar{v})^2}{2 \bar{v}^2} - \dfrac{(v - \bar{v})^3}{3 \bar{v}^3} \right] \\[4mm]
			\displaystyle \le (u - \bar{u}) \left[ \dfrac{R}{v^*} (\theta - \bar{\theta}) - \dfrac{\hat{p}^*}{v^*} (v - \bar{v}) \right] - \sigma \left[ \dfrac{\hat{p}^*}{2 v^*} (v - \bar{v})^2 + \dfrac{R}{2 (\gamma - 1) \theta^* } (\theta - \bar{\theta})^2 + \dfrac{(u - \bar{u})^2}{2} \right] \\[4mm]
			\displaystyle \qquad + C \Big(|v - \bar{v}| + |\bar{v} - v^*| + |\bar{\theta} - \theta^*|\Big) |(v - \bar{v}, \, u - \bar{u}, \, \theta - \bar{\theta})|^2 \\[4mm]
			\displaystyle \le - \dfrac{\sigma^* \hat{p}^*}{2 v^*} \left[ (v - \bar{v}) + \dfrac{1}{\sigma^*} (u - \bar{u}) \right]^2 - \dfrac{ R \sigma^*}{2 (\gamma - 1) \theta^*} \left[ (\theta - \bar{\theta}) - \dfrac{(\gamma - 1) \theta^*}{\sigma^* v^*} (u - \bar{u}) \right]^2 \\[4mm]
			\displaystyle \qquad + C \Big(|v - \bar{v}| + \delta_S + |v^R - v_*| + |\theta^R - \theta_*| + |v^C - v^*| + |\theta^C - \theta^*|\Big) |(v - \bar{v}, \, u - \bar{u}, \, \theta - \bar{\theta})|^2.
		\end{array}
	\end{align*}
	Here the last inequality uses the fact that
	\begin{align*}
		\dfrac{\sigma^* \hat{p}^*}{2 v^*} \cdot \dfrac{1}{(\sigma^*)^2} + \dfrac{ R \sigma^*}{2 (\gamma - 1) \theta^*} \cdot \dfrac{(\gamma - 1)^2 (\theta^*)^2}{(\sigma^*)^2 v^*} = \dfrac{\sigma^*}{2}.
	\end{align*}
	Thus
	\begin{align*}
		\mathbf{B}_2(U) - \mathbf{G}(U) \le - \mathbf{G}_1(U) - \mathbf{G}_2(U) + \mathbf{B}_{\text{new}}(U),
	\end{align*}
	where
	\begin{align} \label{G_12}
		\begin{array}{l}
			\displaystyle \mathbf{G}_1(U) : = \dfrac{\sigma^* \hat{p}^*}{2 v^*} \int_{\mathbb{R}} a_\xi^{-\mathbf{X}} \left[ (v - \bar{v}) + \dfrac{1}{\sigma^*} (u - \bar{u}) \right]^2 d\xi, \\[4mm]
			\displaystyle \mathbf{G}_2(U) : = \dfrac{ R \sigma^*}{2 (\gamma - 1) \theta^*} \int_{\mathbb{R}} a_\xi^{-\mathbf{X}} \left[ (\theta - \bar{\theta}) - \dfrac{(\gamma - 1) \theta^*}{\sigma^* v^*} (u - \bar{u}) \right]^2 d\xi, \\[4mm]
			\displaystyle \mathbf{B}_{\text{new}}(U) : = C \delta_S \int_{\mathbb{R}} a_\xi^{-\mathbf{X}} |(v - \bar{v}, \, u - \bar{u}, \, \theta - \bar{\theta})|^2 d\xi \\[4mm]
			\displaystyle \qquad \qquad \qquad + C \int_{\mathbb{R}} a_\xi^{-\mathbf{X}} |(v - \bar{v}, \, u - \bar{u}, \, \theta - \bar{\theta})|^3 d\xi \\[4mm]
			\displaystyle \qquad \qquad \qquad + C \int_{\mathbb{R}} a_\xi^{-\mathbf{X}} \Big( |v^R - v_*| + |\theta^R - \theta_*| + |v^C - v^*| + |\theta^C - \theta^*| \Big) \\[4mm]
			\displaystyle \qquad \qquad \qquad \qquad \qquad \qquad \times |(v - \bar{v}, \, u - \bar{u}, \, \theta - \bar{\theta})|^2 d\xi.
		\end{array}
	\end{align}
	From above, the bad term $\mathbf{B}_2$ is controlled by the good term $\mathbf{G}$, we could extract the new good terms $\mathbf{G}_1, \mathbf{G}_2$, which play an important role in the following decoupling process.
	In what follows, we will give the estimates on $\mathbf{B}_{new}$. First of all, 
	\begin{align*}
		\begin{array}{l}
			\displaystyle C \delta_S \int_{\mathbb{R}} a_\xi^{-\mathbf{X}} |(v - \bar{v}, \, u - \bar{u}, \, \theta - \bar{\theta})|^2 d\xi \le C \lambda \int_{\mathbb{R}} v_\xi^S |(v - \bar{v}, \, u - \bar{u}, \, \theta - \bar{\theta})|^2 d\xi = C \lambda \mathcal{G}^S.
		\end{array}
	\end{align*}
	By $v^S_\xi > 0$ in Lemma \ref{VS-pro}, \reff{Ass-Inf} and \reff{lam-range}, 
	\begin{align*}
		\begin{array}{l}
			\displaystyle C \int_{\mathbb{R}} a_\xi^{-\mathbf{X}} |(v - \bar{v}, \, u - \bar{u}, \, \theta - \bar{\theta})|^3 d\xi \\[4mm]
			\displaystyle \le C \int_{\mathbb{R}} a_\xi^{-\mathbf{X}} \left[ \left| (v - \bar{v}) + \dfrac{1}{\sigma^*} (u - \bar{u}) \right|^3 + \left| (\theta - \bar{\theta}) - \dfrac{(\gamma - 1) \theta^*}{\sigma^* v^*} (u - \bar{u}) \right|^3 \right] d\xi + C \int_{\mathbb{R}} a_\xi^{-\mathbf{X}} |u - \bar{u}|^3 d\xi \\[4mm]
			\displaystyle \le C \varepsilon_1 (\mathbf{G}_1 + \mathbf{G}_2) + C \dfrac{\lambda}{\delta_S} \| u - \bar{u} \|_{L^\infty(\mathbb{R})}^2 \int_{\mathbb{R}} v_\xi^S  |u - \bar{u}| d\xi \\[4mm]
			\displaystyle \le C \varepsilon_1 (\mathbf{G}_1 + \mathbf{G}_2) + C \dfrac{\lambda}{\delta_S} \| u - \bar{u} \|_{L^2(\mathbb{R})} \cdot \| (u - \bar{u})_\xi \|_{L^2(\mathbb{R})} \left( \int_{\mathbb{R}} v_\xi^S (u - \bar{u})^2 d\xi \right)^{1/2} \left( \int_{\mathbb{R}} v_\xi^S d\xi \right)^{1/2} \\[4mm]
			\displaystyle \le C \varepsilon_1 (\mathbf{G}_1 + \mathbf{G}_2) + C \dfrac{\lambda}{\sqrt{\delta_S}} \varepsilon_1 \sqrt{\mathbf{D}} \sqrt{\mathcal{G}^S} \\[4mm]
			\displaystyle \le C \varepsilon_1 (\mathbf{G}_1 + \mathbf{G}_2 + \mathbf{D} + \mathcal{G}^S).
		\end{array}
	\end{align*}
	Finally, using Lemma \ref{WI-est} and \reff{lam-range},
	\begin{align*}
		\begin{array}{l}
			\displaystyle C \int_{\mathbb{R}} a_\xi^{-\mathbf{X}} \left( |v^R - v_*| + |\theta^R - \theta_*| + |v^C - v^*| + |\theta^C - \theta^*| \right) |(v - \bar{v}, \, u - \bar{u}, \, \theta - \bar{\theta})|^2 d\xi \\[4mm]
			\displaystyle \le C (\delta_R + \delta_C) (\mathbf{G}_1 + \mathbf{G}_2) + C \dfrac{\lambda}{\delta_S} \| u - \bar{u} \|_{L^\infty(\mathbb{R})} \cdot \| u - \bar{u} \|_{L^2(\mathbb{R})} \cdot \left( \| |v_\xi^S| \, |(v^R - v_*, \, \theta^R - \theta_*)| \|_{L^2(\mathbb{R})} \right. \\[4mm]
			\displaystyle \qquad \left. + \| |v_\xi^S| \cdot |(v^C - v^*, \, \theta^C - \theta^*)| \|_{L^2(\mathbb{R})} \right) \\[2mm]
			\displaystyle \le C (\delta_R + \delta_C) (\mathbf{G}_1 + \mathbf{G}_2) + C \dfrac{\lambda}{\delta_S} \| u - \bar{u} \|_{L^2(\mathbb{R})}^{3/2} \| (u - \bar{u})_\xi \|_{L^2(\mathbb{R})}^{1/2} \delta_S^{3/2} (\delta_R + \delta_C) e^{-C \delta_S t} \\[4mm]
			\displaystyle \le C (\delta_R + \delta_C) (\mathbf{G}_1 + \mathbf{G}_2) + C \varepsilon_1^{3/2} \mathbf{D}^{1/4} \delta_S (\delta_R + \delta_C) e^{-C \delta_S t} \\[4mm]
			\displaystyle \le C (\delta_R + \delta_C) (\mathbf{G}_1 + \mathbf{G}_2) + C \varepsilon_1 \mathbf{D} + C \delta_S^{4/3} (\delta_R + \delta_C)^{4/3} e^{-C \delta_S t}.
		\end{array}
	\end{align*}
	Therefore,
	\begin{align} \label{LO-B2G}
		\mathbf{B}_2 - \mathbf{G} \le - \frac{3}{4}(\mathbf{G}_1 + \mathbf{G}_2) + C (\lambda + \varepsilon_1 + \delta_0) ( \mathbf{D} + \mathcal{G}^S) + C \delta_S^{4/3} (\delta_R + \delta_C)^{4/3} e^{-C \delta_S t}.
	\end{align}
	\hspace*{\fill}\\

	\noindent $\bullet$ {\bf Step 2: Estimates on $-\dfrac{\delta_S}{2M} |\mathbf{\dot{X}}|^2$.}

	From Lemma \ref{REM-lem}, we have 
	\begin{align*}
		\begin{array}{l}
			\displaystyle \mathbf{Y}_1(U) = \int_{\mathbb{R}} a^{-\mathbf{X}} \left[ \dfrac{\hat{p}(\bar{v}, \bar{\theta}) v_\xi^S}{\bar{v}} (v - \bar{v}) + u_\xi^S (u - \bar{u}) + \dfrac{R \theta_\xi^S}{(\gamma - 1) \bar{\theta}} (\theta - \bar{\theta}) \right] d\xi \\[4mm]
			\displaystyle \qquad =: \mathbf{Y}_{11}(U) + \mathbf{Y}_{12}(U) + \mathbf{Y}_{13}(U).
		\end{array}
	\end{align*}
	For
	\begin{align*}
		\begin{array}{l}
			\displaystyle \mathbf{Y}_{11} = - \int_{\mathbb{R}} a^{-\mathbf{X}} \dfrac{\hat{p}(\bar{v}, \bar{\theta})}{\sigma^* \bar{v}} v_\xi^S (u - \bar{u}) d\xi + \int_{\mathbb{R}} a^{-\mathbf{X}} \dfrac{\hat{p}(\bar{v}, \bar{\theta})}{\bar{v}} v_\xi^S \left[ (v - \bar{v}) + \dfrac{1}{\sigma^*} (u - \bar{u}) \right] d\xi,
		\end{array}
	\end{align*}
	using $|a^{-\mathbf{X}} - 1| \le \lambda$ from \reff{a-range}, $\left| \dfrac{\hat{p}(\bar{v}, \bar{\theta})}{\bar{v}} - \dfrac{\hat{p}^*}{v^*} \right| \le C \delta_0$, \reff{wy-def} and \reff{par-y}, we have 
	\begin{align*}
		\begin{array}{l}
			\displaystyle \left| \mathbf{Y}_{11} + \dfrac{\hat{p}^* \delta_S}{\sigma^* v^*} \int_0^1 w dy \right| \le C \delta_S (\lambda + \delta_0) \int_0^1 |w| dy + C \dfrac{\delta_S}{\lambda} \int_{\mathbb{R}} a_\xi^{-\mathbf{X}} \left| (v - \bar{v}) + \dfrac{1}{\sigma^*} (u - \bar{u}) \right| d\xi.
		\end{array}
	\end{align*}
	Similarly, using \reff{eq:shock_a} for 
	\begin{align*}
		\mathbf{Y}_{12} = - \sigma \int_{\mathbb{R}} a^{-\mathbf{X}} v^S_\xi (u - \bar{u}) \, d\xi,
	\end{align*} 
	due to \reff{sig_sig-}, 
	\begin{align*}
		\left| \mathbf{Y}_{12} + \sigma^* \delta_S \int_0^1 w \, dy \right| \le C \delta_S (\lambda + \delta_S) \int_0^1 |w| \, dy.
	\end{align*}
	Finally,
	\begin{align*}
		\displaystyle \mathbf{Y}_{13} = \int_{\mathbb{R}} a^{-\mathbf{X}} \dfrac{R \theta_\xi^S}{(\gamma - 1) \bar{\theta}} \dfrac{(\gamma - 1) \theta^*}{\sigma^* v^*} (u - \bar{u}) \, d\xi + \int_{\mathbb{R}} a^{-\mathbf{X}} \dfrac{R \theta_\xi^S}{(\gamma - 1) \bar{\theta}} \left[ (\theta - \bar{\theta}) - \dfrac{(\gamma - 1) \theta^*}{\sigma^* v^*} (u - \bar{u}) \right] d\xi.
	\end{align*}
	Since \reff{exist_theta_v} in Lemma \ref{VS-pro} and $|\bar{\theta} - \theta^*| \le C \delta_0$, it yields 
	\begin{align*}
		\begin{array}{l}
			\displaystyle \left| \mathbf{Y}_{13} + \dfrac{(\gamma - 1) p^*}{\sigma^* v^*} \delta_S \int_0^1 w \, dy \right| \le C \delta_S (\lambda + \delta_0) \int_0^1 |w| \, dy \\[4mm]
			\displaystyle \qquad \qquad \qquad \qquad \qquad \qquad \qquad + C \dfrac{\delta_S}{\lambda} \int_{\mathbb{R}} a_\xi^{-\mathbf{X}} \left| (\theta - \bar{\theta}) - \dfrac{(\gamma - 1) \theta^*}{\sigma^* v^*} (u - \bar{u}) \right| \, d\xi.
		\end{array}	
	\end{align*}

	Therefore, using the above computation with \reff{left-p-sig}, more precisely, $\sigma^* = \sqrt{\dfrac{\gamma p^*}{v^*} + \dfrac{1}{(v^*)^2}},$  we obtain 
	\begin{align*}
		\begin{array}{l}
			\displaystyle \left| \mathbf{Y}_1 + 2 \sigma^* \delta_S \int_0^1 w \, dy \right| \le C \delta_S (\lambda + \delta_0) \int_0^1 |w| \, dy \\[4mm]
			\displaystyle \qquad + C \dfrac{\delta_S}{\lambda} \int_{\mathbb{R}} a_\xi^{-\mathbf{X}} \left| (v - \bar{v}) + \dfrac{1}{\sigma^*} (u - \bar{u}) \right|^2 d\xi + C \dfrac{\delta_S}{\lambda} \int_{\mathbb{R}} a_\xi^{-\mathbf{X}} \left| (\theta - \bar{\theta}) - \dfrac{(\gamma - 1) \theta^*}{\sigma^* v^*} (u - \bar{u}) \right|^2 d\xi,
		\end{array}
	\end{align*}
	which together with \reff{Xt-def} yields 
	\begin{align*}
		\begin{array}{l}
			\displaystyle \left| 2 M \sigma^* \int_0^1 w \, dy \right| - |\dot{\mathbf{X}}(t)| \le \left| \dot{\mathbf{X}}(t) - 2 M \sigma^* \int_0^1 w \, dy \right| \le C (\lambda + \delta_0) \int_0^1 |w| \, dy \\[4mm]
			\displaystyle \qquad + C \dfrac{1}{\lambda} \int_{\mathbb{R}} a_\xi^{-\mathbf{X}} \left( \left| (v - \bar{v}) + \dfrac{1}{\sigma^*} (u - \bar{u}) \right| + \left| (\theta - \bar{\theta}) - \dfrac{(\gamma - 1) \theta^*}{\sigma^* v^*} (u - \bar{u}) \right| \right) d\xi.
		\end{array}
	\end{align*}
	By Hölder's inequality,
	\begin{align*}
		\begin{array}{l}
			\displaystyle \left( \left| 2 M \sigma^* \int_0^1 w \, dy \right| - |\dot{\mathbf{X}}(t)| \right)^2 \le C (\lambda + \delta_0)^2 \int_0^1 w^2 \, dy + C \dfrac{1}{\lambda^2} (\mathbf{G}_1 + \mathbf{G}_2) \int_{\mathbb{R}} a_\xi^{-\mathbf{X}} d\xi \\[4mm]
			\displaystyle \qquad \qquad \qquad \qquad \qquad \qquad \quad \le C (\lambda + \delta_0)^2 \int_0^1 w^2 \, dy + C \dfrac{1}{\lambda} (\mathbf{G}_1 + \mathbf{G}_2).
		\end{array}
	\end{align*}
	This and the algebraic inequality $\dfrac{p^2}{2} - q^2 \le (p - q)^2$ if $p, q \ge 0$ imply 
	\begin{align} \label{LO-Xt}
		\displaystyle - \dfrac{\delta_S}{2M} |\dot{\mathbf{X}}(t)|^2 \le - M (\sigma^*)^2 \delta_S \left( \int_0^1 w \, dy \right)^2 + C \delta_S (\lambda + \delta_0)^2 \int_0^1 w^2 \, dy + C \dfrac{\delta_S}{\lambda} (\mathbf{G}_1 + \mathbf{G}_2).
	\end{align}
	\hspace*{\fill}\\

	\noindent $\bullet$ {\bf Step 3: Estimates on $\mathbf{B}_1$.}

	In what follows, we first will rewrite $\mathbf{B}_1$ into the main term of leading order and the other tractable error terms. First of all, taking the Taylor expansion for $\Psi(z) = z - 1 - \ln z$ near the point $z=1$, we have
	\begin{align*}
		\begin{array}{l}
			\Psi(z) = \Psi(1) + \Psi'(1) (z - 1) + \dfrac{\Psi''(1)}{2} (z - 1)^2 + O(|z - 1|^3) \\[4mm]
			\qquad = \dfrac{1}{2} (z - 1)^2 + O(|z - 1|^3),
		\end{array}
	\end{align*}
	which together with \reff{Ass-Inf} means 
	\begin{align} \label{st3-v}
		\Psi\left( \dfrac{v}{\bar{v}} \right) = \dfrac{(v - \bar{v})^2}{2 \bar{v}^2} + O(|v - \bar{v}|^3),
	\end{align}
	and
	\begin{align} \label{st3-theta}
		\Psi\left( \dfrac{\theta}{\bar{\theta}} \right) = \dfrac{(\theta - \bar{\theta})^2}{2 \bar{\theta}^2} + O(|\theta - \bar{\theta}|^3).
	\end{align}

	Let $J_1, J_2$ denote the first and second bracketed terms below, respectively. From Lemma \ref{REM-lem}, we have 
	\begin{align*} 
		\begin{array}{l}
			\displaystyle \mathbf{B}_1(U):= \int_{\mathbb{R}} a^{-\mathbf{X}} \Bigg\{ \Big[ R \left( (\theta_t^R - \sigma \theta_\xi^R) + (\theta_t^C - \sigma \theta_\xi^C) - \sigma \theta_\xi^S \right) \Psi\left( \dfrac{v}{\bar{v}} \right) - \dfrac{p(\bar{v}, \bar{\theta}) \bar{u}_\xi}{v \bar{v}} (v - \bar{v})^2\\[3mm] 
			\displaystyle \qquad\qquad\qquad\qquad -\left( (v^R_t - \sigma v^R_\xi ) + ( v^C_t - \sigma v^C_\xi ) - \sigma v^S_\xi \right)  \dfrac{(v - \bar{v})^2}{\bar{v}^3} \Big]  \\[4mm]
			\displaystyle \qquad\qquad\qquad\,\,  + \Big[ \dfrac{R}{\gamma - 1} \left( (\theta_t^R - \sigma \theta_\xi^R) + (\theta_t^C - \sigma \theta_\xi^C) - \sigma \theta_\xi^S \right) \Psi\left( \dfrac{\theta}{\bar{\theta}} \right) \\[4mm]
			\displaystyle \qquad\qquad\qquad\qquad - \dfrac{\bar{u}_\xi}{\theta} (\theta - \bar{\theta}) (p(v, \theta) - p(\bar{v}, \bar{\theta})) + \dfrac{p(\bar{v}, \bar{\theta}) \bar{u}_\xi}{\theta \bar{\theta}} (\theta - \bar{\theta})^2 \Big] \Bigg\} \, d\xi, \\[4mm]
			\displaystyle \qquad\quad\,\, =: \int_{\mathbb{R}} a^{-\mathbf{X}} (J_1+J_2) \, d\xi.
		\end{array}
	\end{align*}
	Next, we will handle the terms involving $J_1, J_2$. From $\reff{eq:ARV-change}_3$ and $\reff{eq:ARV-change}_1$, it follows 
	\begin{align*}
		\theta_t^R - \sigma \theta_\xi^R = - \dfrac{(\gamma - 1)}{R} p(v^R, \theta^R) u_\xi^R,
	\end{align*}
	and 
	\begin{align*}
		v_t^R - \sigma v_\xi^R = u_\xi^R,
	\end{align*}
	which together with $\bar{u}_\xi = u_\xi^R + u_\xi^C + u_\xi^S$ implies
	\begin{align*}
		\begin{array}{l}
			\displaystyle J_1 = - u_\xi^R \left[ (\gamma - 1) p(v^R, \theta^R) \Psi\left( \dfrac{v}{\bar{v}} \right) + \dfrac{p(\bar{v}, \bar{\theta})}{\bar{v} v} (v - \bar{v})^2 + \dfrac{(v - \bar{v})^2}{\bar{v}^3}  \right] \\[4mm]
			\displaystyle \qquad + \left[ R (\theta_t^C - \sigma \theta_\xi^C) \Psi\left( \dfrac{v}{\bar{v}} \right) - \dfrac{p(\bar{v}, \bar{\theta}) u_\xi^C}{\bar{v} v} (v - \bar{v})^2 - (v_t^C - \sigma v_\xi^C) \dfrac{(v - \bar{v})^2}{\bar{v}^3}  \right] \\[4mm]
			\displaystyle \qquad + \left[ - \sigma R \theta_\xi^S \, \Psi\left( \dfrac{v}{\bar{v}} \right) - \dfrac{p(\bar{v}, \bar{\theta}) u_\xi^S}{\bar{v} v} (v - \bar{v})^2 + \sigma v_\xi^S \dfrac{(v - \bar{v})^2}{\bar{v}^3} \right] \\[4mm]
			\displaystyle \qquad =: J_{11} + J_{12} + J_{13}.
		\end{array}
	\end{align*}
	By \reff{st3-v} and 
	\begin{align*}
		\begin{array}{l}
			\displaystyle |p(v^R, \theta^R) - p(\bar{v}, \bar{\theta})| \le C (|\bar{v} - v^R| + |\bar{\theta} - \theta^R|) \\[4mm]
			\displaystyle \qquad \le C (|(v^C - v_*, \, \theta^C - \theta_*)| + |(v^S - v^*, \, \theta^S - \theta^*)|) \\[4mm]
			\displaystyle \qquad \le C (\delta_C + \delta_S) \le C \delta_0,
		\end{array}
	\end{align*}
	we obtain
	\begin{align*}
		\displaystyle J_{11} \le - u_\xi^R \, \Big( \dfrac{(\gamma + 1) p(\bar{v}, \bar{\theta})}{2} + \dfrac{1}{\bar{v}} \Big) \dfrac{1}{\bar{v}^2} \, (v - \bar{v})^2 + C (\delta_0 + \varepsilon_1) |u_\xi^R| \cdot (v - \bar{v})^2.
	\end{align*}
	Since $\reff{eq:VCD-change}_1$, $\reff{eq:VCD-change}_3$ and Lemma \ref{VCD}, we have 
	\begin{align*}
		\displaystyle J_{12} \le C (|u_\xi^C| + |\theta_{\xi\xi}^C| + |\theta_\xi^C| \cdot |v_\xi^C| ) (v - \bar{v})^2 \le C \delta_C (1+t)^{-1} e^{-\frac{C_1 |\xi + \sigma t|^2}{1+t}} (v - \bar{v})^2.
	\end{align*}
	By \reff{st3-v}, \reff{exist_uv} and \reff{exist_theta_v} in Lemma \ref{VS-pro}, it follows 
	\begin{align*}
		\displaystyle J_{13} \le \Big( \dfrac{(\gamma + 1) p^*}{2} + \dfrac{1}{v^*} \Big) \dfrac{\sigma^*}{(v^*)^2} v_\xi^S (v - \bar{v})^2 + C (\delta_0 + \varepsilon_1) |v_\xi^S| (v - \bar{v})^2.
	\end{align*}
	Therefore, we have 
	\begin{align*}
		\begin{array}{l}
			\displaystyle J_1 \le  \underset{=:J_{1R}}{\underbrace{- u_\xi^R \, \Big( \dfrac{(\gamma + 1) p(\bar{v}, \bar{\theta})}{2} + \dfrac{1}{\bar{v}} \Big) \dfrac{1}{\bar{v}^2} \, (v - \bar{v})^2} }   + \Big( \dfrac{(\gamma + 1) p^*}{2} + \dfrac{1}{v^*} \Big) \dfrac{\sigma^*}{(v^*)^2} v_\xi^S (v - \bar{v})^2 \\[4mm]
			\displaystyle \qquad \qquad + C \delta_C (1+t)^{-1} e^{-\frac{C_1 |\xi + \sigma t|^2}{1+t}} (v - \bar{v})^2 + C (\delta_0 + \varepsilon_1) (|u_\xi^R| + |v_\xi^S|) \cdot (v - \bar{v})^2.
		\end{array}
	\end{align*}
	Likewise, 
	\begin{align*}
		\begin{array}{l}
			\displaystyle J_2 = \left[ \dfrac{R}{\gamma - 1} (\theta_t^R - \sigma \theta_\xi^R) \Psi\left( \dfrac{\theta}{\bar{\theta}} \right) - \dfrac{u_\xi^R}{\theta} (\theta - \bar{\theta}) (p(v, \theta) - p(\bar{v}, \bar{\theta})) + \dfrac{p(\bar{v}, \bar{\theta}) u_\xi^R}{\bar{\theta} \theta} (\theta - \bar{\theta})^2 \right] \\[4mm]
			\displaystyle \qquad + \left[ \dfrac{R}{\gamma - 1} (\theta_t^C - \sigma \theta_\xi^C) \Psi\left( \dfrac{\theta}{\bar{\theta}} \right) - \dfrac{u_\xi^C}{\theta} (\theta - \bar{\theta}) (p(v, \theta) - p(\bar{v}, \bar{\theta})) + \dfrac{p(\bar{v}, \bar{\theta}) u_\xi^C}{\bar{\theta} \theta} (\theta - \bar{\theta})^2 \right] \\[4mm]
			\displaystyle \qquad + \left[ -\dfrac{\sigma R}{\gamma - 1} \theta_\xi^S \Psi\left( \dfrac{\theta}{\bar{\theta}} \right) - \dfrac{u_\xi^S}{\theta} (\theta - \bar{\theta}) (p(v, \theta) - p(\bar{v}, \bar{\theta})) + \dfrac{p(\bar{v}, \bar{\theta}) u_\xi^S}{\bar{\theta} \theta} (\theta - \bar{\theta})^2 \right] \\[4mm]
			\displaystyle \le \underset{=:J_{2R}}{\underbrace{u_\xi^R \left[ - \dfrac{R}{2 \bar{v} \bar{\theta}} (\theta - \bar{\theta}) + \dfrac{p(\bar{v}, \bar{\theta})}{\bar{v} \bar{\theta}} (v - \bar{v}) \right] (\theta - \bar{\theta})}} + \left[ \dfrac{\sigma^* p^*}{2 (\theta^*)^2} v_\xi^S (\theta - \bar{\theta})^2 - \dfrac{\sigma^* p^*}{\theta^* v^*} v_\xi^S (v - \bar{v}) (\theta - \bar{\theta}) \right] \\[5mm]
			\displaystyle \qquad + C \delta_C (1+t)^{-1} e^{-\frac{C_1 |\xi + \sigma t|^2}{1+t}} |(v - \bar{v}, \theta - \bar{\theta})|^2  + C (\delta_0 + \varepsilon_1) (|u_\xi^R| + |v_\xi^S|) |(v - \bar{v}, \, \theta - \bar{\theta})|^2.
		\end{array}
	\end{align*}

	Using the Young inequality, it follows 
	\begin{align*}
		\begin{array}{l}
			\displaystyle (v - \bar{v}) (\theta - \bar{\theta}) \le \dfrac{\bar{\theta}}{\bar{v}} (v - \bar{v})^2 + \dfrac{\bar{v}}{4 \bar{\theta}} (\theta - \bar{\theta})^2,
		\end{array}
	\end{align*}
	and $u_\xi^R \sim v_\xi^R > 0$ by Lemma \ref{ARV-pro}, we have
	\begin{align*}
		\begin{array}{l}
			\displaystyle J_{1R} + J_{2R} \le - u_\xi^R \left[ \Big(\dfrac{(\gamma - 1) p(\bar{v}, \bar{\theta})}{2} + \dfrac{1}{\bar{v}} \Big) \dfrac{1}{\bar{v}^2} (v - \bar{v})^2 + \dfrac{R}{4 \bar{v} \bar{\theta}} (\theta - \bar{\theta})^2 \right] \\[4mm]
			\displaystyle \qquad \quad \qquad \le - C |v_\xi^R| \cdot |(v - \bar{v}, \, \theta - \bar{\theta})|^2 + C \delta_0 |v^R_\xi|\cdot |(v - \bar{v}, \, \theta - \bar{\theta})|^2.
		\end{array}
	\end{align*}
	Therefore, we have
	\begin{align} \label{LO-B1-middle}
		\begin{array}{l}
			\displaystyle \mathbf{B}_1(U) \le \mathbf{B}_1^S(U) - C \mathcal{G}^R(U) + C \delta_C (1+t)^{-1} \int_{\mathbb{R}} e^{-\frac{C_1 |\xi + \sigma t|^2}{1+t}} |(v - \bar{v}, \, \theta - \bar{\theta})|^2 d\xi \\[4mm]
			\displaystyle \qquad \qquad + C (\delta_0 + \varepsilon_1) \mathcal{G}^S(U),
		\end{array}
	\end{align}
	with $\mathcal{G}^R, \mathcal{G}^S$ defined in \reff{LO-G_SR} and 
	\begin{align*}
		\begin{array}{l}
			\displaystyle \mathbf{B}_1^S(U) := \int_{\mathbb{R}} a^{-\mathbf{X}} \left[ \left( \dfrac{(\gamma + 1) p^*}{2} + \dfrac{1}{v^*} \right) \dfrac{\sigma^*}{(v^*)^2} v_\xi^S (v - \bar{v})^2 + \dfrac{\sigma^* p^*}{2 (\theta^*)^2} v_\xi^S (\theta - \bar{\theta})^2 \right. \\[4mm]
			\displaystyle \qquad \qquad \qquad \qquad \,\, \left. - \dfrac{\sigma^* p^*}{\theta^* v^*} v_\xi^S (v - \bar{v}) (\theta - \bar{\theta}) \right] d\xi, \\[4mm]
			\displaystyle \qquad \quad\,\, =: \mathbf{B}_{11}^S(U) + \mathbf{B}_{12}^S(U) + \mathbf{B}_{13}^S(U).
		\end{array}
	\end{align*}
	
	Next, we will focus on estimating $\mathbf{B}_{1i}^S(U), i=1,2,3$.
	Using \reff{a-def}, \reff{par-a}, \reff{wy-def}, and the Young inequality, it holds that 
	\begin{align*}
		\begin{array}{l}
			\displaystyle \mathbf{B}_{11}^S = \left( \dfrac{(\gamma + 1) p^*}{2} + \dfrac{1}{v^*} \right) \dfrac{\sigma^*}{(v^*)^2} \int_{\mathbb{R}} a^{-\mathbf{X}} v_\xi^S \left| \left[ (v - \bar{v}) + \dfrac{1}{\sigma^*} (u - \bar{u}) \right] - \dfrac{1}{\sigma^*} (u - \bar{u}) \right|^2 d\xi \\[4mm]
			\displaystyle \qquad \le \left( \dfrac{(\gamma + 1) p^*}{2} + \dfrac{1}{v^*} \right) \dfrac{\sigma^*}{(v^*)^2} \left( 1 + \lambda + \left( \dfrac{\delta_S}{\lambda} \right)^{1/2} \right) \int_{\mathbb{R}} v_\xi^S \left| \dfrac{1}{\sigma^*} (u - \bar{u}) \right|^2 d\xi \\[4mm]
			\displaystyle \qquad \qquad + C \left( \dfrac{\delta_S}{\lambda} \right)^{-1/2} \int_{\mathbb{R}} v_\xi^S \left| (v - \bar{v}) + \dfrac{1}{\sigma^*} (u - \bar{u}) \right|^2 d\xi \\[4mm]
			\displaystyle \qquad \le \left( \dfrac{(\gamma + 1) p^*}{2} + \dfrac{1}{v^*} \right) \dfrac{1}{\sigma^* (v^*)^2} \left( 1 + \lambda + \left( \dfrac{\delta_S}{\lambda} \right)^{1/2} \right) \delta_S \int_0^1 w^2 dy + C \left( \dfrac{\delta_S}{\lambda} \right)^{1/2} \mathbf{G}_1.
		\end{array}
	\end{align*}
	Similarly,
	\begin{align*}
		\begin{array}{l}
			\displaystyle \mathbf{B}_{12}^S = \dfrac{\sigma^* p^*}{2 (\theta^*)^2} \int_{\mathbb{R}} a^{-\mathbf{X}} v_\xi^S \left|\left[ (\theta - \bar{\theta}) - \dfrac{(\gamma - 1) \theta^*}{\sigma^* v^*} (u - \bar{u}) \right]+ \dfrac{(\gamma - 1) \theta^*}{\sigma^* v^*} (u - \bar{u}) \right|^2 d\xi \\[4mm]
			\displaystyle \qquad \le \dfrac{\sigma^* p^*}{2 (\theta^*)^2} \left( 1 + \lambda + \left( \dfrac{\delta_S}{\lambda} \right)^{1/2} \right) \int_{\mathbb{R}} v_\xi^S \left| \dfrac{(\gamma - 1) \theta^*}{\sigma^* v^*} (u - \bar{u}) \right|^2 d\xi \\[4mm]
			\displaystyle \qquad \qquad + C \left( \dfrac{\delta_S}{\lambda} \right)^{-1/2} \int_{\mathbb{R}} v_\xi^S \left| (\theta - \bar{\theta}) - \dfrac{(\gamma - 1) \theta^*}{\sigma^* v^*} (u - \bar{u}) \right|^2 d\xi \\[4mm]
			\displaystyle \qquad \le \dfrac{(\gamma - 1)^2 p^*}{2 \sigma^* (v^*)^2} \left( 1 + \lambda + \left( \dfrac{\delta_S}{\lambda} \right)^{1/2} \right) \delta_S \int_0^1 w^2 dy + C \left( \dfrac{\delta_S}{\lambda} \right)^{1/2} \mathbf{G}_2.
		\end{array}
	\end{align*}
	Lastly,
	\begin{align*}
		\begin{array}{l}
			\displaystyle \mathbf{B}_{13}^S = - \dfrac{\sigma^* p^*}{\theta^* v^*} \int_{\mathbb{R}} a^{-\mathbf{X}} v_\xi^S \left[ (v - \bar{v}) + \dfrac{1}{\sigma^*} (u - \bar{u}) - \dfrac{1}{\sigma^*} (u - \bar{u}) \right] \\[4mm]
			\displaystyle \qquad \qquad \qquad \qquad \times \left[ (\theta - \bar{\theta}) - \dfrac{(\gamma - 1) \theta^*}{\sigma^* v^*} (u - \bar{u}) + \dfrac{(\gamma - 1) \theta^*}{\sigma^* v^*} (u - \bar{u}) \right] d\xi \\[4mm]
			\displaystyle \qquad \le \dfrac{(\gamma - 1) p^*}{\sigma^* (v^*)^2} \left( 1 + \lambda + \left( \dfrac{\delta_S}{\lambda} \right)^{1/2} \right) \int_{\mathbb{R}} v_\xi^S (u - \bar{u})^2 d\xi \\[4mm]
			\displaystyle \qquad \qquad + C \left( \dfrac{\delta_S}{\lambda} \right)^{-1/2} \int_{\mathbb{R}} v_\xi^S \left[ \left| (v - \bar{v}) + \dfrac{1}{\sigma^*} (u - \bar{u}) \right|^2 + \left| (\theta - \bar{\theta}) - \dfrac{(\gamma - 1) \theta^*}{\sigma^* v^*} (u - \bar{u}) \right|^2 \right] d\xi \\[4mm]
			\displaystyle \qquad \le \dfrac{(\gamma - 1) p^*}{\sigma^* (v^*)^2} \left( 1 + \lambda + \left( \dfrac{\delta_S}{\lambda} \right)^{1/2} \right) \delta_S \int_0^1 w^2 dy + C \left( \dfrac{\delta_S}{\lambda} \right)^{1/2} (\mathbf{G}_1 + \mathbf{G}_2).
		\end{array}
	\end{align*}
	Let $\alpha^*$ be the $O(1)$-constant defined by
	\begin{align*}
		\begin{array}{l}
			\displaystyle \alpha^* := \left( \dfrac{(\gamma + 1) p^*}{2} + \dfrac{1}{v^*} \right) \dfrac{1}{\sigma^* (v^*)^2} + \dfrac{(\gamma - 1)^2 p^*}{2 \sigma^* (v^*)^2} + \dfrac{(\gamma - 1) p^*}{\sigma^* (v^*)^2} \\[4mm]
			\displaystyle \qquad = \left( \dfrac{\gamma (\gamma + 1) p^*}{2} + \dfrac{1}{v^*} \right) \dfrac{1}{\sigma^* (v^*)^2},
		\end{array}
	\end{align*}
	Thus
	\begin{align*}
		\mathbf{B}_1^S \le \alpha^* \left( 1 + \lambda + \left( \dfrac{\delta_S}{\lambda} \right)^{1/2} \right) \delta_S \int_0^1 w^2 dy + C \left( \dfrac{\delta_S}{\lambda} \right)^{1/2} (\mathbf{G}_1 + \mathbf{G}_2),
	\end{align*}
	which together with \reff{LO-B1-middle} yields
	\begin{align} \label{LO-B1}
		\begin{array}{l}
			\displaystyle \mathbf{B}_1 \le \alpha^* \left( 1 + \lambda + \left( \dfrac{\delta_S}{\lambda} \right)^{1/2} \right) \delta_S \int_0^1 w^2 dy - C \mathcal{G}^R \\[4mm]
			\displaystyle \qquad + C \delta_C (1+t)^{-1} \int_{\mathbb{R}} e^{-\frac{C_1 |\xi + \sigma t|^2}{1+t}} |(v - \bar{v}, \, \theta - \bar{\theta})|^2 d\xi \\[3mm]
			\displaystyle \qquad + C (\delta_0 + \varepsilon_1) \mathcal{G}^S + C \left( \dfrac{\delta_S}{\lambda} \right)^{1/2} (\mathbf{G}_1 + \mathbf{G}_2).
		\end{array}
	\end{align}
	\hspace*{\fill}\\

	\noindent $\bullet$ {\bf Step 4: Estimates on $\mathbf{D}$.}
	
	By Lemma \reff{REM-lem},
	\begin{align*}
		\mathbf{D}(U) = \int_{\mathbb{R}} a^{-\mathbf{X}} \left( \dfrac{\mu}{v} |(u - \bar{u})_\xi|^2 + \dfrac{k}{v \theta} |(\theta - \bar{\theta})_\xi|^2 \right) d\xi =: \mathbf{D}_1(U) + \mathbf{D}_2(U).
	\end{align*}
	From \reff{a-def} and \reff{par-y}, it immediately follows 
	\begin{align} \label{LO-D1-middle}
		\begin{array}{l}
			\displaystyle \mathbf{D}_1 \ge \int_{\mathbb{R}} \dfrac{\mu}{v} |(u - \bar{u})_y|^2 \left| \dfrac{dy}{d\xi} \right|^2 d\xi = \int_0^1 \dfrac{\mu v_\xi^S}{\delta_S v} |w_y|^2 dy \\[4mm]
			\displaystyle \qquad \ge \dfrac{1 - C(\delta_0 + \varepsilon_1)}{\delta_S} \int_0^1 \dfrac{\mu v_\xi^S}{v^S} |w_y|^2 dy \\[4mm]
			\displaystyle \qquad = \dfrac{1 - C(\delta_0 + \varepsilon_1)}{\delta_S} \int_0^1 \left( \dfrac{\mu v_\xi^S}{v^S} - \dfrac{1}{\sigma} \Phi(v^S, \phi^S) + \dfrac{1}{\sigma} \Phi(v^S, \phi^S) \right) |w_y|^2 dy \\[4mm]
			\displaystyle \qquad \ge \dfrac{1 - C(\delta_0 + \varepsilon_1)}{\delta_S} \int_0^1 \left( \dfrac{\mu v_\xi^S}{v^S} - \dfrac{1}{\sigma} \Phi(v^S, \phi^S) \right) |w_y|^2 dy,
		\end{array}
	\end{align}
	where the last inequality is derived from the fact that
	\begin{align*}
		|\Phi(v^S, \phi^S)| \le C \delta_S |v_\xi^S|  \le C \delta_S \left| \dfrac{\mu v_\xi^S}{v^S} - \dfrac{1}{\sigma} \Phi(v^S, \phi^S) \right|.
	\end{align*}
	This means that $\Phi(v^S, \phi^S)$ consists of the higher-order terms than $v_\xi^S$.
	Using \reff{eq:shock_a}, \reff{eq:shock_b} and \reff{ele-add}, we have 
	\begin{align*}
		\begin{array}{l}
			- \sigma v_\xi^S - u_\xi^S = 0, \\
			- \sigma u_\xi^S + \hat{p}(v^S, \theta^S)_\xi = \mu \left( \dfrac{u_\xi^S}{v^S} \right)_\xi + \Phi(v^S, \phi^S)_\xi.
		\end{array}
	\end{align*}
	Integrating the space variable over $(-\infty, \xi]$ yields
	\begin{align*}
		\begin{array}{l}
			- \sigma (v^S - v^*) - (u^S - u^*) = 0, \\
			- \sigma (u^S - u^*) + \hat{p}(v^S, \theta^S) - \hat{p}^* = \mu \dfrac{u_\xi^S}{v^S} + \Phi(v^S, \phi^S) = - \sigma \dfrac{\mu v_\xi^S}{v^S} + \Phi(v^S, \phi^S).
		\end{array}
	\end{align*}
	Since  
	\begin{align*}
		\sigma = \sqrt{-\dfrac{\hat{p}_+ - \hat{p}^*}{v_+ - v^*}}
	\end{align*} 
	by \reff{sig-def}, \reff{left-p-sig} and $\hat{p}_+ : = \hat{p}(v_+, \theta_+)$, we have 
	\begin{align*}
		\begin{array}{l}
			\displaystyle \dfrac{\mu v_\xi^S}{v^S} - \dfrac{1}{\sigma} \Phi(v^S, \phi^S) = - \dfrac{1}{\sigma} \left[ \sigma^2 (v^S - v^*) + (\hat{p}(v^S, \theta^S) - \hat{p}^*) \right] \\[4mm]
			\displaystyle = - \dfrac{1}{\sigma} \left[ - \dfrac{\hat{p}_+ - \hat{p}^*}{v_+ - v^*} (v^S - v^*) + (\hat{p}(v^S, \theta^S) - \hat{p}^*) \right] \\[4mm]
			\displaystyle = - \dfrac{1}{\sigma (v_+ - v^*)} \left[ (\hat{p}^* - \hat{p}_+) (v^S - v^*) + (\hat{p}(v^S, \theta^S) - \hat{p}^*) (v_+ - v^*) \right] \\[4mm]
			\displaystyle = - \dfrac{1}{\sigma (v_+ - v^*)} \left[ (\hat{p}^* - \hat{p}(v^S, \theta^S)) (v^S - v^*) + (\hat{p}(v^S, \theta^S) - \hat{p}_+) (v^S - v^*) \right. \\[4mm]
			\displaystyle \qquad\qquad\qquad\qquad \left. + (\hat{p}(v^S, \theta^S) - \hat{p}^*) (v^S - v^*) + (\hat{p}(v^S, \theta^S) - \hat{p}^*) (v_+ - v^S) \right] \\[4mm]
			\displaystyle = - \dfrac{1}{\sigma (v_+ - v^*)} \left[ (\hat{p}(v^S, \theta^S) - \hat{p}_+) (v^S - v^*) + (\hat{p}(v^S, \theta^S) - \hat{p}^*) (v_+ - v^S) \right] \\[4mm]
			\displaystyle = \dfrac{(v^S - v^*) (v_+ - v^S)}{\sigma (v_+ - v^*)} \left[ \dfrac{\hat{p}(v^S, \theta^S) - \hat{p}_+}{v^S - v_+} - \dfrac{\hat{p}(v^S, \theta^S) - \hat{p}^*}{v^S - v^*} \right] \\[4mm]
			\displaystyle = \dfrac{y (1 - y) \delta_S}{\sigma} \left[ \dfrac{\hat{p}(v^S, \theta^S) - \hat{p}_+}{v^S - v_+} - \dfrac{\hat{p}(v^S, \theta^S) - \hat{p}^*}{v^S - v^*} \right],
		\end{array}
	\end{align*}
	where the last equality is due to $y = \dfrac{v^S - v^*}{\delta_S}$, $1 - y = \dfrac{v_+ - v^S}{\delta_S}$ and $\delta_S := v_+ - v^* > 0$. By Lemma \ref{Crucial}, we have
	\begin{align*}
		\left| \dfrac{1}{y(1-y)} \left( \dfrac{\mu v_\xi^S}{v^S} - \dfrac{1}{\sigma} \Phi(v^S, \phi^S) \right) - \alpha^* \beta^* \delta_S^2 \right| \le C \delta_S^3,
	\end{align*}
	which together with \reff{LO-D1-middle} implies 
	\begin{align} \label{LO-D1}
		\mathbf{D}_1 \ge (1 - C(\delta_0 + \varepsilon_1)) \alpha^* \beta^* \delta_S \int_0^1 y(1-y) |w_y|^2 dy. 
	\end{align}
	Similarly,
	\begin{align} \label{LO-D2}
		\mathbf{D}_2 \ge (1 - C(\delta_0 + \varepsilon_1)) \dfrac{\kappa}{\mu \theta^*} \alpha^* \beta^* \delta_S \int_0^1 y(1-y) |(\theta - \bar{\theta})_y|^2 dy. 
	\end{align}
	
	Our goal is to employ the Poincaré-type Lemma \reff{Poincare} to absorb the bad term $\mathbf{B}_1$ by the diffusion term $\mathbf{D}$ and $-\dfrac{\delta_S}{2M} |\dot{\mathbf{X}}(t)|^2$. However, it follows from \reff{LO-B1} and \reff{LO-D1} that $\mathbf{D}_1$ can not fully control the bad term $\mathbf{B}_1$, since ($\beta^*$ defined in Lemma \ref{Crucial})
	\begin{align*}
		1 > \beta^* = \dfrac{R \mu (\gamma p^* + \dfrac{1}{v^*})}{k(\gamma - 1)^2 p^* + R \mu (\gamma p^* + \dfrac{1}{v^*})} \to 0 \quad \text{as } k \to +\infty.
	\end{align*}
	Hence, we need to extract the good terms from the diffusion term $\mathbf{D}_2$ to couple with this difficulty. First, using the Poincaré-type Lemma \ref{Poincare} and 
	\begin{align*}
		\displaystyle \int_{0}^{1} \left|f-\int_{0}^{1} f \, dy \right|^2 dy =\int_{0}^{1} f^2 \, dy - \left(\int_{0}^{1} f \, dy \right)^2,
	\end{align*}
	it follows from \reff{LO-D1} and \reff{LO-D2} that 
	\begin{align} \label{LO-D-middle}
		\begin{array}{l}
			\displaystyle \mathbf{D} \ge 2 (1 - C(\delta_0 + \varepsilon_1)) \alpha^* \beta^* \delta_S \left[ \int_0^1 w^2 dy - \left( \int_0^1 w dy \right)^2 \right. \\[4mm]
			\displaystyle \qquad \left. + \dfrac{\kappa}{\mu \theta^*} \left( \int_0^1 (\theta - \bar{\theta})^2 dy - \left( \int_0^1 (\theta - \bar{\theta}) dy \right)^2 \right) \right].
		\end{array}
	\end{align}
	Using the Young inequality, 
	\begin{align*}
		\begin{array}{l}
			\displaystyle \int_0^1 (\theta - \bar{\theta})^2 dy = \int_0^1 \left| \left[ (\theta - \bar{\theta}) - \dfrac{(\gamma - 1) \theta^*}{\sigma^* v^*} (u - \bar{u}) \right] + \dfrac{(\gamma - 1) \theta^*}{\sigma^* v^*} (u - \bar{u}) \right|^2 dy \\[4mm]
			\displaystyle \ge \left( \dfrac{(\gamma - 1) \theta^*}{\sigma^* v^*} \right)^2 \left( 1 - \left( \dfrac{\delta_S}{\lambda} \right)^{1/2} \right) \int_0^1 (u - \bar{u})^2 dy \\[4mm]
			\displaystyle \qquad \qquad \qquad - C \left( \dfrac{\delta_S}{\lambda} \right)^{-1/2} \int_0^1 \left| (\theta - \bar{\theta}) - \dfrac{(\gamma - 1) \theta^*}{\sigma^* v^*} (u - \bar{u}) \right|^2 dy \\[4mm]
			\displaystyle = \left( \dfrac{(\gamma - 1) \theta^*}{\sigma^* v^*} \right)^2 \left( 1 - \left( \dfrac{\delta_S}{\lambda} \right)^{1/2} \right) \int_0^1 w^2 dy - \dfrac{C}{\delta_S} \left( \dfrac{\delta_S}{\lambda} \right)^{1/2} \mathbf{G}_2,
		\end{array}
	\end{align*}
	and by the H\"older inequality,
	\begin{align*}
		\begin{array}{l}
			\displaystyle \left( \int_0^1 (\theta - \bar{\theta}) dy \right)^2 = \left( \int_0^1 \left\{ \left[ (\theta - \bar{\theta}) - \dfrac{(\gamma - 1) \theta^*}{\sigma^* v^*} (u - \bar{u}) \right] + \dfrac{(\gamma - 1) \theta^*}{\sigma^* v^*} (u - \bar{u}) \right\} dy \right)^2 \\[4mm]
			\displaystyle \le 2 \left( \dfrac{(\gamma - 1) \theta^*}{\sigma^* v^*} \right)^2 \left( \int_0^1 (u - \bar{u}) dy \right)^2 + 2 \left( \int_0^1 \left[ (\theta - \bar{\theta}) - \dfrac{(\gamma - 1) \theta^*}{\sigma^* v^*} (u - \bar{u}) \right] dy \right)^2 \\[4mm]
			\displaystyle = 2 \left( \dfrac{(\gamma - 1) \theta^*}{\sigma^* v^*} \right)^2 \left( \int_0^1 w dy \right)^2 + \dfrac{C}{\lambda} \mathbf{G}_2.
		\end{array}
	\end{align*}
	Since 
	\begin{align*}
		1 + \dfrac{\kappa}{\mu \theta^*} \dfrac{(\gamma - 1)^2 (\theta^*)^2}{(\sigma^*)^2 (v^*)^2} = \dfrac{1}{\beta^*}
	\end{align*}
	with $\beta^*$ defined in Lemma \ref{Crucial} and
	\begin{align*}
		1 + \dfrac{\kappa}{\mu \theta^*} \dfrac{2 (\gamma - 1)^2 (\theta^*)^2}{(\sigma^*)^2 (v^*)^2} = 1 + \dfrac{2 k (\gamma - 1)^2 p^*}{R \mu (\gamma p^* + \frac{1}{v^*})}
	\end{align*}
	with $\sigma^* = \sqrt{\dfrac{\gamma p^*}{v^*} + \dfrac{1}{(v^*)^2}}$ by \reff{left-p-sig}, it follows from \reff{LO-D-middle} that 
	\begin{align} \label{LO-D}
		\begin{array}{l}
			\displaystyle \mathbf{D} \ge 2 \left( 1 - C \left( \delta_0 + \varepsilon_1 + \left( \dfrac{\delta_S}{\lambda} \right)^{1/2} \right) \right) \alpha^* \delta_S \int_0^1 w^2 dy \\[4mm]
			\displaystyle \qquad - 2 \alpha^* \beta^* \left( 1 + \dfrac{2 \kappa (\gamma - 1)^2 p^*}{R \mu (\gamma p^* + \frac{1}{v^*})} \right) \delta_S \left( \int_0^1 w dy \right)^2 - C \left( \dfrac{\delta_S}{\lambda} \right)^{1/2} \mathbf{G}_2.
		\end{array}
	\end{align}
	\hspace*{\fill}\\

	\noindent $\bullet$ {\bf Step 5: Conclusion.}

	Using \reff{LO-B2G}, \reff{LO-Xt}, \reff{LO-B1} and \reff{LO-D}, it follows that 
	\begin{align*}
		\begin{array}{l}
			\displaystyle - \dfrac{\delta_S}{2M} |\dot{\mathbf{X}}(t)|^2 + \mathbf{B}_1 + \mathbf{B}_2 - \mathbf{G} - \dfrac{3}{4} \mathbf{D} \\[4mm]
			\displaystyle \le - \dfrac{1}{4} \alpha^* \delta_S \int_0^1 w^2 dy - \dfrac{1}{2} (\mathbf{G}_1 + \mathbf{G}_2) - C \mathcal{G}^R \\[4mm]
			\displaystyle \qquad + \left[ \dfrac{3}{2} \alpha^* \beta^* \left( 1 + \dfrac{2 k (\gamma - 1)^2 p^*}{R \mu (\gamma p^* + \frac{1}{v^*})} \right) - M (\sigma^*)^2 \right] \delta_S \left( \int_0^1 w dy \right)^2 \\[4mm]
			\displaystyle \qquad + C (\lambda + \delta_0 + \varepsilon_1) \mathcal{G}^S + C \delta_C (1+t)^{-1} \int_{\mathbb{R}} e^{-\frac{C_1 |\xi + \sigma t|^2}{1+t}} |(v - \bar{v}, \, \theta - \bar{\theta})|^2 d\xi \\[4mm]
			\displaystyle \qquad + C \delta_S^{4/3} (\delta_R + \delta_C)^{4/3} e^{-C \delta_S t}.
		\end{array}
	\end{align*}
	Take the specific $O(1)$-constant $M$ as below \reff{Xt-def}, that is,
	\begin{align} \label{M-def}
		M := \dfrac{3 \alpha^* \beta^*}{2 (\sigma^*)^2} \left( 1 + \dfrac{2 k (\gamma - 1)^2 p^*}{R \mu (\gamma p^* + \frac{1}{v^*})} \right).
	\end{align}
	Note that
	\begin{align*}
		\delta_S \int_0^1 w^2 dy = \int_{\mathbb{R}} v_\xi^S (u - \bar{u})^2 d\xi,
	\end{align*}
	\begin{align*}
		\begin{array}{l}
			\displaystyle \int_{\mathbb{R}} v_\xi^S (v - \bar{v})^2 d\xi = \int_{\mathbb{R}} v_\xi^S \left[ (v - \bar{v}) + \dfrac{1}{\sigma^*} (u - \bar{u}) - \dfrac{1}{\sigma^*} (u - \bar{u}) \right]^2 d\xi \\[4mm]
			\displaystyle \qquad \le C \int_{\mathbb{R}} v_\xi^S (u - \bar{u})^2 d\xi + C \int_{\mathbb{R}} v_\xi^S \left| (v - \bar{v}) + \dfrac{1}{\sigma^*} (u - \bar{u}) \right|^2 d\xi \\[4mm]
			\displaystyle \qquad \le C \int_{\mathbb{R}} v_\xi^S (u - \bar{u})^2 d\xi + C \dfrac{\delta_S}{\lambda} \mathbf{G}_1,
		\end{array}
	\end{align*}
	and
	\begin{align*}
		\begin{array}{l}
			\displaystyle \int_{\mathbb{R}} v_\xi^S (\theta - \bar{\theta})^2 d\xi = \int_{\mathbb{R}} v_\xi^S \left| (\theta - \bar{\theta}) - \dfrac{(\gamma - 1) \theta^*}{\sigma^* v^*} (u - \bar{u}) + \dfrac{(\gamma - 1) \theta^*}{\sigma^* v^*} (u - \bar{u}) \right|^2 d\xi \\[4mm]
			\displaystyle \qquad \le C \int_{\mathbb{R}} v_\xi^S (u - \bar{u})^2 d\xi + C \dfrac{\delta_S}{\lambda} \mathbf{G}_2.
		\end{array}
	\end{align*}
	Therefore, we have
	\begin{align*}
		\begin{array}{l}
			\displaystyle - \dfrac{\delta_S}{2M} |\dot{\mathbf{X}}(t)|^2 + \mathbf{B}_1 + \mathbf{B}_2 - \mathbf{G} - \dfrac{3}{4} \mathbf{D} \\[4mm]
			\displaystyle \le - C \mathcal{G}^R - C \mathcal{G}^S + C \delta_S^{4/3} (\delta_R + \delta_C)^{4/3} e^{-C \delta_S t} \\[4mm]
			\displaystyle \qquad + C \delta_C (1+t)^{-1} \int_{\mathbb{R}} e^{-\frac{C_1 |\xi + \sigma t|^2}{1+t}} |(v - \bar{v}, \, \theta - \bar{\theta})|^2 d\xi.
		\end{array}
	\end{align*}
\end{proof}
\hspace*{\fill}

\subsection{Estimates on \texorpdfstring{$\mathbf{Y}_i$ $(i=2,3,4)$, $\mathbf{E}_1$ and $\mathbf{E}_2$}{Y\_i (i=2,3,4), E\_1 and E\_2}}
\hspace*{\fill}\\ 

Substituting the estimates in Lemma \ref{LO-lem} into the inequality \reff{RE-mid1}, we have 
\begin{align} \label{RE-mid2}
	\begin{array}{l}
		\displaystyle \dfrac{d}{dt} \int_{\mathbb{R}} a^{-\mathbf{X}}(\xi) \tilde{\eta}((U, \phi) | (\bar{U}, \bar{\phi})) \, d\xi \le - C (\mathcal{G}^S + \mathcal{G}^R) + C \delta_S^{4/3} (\delta_R + \delta_C)^{4/3} e^{-C \delta_S t} \\[4mm]
		\displaystyle \qquad + C \delta_C (1+t)^{-1} \int_{\mathbb{R}} e^{-\frac{C_1 |\xi + \sigma t|^2}{1+t}} |(v - \bar{v}, \, \theta - \bar{\theta})|^2 d\xi \\[4mm]
		\displaystyle \qquad - \dfrac{\delta_S}{2M} |\dot{\mathbf{X}}(t)|^2 + \dot{\mathbf{X}}(t) \sum_{i=2}^4 \mathbf{Y}_i - \dfrac{1}{4} \mathbf{D} + \mathbf{E}_1 + \mathbf{E}_2.
	\end{array}
\end{align}
In what follows, we will estimate the rest error terms $Y_i (i=2,3,4)$, $\mathbf{E}_1$ and $\mathbf{E}_2$.
\\

\noindent $\bullet$ {\bf Estimates on $\mathbf{Y}_i (i=2,3,4)$.}

\begin{align*} 
	\dot{\mathbf{X}}(t) \sum_{i=2}^4 \mathbf{Y}_i \le \dfrac{\delta_S}{8M} |\dot{\mathbf{X}}(t)|^2 + \dfrac{C}{\delta_S} \sum_{i=2}^4 |\mathbf{Y}_i|^2.
\end{align*}

For $\mathbf{Y}_2$ defined in Lemma \ref{REM-lem}, due to $|\theta_\xi^S| \sim |v_\xi^S|$ in Lemma \ref{VS-pro},
\begin{align*} 
	\displaystyle \dfrac{|\mathbf{Y}_2|^2}{\delta_S} \le \dfrac{C}{\delta_S} \| v_\xi^S \|_{L^\infty(\mathbb{R})} \cdot \| (v - \bar{v}, \, \theta - \bar{\theta}) \|_{L^2(\mathbb{R})}^2 \mathcal{G}^S  \le C \varepsilon_1 \mathcal{G}^S.
\end{align*}
For
\begin{align*}
	\begin{array}{l}
		\displaystyle \mathbf{Y}_3 = \int_{\mathbb{R}} a^{-\mathbf{X}}_\xi \left[ \bar{\theta} \eta(U|\bar{U}) + \dfrac{(v - \bar{v})^2}{2 \bar{v}^2} - \dfrac{(v - \bar{v})^3}{3 \bar{v}^3} \right] d\xi \\[4mm]
		\displaystyle \qquad \qquad - \int_{\mathbb{R}} a^{-\mathbf{X}}_\xi \left[ - \dfrac{\nu^2}{\bar{v}^2} (v - \bar{v}) \tilde{\phi}_{\xi\xi} + \dfrac{\nu^4}{2 \bar{v}^3 e^{\bar{\phi}}} \tilde{\phi}_{\xi\xi}^2 + \dfrac{\nu^2}{2 \bar{v}} \tilde{\phi}_\xi^2 \right] d\xi \\[4mm]
		\displaystyle \qquad =: \mathbf{Y}_{31} + \mathbf{Y}_{32}.
	\end{array}
\end{align*}
By $\| v_\xi^S \|_{L^\infty(\mathbb{R})} \le C \delta_S^2$ in Lemma \ref{VS-pro} and \reff{lam-range},
\begin{align*}
	\begin{array}{l}
		\displaystyle \dfrac{|\mathbf{Y}_{31}|^2}{\delta_S} \le \dfrac{C \lambda^2}{\delta_S^3} \left( \int_{\mathbb{R}} |v_\xi^S| \cdot |(v - \bar{v}, \, \theta - \bar{\theta})|^2 d\xi \right)^2 \\[4mm]
		\displaystyle \qquad \quad \le \dfrac{C \lambda^2}{\delta_S^3} \| v_\xi^S \|_{L^\infty(\mathbb{R})} \| (v - \bar{v}, \, \theta - \bar{\theta}) \|_{L^2(\mathbb{R})}^2 \cdot \mathcal{G}^S \\[4mm]
		\displaystyle \qquad \quad \le \dfrac{C \lambda^2 \varepsilon_1^2}{\delta_S} \mathcal{G}^S \le C \varepsilon_1 \mathcal{G}^S,
	\end{array}
\end{align*}
and using 
\begin{align} \label{phi-H2}
	\| \tilde{\phi} \|_{H^2(\mathbb{R})} \le C (\| v-\bar{v} \|_{H^1(\mathbb{R})} + \delta_0 ) \le C(\delta_0 + \varepsilon_1)
\end{align} 
from Lemma \ref{Basic estimates},
\begin{align*}
	\begin{array}{l}
		\displaystyle \dfrac{|\mathbf{Y}_{32}|^2}{\delta_S} \le \dfrac{C \lambda^2}{\delta_S^3} \left( \int_{\mathbb{R}} |v_\xi^S| \big(|v - \bar{v}|^2 + |\tilde{\phi}_{\xi\xi}|^2 + |\tilde{\phi}_\xi|^2\big) d\xi \right)^2 \\[4mm]
		\displaystyle \qquad\quad \le \dfrac{C \lambda^2}{\delta_S^3} \| v_\xi^S \|_{L^\infty(\mathbb{R})}^2 \| v - \bar{v} \|_{L^2(\mathbb{R})}^2 \mathcal{G}^S + \dfrac{C \lambda^2}{\delta_S^3} \| v_\xi^S \|_{L^\infty(\mathbb{R})}^2 \| \tilde{\phi}_\xi \|_{H^1(\mathbb{R})}^4 \\[4mm]
		\displaystyle \qquad\quad \le \dfrac{C \lambda^2}{\delta_S^3} \delta_S^4 \cdot \varepsilon_1^2 \mathcal{G}^S + \dfrac{C \lambda^2}{\delta_S^3} \delta_S^4 \cdot (\delta_0 + \varepsilon_1)^2 \cdot \| \tilde{\phi}_\xi \|_{H^1(\mathbb{R})}^2 \\[4mm]
		\displaystyle \qquad\quad \le C (\delta_0 + \varepsilon_1) (\mathcal{G}^S + \| \tilde{\phi}_\xi \|_{H^1(\mathbb{R})}^2),
	\end{array}
\end{align*}
Thus
\begin{align*} 
	\dfrac{|\mathbf{Y}_3|^2}{\delta_S} \le C (\varepsilon_1 + \delta_0) (\mathcal{G}^S + \| \tilde{\phi}_\xi \|_{H^1(\mathbb{R})}^2).
\end{align*}
It follows from $|\phi_\xi^S| \sim |v_\xi^S|$ in Lemma \ref{VS-pro}, \reff{Ass-Inf} and \reff{phi-H2} that
\begin{align*}
	\begin{array}{l}
		\displaystyle \dfrac{|\mathbf{Y}_4|^2}{\delta_S} \le \dfrac{C}{\delta_S} \left( \int_{\mathbb{R}} \left[ |v_\xi^S| \cdot |v - \bar{v}| \cdot |\tilde{\phi}_{\xi\xi}| + |v_\xi^S| \cdot |\tilde{\phi}_{\xi\xi}| + |v_\xi^S| \cdot |\tilde{\phi}_{\xi\xi}|^2 + |v_\xi^S| \cdot |\tilde{\phi}_\xi|^2 \right] d\xi \right)^2 \\[4mm]
		\displaystyle \qquad\quad \le \dfrac{C}{\delta_S} \left( \| v_\xi^S \|_{L^\infty(\mathbb{R})} \int_{\mathbb{R}} v_\xi^S d\xi \cdot \| \tilde{\phi}_{\xi\xi} \|_{L^2(\mathbb{R})}^2 + \| v_\xi^S \|_{L^\infty(\mathbb{R})}^2 \| \tilde{\phi}_\xi \|_{H^1(\mathbb{R})}^4 \right) \\[4mm]
		\displaystyle \qquad\quad \le C \delta_S \| \tilde{\phi}_\xi \|_{H^1(\mathbb{R})}^2.
	\end{array}
\end{align*}

Therefore, we have 
\begin{align}  \label{Y-Est}
	\dot{\mathbf{X}}(t) \sum_{i=2}^4 \mathbf{Y}_i \le \dfrac{\delta_S}{8M} |\dot{\mathbf{X}}(t)|^2 + C (\delta_0 + \varepsilon_1) \left( \mathcal{G}^S + \| \tilde{\phi}_\xi \|_{H^1(\mathbb{R})}^2 \right).
\end{align}
\hspace*{\fill}\\

\noindent $\bullet$ {\bf Estimates on $\mathbf{E}_1$.}

From Lemma \ref{REM-lem}, we have 
\begin{align*}
	\begin{array}{l}
		\displaystyle \mathbf{E}_1(U):= \int_{\mathbb{R}} a^{-\mathbf{X}} \left\{ \left[ \left( v_t^R + v_t^C - \sigma \bar{v}_\xi \right) \dfrac{(v - \bar{v})^3}{\bar{v}^4} - \dfrac{(v - \bar{v})^3}{\bar{v}^3 v} (u - \bar{u})_\xi \right] \right. \\[3mm]
		\displaystyle \qquad + \left[ \kappa \dfrac{\theta - \bar{\theta}}{\theta} \Big( \dfrac{\theta_\xi}{v} - \dfrac{\bar{\theta}_\xi}{\bar{v}} \Big) + \mu (u - \bar{u}) \left( \dfrac{u_\xi}{v} - \dfrac{\bar{u}_\xi}{\bar{v}} \right) \right]_\xi + \Bigg[ -\kappa \dfrac{(\theta - \bar{\theta})_\xi}{\theta}  \Big( \dfrac{1}{v} - \dfrac{1}{\bar{v}} \Big) \bar{\theta}_\xi \\[3mm]
		\displaystyle \qquad \qquad + \kappa \dfrac{(\theta - \bar{\theta}) \theta_\xi}{\theta^2} \Big( \dfrac{\theta_\xi}{v} - \dfrac{\bar{\theta}_\xi}{\bar{v}} \Big) + \mu \dfrac{\theta - \bar{\theta}}{\theta} \Big( \dfrac{u_\xi^2}{v} - \dfrac{\bar{u}_\xi^2}{\bar{v}} \Big) - \dfrac{(\theta - \bar{\theta})^2}{\theta \bar{\theta}} \Big( \kappa \big( \dfrac{\bar{\theta}_\xi}{\bar{v}} \big)_\xi + \mu \dfrac{\bar{u}_\xi^2}{\bar{v}} \Big) \Bigg] \\[3mm]
		\displaystyle \qquad \left. - \left[ (u - \bar{u}) Q_1 + \dfrac{\theta - \bar{\theta}}{\theta} Q_2 \right] \right\} d\xi  \\[3mm] 
		\displaystyle \qquad\quad =: \sum_{i=1}^4 \mathbf{E}_{1i}(U). 
	\end{array}
\end{align*}
Using \reff{Ass-Inf}, \reff{ProAss} and the Sobolev inequality, it yields 
\begin{align*}
	\begin{array}{l}
		\displaystyle |\mathbf{E}_{11}| \le C \| v - \bar{v} \|_{L^\infty(\mathbb{R})} (\mathcal{G}^R + \mathcal{G}^S) + C \| v - \bar{v} \|_{L^\infty(\mathbb{R})} (1+t)^{-1} \int_{\mathbb{R}} e^{-\frac{C_1 |\xi + \sigma t|^2}{1+t}} (v - \bar{v})^2 d\xi \\[4mm]
		\displaystyle \qquad \qquad + C \| v - \bar{v} \|_{L^\infty(\mathbb{R})}^2 \cdot \| v - \bar{v} \|_{L^2(\mathbb{R})} \cdot \| (u - \bar{u})_\xi \|_{L^2(\mathbb{R})} \\[4mm]
		\displaystyle \qquad \le C \varepsilon_1 \left( \mathcal{G}^R + \mathcal{G}^S + (1+t)^{-1} \int_{\mathbb{R}} e^{-\frac{C_1 |\xi + \sigma t|^2}{1+t}} (v - \bar{v})^2 d\xi \right) \\[4mm]
		\displaystyle \qquad \qquad + C \| v - \bar{v} \|_{L^2}^2 \cdot \| (v - \bar{v})_\xi \|_{L^2(\mathbb{R})} \| (u - \bar{u})_\xi \|_{L^2(\mathbb{R})} \\[4mm]
		\displaystyle \qquad \le C \varepsilon_1 \left( \mathcal{G}^R + \mathcal{G}^S + \mathbf{D} + (1+t)^{-1} \int_{\mathbb{R}} e^{-\frac{C_1 |\xi + \sigma t|^2}{1+t}} (v - \bar{v})^2 d\xi + \| ( v-\bar{v} )_\xi \|_{L^2(\mathbb{R})}^2 \right).
	\end{array}
\end{align*}
By Lemma \ref{ARV-pro}, \ref{VCD} and \ref{VS-pro}, we have 
\begin{align*}
	|\bar{\theta}_\xi| + |\bar{u}_\xi| \le C (|v_\xi^R| + |v_\xi^C| + |v_\xi^S|),
\end{align*}
which together with \reff{par-a} yields 
\begin{align*}
	\begin{array}{l}
		\displaystyle |\mathbf{E}_{12}| \le C \dfrac{\lambda}{\delta_S} \int_{\mathbb{R}} |v_\xi^S| \cdot \Big( |\theta - \bar{\theta}| \cdot |(\theta - \bar{\theta})_\xi| + (|v_\xi^R| + |v_\xi^C| + |v_\xi^S|) |\theta - \bar{\theta}| \cdot |v - \bar{v}| \\[4mm]
		\displaystyle \qquad \qquad \qquad \qquad \qquad + |u - \bar{u}| \cdot |(u - \bar{u})_\xi| + (|v_\xi^R| + |v_\xi^C| + |v_\xi^S|) |u - \bar{u}| \cdot |v - \bar{v}| \Big) d\xi \\[4mm]
		\displaystyle \qquad \le C \dfrac{\lambda}{\delta_S} \left( \| v_\xi^S \|_{L^\infty(\mathbb{R})}^{1/2} \sqrt{\mathcal{G}^S} \sqrt{\mathbf{D}} + \| v_\xi^S \|_{L^\infty(\mathbb{R})} \big( \mathcal{G}^R + \mathcal{G}^S \right. \\[4mm]
		\displaystyle \qquad \qquad \quad + (1+t)^{-1} \int_{\mathbb{R}} e^{-\frac{C_1 |\xi + \sigma t|^2}{1+t}} |(v - \bar{v}, \, \theta - \bar{\theta})|^2 d\xi \big) + \| v_\xi^S \|_{L^\infty(\mathbb{R})}^{1/2} \| v_\xi^R \|_{L^\infty(\mathbb{R})}^{1/2} \sqrt{\mathcal{G}^S} \sqrt{\mathcal{G}^R} \\[4mm]
		\displaystyle \qquad \qquad \left. + \| v_\xi^S \|_{L^\infty(\mathbb{R})}^{1/2} \sqrt{\mathcal{G}^S} (1+t)^{-1} \int_{\mathbb{R}} e^{-\frac{C_1 |\xi + \sigma t|^2}{1+t}} |v - \bar{v}|^2 d\xi \right) \\[4mm]
		\displaystyle \qquad \le C \lambda \left( \mathcal{G}^S + \mathcal{G}^R + \mathbf{D} + (1+t)^{-1} \int_{\mathbb{R}} e^{-\frac{C_1 |\xi + \sigma t|^2}{1+t}} |(v - \bar{v}, \, \theta - \bar{\theta})|^2 d\xi \right).
	\end{array}
\end{align*}
Similarly, 
\begin{align*}
	\begin{array}{l}
		\displaystyle |\mathbf{E}_{13}| \le C \int_{\mathbb{R}} \Big( |(\theta - \bar{\theta})_\xi| \cdot |v - \bar{v}| \cdot |\bar{\theta}_\xi| + |\theta - \bar{\theta}| \cdot (|(\theta - \bar{\theta})_\xi| + |\bar{\theta}_\xi|) (|(\theta - \bar{\theta})_\xi| + |\bar{\theta}_\xi| \cdot |v - \bar{v}|) \\[4mm]
		\displaystyle \qquad \qquad + |\theta - \bar{\theta}| \cdot \left( |(u - \bar{u})_\xi|^2 + |(u - \bar{u})_\xi| \cdot |\bar{u}_\xi| + |\bar{u}_\xi|^2 \cdot |v - \bar{v}| \right)  + |\theta - \bar{\theta}|^2 \cdot |\bar{u}_\xi|^2  \Big) d\xi \\[4mm]
		\displaystyle \qquad \qquad + \kappa \int_{\mathbb{R}} \left( \dfrac{(\theta - \bar{\theta})^2}{\theta \bar{\theta}} \right)_\xi \cdot \left( \dfrac{\bar{\theta}_\xi}{\bar{v}} \right) d\xi \\[4mm]
		\displaystyle \qquad \le C (\delta_0 + \varepsilon_1) \left( \mathcal{G}^R + \mathcal{G}^S + \mathbf{D} + (1+t)^{-1} \int_{\mathbb{R}} e^{-\frac{C_1 |\xi + \sigma t|^2}{1+t}} |(v - \bar{v}, \, \theta - \bar{\theta})|^2 d\xi \right).
	\end{array}
\end{align*}

Finally, using $\reff{eq:ARV-change}_4$ and $\reff{hat-p-Phi}_2$ for \reff{Q123_R} $i=1,2$, and Lemma \ref{ARV-pro}, we have 
\begin{align}  \label{Q12R-upbound}
	\begin{array}{l}
		\displaystyle |(Q_1^R, Q_2^R)| \le C (|v_{\xi\xi\xi}^R| + |v_{\xi\xi}^R| \cdot |v_\xi^R| + |v_\xi^R|^3 + |(u_{\xi\xi}^R, \theta_{\xi\xi}^R)| + |(v_\xi^R, u_\xi^R, \theta_\xi^R)|^2) \\[4mm]
		\displaystyle \qquad \le C (|(v_{\xi\xi}^S, u_{\xi\xi}^S, \theta_{\xi\xi}^S)| + |(v_\xi^S, u_\xi^S, \theta_\xi^S)|^2),
	\end{array}
\end{align}
which yields 
\begin{align*}
	\| (Q_1^R, \, Q_2^R) \|_{L^1(\mathbb{R})} \le C ( \| (v_{\xi\xi}^S, \, u_{\xi\xi}^S, \, \theta_{\xi\xi}^S) \|_{L^1(\mathbb{R})} + \| (v_\xi^S, \, u_\xi^S, \, \theta_\xi^S) \|_{L^2(\mathbb{R})}^2 ).
\end{align*}
Note that from Lemma \ref{ARV-pro},
\begin{align*}
	\| (v_{\xi\xi}^R, \, u_{\xi\xi}^R, \, \theta_{\xi\xi}^R) \|_{L^1(\mathbb{R})} \le 
	\begin{cases}
		C \delta_R, & \text{if } 1+t \le \delta_R^{-1} \\[2mm]
		C (1+t)^{-1}, & \text{if } 1+t \ge \delta_R^{-1},
	\end{cases}
\end{align*}
and
\begin{align*}
	\| (v_\xi^R, \, u_\xi^R, \, \theta_\xi^R) \|_{L^2(\mathbb{R})} \le 
	\begin{cases}
		C \delta_R, & \text{if } 1+t \le \delta_R^{-1} \\[2mm]
		C \delta_R^{1/2} (1+t)^{-1/2}, & \text{if } 1+t \ge \delta_R^{-1}.
	\end{cases}
\end{align*}
Thus
\begin{align} \label{Q12R-Est}
	\int_0^\infty \| (Q_1^R, \, Q_2^R) \|_{L^2(\mathbb{R})}^{4/3} dt \le C \delta_R^{1/3}.
\end{align}
By \reff{Q123_C} and Lemma \ref{VCD}, it follows that
\begin{align}  \label{Q12C-Est}
	\| Q_1^C \|_{L^1(\mathbb{R})} \le C \delta_C (1+t)^{-5/4}, \quad \| Q_2^C \|_{L^1(\mathbb{R})} \le C \delta_C (1+t)^{-7/4}.
\end{align}
Together with Lemma \ref{WI-est}, that is, 
\begin{align}  \label{Q12I-Est}
	\| (Q_1^I, \, Q_2^I) \|_{L^2(\mathbb{R})} \le C \delta_S (\delta_R + \delta_C) e^{-C \delta_S t} + C \delta_R \delta_C e^{-Ct},
\end{align}
and using $Q_i := Q_i^I + Q_i^R + Q_i^C, \, i=1,2$ from \reff{Q123},
we have 
\begin{align*}
	\begin{array}{l}
		\displaystyle |\mathbf{E}_{14}| \le C \int_{\mathbb{R}} \left[ (|Q_1^I| + |Q_1^R| + |Q_1^C|) (u - \bar{u}) + (|Q_2^I| + |Q_2^R| + |Q_2^C|) (\theta - \bar{\theta}) \right] d\xi \\[4mm]
		\displaystyle \qquad \le C \left( \| Q_1^I \|_{L^2(\mathbb{R})} + \| Q_1^C \|_{L^2(\mathbb{R})} \right) \| u - \bar{u} \|_{L^2(\mathbb{R})} \\[4mm]
		\displaystyle \qquad \qquad + C \left( \| Q_2^I \|_{L^2(\mathbb{R})} + \| Q_2^C \|_{L^2(\mathbb{R})} \right) \| \theta - \bar{\theta} \|_{L^2(\mathbb{R})} \\[4mm]
		\displaystyle \qquad \qquad + C \| Q_1^R \|_{L^1(\mathbb{R})} \| u - \bar{u} \|_{L^\infty(\mathbb{R})} + C \| Q_2^R \|_{L^1(\mathbb{R})} \| \theta - \bar{\theta} \|_{L^\infty(\mathbb{R})} \\[4mm]
		\displaystyle \qquad \le C \left[ \delta_S (\delta_R + \delta_C) e^{-C \delta_S t} + \delta_R \delta_C e^{-Ct} + \delta_C (1+t)^{-5/4} \right] ( \| u - \bar{u} \|_{L^2(\mathbb{R})} + \| \theta - \bar{\theta} \|_{L^2(\mathbb{R})} ) \\[4mm]
		\displaystyle \qquad \qquad + C \| Q_1^R \|_{L^1(\mathbb{R})} \| u - \bar{u} \|_{L^2(\mathbb{R})}^{1/2} \| (u - \bar{u})_\xi \|_{L^2(\mathbb{R})}^{1/2} \\[4mm]
		\displaystyle \qquad \qquad + C \| Q_2^R \|_{L^1(\mathbb{R})} \| \theta - \bar{\theta} \|_{L^2(\mathbb{R})}^{1/2} \| (\theta - \bar{\theta})_\xi \|_{L^2(\mathbb{R})}^{1/2} \\[4mm]
		\displaystyle \qquad \le C \left[ \delta_S (\delta_R + \delta_C) e^{-C \delta_S t} + \delta_R \delta_C e^{-Ct} + \delta_C (1+t)^{-5/4} \right] ( \| u - \bar{u} \|_{L^2(\mathbb{R})} + \| \theta - \bar{\theta} \|_{L^2(\mathbb{R})} ) \\[4mm]
		\displaystyle \qquad \qquad + C \left( \| (Q_1^R, \, Q_2^R) \|_{L^1(\mathbb{R})}^{4/3} \right) \| (u - \bar{u}, \, \theta - \bar{\theta}) \|_{L^2(\mathbb{R})}^{2/3}  + \dfrac{1}{40} \mathbf{D}.
	\end{array}
\end{align*}

Thus,
\begin{align} \label{E1-Est}
	\begin{array}{l}
		\displaystyle |\mathbf{E}_1| \le C \left[ \delta_S (\delta_R + \delta_C) e^{-C \delta_S t} + \delta_R \delta_C e^{-Ct} + \delta_C (1+t)^{-5/4} \right] \| (u - \bar{u}, \, \theta - \bar{\theta}) \|_{L^\infty(0, T; L^2(\mathbb{R}))} \\[4mm]
		\displaystyle \qquad + C \| (Q_1^R, \, Q_2^R) \|_{L^1(\mathbb{R})}^{4/3} \| (u - \bar{u}, \, \theta - \bar{\theta}) \|_{L^\infty(0, T; L^2(\mathbb{R}))}^{2/3} + \dfrac{1}{20} \mathbf{D}\\[4mm]
		\displaystyle \qquad + (\lambda + \delta_0 + \varepsilon_1) \left( \mathcal{G}^S + \mathcal{G}^R + (1+t)^{-1} \int_{\mathbb{R}} e^{-\frac{C_1 |\xi + \sigma t|^2}{1+t}} |(v - \bar{v}, \, \theta - \bar{\theta})|^2 d\xi + \| (v-\bar{v})_\xi \|_{L^2(\mathbb{R})}^2 \right).
	\end{array}
\end{align}
\hspace*{\fill}\\

\noindent $\bullet$ {\bf Estimates on $\mathbf{E}_2$.}

From Lemma \ref{REM-lem}, we have 
\begin{align*}
	\begin{array}{l}
		\displaystyle \mathbf{E}_2(U) := \int_{\mathbb{R}} \Bigg\{ - (a^{-\mathbf{X}} (u - \bar{u}))_\xi \left[ \dfrac{\nu^2}{2} \Big( \big( \dfrac{\phi_\xi}{v} \big)^2 - \big( \dfrac{\bar{\phi}_\xi}{\bar{v}} \big)^2 \Big) + \nu^2 \Big( \dfrac{\phi_\xi v_\xi}{v^3} - \dfrac{\bar{\phi}_\xi \bar{v}_\xi}{\bar{v}^3} \Big) - \nu^2 \Big( \dfrac{1}{v^2} - \dfrac{1}{\bar{v}^2} \Big) \phi_{\xi\xi} \right] \\[3mm]
		\displaystyle \qquad + \left[\nu^2 a^{-\mathbf{X}}_\xi (u - \bar{u}) \dfrac{\tilde{\phi}_{\xi\xi}}{\bar{v}^2} + 2 \nu^2 a^{-\mathbf{X}} \left( v_t^R + v_t^C \right) (v - \bar{v}) \dfrac{\tilde{\phi}_{\xi\xi}}{\bar{v}^3} + \nu^2 \left( \dfrac{a^{-\mathbf{X}}}{\bar{v}^2} \right)_\xi (v - \bar{v}) \tilde{\phi}_{\xi t} \right]\\[4mm]
		\displaystyle \qquad + \Bigg[ a^{-\mathbf{X}} \dfrac{\nu^4}{2} \bigg( - \dfrac{3 (v_t^R + v_t^C)}{\bar{v}^4 e^{\bar{\phi}}} - \dfrac{\phi_t^R + \phi_t^C}{\bar{v}^3 e^{\bar{\phi}}} \bigg) \tilde{\phi}_{\xi\xi}^2 - \nu^2 a^{-\mathbf{X}} \dfrac{v_t^R + v_t^C}{2 \bar{v}^2} \tilde{\phi}_\xi^2  \\[3mm]
		\displaystyle \qquad \qquad - \dfrac{\nu^4}{\bar{v} e^{\bar{\phi}}} \Big( \dfrac{a^{-\mathbf{X}}}{\bar{v}^2} \Big)_\xi \tilde{\phi}_{\xi t} \cdot \tilde{\phi}_{\xi\xi} + \Big( \dfrac{\nu^4 a^{-\mathbf{X}}}{\bar{v}^3 e^{\bar{\phi}}} \Big)_\xi \dfrac{\sigma \tilde{\phi}_{\xi\xi}^2}{2} - \dfrac{\sigma \nu^4 a^{-\mathbf{X}}}{\bar{v}^2} \Big( \dfrac{1}{\bar{v} e^{\bar{\phi}}} \Big)_\xi \tilde{\phi}_{\xi\xi}^2 - \Big(\dfrac{a^{-\mathbf{X}}}{\bar{v}} \Big)_\xi \dfrac{\sigma \nu^2 \tilde{\phi}_\xi^2}{2} \Bigg] \\[3mm]
		\displaystyle \qquad + \nu^2 \dfrac{a^{-\mathbf{X}}}{\bar{v}^2} \left( \tilde{\phi}_{\xi t} - \sigma \tilde{\phi}_{\xi\xi} \right) \left[ \Big( - \nu^2 e^{\bar{\phi}} \tilde{\phi}_\xi \dfrac{\bar{v}_\xi}{\bar{v}^2} - \nu^2 e^{-\bar{\phi}} \Big( \dfrac{v - \bar{v}}{\bar{v} v} \phi_\xi \Big)_\xi + N_p \Big)_\xi -\bar{v}_\xi \tilde{\phi} \right] \\[3mm]
		\displaystyle \qquad - \nu^2 \dfrac{a^{-\mathbf{X}}}{\bar{v}^2} \left( \tilde{\phi}_{\xi t} - \sigma \tilde{\phi}_{\xi\xi} \right) (e^{-\bar{\phi}} Q_3)_\xi \Bigg\} d\xi\\[3mm] 
		\displaystyle \qquad\quad =: \sum_{i=1}^{5} \mathbf{E}_{2i}(U),
	\end{array}
\end{align*}

For 
\begin{align*}
	\begin{array}{l}
		\displaystyle |\mathbf{E}_{21}| \le C \int_{\mathbb{R}} \left( \dfrac{\lambda}{\delta_S} |v_\xi^S| \cdot |u - \bar{u}| + |(u - \bar{u})_\xi| \right) \left( |\tilde{\phi}_\xi|^2 + |\tilde{\phi}_\xi| \cdot |\bar{\phi}_\xi| + |\bar{\phi}_\xi|^2 \cdot |v - \bar{v}| \right. \\[4mm]
		\displaystyle \qquad \qquad + |\tilde{\phi}_\xi| \cdot (|(v - \bar{v})_\xi| + |\bar{v}_\xi|) + |\bar{\phi}_\xi| \cdot |(v - \bar{v})_\xi| + |\bar{\phi}_\xi| \cdot |\bar{v}_\xi| \cdot |v - \bar{v}| \\[4mm]
		\displaystyle \qquad \qquad \left. + |v - \bar{v}| \cdot (|\tilde{\phi}_{\xi\xi}| + |\bar{\phi}_{\xi\xi}|) \right) d\xi \\[4mm]
		\displaystyle \qquad \le C (\lambda + \delta_0 + \varepsilon_1) \left( \mathcal{G}^S + \mathcal{G}^R + \mathbf{D} + \| \tilde{\phi}_\xi \|_{H^1(\mathbb{R})}^2 + (1+t)^{-1} \int_{\mathbb{R}} e^{-\frac{C_1 |\xi + \sigma t|^2}{1+t}} |v - \bar{v}|^2 d\xi \right).
	\end{array}
\end{align*}
Similarly,
\begin{align*}
	\begin{array}{l}
		\displaystyle |\mathbf{E}_{22}| \le C \int_{\mathbb{R}} \left( \dfrac{\lambda}{\delta_S} |v_\xi^S| \cdot |u - \bar{u}| \cdot |\tilde{\phi}_{\xi\xi}| + (|v_\xi^S| + |v_\xi^R| + |v_\xi^C|) |v - \bar{v}| \cdot |\tilde{\phi}_{\xi\xi}|  \right. \\[4mm]
		\displaystyle \qquad \qquad \qquad \qquad \left. + (|v_\xi^R| + |u_\xi^R| + |v_\xi^C| + |u_\xi^C|) |v - \bar{v}| \cdot |\tilde{\phi}_{\xi t}| \right) d\xi \\[4mm]
		\displaystyle \qquad \le C (\lambda + \delta_0^{1/2}) \left( \mathcal{G}^S + \mathcal{G}^R + \| \tilde{\phi}_{\xi\xi} \|_{L^2(\mathbb{R})}^2 + \| \tilde{\phi}_{\xi t} \|_{L^2(\mathbb{R})}^2 + (1+t)^{-1} \int_{\mathbb{R}} e^{-\frac{C_1 |\xi + \sigma t|^2}{1+t}} | v - \bar{v} |^2 d\xi \right),
	\end{array}
\end{align*}
and 
\begin{align*}
	\begin{array}{l}
		\displaystyle |\mathbf{E}_{23}| \le C \int_{\mathbb{R}} \left[ (|v_t^R| + |v_t^C| + |\phi_t^R| + |\phi_t^C|) \cdot |\tilde{\phi}_{\xi\xi}|^2 + (|v_t^R| + |v_t^C|) |\tilde{\phi}_\xi|^2 \right. \\[4mm]
		\displaystyle \qquad \qquad \left. + |\bar{v}_\xi| \cdot |\tilde{\phi}_{\xi t}| \cdot |\tilde{\phi}_{\xi\xi}| + (|\bar{v}_\xi| + |\bar{\phi}_\xi|) \cdot |\tilde{\phi}_{\xi\xi}|^2 + |\bar{v}_\xi| \cdot |\tilde{\phi}_\xi|^2 \right] d\xi \\[4mm]
		\displaystyle \qquad \le C \delta_0 \left( \| \tilde{\phi}_\xi \|_{H^1(\mathbb{R})}^2 + \| \tilde{\phi}_{\xi t} \|_{L^2(\mathbb{R})}^2 \right).
	\end{array}
\end{align*}
By the Taylor expansion $|1 + \tilde{\phi} - e^{\tilde{\phi}}| \le C |\tilde{\phi}|^2$, $|1 - e^{\tilde{\phi}}| \le C |\tilde{\phi}|$, then
\begin{align*}
	\begin{array}{l}
		\displaystyle |\mathbf{E}_{24}| \le C \int_{\mathbb{R}} \bigg[ (|\tilde{\phi}_{\xi t}| + |\tilde{\phi}_{\xi\xi}|) \cdot \Big(|\bar{v}_\xi|^2 \cdot \tilde{\phi}_\xi + |\bar{v}_{\xi\xi}| |\tilde{\phi}_\xi| + |\bar{\phi}_\xi| \cdot |\bar{v}_\xi| \cdot |\tilde{\phi}_\xi| + |\bar{v}_\xi| \cdot |\tilde{\phi}_{\xi\xi}|  \\[4mm]
		\displaystyle \qquad \qquad + |\bar{v}_\xi| \cdot |\tilde{\phi}|^2 + |\bar{\phi}| \cdot |\tilde{\phi}_\xi| + |(v - \bar{v})_\xi| \cdot |\tilde{\phi}| + |v - \bar{v}| \cdot |\tilde{\phi}_\xi| + |\bar{v}_\xi| |\tilde{\phi}| \Big) \\[4mm]
		\displaystyle \qquad \qquad + \left(|\bar{v}_\xi|(|\tilde{\phi}_{\xi t}| + |\tilde{\phi}_{\xi\xi}|) + |\tilde{\phi}_{\xi\xi t}| + |\tilde{\phi}_{\xi\xi\xi}| \right) \Big( ( |(v - \bar{v})_\xi| + |\bar{v}_\xi| \cdot |v - \bar{v}| + |(v - \bar{v})_\xi| \cdot |v - \bar{v}| ) \\[4mm]
		\displaystyle \qquad \qquad \times (|\tilde{\phi}_\xi| + |\bar{\phi}_\xi|) + |v - \bar{v}| (|\tilde{\phi}_{\xi\xi}| + |\bar{\phi}_{\xi\xi}|) \Big) \bigg] d\xi \\[4mm]
		\displaystyle \le C (\delta_0^{1/2} + \varepsilon_1) \left( \mathcal{G}^S + \mathcal{G}^R + \mathbf{D} + (1+t)^{-1} \int_{\mathbb{R}} e^{-\frac{C_1 |\xi + \sigma t|^2}{1+t}} | v - \bar{v} |^2 d\xi  \right. \\[4mm]
		\displaystyle \left. \qquad + \| (|v_\xi^R|^{1/2}, \, (1+t)^{-1/2} e^{-\frac{C_1 |\xi + \sigma t|^2}{2(1+t)}}, \, |v_\xi^S|^{1/2}) \tilde{\phi} \|_{L^2(\mathbb{R})}^2 + \| (\tilde{\phi}_\xi, \, \tilde{\phi}_{\xi t}, \, \tilde{\phi}_{\xi\xi}, \, \tilde{\phi}_{\xi\xi t}, \, \tilde{\phi}_{\xi\xi\xi}) \|_{L^2(\mathbb{R})}^2   \right).
	\end{array}
\end{align*}

Finally, 
\begin{align*}
	\begin{array}{l}
		\displaystyle |\mathbf{E}_{25}| \le \displaystyle \int_{\mathbb{R}} (|\tilde{\phi}_{\xi t}| + |\tilde{\phi}_{\xi\xi}|) (|\tilde{\phi}_\xi| |Q_3| + |(Q_{3})_\xi|) d\xi \\[4mm]
		\displaystyle \qquad \quad \le C ( \| \tilde{\phi}_{\xi t} \|_{L^2(\mathbb{R})} + \| \tilde{\phi}_{\xi\xi} \|_{L^2(\mathbb{R})} ) \cdot ( \| \tilde{\phi}_\xi Q_3 \|_{L^2(\mathbb{R})} + \| (Q_{3})_\xi \|_{L^2(\mathbb{R})}).
	\end{array}
\end{align*}
From $Q_3 = Q_3^I + Q_3^R + Q_3^C$ in \reff{Q123}, by lemma \ref{WI-est},
\begin{align*}
	\| ((Q_3^I)_\xi, \bar{v}_\xi Q_3^I) \|_{L^2(\mathbb{R})} \le C\delta_S (\delta_R + \delta_C) e^{-C \delta_S t} + C \delta_R \delta_C  e^{-C t}.
\end{align*}
By $\reff{Q123_R}_3$, it holds that $|(Q_3^R, (Q_3^R)_\xi)| \le C (|v_{\xi\xi}^R| + |v_\xi^R|^2)$ and so 
\begin{align*}
	\| (Q_3^R, (Q_3^R)_\xi) \|_{L^2(\mathbb{R})}^2 \le C \left( \| v_{\xi\xi}^R \|_{L^2(\mathbb{R})}^2 + \| v_\xi^R \|_{L^4(\mathbb{R})}^4 \right).
\end{align*}
Note that by Lemma \ref{ARV-pro},
\begin{align*}
	\| v_{\xi\xi}^R \|_{L^2(\mathbb{R})} \le 
	\begin{cases}
		C \delta_R, & \text{if } 1+t \le \delta_R^{-1} \\
		C (1+t)^{-1}, & \text{if } 1+t \ge \delta_R^{-1}
	\end{cases}
\end{align*}
and
\begin{align*}
	\| v_\xi^R \|_{L^4(\mathbb{R})} \le 
	\begin{cases}
		C \delta_R, & \text{if } 1+t \le \delta_R^{-1} \\
		C \delta_R^{1/4} (1+t)^{-3/4}, & \text{if } 1+t \ge \delta_R^{-1}
	\end{cases}
\end{align*}
Thus,
\begin{align}  \label{Q3R_Est}
	\int_0^{+\infty} \| (Q_3^R, (Q_3^R)_\xi) \|_{L^2(\mathbb{R})}^2 dt \le C \delta_R.
\end{align}
By \reff{Q123_C} and Lemma \ref{VCD}, it follows that
\begin{align}  \label{Q3C_Est}
	\| Q_3^C \|_{L^2(\mathbb{R})} \le C \delta_C (1+t)^{-3/4}, \quad \| (Q_3^C)_\xi \|_{L^2(\mathbb{R})} \le C \delta_C (1+t)^{-5/4}.
\end{align}
Thus by the Young inequality,
\begin{align*}
	\begin{array}{l}
		\displaystyle |\mathbf{E}_{25}|  \le C \delta_0^{1/3} \left( \| \tilde{\phi}_{\xi t} \|_{L^2(\mathbb{R})}^2 + \| \tilde{\phi}_{\xi\xi} \|_{L^2(\mathbb{R})}^2 \right) + C \delta_0^{5/3} \left( \| Q_3^R \|_{L^2(\mathbb{R})}^2 + \delta_C^2 (1+t)^{-3/2} \right) \\[4mm]
		\displaystyle \qquad \qquad + C \delta_0^{-1/3} \left( \delta_S^2 (\delta_R + \delta_C)^2 e^{-C \delta_S t} + \delta_R^2 \delta_C^2 e^{-Ct} + \| (Q_3^R)_\xi \|_{L^2(\mathbb{R})}^2 + \delta_C^2 (1+t)^{-5/2} \right).
	\end{array}
\end{align*}
Therefore,
\begin{align}  \label{E2-Est}
	\begin{array}{l}
		\displaystyle |\mathbf{E}_2| \le C (\lambda + \delta_0^{1/3} + \varepsilon_1) \left( \mathcal{G}^S + \mathcal{G}^R + \mathbf{D} + (1+t)^{-1} \int_{\mathbb{R}} e^{-\frac{C_1 |\xi + \sigma t|^2}{1+t}} |v - \bar{v}|^2 d\xi \right. \\[4mm]
		\displaystyle \qquad \left. + \| (|v_\xi^R|^{1/2}, \, (1+t)^{-1/2} e^{-\frac{C_1 |\xi + \sigma t|^2}{2(1+t)}}, \, |v_\xi^S|^{1/2}) \tilde{\phi} \|_{L^2(\mathbb{R})}^2 + \| (\tilde{\phi}_\xi, \, \tilde{\phi}_{\xi\xi}, \, \tilde{\phi}_{\xi t}, \, \tilde{\phi}_{\xi\xi t}, \, \tilde{\phi}_{\xi\xi\xi}) \|_{L^2(\mathbb{R})}^2 \right) \\[4mm]
		\displaystyle \qquad + C \delta_0^{5/3} \left( \| Q_3^R \|_{L^2(\mathbb{R})}^2 + \delta_C^2 (1+t)^{-3/2} \right) \\[4mm]
		\displaystyle \qquad + C \delta_0^{-1/3} \left( \delta_S^2 (\delta_R + \delta_C)^2 e^{-C \delta_S t} + C \delta_R^2 \delta_C^2 e^{-Ct} + \| (Q_3^R)_\xi \|_{L^2(\mathbb{R})}^2 + \delta_C^2 (1+t)^{-5/2} \right).
	\end{array}
\end{align}
\hspace*{\fill}\\

\noindent $\bullet$ {\bf Conclusion.}

Substituting \reff{Y-Est}, \reff{E1-Est} and \reff{E2-Est} into \reff{RE-mid2}, it follows that 
\begin{align*}
	\begin{array}{l}
		\displaystyle \dfrac{d}{dt} \int_{\mathbb{R}} a^{-\mathbf{X}} \tilde{\eta}((U, \phi) | (\bar{U}, \bar{\phi})) \, d\xi \le - C (\mathcal{G}^S + \mathcal{G}^R) - \dfrac{\delta_S}{4M} |\dot{\mathbf{X}}(t)|^2 - \dfrac{1}{8} \mathbf{D} + C \delta_S^{4/3} (\delta_R^{4/3} + \delta_C^{4/3}) e^{-C \delta_S t}\\[4mm]
		\displaystyle \qquad + C \left[ \delta_S (\delta_R + \delta_C) e^{-C \delta_S t} + \delta_R \delta_C e^{-Ct} + \delta_C (1+t)^{-5/4} \right] \| (u - \bar{u}, \, \theta - \bar{\theta}) \|_{L^\infty(0, T; L^2(\mathbb{R}))}  \\[4mm]
		\displaystyle \qquad + C \| (Q_1^R, \, Q_2^R) \|_{L^1(\mathbb{R})}^{4/3} \| (u - \bar{u}, \, \theta - \bar{\theta}) \|_{L^\infty(0, T; L^2(\mathbb{R}))}^{2/3} + C \delta_0^{5/3} \left( \| Q_3^R \|_{L^2(\mathbb{R})}^2 + \delta_C^2 (1+t)^{-3/2} \right) \\[4mm]
		\displaystyle \qquad + C \delta_0^{-1/3} \left( \delta_S^2 (\delta_R + \delta_C)^2 e^{-C \delta_S t} + \delta_R^2 \delta_C^2 e^{-Ct} + \| (Q_3^R)_\xi \|_{L^2(\mathbb{R})}^2 + \delta_C^2 (1+t)^{-5/2} \right) \\[4mm]
		\displaystyle \qquad + C (\lambda + \delta_0^{1/3} + \varepsilon_1) \left( (1+t)^{-1} \int_{\mathbb{R}} e^{-\frac{C_1 |\xi + \sigma t|^2}{1+t}} |(v - \bar{v}, \, \theta - \bar{\theta})|^2 d\xi + \| (v-\bar{v})_\xi \|_{L^2(\mathbb{R})}^2  \right. \\[4mm]
		\displaystyle \qquad \quad \left. + \| (|v_\xi^R|^{1/2}, \, (1+t)^{-1/2} e^{-\frac{C_1 |\xi + \sigma t|^2}{2(1+t)}}, \, |v_\xi^S|^{1/2}) \tilde{\phi} \|_{L^2(\mathbb{R})}^2 + \| (\tilde{\phi}_\xi, \, \tilde{\phi}_{\xi t}, \, \tilde{\phi}_{\xi\xi}, \, \tilde{\phi}_{\xi\xi t}, \, \tilde{\phi}_{\xi\xi\xi}) \|_{L^2(\mathbb{R})}^2 \right).
	\end{array}
\end{align*}

From the Poisson equation $\reff{RCDS}_4$ and Lemma \ref{ARV-pro}-\ref{VS-pro}, it follows 
\begin{align*}
	|\bar{v} e^{\bar{\phi}} - 1| \le |Q_3 + \nu^2 \left( \dfrac{\bar{\phi}_\xi}{\bar{v}} \right)_\xi| \le C \delta_0,
\end{align*}
which together with \reff{Ass-Inf} and \reff{st3-v} yields
\begin{align*}
	\begin{array}{l}
		\displaystyle R \bar{\theta} \, \Phi\left( \dfrac{v}{\bar{v}} \right) + \dfrac{(v - \bar{v})^2}{2 \bar{v}^2} - \dfrac{(v - \bar{v})^3}{3 \bar{v}^3} - \dfrac{\nu^2}{\bar{v}^2} (v - \bar{v}) \tilde{\phi}_{\xi\xi} + \dfrac{\nu^4}{2 \bar{v}^3 e^{\bar{\phi}}} \tilde{\phi}_{\xi\xi}^2 \\[4mm]
		\displaystyle \ge (1 - C \varepsilon_1) R \bar{\theta} \dfrac{(v - \bar{v})^2}{2 \bar{v}^2} + \dfrac{(v - \bar{v})^2}{2 \bar{v}^2} - \dfrac{\nu^2}{\bar{v}^2} (v - \bar{v}) \tilde{\phi}_{\xi\xi} + (1 - C \delta_0) \dfrac{\nu^4}{2 \bar{v}^2} \tilde{\phi}_{\xi\xi}^2 \\[4mm]
		\displaystyle \ge C (|v - \bar{v}|^2 + |\tilde{\phi}_{\xi\xi}|^2),
	\end{array}
\end{align*}
which yields 
\begin{align*}
	\int_{\mathbb{R}} a^{-\mathbf{X}} \tilde{\eta}((U, \, \phi) | (\bar{U}, \, \bar{\phi})) \, d\xi \sim \| (U - \bar{U}, \, \tilde{\phi}_\xi, \, \tilde{\phi}_{\xi\xi}) \|_{L^2(\mathbb{R})}^2.
\end{align*}
Integrating the above inequality over $[0, t]$ for all $t \le T$, and using \reff{Q12R-Est} and \reff{Q3R_Est} and the Young inequality yields 
\begin{align}  \label{RE-mid3}
	\begin{array}{l}
		\displaystyle \sup_{t \in [0, T]} \| (U - \bar{U}, \, \tilde{\phi}_\xi, \, \tilde{\phi}_{\xi\xi}) \|_{L^2(\mathbb{R})}^2 + \int_0^t (|\dot{\mathbf{X}}(\tau)|^2 + \mathcal{G}^S + \mathcal{G}^R + \mathbf{D}) d\tau \\[4mm]
		\displaystyle \le C \| (U - \bar{U}, \, \tilde{\phi}_\xi, \, \tilde{\phi}_{\xi\xi})(t=0) \|_{L^2(\mathbb{R})}^2 + C \delta_0^{1/2} \\[4mm]
		\displaystyle + C (\lambda + \delta_0^{1/3} + \varepsilon_1) \left( \int_{0}^{t} (1+\tau)^{-1} \int_{\mathbb{R}} e^{-\frac{C_1 |\xi + \sigma \tau|^2}{ 1 + \tau}} |(v - \bar{v}, \, \theta - \bar{\theta})|^2 d\xi d\tau + \int_{0}^{t}\| (v-\bar{v})_\xi \|_{L^2(\mathbb{R})}^2 d\tau \right. \\[4mm]
		\displaystyle \left. + \int_{0}^{t} \| (|v_\xi^R|^{1/2}, (1+\tau)^{-1/2} e^{-\frac{C_1 |\xi + \sigma \tau|^2}{2( 1 + \tau)}},  |v_\xi^S|^{1/2}) \tilde{\phi} \|_{L^2(\mathbb{R})}^2 d\tau + \int_{0}^{t} \| (\tilde{\phi}_\xi,  \tilde{\phi}_{\xi \tau},  \tilde{\phi}_{\xi\xi},  \tilde{\phi}_{\xi\xi \tau},  \tilde{\phi}_{\xi\xi\xi}) \|_{L^2(\mathbb{R})}^2 d\tau \right).
	\end{array}
\end{align}
with $U, \bar{U}$ defined in \reff{U-nota}.
\\

\begin{lemma} \label{Heat-lem}
	The following estimates hold.
	\begin{align*}
		\begin{array}{l}
			\displaystyle \int_{0}^{t} \mathcal{G}^C d\tau : = \int_0^t (1+\tau)^{-1} \int_{\mathbb{R}} e^{-\frac{C_1 |\xi + \sigma \tau|^2}{1+\tau}} |(v - \bar{v}, \, u - \bar{u}, \, \theta - \bar{\theta})|^2 d\xi d\tau \\[4mm]
			\displaystyle \le C \sup_{t \in [0, T]} \| U - \bar{U} \|_{L^2(\mathbb{R})}^2 + C \int_0^t (\delta_S|\dot{\mathbf{X}}(\tau)|^2 + \mathcal{G}^S + \mathcal{G}^R + \mathbf{D}) d\tau \\[4mm]
			\displaystyle \qquad + C\int_{\mathbb{R}} \| \tilde{\phi}_{\xi}\|_{H^1(\mathbb{R})}^2 d\tau + C \int_0^t \| (v - \bar{v})_\xi \|_{L^2(\mathbb{R})}^2 d\tau + C \delta_0^{1/3}.
		\end{array}
	\end{align*}
\end{lemma}
\hspace*{\fill}

\begin{remark}
	Note that compared with the proof of this heat kernel type lemma in \cite{Kang25}, the only difference for the subsystem $\reff{PS}_i (i=1,2,3)$ lies in the additional higher-order term $\Phi(v, \phi)$ and the slightly different term $Q_1$. 
	We expect that these terms do not affect the structure of the final estimate, except for yielding the extra space-time estimates for $\tilde{\phi}_\xi, \tilde{\phi}_{\xi\xi}$ on the right-hand side.
\end{remark}
\hspace*{\fill} \\

\subsection{Elliptic estimates}
\hspace*{\fill}\\ 

This subsection is devoted to proving estimates for $\tilde{\phi}$ and its derivations appearing on the right-hand side of \reff{RE-mid3}. 
The proof is mainly based on the following Poisson equation. For convenience, set 
\begin{align*}
	\tilde{v}:=v-\bar{v}, \qquad \tilde{u}:=u-\bar{u},
\end{align*} 
and then $\reff{PS}_4$ has the form 
\begin{align} \label{poisson-sing}
	-\nu^2 \left( \dfrac{\phi_\xi}{v} - \dfrac{\bar{\phi}_\xi}{\bar{v}} \right)_\xi = - e^{\bar{\phi}} (\tilde{v} + \bar{v} \tilde{\phi}) + e^{\bar{\phi}} N_p - Q_3
\end{align}
with  $Q_3$ defined in \reff{Q123} and 
\begin{align} \label{Np-sig}
	\displaystyle N_p := \bar{v} (1 + \tilde{\phi} - e^{\tilde{\phi}}) + \tilde{v} (1 - e^{\tilde{\phi}}).
\end{align}
\hspace*{\fill}

\begin{lemma}  \label{phi-low}
	There exists a generic constant $C > 0$ such that
	\begin{align*}
		\begin{array}{l}
			\displaystyle \int_0^t \| (|v_\xi^R|^{1/2}, (1+\tau)^{-1/2} e^{-\frac{C_1 |\xi + \sigma \tau|^2}{2( 1 + \tau)}},  |v_\xi^S|^{1/2}) \tilde{\phi} \|_{L^2(\mathbb{R})}^2 d\tau  \\[4mm]
			\displaystyle \le C \delta_0^2 + C \delta_0 \int_0^t \| \tilde{\phi}_\xi \|_{L^2(\mathbb{R})}^2 d\tau + C \int_0^t \left( \mathcal{G}^S + \mathcal{G}^R + \mathcal{G}^C \right) d\tau.
		\end{array}
	\end{align*}
\end{lemma}

\begin{proof}[Proof]
	From the above equation \reff{poisson-sing}, we have
	\begin{align} \label{for basic}
		\bar{v} e^{\bar{\phi}} \tilde{\phi} - \nu^2 \left[ \dfrac{\tilde{\phi}_\xi}{v} + \bar{\phi}_\xi \left( \dfrac{1}{\bar{v}} - \dfrac{1}{v} \right) \right]_\xi = - e^{\bar{\phi}} \tilde{v} + e^{\bar{\phi}} N_p - Q_3.
	\end{align}
	Multiplying by $v_\xi^S \tilde{\phi}$ and integrating over $[0, t] \times \mathbb{R}$ yields
	\begin{align*}
		\begin{array}{l}
			\displaystyle \int_0^t \int_{\mathbb{R}} \left( \bar{v} e^{\bar{\phi}} v_\xi^S |\tilde{\phi}|^2 + \nu^2 v_\xi^S \dfrac{|\tilde{\phi}_\xi|^2}{v} \right) d\xi d\tau \\[4mm]
			\displaystyle = \int_0^t \int_{\mathbb{R}} \left[ - \nu^2 \left( \dfrac{v_{\xi\xi}^S \tilde{\phi} \tilde{\phi}_\xi}{v} + (v_\xi^S \tilde{\phi})_\xi \bar{\phi}_\xi \left( \dfrac{1}{v} - \dfrac{1}{\bar{v}} \right) \right) + v_\xi^S \tilde{\phi} (- e^{\bar{\phi}} \tilde{v} + e^{\bar{\phi}} N_p - Q_3) \right] d\xi d\tau \\[4mm]
			\displaystyle \le C \int_0^t \int_{\mathbb{R}} \left[ \delta_S |v_\xi^S| \cdot |\tilde{\phi}| \cdot |\tilde{\phi}_\xi| + (\delta_S |v_\xi^S| \cdot |\tilde{\phi}| + |v_\xi^S| \cdot |\tilde{\phi}_\xi|) |\bar{\phi}_\xi| \cdot |v - \bar{v}| \right. \\[4mm]
			\displaystyle \qquad \left. + |v_\xi^S| \cdot |\tilde{\phi}| (|\tilde{v}| + |\tilde{\phi}|^2 + |\tilde{\phi}| \cdot |\tilde{v}| + |Q_3|) \right] d\xi d\tau \\[4mm]
			\displaystyle \le C \delta_S \int_0^t \int_{\mathbb{R}} \nu^2 v_\xi^S \dfrac{|\tilde{\phi}_\xi|^2}{v} d\xi d\tau + \dfrac{1}{2} \int_0^t \int_{\mathbb{R}} \bar{v} e^{\bar{\phi}} v_\xi^S |\tilde{\phi}|^2 d\xi d\tau \\[4mm]
			\displaystyle \qquad + C \int_0^t \left( \mathcal{G}^S + \mathcal{G}^R + \mathcal{G}^C \right) d\tau + C \delta_S \int_0^t \| Q_3 \|_{L^2(\mathbb{R})}^2 d\tau.
		\end{array}
	\end{align*}
	Using Lemma \ref{WI-est}, \reff{Q3R_Est} and \reff{Q3C_Est}, we have 
	\begin{align*}  
		\int_0^t \| Q_3 \|_{L^2(\mathbb{R})}^2 d\tau \le C \delta_0,
	\end{align*} 
	which yields
	\begin{align*}
		\begin{array}{l}
			\displaystyle \int_0^t \| (v_\xi^S)^{1/2} \tilde{\phi} \|_{L^2(\mathbb{R})}^2 d\tau \le C \delta_0^2 + C \int_0^t \left( \mathcal{G}^S + \mathcal{G}^R + \mathcal{G}^C \right) d\tau.
		\end{array}
	\end{align*}

	Similarly, it holds 
	\begin{align*}
		\displaystyle \int_0^t \| (v_\xi^R)^{1/2} \tilde{\phi} \|_{L^2(\mathbb{R})}^2 d\tau \le C \delta_0^2 + C \delta_R \int_0^t \| \tilde{\phi}_\xi \|_{L^2(\mathbb{R})}^2 d\tau + C \int_0^t (\mathcal{G}^S + \mathcal{G}^R + \mathcal{G}^C ) d\tau,
	\end{align*}
	where the second term on the right hand is due to the term
	\begin{align*}
		-\nu^2 \int_0^t \int_{\mathbb{R}} \dfrac{v_{\xi\xi}^R \tilde{\phi} \tilde{\phi}_\xi}{v} d\xi d\tau \le \dfrac{1}{4} \int_0^t \int_{\mathbb{R}} \bar{v} e^{\bar{\phi}} v_\xi^R |\tilde{\phi}|^2 d\xi d\tau + C \delta_R \int_0^t \| \tilde{\phi}_\xi \|_{L^2(\mathbb{R})}^2 d\tau.
	\end{align*}
	
	Likewise, we have
	\begin{align*}
		\displaystyle \int_0^t \| (1+\tau)^{-1/2} e^{-\frac{C_1 |\xi + \sigma \tau|^2}{2( 1 + \tau)}} \tilde{\phi} \|_{L^2(\mathbb{R})}^2 d\tau \le C \delta_0^2 + C \delta_C \int_0^t \| \tilde{\phi}_\xi \|_{L^2(\mathbb{R})}^2 d\tau + C \int_0^t (\mathcal{G}^S + \mathcal{G}^R + \mathcal{G}^C ) d\tau,
	\end{align*}
	where the second term on the right-hand side is attributed to
	\begin{align*}
		-\nu^2 \int_0^t \int_{\mathbb{R}} \dfrac{(h^2)_\xi \tilde{\phi} \tilde{\phi}_\xi}{v} d\xi d\tau \le \dfrac{1}{4} \int_0^t \int_{\mathbb{R}} \bar{v} e^{\bar{\phi}} h^2 |\tilde{\phi}|^2 d\xi d\tau + C \delta_C \int_0^t \| \tilde{\phi}_\xi \|_{L^2(\mathbb{R})}^2 d\tau
	\end{align*}
	with $h:=(1+\tau)^{-1/2} e^{-\frac{C_1 |\xi + \sigma \tau|^2}{2( 1 + \tau)}}$.
	Combining this with the computation above, we conclude the proof of the lemma.
\end{proof}
\hspace*{\fill}

Subsequently, we will give the higher-order estimates for $\tilde{\phi}$, more precisely, i.e., the last term on the right-hand side of \reff{RE-mid3}.
\\ 

\begin{lemma}  \label{phi-hig}
	There exists a generic constant $C > 0$ such that
	\begin{align*}
		\begin{array}{l}
			\displaystyle \int_0^t \| (\tilde{\phi}_\tau, \, \tilde{\phi}_\xi, \, \tilde{\phi}_{\xi \tau}, \, \tilde{\phi}_{\xi\xi}, \, \tilde{\phi}_{\xi\xi \tau}, \, \tilde{\phi}_{\xi\xi\xi}) \|_{L^2(\mathbb{R})}^2 d\tau  \\[4mm]
			\displaystyle \le C \delta_0 + C \int_0^t \left( \delta_S |\dot{\mathbf{X}}(\tau)|^2 + \mathcal{G}^S + \mathcal{G}^R + \delta_C\mathcal{G}^C + \mathbf{D} + \| \tilde{v}_\xi \|_{H^1(\mathbb{R})}^2 + \| \tilde{u}_{\xi\xi} \|_{L^2(\mathbb{R})}^2 \right) d\tau.
		\end{array}
	\end{align*}
\end{lemma}

\begin{proof} [Proof]
	Differentiating the equation \reff{poisson-sing} with $\xi$, multiplying by $\tilde{\phi}_\xi$ and integrating over $[0, t] \times \mathbb{R}$, it follows that
	\begin{align*}
		\begin{array}{l}
			\displaystyle \int_0^t \int_{\mathbb{R}} \left( \bar{v} e^{\bar{\phi}} |\tilde{\phi}_\xi|^2 + \nu^2 \dfrac{|\tilde{\phi}_{\xi\xi}|^2}{v} \right) d\xi d\tau \\[4mm]
			\displaystyle = \int_0^t \int_{\mathbb{R}} \bigg[ \nu^2 \left( \dfrac{v_\xi \tilde{\phi}_\xi}{v^2} + \dfrac{\tilde{v} \bar{\phi}_{\xi\xi}}{\bar{v} v} + \left( \dfrac{\tilde{v}}{\bar{v} v} \right)_\xi \bar{\phi}_\xi \right) \tilde{\phi}_{\xi\xi} - \Big( (e^{\bar{\phi}})_\xi \tilde{v} + e^{\bar{\phi}} \tilde{v}_\xi + (\bar{v} e^{\bar{\phi}})_\xi \tilde{\phi}_\xi \\[4mm]
			\displaystyle \qquad \qquad \qquad - (e^{\bar{\phi}})_\xi N_p - e^{\bar{\phi}} (N_p)_\xi + (Q_3)_\xi \Big) \tilde{\phi}_\xi \bigg] d\xi d\tau \\[4mm]
			\displaystyle \le C (\delta_0 + \varepsilon_1) \int_0^t \int_{\mathbb{R}} \left( \bar{v} e^{\bar{\phi}} |\tilde{\phi}_\xi|^2 + \nu^2 \dfrac{|\tilde{\phi}_{\xi\xi}|^2}{v} \right) d\xi d\tau \\[4mm]
			\displaystyle \qquad + C \delta_0 \int_0^t \| (|v_\xi^R|^{1/2}, (1+\tau)^{-1/2} e^{-\frac{C_1 |\xi + \sigma \tau|^2}{2( 1 + \tau)}},  |v_\xi^S|^{1/2}) \tilde{\phi} \|_{L^2(\mathbb{R})}^2 d\tau \\[4mm]
			\displaystyle \qquad + C \int_0^t (\mathcal{G}^S + \mathcal{G}^R +  \delta_C\mathcal{G}^C + \| \tilde{v}_\xi \|_{L^2(\mathbb{R})}^2 ) d\tau + C \int_0^t \| (Q_3)_\xi \|_{L^2(\mathbb{R})}^2 d\tau.
		\end{array}
	\end{align*}
	Using Lemma \ref{WI-est}, \reff{Q3R_Est} and \reff{Q3C_Est}, it holds that
	\begin{align*}
		\int_0^t \| (Q_3)_\xi \|_{L^2(\mathbb{R})}^2 d\tau \le C \delta_0,
	\end{align*}
	which together with Lemma \ref{phi-low} yields
	\begin{align} \label{phi_H2E}
		\int_0^t \| \tilde{\phi}_\xi \|_{H^1(\mathbb{R})}^2 d\tau \le C \delta_0 + C \int_0^t (\mathcal{G}^S + \mathcal{G}^R + \delta_C\mathcal{G}^C + \| \tilde{v}_\xi \|_{L^2(\mathbb{R})}^2) d\tau.
	\end{align}
	Similarly, differentiating the equation \reff{poisson-sing} with $\xi$ twice, multiply by $\tilde{\phi}_{\xi\xi}$ and integrating over $[0, t] \times \mathbb{R}$, it implies that
	\begin{align*}
		\int_0^t \| \tilde{\phi}_{\xi\xi} \|_{H^1(\mathbb{R})}^2 d\tau \le C \delta_0 + C \int_0^t (\mathcal{G}^S + \mathcal{G}^R + \delta_C\mathcal{G}^C + \| \tilde{v}_\xi \|_{H^1(\mathbb{R})}^2) d\tau.
	\end{align*}

	Differentiating the equation w.r.t $t$, multiplying by $\tilde{\phi}_t$ and integrating over $[0, t] \times \mathbb{R}$, it implies that
	\begin{align*}
		\begin{array}{l}
			\displaystyle \int_0^t \int_{\mathbb{R}} \left( \bar{v} e^{\bar{\phi}} |\tilde{\phi}_\tau|^2 + \nu^2 \dfrac{|\tilde{\phi}_{\xi \tau}|^2}{v} \right) d\xi d\tau \\[4mm]
			\displaystyle = \int_0^t \int_{\mathbb{R}} \bigg[ \nu^2 \left( \dfrac{v_\tau \tilde{\phi}_{\xi}}{v^2} + \left( \dfrac{\bar{\phi}_\xi}{\bar{v}} \right)_{\tau} \dfrac{\tilde{v}}{v} + \dfrac{\bar{\phi}_\xi}{\bar{v}} \left( \dfrac{\tilde{v}}{v} \right)_\tau \right) \tilde{\phi}_{\xi \tau} \\[4mm]
			\displaystyle \qquad - \Big( (e^{\bar{\phi}})_\tau \tilde{v} + e^{\bar{\phi}} \tilde{v}_\tau + (\bar{v} e^{\bar{\phi}})_\tau \tilde{\phi} - (e^{\bar{\phi}})_\tau N_p - e^{\bar{\phi}} (N_p)_\tau + (Q_3)_\tau \Big) \tilde{\phi}_\tau \bigg] d\xi d\tau.
		\end{array}
	\end{align*}
	Plugging $\reff{eq:CNSP-change}_1, \reff{RCDS}_1$ and $\reff{PS}_1$ into back the above equality, together with Lemma \ref{phi-low} and \reff{phi_H2E}, we have
	\begin{align*}
		\displaystyle \int_0^t \| \tilde{\phi}_\tau \|_{H^1(\mathbb{R})}^2 d\tau \le C \delta_0 + C \int_0^t \left( \delta_S |\dot{\mathbf{X}}(\tau)|^2 + \mathcal{G}^S + \mathcal{G}^R + \delta_C\mathcal{G}^C + \mathbf{D} + \| \tilde{v}_\xi \|_{L^2(\mathbb{R})}^2 \right) d\tau.
	\end{align*}
	Similarly, it follows that
	\begin{align*}
		\displaystyle \int_0^t \| \tilde{\phi}_{\xi \tau} \|_{H^1(\mathbb{R})}^2 d\tau \le C \delta_0 + C \int_0^t \left( \delta_S |\dot{\mathbf{X}}(\tau)|^2 + \mathcal{G}^S + \mathcal{G}^R + \delta_C\mathcal{G}^C + \mathbf{D} + \| \tilde{v}_\xi \|_{H^1(\mathbb{R})}^2 + \| \tilde{u}_{\xi\xi} \|_{L^2(\mathbb{R})}^2 \right) d\tau.
	\end{align*}
	Therefore, we complete the proof of this lemma.
\end{proof}
\hspace*{\fill}

Finally, we give the proof of Lemma \ref{Low-Est-lem}. First of all, plugging \reff{phi_H2E} into Lemma \ref{Heat-lem}, we have  
\begin{align*}
	\begin{array}{l}
		\displaystyle \int_{0}^{t} \mathcal{G}^C d\tau \le C \sup_{t \in [0, T]} \| U - \bar{U} \|_{L^2(\mathbb{R})}^2 + C \int_0^t (\delta_S|\dot{\mathbf{X}}(\tau)|^2 + \mathcal{G}^S + \mathcal{G}^R + \mathbf{D}) d\tau \\[4mm]
		\displaystyle \qquad \qquad + C \int_0^t \| (v - \bar{v})_\xi \|_{L^2(\mathbb{R})}^2 d\tau + C \delta_0^{1/3},
	\end{array}
\end{align*}
which together with Lemma \ref{phi-low} and \ref{phi-hig}, substituted into \reff{RE-mid3}, yields Lemma \ref{Low-Est-lem}.
\hspace*{\fill} \\

\section{Proof of Proposition \ref{A prior estimate}}
\hspace*{\fill}\\ 

In order to close the a prior $H^2$-assumption \reff{ProAss} for $U-\bar{U}$, we need the following lemmas for $L^2$ estimates of the higer order derivatives to complete the proof of Proposition \ref{A prior estimate}. 
\\ 

\begin{lemma}
	There exists a generic constant $C>0$ such that 
	\begin{align*}
		\begin{array}{l}
			\displaystyle \sup_{t \in [0, T]} \left( \| (U - \bar{U}) \|_{H^1(\mathbb{R})}^2 + \| \tilde{\phi} \|_{H^2(\mathbb{R})}^2 \right) \\[4mm]
			\displaystyle + \int_0^t \left( \delta_S |\dot{\mathbf{X}}(\tau)|^2 + \mathcal{G}^S + \mathcal{G}^R + \mathcal{G}^C + \| (v - \bar{v})_\xi \|_{L^2(\mathbb{R})}^2 + \| (u - \bar{u}, \, \theta - \bar{\theta})_\xi \|_{H^1(\mathbb{R})}^2 \right) d\tau \\[4mm]
			\displaystyle \le C \| (U - \bar{U})(t=0) \|_{H^1(\mathbb{R})}^2 + C \delta_0^{1/3} + C (\lambda + \delta_0^{1/3} + \varepsilon_1) \int_0^t \| (v - \bar{v})_{\xi\xi} \|_{L^2(\mathbb{R})}^2 d\tau.
		\end{array}
	\end{align*}
\end{lemma}

\begin{proof} [Proof]
	In this proof, we will give the $H^1$ estimates on $U-\bar{U}$.
	\\ 

	\noindent $\bullet$ {\bf $H^1$ estimates on $v-\bar{v}$.}
	\hspace*{\fill}\\ 

	Differentiating $\reff{PS}_1$ w.r.t $\xi$ and multiplying by $\mu (v - \bar{v})_\xi$ yields
	\begin{align} \label{v-xi}
		\dfrac{\mu}{2} \left( |(v - \bar{v})_\xi|^2 \right)_t - \dfrac{\sigma \mu}{2} \left( |(v - \bar{v})_\xi|^2 \right)_\xi - \dot{\mathbf{X}}(t) v_{\xi\xi}^S \mu (v - \bar{v})_\xi = \mu (u - \bar{u})_{\xi\xi} (v - \bar{v})_\xi. 
	\end{align}
	On the other hand, it follows from \reff{ele-add}, $\reff{hat-p-Phi}_1$, and $\reff{PS}_2$ that
	\begin{align}  \label{PS-u-form}
		\begin{array}{l}
			\displaystyle (u - \bar{u})_t - \sigma (u - \bar{u})_\xi - \dot{\mathbf{X}}(t) u_\xi^S + (p(v, \theta) - p(\bar{v}, \bar{\theta}))_\xi \\[4mm]
			\displaystyle = \mu \left( \dfrac{u_\xi}{v} - \dfrac{\bar{u}_\xi}{\bar{v}} \right)_\xi - \left( \dfrac{\phi_\xi}{v} - \dfrac{\bar{\phi}_\xi}{\bar{v}} \right) - Q_1.
		\end{array}
	\end{align}
	Multiplying the above equation by $- v (v - \bar{v})_\xi$, we have
	\begin{align}  \label{v-xi1}
		\begin{array}{l}
			\displaystyle - v (v - \bar{v})_\xi (u - \bar{u})_t + \sigma v (v - \bar{v})_\xi (u - \bar{u})_\xi + \dot{\mathbf{X}}(t) u_\xi^S v (v - \bar{v})_\xi - (p(v, \theta) - p(\bar{v}, \bar{\theta}))_\xi v (v - \bar{v})_\xi \\[4mm]
			\displaystyle = - \mu \left( \dfrac{u_\xi}{v} - \dfrac{\bar{u}_\xi}{\bar{v}} \right)_\xi v (v - \bar{v})_\xi - \left( \dfrac{\phi_\xi}{v} - \dfrac{\bar{\phi}_\xi}{\bar{v}} \right) v (v - \bar{v})_\xi - v (v - \bar{v})_\xi Q_1.
		\end{array}
	\end{align}
	First,
	\begin{align*}
		\begin{array}{l}
			\displaystyle - v (v - \bar{v})_\xi (u - \bar{u})_t + \sigma v (v - \bar{v})_\xi (u - \bar{u})_\xi \\[4mm]
			\displaystyle = - \partial_t \left[ v (v - \bar{v})_\xi (u - \bar{u}) \right] + \partial_\xi \left[ v (v - \bar{v})_t (u - \bar{u}) \right] \\[4mm]
			\displaystyle \qquad + (v_t - \sigma v_\xi) (v - \bar{v})_\xi (u - \bar{u}) - (v_\xi (u - \bar{u}) + v (u - \bar{u})_\xi) \left( (v - \bar{v})_t - \sigma (v - \bar{v})_\xi \right).
		\end{array}
	\end{align*}
	Second,
	\begin{align*}
		\begin{array}{l}
			\displaystyle - (p(v, \theta) - p(\bar{v}, \bar{\theta}))_\xi v (v - \bar{v})_\xi \\[4mm]
			\displaystyle = p(v, \theta) |(v - \bar{v})_\xi|^2 - R (\theta - \bar{\theta})_\xi (v - \bar{v})_\xi - R v \bar{\theta}_\xi \left( \dfrac{1}{v} - \dfrac{1}{\bar{v}} \right) (v - \bar{v})_\xi \\[4mm]
			\displaystyle \qquad + v \bar{v}_\xi \left( \dfrac{p(v, \theta)}{v} - \dfrac{p(\bar{v}, \bar{\theta})}{\bar{v}} \right) (v - \bar{v})_\xi,
		\end{array}
	\end{align*}
	and
	\begin{align*}
		\begin{array}{l}
			\displaystyle - \mu \left( \dfrac{u_\xi}{v} - \dfrac{\bar{u}_\xi}{\bar{v}} \right)_\xi v (v - \bar{v})_\xi \\[4mm]
			\displaystyle = - \mu (u - \bar{u})_{\xi\xi} (v - \bar{v})_\xi - \mu \bar{u}_{\xi\xi} \left( \dfrac{1}{v} - \dfrac{1}{\bar{v}} \right) v (v - \bar{v})_\xi \\[4mm]
			\displaystyle \qquad + \mu \dfrac{v_\xi}{v} (u - \bar{u})_\xi (v - \bar{v})_\xi + \mu \bar{u}_\xi \left( \dfrac{v_\xi}{v^2} - \dfrac{\bar{v}_\xi}{\bar{v}^2} \right) v (v - \bar{v})_\xi.
		\end{array}
	\end{align*}
	Substituting the above computation into \reff{v-xi1}, adding with \reff{v-xi}, and integrating over $\mathbb{R}$, we have
	\begin{align}  \label{v-xi2}
		\begin{array}{l}
			\displaystyle \dfrac{d}{dt} \int_{\mathbb{R}} \left[ \dfrac{\mu}{2} |(v - \bar{v})_\xi|^2 - v (v - \bar{v})_\xi (u - \bar{u}) \right] d\xi + \int_{\mathbb{R}} p(v, \theta) |(v - \bar{v})_\xi|^2 d\xi \\[4mm]
			\displaystyle = \int_{\mathbb{R}} \left[ R (\theta - \bar{\theta})_\xi (v - \bar{v})_\xi + \dot{\mathbf{X}}(t) (\mu v_{\xi\xi}^S - v u_\xi^S) (v - \bar{v})_\xi \right. \\[4mm]
			\displaystyle \qquad - (v_t - \sigma v_\xi) (v - \bar{v})_\xi (u - \bar{u}) + (v_\xi (u - \bar{u}) + v (u - \bar{u})_\xi) ((v - \bar{v})_t - \sigma (v - \bar{v})_\xi) \\[4mm]
			\displaystyle \qquad + R v \bar{\theta}_\xi \left( \dfrac{1}{v} - \dfrac{1}{\bar{v}} \right) (v - \bar{v})_\xi - v \bar{v}_\xi \left( \dfrac{p(v, \theta)}{v} - \dfrac{p(\bar{v}, \bar{\theta})}{\bar{v}} \right) (v - \bar{v})_\xi \\[4mm]
			\displaystyle \qquad - \mu \bar{u}_{\xi\xi} \left( \dfrac{1}{v} - \dfrac{1}{\bar{v}} \right) v (v - \bar{v})_\xi + \mu \dfrac{v_\xi}{v} (u - \bar{u})_\xi (v - \bar{v})_\xi \\[4mm]
			\displaystyle \qquad \left. + \mu \bar{u}_\xi \left( \dfrac{v_\xi}{v^2} - \dfrac{\bar{v}_\xi}{\bar{v}^2} \right) v (v - \bar{v})_\xi - \left( \dfrac{\phi_\xi}{v} - \dfrac{\bar{\phi}_\xi}{\bar{v}} \right) v (v - \bar{v})_\xi - v (v - \bar{v})_\xi Q_1 \right] d\xi.
		\end{array}
	\end{align}
	
	Next, we only estimate some special terms and the other terms can be handled easily. (see \cite{Kang25}). First of all,
	\begin{align*}
		\int_{\mathbb{R}} R (\theta - \bar{\theta})_\xi (v - \bar{v})_\xi d\xi \le \dfrac{1}{4} \int_{\mathbb{R}} p(v, \theta) |(v - \bar{v})_\xi|^2 d\xi + C \mathbf{D}.
	\end{align*}
	From $\reff{eq:CNSP-change}_1$, it follows $v_t - \sigma v_\xi = u_\xi = (u - \bar{u})_\xi + u_\xi^R + u_\xi^C + u_\xi^S$. Together with \reff{ProAss} and \reff{Ass-Inf}, we have
	\begin{align*}
		\begin{array}{l}
			\displaystyle - \int_{\mathbb{R}} (v_t - \sigma v_\xi) (u - \bar{u}) (v - \bar{v})_\xi d\xi \\[4mm]
			\displaystyle = - \int_{\mathbb{R}} \left[ (u - \bar{u})_\xi + u_\xi^R + u_\xi^C + u_\xi^S \right] (u - \bar{u}) (v - \bar{v})_\xi d\xi \\[4mm]
			\displaystyle \le C \left( \| u - \bar{u} \|_{L^\infty(\mathbb{R})} \cdot \sqrt{\mathbf{D}} + \| u - \bar{u} \|_{L^2(\mathbb{R})}^{1/2} \| (u - \bar{u})_\xi \|_{L^2(\mathbb{R})}^{1/2} \| u_\xi^R \|_{L^2(\mathbb{R})} \right. \\[4mm]
			\displaystyle \qquad \left. + \delta_C \sqrt{\mathcal{G}^C} + \delta_S \sqrt{\mathcal{G}^S} \right) \| (v - \bar{v})_\xi \|_{L^2(\mathbb{R})} \\[3mm]
			\displaystyle \le \dfrac{1}{8} \| p(v, \theta)^{1/2} (v - \bar{v})_\xi \|_{L^2(\mathbb{R})}^2 + C (\delta_0 + \varepsilon_1) (\mathcal{G}^S + \mathcal{G}^C + \mathbf{D}) + C \varepsilon_1 \| u_\xi^R \|_{L^2(\mathbb{R})}^4.
		\end{array}
	\end{align*}
	Note that from Lemma \ref{ARV-pro},
	\begin{align*}
		\| (v_\xi^R, \, u_\xi^R) \|_{L^2(\mathbb{R})} \le 
		\begin{cases}
			C \delta_R, & \text{if } 1+t \le \delta_R^{-1}, \\[2mm]
			C \delta_R^{1/2} (1+t)^{-1/2}, & \text{if } 1+t \ge \delta_R^{-1},
		\end{cases}
	\end{align*}
	which yields
	\begin{align*} 
		\int_0^{+\infty} \| (v_\xi^R, \, u_\xi^R) \|_{L^2(\mathbb{R})}^4 dt \le C \delta_R^3.
	\end{align*}
	
	Second, as for 
	\begin{align*}
		\begin{array}{l}
			\displaystyle - \int_{\mathbb{R}} \left( \dfrac{\phi_\xi}{v} - \dfrac{\bar{\phi}_\xi}{\bar{v}} \right) v (v - \bar{v})_\xi d\xi = - \int_{\mathbb{R}} \left[ \dfrac{\tilde{\phi}_\xi}{v} + \bar{\phi}_\xi \left( \dfrac{1}{v} - \dfrac{1}{\bar{v}} \right) \right] v (v - \bar{v})_\xi d\xi \\[4mm]
			\displaystyle \le C \int_{\mathbb{R}} (|\tilde{\phi}_\xi| + |\bar{\phi}_\xi| \cdot |v - \bar{v}|) |(v - \bar{v})_\xi| d\xi \\[4mm]
			\displaystyle \le \dfrac{1}{8} \| p(v, \theta)^{1/2} (v - \bar{v})_\xi \|_{L^2(\mathbb{R})}^2 + C \| \tilde{\phi}_\xi \|_{L^2(\mathbb{R})}^2 + C \delta_0^{1/2} (\mathcal{G}^S + \mathcal{G}^R + \mathcal{G}^C).
		\end{array}
	\end{align*}
	Finally, for
	\begin{align*}
		- \int_{\mathbb{R}} v (v - \bar{v})_\xi Q_1 d\xi \le \dfrac{1}{8} \| p(v, \theta)^{1/2} (v - \bar{v})_\xi \|_{L^2(\mathbb{R})}^2 + C \| Q_1 \|_{L^2(\mathbb{R})}^2
	\end{align*}
	
	Since $Q_i = Q_i^I + Q_i^R + Q_i^C, i=1,2$ and \reff{Q12R-upbound},
	\begin{align*}
		\| (Q_1^R, \, Q_2^R) \|_{L^2(\mathbb{R})}^2 \le C \left( \| (v_{\xi\xi}^R, \, u_{\xi\xi}^R, \, \theta_{\xi\xi}^R) \|_{L^2(\mathbb{R})}^2 + \| (v_\xi^R, \, u_\xi^R, \, \theta_\xi^R) \|_{L^4(\mathbb{R})}^4 \right).
	\end{align*}
	By Lemma \ref{ARV-pro},
	\begin{align*}
		\| (v_{\xi\xi}^R, \, u_{\xi\xi}^R, \, \theta_{\xi\xi}^R) \|_{L^2(\mathbb{R})} \le 
		\begin{cases}
			C \delta_R, & \text{if } 1+t \le \delta_R^{-1}, \\[2mm]
			C (1+t)^{-1}, & \text{if } 1+t \ge \delta_R^{-1},
		\end{cases}
	\end{align*}
	and
	\begin{align*}
		\| (v_\xi^R, \, u_\xi^R, \, \theta_\xi^R) \|_{L^4(\mathbb{R})} \le 
		\begin{cases}
			C \delta_R, & \text{if } 1+t \le \delta_R^{-1}, \\[2mm]
			C \delta_R^{1/4} (1+t)^{-3/4}, & \text{if } 1+t \ge \delta_R^{-1}.
		\end{cases}
	\end{align*}
	Then 
	\begin{align*}  
		\int_{0}^{+\infty} \| (Q_1^R, \, Q_2^R) \|_{L^2(\mathbb{R})}^2 dt \le C \delta_R.
	\end{align*}
	Together with \reff{Q12C-Est} and \reff{Q12I-Est} defined in Lemma \ref{WI-est},
	\begin{align*}
		\int_0^{+\infty} \| (Q_1^C, Q_2^C) \|_{L^2(\mathbb{R})}^2 dt \le C \delta_C^2, \quad \int_0^{+\infty} \| (Q_1^I, Q_2^I) \|_{L^2(\mathbb{R})}^2 dt \le C \delta_0^3,
	\end{align*}
	which yields
	\begin{align}  \label{Q12-Est-all}
		\int_0^{+\infty} \| (Q_1, Q_2) \|_{L^2(\mathbb{R})}^2 dt \le C \delta_0.
	\end{align}
	Integrating \reff{v-xi2} over $[0, t]$ and using the above computation, we have
	\begin{align} \label{v-H1-est}
		\begin{array}{l}
			\displaystyle \| (v - \bar{v})_\xi \|_{L^2(\mathbb{R})}^2 + \int_0^t \| (v - \bar{v})_\xi \|_{L^2(\mathbb{R})}^2 d\tau \\[4mm]
			\displaystyle \le C \left( \| (v - \bar{v})_\xi(t=0) \|_{L^2(\mathbb{R})}^2 + \| (u - \bar{u})(t=0) \|_{L^2(\mathbb{R})}^2 \right) + C \left( \| u - \bar{u} \|_{L^2(\mathbb{R})}^2 + \int_0^t \mathbf{D} d\tau \right) \\[4mm]
			\displaystyle \qquad + C (\delta_0^{1/2} + \varepsilon_1) \int_0^t (\delta_S |\dot{\mathbf{X}}(\tau)|^2 + \mathcal{G}^S + \mathcal{G}^R + \mathcal{G}^C + \| \tilde{\phi}_\xi \|_{L^2(\mathbb{R})}^2) d\tau + C \delta_0. 
		\end{array}
	\end{align}
	\hspace*{\fill}\\

	\noindent $\bullet$ {\bf $H^1$ estimates on $u-\bar{u}$.}
	\hspace*{\fill}\\ 

	Multiplying $\reff{PS}_2$ by $-(u - \bar{u})_{\xi\xi}$ and integrating by parts over $\mathbb{R}$ implies that
	\begin{align}  \label{u-par-est}
		\begin{array}{l}
			\displaystyle \dfrac{d}{dt} \int_{\mathbb{R}} \dfrac{1}{2}|(u - \bar{u})_\xi|^2 d\xi + \mu \int_{\mathbb{R}} \dfrac{|(u - \bar{u})_{\xi\xi}|^2}{v} d\xi \\[4mm]
			\displaystyle = \int_{\mathbb{R}} \bigg\{ - \dot{\mathbf{X}}(t) u_\xi^S (u - \bar{u})_{\xi\xi} + (\hat{p}(v, \theta) - \hat{p}(\bar{v}, \bar{\theta}))_\xi (u - \bar{u})_{\xi\xi} \\[4mm]
			\displaystyle \qquad \quad + \dfrac{(u - \bar{u})_\xi v_\xi (u - \bar{u})_{\xi\xi}}{v^2} - \left[ \bar{u}_\xi \left( \dfrac{1}{v} - \dfrac{1}{\bar{v}} \right) \right]_\xi (u - \bar{u})_{\xi\xi} \\[4mm]
			\displaystyle \qquad \quad - (\Phi(v, \phi) - \Phi(\bar{v}, \bar{\phi})) (u - \bar{u})_{\xi\xi} + (u - \bar{u})_{\xi\xi} Q_1 \bigg\} d\xi.
		\end{array}
	\end{align}
	Next, we only estimate some specific terms. Using the Young inequality and Lemma \ref{VS-pro},
	\begin{align*}
		\begin{array}{l}
			\displaystyle \int_{\mathbb{R}} \dot{\mathbf{X}}(t) u_\xi^S (u - \bar{u})_{\xi\xi} d\xi \le C |\dot{\mathbf{X}}(t)| \int_{\mathbb{R}} |v_\xi^S| \cdot |(u - \bar{u})_{\xi\xi}| d\xi \\[4mm]
			\displaystyle \le C \delta_S^2 |\dot{\mathbf{X}}(t)|^2 + \dfrac{C}{\delta_S^2} \left( \int_{\mathbb{R}} |v_\xi^S| \cdot |(u - \bar{u})_{\xi\xi}| d\xi \right)^2 \\[4mm]
			\displaystyle \le C \delta_S^2 |\dot{\mathbf{X}}(t)|^2 + \dfrac{C}{\delta_S^2} \| (u - \bar{u})_{\xi\xi} \|_{L^2(\mathbb{R})}^2 \cdot \| v_\xi^S \|_{L^\infty(\mathbb{R})} \cdot \| v_\xi^S \|_{L^1(\mathbb{R})} \\[4mm]
			\displaystyle \le C \delta_S \left( \delta_S |\dot{\mathbf{X}}(t)|^2 + \| (u - \bar{u})_{\xi\xi} \|_{L^2(\mathbb{R})}^2 \right).
		\end{array}
	\end{align*}
	Similarly, 
	\begin{align*}
		\begin{array}{l}
			\displaystyle \int_{\mathbb{R}} (\hat{p}(v, \theta) - \hat{p}(\bar{v}, \bar{\theta}))_\xi (u - \bar{u})_{\xi\xi} \, d\xi \\[4mm]
			\displaystyle = \int_{\mathbb{R}} \left[ \dfrac{R(\theta - \bar{\theta})_\xi}{v} + R \bar{\theta}_\xi \left( \dfrac{1}{v} - \dfrac{1}{\bar{v}} \right) - \hat{p}(v, \theta) \dfrac{(v - \bar{v})_\xi}{v} - \bar{v}_\xi \left( \dfrac{\hat{p}(v, \theta)}{v} - \dfrac{\hat{p}(\bar{v}, \bar{\theta})}{\bar{v}} \right) \right] (u - \bar{u})_{\xi\xi} \, d\xi \\[4mm]
			\displaystyle \le \dfrac{\mu}{4} \int_{\mathbb{R}} \dfrac{|(u - \bar{u})_{\xi\xi}|^2}{v} \, d\xi + C \mathbf{D} + C \delta_0^{1/2} (\mathcal{G}^S + \mathcal{G}^R + \mathcal{G}^C).
		\end{array}
	\end{align*}
	Likewise,
	\begin{align*}
		\begin{array}{l}
			\displaystyle - \int_{\mathbb{R}} (\Phi(v, \phi) - \Phi(\bar{v}, \bar{\phi})) (u - \bar{u})_{\xi\xi} d\xi \\[4mm]
			\displaystyle = - \int_{\mathbb{R}} \bigg[ \dfrac{\nu^2}{2} \left( \dfrac{\phi_\xi^2 - \bar{\phi}_\xi^2}{v^2} + \bar{\phi}_\xi^2 \left( \dfrac{1}{v^2} - \dfrac{1}{\bar{v}^2} \right) \right) -\nu^2 \Big( \dfrac{\bar{\phi}_{\xi\xi}}{v^2} + \bar{\phi}_{\xi\xi} \big( \dfrac{1}{v^2} - \dfrac{1}{\bar{v}^2} \big) \\[4mm]
			\displaystyle \qquad \qquad \qquad - \tilde{\phi}_\xi \cdot \dfrac{v_\xi}{v^3} - \bar{\phi}_\xi \big( \dfrac{v_\xi}{v^3} - \dfrac{\bar{v}_\xi}{\bar{v}^3} \big) \Big) \bigg] (u - \bar{u})_{\xi\xi} d\xi \\[4mm]
			\displaystyle \le C \int_{\mathbb{R}} \bigg[ |\tilde{\phi}_\xi|^2 + |\bar{\phi}_\xi| \cdot |\tilde{\phi}_\xi| + |\bar{\phi}_\xi|^2 \cdot |v - \bar{v}| + |\tilde{\phi}_{\xi\xi}| + |\bar{\phi}_{\xi\xi}| \cdot |v - \bar{v}| \\[4mm]
			\displaystyle \qquad + |\tilde{\phi}_\xi| \cdot |v_\xi| + |\bar{\phi}_\xi| (|(v - \bar{v})_\xi| + |\bar{v}_\xi| \cdot |v - \bar{v}|) \bigg] |(u - \bar{u})_{\xi\xi}| d\xi \\[4mm]
			\displaystyle \le C (\delta_0^{1/2} + \varepsilon_1) \left( \| (u - \bar{u})_{\xi\xi} \|_{L^2(\mathbb{R})}^2 + \| \tilde{\phi}_\xi \|_{H^1(\mathbb{R})}^2 + \mathcal{G}^S + \mathcal{G}^R + \mathcal{G}^C + \mathbf{D} \right) d\xi,
		\end{array}
	\end{align*}
	and 
	\begin{align*}
		\displaystyle \int_{\mathbb{R}} Q_1 (u - \bar{u})_{\xi\xi} d\xi \le \dfrac{\mu}{8} \int_{\mathbb{R}} \dfrac{|(u - \bar{u})_{\xi\xi}|^2}{v} d\xi + C \| Q_1 \|_{L^2(\mathbb{R})}^2,
	\end{align*}
	together with \reff{Q12-Est-all}, and combining with \reff{u-par-est}, we have 
	\begin{align} \label{u-H1-est}
		\begin{array}{l}
			\displaystyle \| (u - \bar{u})_\xi \|_{L^2(\mathbb{R})}^2 + \int_0^t \| (u - \bar{u})_{\xi\xi} \|_{L^2(\mathbb{R})}^2 d\tau \\[4mm]
			\displaystyle \le \| (u - \bar{u})_\xi(t=0) \|_{L^2(\mathbb{R})}^2 + C \int_0^t \mathbf{D} d\tau \\[4mm]
			\displaystyle \quad + C (\delta_0^{\frac{1}{2}} + \varepsilon_1) \int_0^t (\delta_S |\dot{\mathbf{X}}(t)|^2 + \mathcal{G}^S + \mathcal{G}^C + \mathcal{G}^R + \| \tilde{\phi}_\xi \|_{H^1(\mathbb{R})}^2 + \| (v-\bar{v})_\xi \|_{L^2(\mathbb{R})}^2) d\tau + C \delta_0.
		\end{array}
	\end{align}
	\hspace*{\fill}\\

	\noindent $\bullet$ {\bf $H^1$ estimates on $\theta-\bar{\theta}$.}
	\hspace*{\fill}\\ 

	Similarly, multiplying $\reff{PS}_3$ by  $-( \theta - \bar{\theta} )_{\xi\xi}$ and integrating by parts over $\mathbb{R}$, it follows that
	\begin{align*}
		\begin{array}{l}
			\displaystyle \dfrac{d}{dt} \int_{\mathbb{R}} \dfrac{R}{2(\gamma - 1)}|(\theta - \bar{\theta})_\xi|^2 d\xi + \kappa \int_{\mathbb{R}} \dfrac{|(\theta - \bar{\theta})_{\xi\xi}|^2}{v} d\xi \\[4mm]
			\displaystyle = \int_{\mathbb{R}} \bigg\{ - \dot{\mathbf{X}}(t) \dfrac{R}{\gamma - 1} \theta_\xi^S (\theta - \bar{\theta})_{\xi\xi} + (p(v, \theta) u_\xi - p(\bar{v}, \bar{\theta}) \bar{u}_\xi) (\theta - \bar{\theta})_{\xi\xi} \\[4mm]
			\displaystyle \qquad \qquad + \dfrac{(\theta - \bar{\theta})_\xi v_\xi}{v^2} (\theta - \bar{\theta})_{\xi\xi} - \left[ \bar{\theta}_\xi \left( \dfrac{1}{v} - \dfrac{1}{\bar{v}} \right) \right]_\xi (\theta - \bar{\theta})_{\xi\xi} \\[4mm]
			\displaystyle \qquad \qquad - \mu \left( \dfrac{u_\xi^2}{v} - \dfrac{\bar{u}_\xi^2}{\bar{v}} \right) (\theta - \bar{\theta})_{\xi\xi} + (\theta - \bar{\theta})_{\xi\xi} Q_2\bigg\} d\xi.
		\end{array}
	\end{align*}
	All the terms on the right-hand side can be estimated in a way similar to the above $H^1$ estimate for $u - \bar{u}$.
	Hence we have
	\begin{align} \label{theta-H1-est}
		\begin{array}{l}
			\displaystyle \| (\theta - \bar{\theta})_\xi \|_{L^2(\mathbb{R})}^2 + \int_0^t \| (\theta - \bar{\theta})_{\xi\xi} \|_{L^2(\mathbb{R})}^2 d\tau \\[4mm]
			\displaystyle \le \| (\theta - \bar{\theta})_\xi(t=0) \|_{L^2(\mathbb{R})}^2 + C \int_0^t \mathbf{D} d\tau + C (\delta_0^{1/2} + \varepsilon_1) \int_0^t (\delta_S |\dot{\mathbf{X}}(t)|^2 + \mathcal{G}^S + \mathcal{G}^R + \mathcal{G}^C) d\tau + C \delta_0.
		\end{array}
	\end{align}
	combining with \reff{v-H1-est}, \reff{u-H1-est} and \reff{theta-H1-est}, and using \reff{phi_H2E}, Lemma \ref{Low-Est-lem} and the below Lemma \ref{Basic estimates}, we complete the proof of this lemma.
\end{proof}
\hspace*{\fill}

\begin{lemma} \label{Basic estimates}
	Under the assumption of Proposition \ref{A prior estimate}, there exists a generic constant $C > 0$ such that
	\begin{align*}
		\| \tilde{\phi} \|_{H^k(\mathbb{R})} \le C \| \tilde{v} \|_{H^{k-1}(\mathbb{R})} + C \delta_0, \quad k = 1, 2, 3,
	\end{align*}
	where $\tilde{v} := v - \bar{v}$ and $\tilde{\phi} := \phi - \bar{\phi}$.
\end{lemma}

\begin{proof}
	For coherence of the exposition, the proof is given in the Appendix C.
\end{proof}
\hspace*{\fill}

The following lemma is mainly used to prove the $H^2$ estimates for $(U-\bar{U})$.
Compared with the above proof of $H^1$ estimates for $U - \bar{U}$, we believe that the corresponding $H^2$ estimates can be obtained by the same method, since this method does not introduce any additional bad terms. 
Therefore, we will present the results without giving the proof. Meanwhile, the following lemma also indicates that Proposition \ref{A prior estimate} holds.
\\ 

\begin{lemma}
	There exists a generic constant $C>0$ such that 
	\begin{align*}
		\begin{array}{l}
			\displaystyle \sup_{t \in [0, T]} \left( \| (U - \bar{U}) \|_{H^2(\mathbb{R})}^2 + \| \tilde{\phi} \|_{H^3(\mathbb{R})}^2 \right) \\[4mm]
			\displaystyle + \int_0^t (\delta_S |\dot{\mathbf{X}}(\tau)|^2 + \mathcal{G}^S + \mathcal{G}^R + \mathcal{G}^C + \| (v - \bar{v})_\xi \|_{H^1(\mathbb{R})}^2 + \| (u - \bar{u}, \, \theta - \bar{\theta})_\xi \|_{H^2(\mathbb{R})}^2) d\tau \\[4mm]
			\displaystyle \le C \| (U - \bar{U})(t=0) \|_{L^2(\mathbb{R})}^2  + C \delta_0^{1/3}.
		\end{array}
	\end{align*}
\end{lemma}
\hspace*{\fill} \\

\section*{Appendix}
\hspace*{\fill}\\ 

\subsection*{Appendix A. Existence of viscous shock profile}

\begin{proof} [\bf Proof of Lemma \ref{VS-pro}]
	\vspace*{1.5em}
	The proof is divided into three steps. Since the solvability of \reff{eq:viscous_shock} is equivalent to that of \reff{eq:v_xi}-\reff{eq:phi_xi_xi}, we first prove the existence of viscous shock waves to the latter.
	And then, in the second and third steps, we give the estimates in \reff{exist_all} and \reff{exist_uv}-\reff{exist_theta_v}, respectively.   
	\\ 

	$\bullet$ \noindent{\bf Step 1 (Existence):} Define the new variables $Z=e^\phi, W=\phi_\xi$, it follows from \reff{eq:v_xi}-\reff{eq:phi_xi_xi} that 
	\begin{align}\label{eq:ZW}
		\left\{
		\begin{array}{l}
			v_\xi = - \dfrac{v}{\mu \sigma} \left[ \sigma^2 (v - v^*) + (p - p(v^*, {\theta^*})) - \dfrac{\nu^2}{2} \dfrac{W^2}{v^2} + Z - Z^* \right], \\[3mm]
			\theta_\xi = - \dfrac{\sigma v}{\kappa} \bigg[ \dfrac{R(\theta - {\theta^*})}{\gamma - 1} - \dfrac{\sigma^2}{2} (v - v^*)^2 + p(v^*, {\theta^*})(v - v^*) + \dfrac{\nu^2}{2} \dfrac{W^2}{v}  \\[1mm]
			\qquad\qquad\qquad- v (Z - Z^*) + (\ln Z - \ln Z^*) \bigg], \\[2mm]
			Z_\xi=ZW, \\[2mm]
			W_{\xi} = - \dfrac{v}{\nu^2} (1 - v Z) - \dfrac{W}{\mu \sigma} \left[ \sigma^2 (v - v^*) + (p - p(v^*, {\theta^*})) - \dfrac{\nu^2}{2} \dfrac{W^2}{v^2} + Z - Z^* \right],
		\end{array}
		\right.
	\end{align}
	with the far fields given by 
	\begin{align*} 
		(v, \theta, Z, W)(-\infty)=(v^*, \theta^*, \frac{1}{v^*}, 0), \qquad 
		(v, \theta, Z, W)(+\infty)=(v_+, \theta_+, \frac{1}{v_+}, 0).
	\end{align*} 
	We introduce $U=(v, \theta, Z, W), U^*=(v^*, \theta^*, \frac{1}{v^*}, 0), U_+=(v_+, \theta_+, \frac{1}{v_+}, 0)$ and rewrite \reff{eq:ZW} in the form $U_\xi=F(U)$, 
	where $F(U)$ denotes the vector fields on the right hand of \reff{eq:ZW}. 
	Motivated by Majda and Duan $et~ al.$, we take $(U, \tau)$ satisfying the following extended ODEs 
	\begin{align}\label{eq:ZW_tau}
		\left\{
		\begin{array}{l}
			v_\xi = - \dfrac{v}{\mu \tau} \left[ \tau^2 (v - v^*) + (p - p(v^*, {\theta^*})) - \dfrac{\nu^2}{2} \dfrac{W^2}{v^2} + Z - Z^* \right], \\[3mm]
			\theta_\xi = - \dfrac{\tau v}{\kappa} \bigg[ \dfrac{R(\theta - {\theta^*})}{\gamma - 1} - \dfrac{\tau^2}{2} (v - v^*)^2 + p(v^*, {\theta^*})(v - v^*) + \dfrac{\nu^2}{2} \dfrac{W^2}{v}  \\[1mm]
			\qquad\qquad\qquad- v (Z - Z^*) + (\ln Z - \ln Z^*) \bigg], \\[2mm]
			Z_\xi=ZW, \\[2mm]
			W_{\xi} = - \dfrac{v}{\nu^2} (1 - v Z) - \dfrac{W}{\mu \tau} \left[ \tau^2 (v - v^*) + (p - p(v^*, {\theta^*})) - \dfrac{\nu^2}{2} \dfrac{W^2}{v^2} + Z - Z^* \right], \\[3mm]
			\,\,\tau_\xi = 0.
		\end{array}
		\right.
	\end{align}
	For fixed $\tau=\sigma$ and from \reff{eq:ZW}, $(U^*, \sigma)$ and $(U_+, \sigma)$ are the only two critical points of \reff{eq:ZW_tau}.
	For every $\tau$, $(U^*, \tau)$ is always a critical point of \reff{eq:ZW_tau}. 
	Later, we will take $(U^*, \sigma^*)$ ($\sigma^*=\sigma_3(v^*, {\theta^*})$ defined in \reff{Eigen}) as a reference state 
	and construct the center manifold of \reff{eq:ZW_tau} near this reference state.
	Compute the Jacobian matrix of \reff{eq:ZW_tau} at $(U^*, \sigma^*)$ as

	\begin{align*}
	J &= \left( 
		\begin{array}{ccccc}
		-\dfrac{(\gamma-1) R {\theta^*} + 1}{\mu \sigma^* v^*} & -\dfrac{R}{\mu \sigma^*} & -\dfrac{v^*}{\mu \sigma^*} & 0 & 0 \\[4mm]
		-\dfrac{ R \sigma^* {\theta^*}}{\kappa} & -\dfrac{ R \sigma^* v^*}{\kappa(\gamma-1)} & 0 & 0 & 0 \\[4mm]
		0 & 0 & 0 & \dfrac{1}{v^*} & 0 \\[4mm]
		\dfrac{1}{\nu^2} & 0 & \dfrac{(v^*)^2}{\nu^2} & 0 & 0 \\[4mm]
		0 & 0 & 0 & 0 & 0
		\end{array} 
		\right).
	\end{align*}

	The matrix $J$ has five eigenvalues, among which $\chi=0$ is a double eigenvalue and the remaining three are not purely imaginary. 	
	Indeed, the characteristic function of the $4 \times 4$ square submatrix $J_{ul}$ in the upper-left corner of $J$ is 
	\begin{align*}
		|J_{ul} - \chi I| = - \frac{v^*}{\chi} \Bigg\{ 
		& -v^* \chi^3 - v^* \left( \frac{(\gamma - 1) R {\theta^*} + 1}{\mu \sigma^* v^*} + \frac{R\sigma^* v^*}{\kappa(\gamma - 1)} \right) \chi^2  \\
		& + \left( - \frac{R v^*}{\mu \kappa (\gamma - 1)} + \frac{(v^*)^2}{\nu^2} \right) \chi  
		  + \frac{(v^*)^2}{\nu^2} \left( \frac{(\gamma - 1) R {\theta^*}}{\mu \sigma^* v^*} + \frac{R\sigma^* v^*}{\kappa (\gamma - 1)} \right) \Bigg\}.
	\end{align*}
	Since the characteristic equation has real coefficients, its complex roots, if they exists, must occur in conjugate pairs.
	Therefore, plugging the conjugated purely imaginary roots into the above equation yields 
	\begin{align*}
		- v^* \left( \frac{(\gamma - 1) R {\theta^*} + 1}{\mu \sigma^* v^*} + \frac{R\sigma^* v^*}{\kappa(\gamma - 1)} \right) \chi^2
		+ \frac{(v^*)^2}{\nu^2} \left( \frac{(\gamma - 1) R {\theta^*}}{\mu \sigma^* v^*} + \frac{R\sigma^* v^*}{\kappa (\gamma - 1)} \right) = 0,
	\end{align*}
	which only has real roots and thus contradicts the hypotheses.

	The two eigenvectors corresponding to the zero eigenvalue is given by 
	\begin{align*}
		\mathcal{R}_1=\left(1, \,\,-\frac{(\gamma-1){\theta^*}}{v^*}, \,\,-\frac{1}{(v^*)^2}, \,\,0, \,\,0\right),
		\quad \mathcal{R}_2=(0, \,0, \,0, \,0, \,1).
	\end{align*}
	Thus, by the $\mathit{Center ~ Manifold ~ Theorem}$ (e.g. Proposition 3.2 in Majda \cite{Majda85}), we obtain a 2-d manifold $M^*$ which is invariant by the flow \reff{eq:ZW_tau} 
	\begin{align} \label{Center_man}
		\left\{
		\begin{array}{l}
			U=U(\eta, \tau)=U^* + \eta \left(1, \,\,-\dfrac{(\gamma-1){\theta^*}}{v^*}, \,\,-\dfrac{1}{(v^*)^2}, \,\,0 \right) + \mathcal{H}(\eta, \tau), \\
			\tau=\tau
		\end{array}
		\right.
	\end{align}
	with $\eta^2+(\tau-\sigma^*)^2\le c^2$ for some small constant $c>0$. The vector-valued function 
	$$ \mathcal{H}(\eta, \tau) = ( \mathcal{H}^1, \mathcal{H}^2, \mathcal{H}^3, \mathcal{H}^4 )(\eta, \tau) $$ 
	satisfies the tangency condition 
	\begin{align*} 
		\mathcal{H}(0, \sigma^*)=\mathcal{H}_\eta(0, \sigma^*)=\mathcal{H}_\tau(0, \sigma^*)=0,
	\end{align*}
	and thus represents the higher-order terms near the fixed point $(0, \sigma^*)$, 
	namely $O(\eta^2+(\tau-\sigma^*)^2)$ by the Taylor expansion for $\eta^2+(\tau-\sigma^*)^2\le c^2$.
	In addition, since the varying critical points $(U^*, \tau)$ of \reff{eq:ZW_tau} lie in $M^*$ (e.g. Proposition 3.2 (3) in Majda), 
	we obtain 
	\begin{align} \label{H_tau}
		\mathcal{H}(\eta, \tau)=O(\eta^2)
	\end{align}
	near the moving points $(0, \tau)$ for $|\tau-\sigma^*|\le \dfrac{c}{2}$.
	On the other hand, there exists a sufficiently small constant $\delta>0$ such that if $\delta_S=|v_+ - v^*|<\delta$ and $|\sigma-\sigma^*|\le \dfrac{c}{2}$,
	then two critical points $(U^*, \sigma)$ and $(U_+, \sigma)$ are located in $M*$. 
	This ensures the existence of $\eta_+>0$ (determined by \reff{v_symbol}) satisfying
	$$ U_+ = U^* + \eta_+ \left(1, \,-\frac{(\gamma-1){\theta^*}}{v^*}, \,-\frac{1}{(v^*)^2}, \,0 \right) + \mathcal{H}(\eta_+, \sigma).$$
	Define the curve 
	\begin{align} \label{U_sol}
		U(\eta) = (v(\eta), \theta(\eta), Z(\eta), W(\eta)):=U(\eta, \tau)|_{\tau=\sigma}.
	\end{align}
	Note that $U(\eta)$ is invariant by the flow \reff{eq:ZW}, since $(U(\eta, \tau), \tau)$ is invariant by the flow \reff{eq:ZW_tau}.
	Therefore, the ODEs \reff{eq:ZW} admits a smooth solution if and only if the associated ODEs 
	\begin{align} \label{eta_xi}
		\left\{
		\begin{array}{l}
			\dfrac{d\eta}{d\xi} = g(\eta), \qquad \forall \xi\in\mathbb{R}, \\[2mm]
			\underset{\xi \rightarrow -\infty}{\lim} \eta(\xi) = \eta^*=0,\quad \underset{\xi \rightarrow +\infty}{\lim} \eta(\xi) = \eta_+ > 0
		\end{array}
		\right.
	\end{align}
	deduces a well-defined trajectory $\eta=\eta(\xi)$ for all $\xi\in\mathbb{R}$, where 
	\begin{align} \label{g_eta}
		g(\eta)U_\eta(\eta) = F(U(\eta)).
	\end{align}

	In what follows, we concentrate on the existence of the solution to \reff{eta_xi}-\reff{g_eta}.
	Since $U^*$ and $U_+$ are the only two critical points of \reff{eq:ZW}, i.e., $F(U^*)=0$ and $F(U_+)=0$, and 
	\begin{align*}
		U_\eta(\eta)|_{\eta=\eta_\pm}= \left(1, \,-\dfrac{(\gamma-1){\theta^*}}{v^*}, \,-\dfrac{1}{(v^*)^2}, \,0\right) + \mathcal{H}_\eta(\eta_\pm, \sigma) \ne 0,
	\end{align*}
	it follows from \reff{g_eta} that $g(\eta)$ vanishes only at two end points of $[0, \eta_+]$.
	The standard ODE theory shows that \reff{eta_xi}-\reff{g_eta} has a smooth solution if and only if $g(\eta)>0$ for all $\eta\in(0,\eta_+)$.
	We then determine the sign of $g(\eta)$ around $\eta^* = 0$. 
	From the first component of \reff{g_eta}, the definition \reff{U_sol} and \reff{Center_man}, and the shock speed \reff{sig-def},
	it holds that
	\begin{align} \label{g_eta_0}
		\begin{array}{l}
			g(\eta)=\dfrac{v^*}{\mu \sigma^*}\left[
			\dfrac{\gamma R {\theta^*} + 1}{ v^* ( v^* + \eta_+ ) ( v^* + \eta ) } (\eta_+ - \eta) + \dfrac{1}{ (v^*)^2 ( v^* + \eta_+ ) } \eta_+
			\right]  \eta + O (\eta^2) > 0, 
		\end{array}
	\end{align}
	where we use \reff{H_tau} near the point $(0, \sigma)$ for $|\sigma-\sigma^*| \le \dfrac{c}{2}$.
	Consequently, \reff{eta_xi}-\reff{g_eta} has a smooth solution. Let 
	\begin{align*}
		({v^S}, u^S, {\theta^S}, \phi^S)(\xi) := \big(v,\,\, (u^*+\sigma v^*)-\sigma v,\,\, \theta,\,\, \ln Z \big)(\eta(\xi)).
	\end{align*}
	Here $v(.), \theta(.)$ and $Z(.)$ is defined in \reff{U_sol}, and the definition of $u^S$ is due to \reff{u_relation}. 
	It is easy to check that $({v^S}, u^S, {\theta^S}, \phi^S)$ solves \reff{eq:viscous_shock}-\reff{far-field}.
	As for the uniqueness,  it is attributed to the uniqueness of the solution to the ODEs \eqref{eq:ZW} with the same far-field conditions and its proof is omitted for brevity.
	\\ 

	$\bullet$ \noindent{\bf Step 2 (Derivation of \reff{exist_all}):}
	First, the monotonicity is given by 
	\begin{align*}
		\begin{array}{l}
			v^S_\xi = \dfrac{dv}{d\eta} \dfrac{d\eta}{d\xi} = (1+\mathcal{H}^1_\eta(\eta(\xi), \sigma))\dfrac{d\eta}{d\xi} \sim \dfrac{d\eta}{d\xi} > 0,
			\\[3mm]
			u^S_\xi = -\sigma v^S_\xi \sim -\dfrac{d\eta}{d\xi} < 0, 
			\\[3mm]
			\theta^S_\xi = \dfrac{d\theta}{d\xi} \dfrac{d\eta}{d\xi} = \left(-\dfrac{{(\gamma-1){\theta^*}}}{v^*} + \mathcal{H}^2_\eta(\eta(\xi), \sigma) \right) \dfrac{d\eta}{d\xi} \sim -\dfrac{d\eta}{d\xi} < 0,
			\\[3mm]
			\phi^S_\xi = \dfrac{1}{Z} \dfrac{dZ}{d\eta} \dfrac{d\eta}{d\xi} = \dfrac{1}{e^\phi}\left( -\dfrac{1}{(v^*)^2} +\mathcal{H}^3_\eta(\eta(\xi),\sigma) \right) \dfrac{d\eta}{d\xi} \sim -\dfrac{d\eta}{d\xi} <0,
		\end{array}
	\end{align*}
	which yields the equivalence 
	\begin{align*}
		|v^S_\xi| \sim |u^S_\xi| \sim |\theta^S_\xi| \sim |\phi^S_\xi|, \quad \forall \xi \in \mathbb{R}.
	\end{align*}
	It follows from \reff{eta_xi} and \reff{g_eta_0} that 
	\begin{align*}
		v^S_\xi(\xi=0)=C\delta_S^2,
	\end{align*} 
	and 
	$|\eta(\xi)| \le C\delta_S e^{-C\delta_S |\xi|}$,
	which implies  
	\begin{align*}
		\left| ({v^S} - v^*, u^S - u^*, {\theta^S} - {\theta^*}, \phi^S - \phi^* )\right| \le C \delta_S e^{-C \delta_S |\xi|}, \qquad  \xi<0.
	\end{align*}
	Similarly, we can also derive 
	\begin{align*}
		\left| ({v^S} - v_+, u^S - u_+, {\theta^S} - \theta_+, \phi^S - \phi_+ )\right| \le C \delta_S e^{-C \delta_S |\xi|}, \qquad  \xi>0.
	\end{align*}
	The remaining estimates for the higher-order derivatives of $({v^S}, u^S, {\theta^S}, \phi^S)$ can be similarly obtained and it is omitted for brevity.
	\\ 

	$\bullet$ \noindent{\bf Step 3 (Derivation of \reff{exist_uv}-\reff{exist_theta_v}):} 
	The estimate \reff{exist_uv} is easily shown by \reff{eq:shock_a} and \reff{sig_sig-},
	since $({v^S}, u^S, {\theta^S}, \phi^S)$ solves \reff{eq:viscous_shock}-\reff{far-field}.
	To prove \reff{exist_theta_v}, from the Poisson equation \reff{eq:shock_d}, we derive the following relation 
	\begin{align*}
		\frac{\phi^S_\xi}{{v^S}} &= - \left[ \frac{\nu^2}{2} \left(\frac{\phi^S_\xi}{{v^S}}\right)^2 - e^{\phi^S} \right]_\xi = - \left[ \frac{\nu^2}{2} \left(\frac{\phi^S_\xi}{{v^S}}\right)^2 - \frac{1}{{v^S}} - \frac{\nu^2}{{v^S}} \left(\frac{\phi^S_\xi}{{v^S}}\right)_\xi \right]_\xi. 
	\end{align*}
	Substituting it and \reff{eq:shock_a} into \reff{eq:shock_b} gives 
	\begin{align*}
		\sigma^2 v^S_\xi + \left(p({v^S}, {\theta^S}) + \frac{1}{{v^S}}\right)_\xi &= \mu \left(\frac{u^S_\xi}{{v^S}}\right)_\xi + \left[ \frac{\nu^2}{2} \left(\frac{\phi^S_\xi}{{v^S}}\right)^2 - \frac{\nu^2}{{v^S}} \left(\frac{\phi^S_\xi}{{v^S}}\right)_\xi \right]_\xi,
	\end{align*}
	which together with the second and sixth estimates of \reff{exist_all}, and \reff{sig_sig-} yields 
	\begin{align*}
		\left|\theta^S_\xi + \dfrac{(\gamma-1){\theta^*}}{v^*} v^S_\xi \right| \le C\delta_S |v^S_\xi|,
	\end{align*}
	i.e., \reff{exist_theta_v}. Therefore, we complete the proof. 
\end{proof}
\hspace*{\fill} 

\subsection*{Appendix B. Sharp estimates for the Diffusion}

\begin{proof} [\bf Proof of Lemma \ref{Crucial}]
	\vspace*{1.5em}
	This proof is separated into two steps. The first step is to get the equivalent form of 
	\begin{align} \label{hat_p_taylor}
		\frac{\hat{p} - \hat{p}_+}{{v^S} - v_+} - \frac{\hat{p} - \hat{p}^*}{{v^S} - v^*}
	\end{align}
	up to an error of order $\delta_S^2$. The second step is to derive the critical term $\dfrac{d^2{\theta^S}}{d(v^S)^2}\bigg|_{v^*}$ in this equivalent form.
	\\ 

	$\bullet$ \noindent{\bf Step 1 (Equivalent form of \reff{hat_p_taylor}):}
	First of all, from \reff{exist_theta_v} and $v^S_\xi>0$ in Lemma \ref{VS-pro}, we obtain
	\begin{align*}
		\Big| ({\theta^S} - \theta_+) 
		&+\frac{(\gamma-1)p(v^*, {\theta^*})}{R} ({v^S} - v_+) \Big| 
		\leq \int_\xi^{+\infty} \left| \theta^S_\xi + \frac{(\gamma-1)p({v^S}, \theta)}{R} v^S_\xi \right| d\xi \\
		&\leq \int_0^{+\infty} C \delta_S |v^S_\xi| d\xi = C \delta_S |{v^S} - v_+|,
	\end{align*}
	which implies
	\begin{align} \label{theta_t+}
		\frac{{\theta^S} - \theta_+}{{v^S} - v_+} = - \frac{(\gamma-1)p(v^*,  {\theta^*})}{R} + O(\delta_S). 
	\end{align}
	Similarly,
	\begin{align} \label{theta_t-}
		\frac{{\theta^S} - {\theta^*}}{{v^S} - v^*} = - \frac{(\gamma-1)p(v^*, {\theta^*})}{R} + O(\delta_S). 
	\end{align}
	
	Then by the Taylor expansion with $\dfrac{\partial^2 \hat{p}}{\partial {\theta^S}^2} = 0$, and \reff{theta_t+} and \reff{theta_t-}, it holds
	\begin{align*}
		\hat{p} - \hat{p}^* 
		&=\frac{\partial \hat{p}}{\partial {v^S}} \bigg|_{(v^*, \theta^*)} ({v^S} - v^*) + \frac{\partial \hat{p}}{\partial {\theta^S}} \bigg|_{(v^*, \theta^*)} ({\theta^S} - \theta^*) 
		 +\frac{1}{2} \frac{\partial^2 \hat{p}}{\partial (v^S)^2} \bigg|_{(v^*, \theta^*)} ({v^S} - v^*)^2 \\
		&\quad + \frac{\partial^2 \hat{p}}{\partial {v^S} \partial {\theta^S}} \bigg|_{(v^*, \theta^*)} ({v^S} - v^*)({\theta^S} - \theta^*) + O(\delta_S^3),
	\end{align*}
	and 
	\begin{align*}
		\hat{p} - \hat{p}_+ 
		&=\frac{\partial \hat{p}}{\partial {v^S}} \bigg|_{(v_+, \theta_+)} ({v^S} - v_+) + \frac{\partial \hat{p}}{\partial {\theta^S}} \bigg|_{(v_+, \theta_+)} ({\theta^S} - \theta_+) 
		 +\frac{1}{2} \frac{\partial^2 \hat{p}}{\partial (v^S)^2} \bigg|_{(v_+, \theta_+)} ({v^S} - v_+)^2 \\
		&\quad + \frac{\partial^2 \hat{p}}{\partial {v^S} \partial {\theta^S}} \bigg|_{(v_+, \theta_+)} ({v^S} - v_+)({\theta^S} - \theta_+) + O(\delta_S^3),
	\end{align*}
	which yields
	\begin{align*}
		\frac{\hat{p} - \hat{p}_+}{{v^S} - v_+} 
		&- \frac{\hat{p} - \hat{p}^*}{{v^S} - v^*} 
		= \left( \frac{\partial \hat{p}}{\partial {v^S}} \bigg|_{(v_+, \theta_+)} - \frac{\partial \hat{p}}{\partial {v^S}} \bigg|_{(v^*, {\theta^*})} \right) \\
		&\quad + \left( \frac{\partial \hat{p}}{\partial {\theta^S}} \bigg|_{(v_+, \theta_+)} \frac{{\theta^S} - \theta_+}{{v^S} - v_+} - \frac{\partial \hat{p}}{\partial {\theta^S}} \bigg|_{(v^*, {\theta^*})} \frac{{\theta^S} - {\theta^*}}{{v^S} - v^*} \right) \\
		&\quad + \frac{1}{2} \left( \frac{\partial^2 \hat{p}}{\partial (v^S)^2} \bigg|_{(v_+, \theta_+)} ({v^S} - v_+) - \frac{\partial^2 \hat{p}}{\partial (v^S)^2} \bigg|_{(v^*, {\theta^*})} ({v^S} - v^*) \right) \\
		&\quad + \left( \frac{\partial^2 \hat{p}}{\partial {v^S} \partial {\theta^S}} \bigg|_{(v_+, \theta_+)} ({\theta^S} - \theta_+) - \frac{\partial^2 \hat{p}}{\partial {v^S} \partial {\theta^S}} \bigg|_{(v^*, {\theta^*})} ({\theta^S} - {\theta^*}) \right) 
		+ O(\delta_S^2) \\[2mm]
		&\quad := \mathcal{J}_1 + \mathcal{J}_2 + \mathcal{J}_3 + \mathcal{J}_4 + O(\delta_S^2).
	\end{align*}
	Since
	\begin{align*}
		\mathcal{J}_1 
		&= \frac{\partial^2 \hat{p}}{\partial (v^S)^2} \bigg|_{(v^*, {\theta^*})} (v_+ - v^*) + \frac{\partial^2 \hat{p}}{\partial {v^S} \partial {\theta^S}} \bigg|_{(v^*, {\theta^*})} (\theta_+ - {\theta^*}) + O(\delta_S^2). \\
		\mathcal{J}_3 
		&= - \frac{1}{2} \frac{\partial^2 \hat{p}}{\partial (v^S)^2} \bigg|_{(v^*, {\theta^*})} (v_+ - v^*) + \frac{1}{2} \left( \frac{\partial^2 \hat{p}}{\partial (v^S)^2} \bigg|_{(v_+, \theta_+)} - \frac{\partial^2 \hat{p}}{\partial (v^S)^2} \bigg|_{(v^*, {\theta^*})} \right) ({v^S} - v_+) \\
		&= - \frac{1}{2} \frac{\partial^2 \hat{p}}{\partial (v^S)^2} \bigg|_{(v^*, {\theta^*})} (v_+ - v^*) + O(\delta_S^2) \\
		\mathcal{J}_4 
		&= - \frac{\partial^2 \hat{p}}{\partial {v^S} \partial {\theta^S}} \bigg|_{(v^*, {\theta^*})} (\theta_+ - {\theta^*}) + \left( \frac{\partial^2 \hat{p}}{\partial {v^S} \partial {\theta^S}} \bigg|_{(v_+, \theta_+)} - \frac{\partial^2 \hat{p}}{\partial {v^S} \partial {\theta^S}} \bigg|_{(v^*, {\theta^*})} \right) ({\theta^S} - \theta_+) \\
		&= - \frac{\partial^2 \hat{p}}{\partial {v^S} \partial {\theta^S}} \bigg|_{(v^*, {\theta^*})} (\theta_+ - {\theta^*}) + O(\delta_S^2),
	\end{align*}
	then
	\begin{align*}
		\frac{\hat{p} - \hat{p}_+}{{v^S} - v_+} - \frac{\hat{p} - \hat{p}^*}{{v^S} - v^*} 
		&= \left( \frac{\partial \hat{p}}{\partial {\theta^S}} \bigg|_{(v_+, \theta_+)} \frac{{\theta^S} - \theta_+}{{v^S} - v_+} - \frac{\partial \hat{p}}{\partial {\theta^S}} \bigg|_{(v^*, {\theta^*})} \frac{{\theta^S} - {\theta^*}}{{v^S} - v^*} \right) \\
		&\quad + \frac{1}{2} \frac{\partial^2 \hat{p}}{\partial (v^S)^2} \bigg|_{(v^*, {\theta^*})} (v_+ - v^*) + O(\delta_S^2).
	\end{align*}
	Together with 
	$$\frac{\partial \hat{p}}{\partial {\theta^S}} = \frac{R}{{v^S}}, \quad \frac{\partial^2 \hat{p}}{\partial (v^S)^2} = \frac{2\hat{p}}{(v^S)^2},$$
	and \reff{theta_t+},
	we have
	\begin{align*}
		&\frac{\hat{p} - \hat{p}_+}{{v^S} - v_+} - \frac{\hat{p} - \hat{p}^*}{{v^S} - v^*} 
		= \left( \frac{R}{v_+} \frac{{\theta^S} - \theta_+}{{v^S} - v_+} - \frac{R}{v^*} \frac{{\theta^S} - {\theta^*}}{{v^S} - v^*} \right) + \frac{\hat{p}^*}{(v^*)^2} (v_+ - v^*) + O(\delta_S^2) \\
		&= \frac{{\theta^S} - \theta_+}{{v^S} - v_+} \left( \frac{R}{v_+} - \frac{R}{v^*} \right) + \frac{R}{v^*} \left( \frac{{\theta^S} - \theta_+}{{v^S} - v_+} - \frac{{\theta^S} - {\theta^*}}{{v^S} - v^*} \right) + \frac{\hat{p}^*}{(v^*)^2} (v_+ - v^*) + O(\delta_S^2) \\
		&= \frac{R}{v^*} \left( \frac{{\theta^S} - \theta_+}{{v^S} - v_+} - \frac{{\theta^S} - {\theta^*}}{{v^S} - v^*} \right) + \left( \gamma p(v^*, {\theta^*}) + \frac{1}{v^*} \right) \frac{1}{(v^*)^2}  (v_+ - v^*) + O(\delta_S^2).
	\end{align*}

	Now it remains to estimate
	$$\frac{{\theta^S} - \theta_+}{{v^S} - v_+} - \frac{{\theta^S} - {\theta^*}}{{v^S} - v^*}.$$
	Since $\theta^S_\xi < 0$ and $v^S_\xi > 0$ in Lemma \ref{VS-pro}, we can regard ${\theta^S}$ as a smooth function of ${v^S}$ by taking ${\theta^S} \circ {v^S}^{-1}: [v^*, v_+] \mapsto [\theta_+, {\theta^*}]$ with $\theta^* = \theta(v^*)$ and $\theta_+ = \theta(v_+)$. 
	Using the Taylor expansion at points $v_+$ and $v^*$, it holds
	$${\theta^S} - \theta^* = \frac{d{\theta^S}}{d{v^S}} \bigg|_{v^*} ({v^S} - v^*) + \frac{1}{2} \frac{d^2{\theta^S}}{d(v^S)^2} \bigg|_{v^*} ({v^S} - v^*)^2 + O(|{v^S} - v^*|^3),$$
	and 
	$${\theta^S} - \theta_+ = \frac{d{\theta^S}}{d{v^S}} \bigg|_{v_+} ({v^S} - v_+) + \frac{1}{2} \frac{d^2{\theta^S}}{d(v^S)^2} \bigg|_{v_+} ({v^S} - v_+)^2 + O(|{v^S} - v_+|^3),$$
	which gives
	\begin{align*}
		&\frac{{\theta^S} - \theta_+}{{v^S} - v_+} - \frac{{\theta^S} - {\theta^*}}{{v^S} - v^*} \\
		&= \left( \frac{d{\theta^S}}{d{v^S}}\bigg|_{v_+} - \frac{d{\theta^S}}{d{v^S}}\bigg|_{v^*} \right) + \frac{1}{2} \left( \frac{d^2{\theta^S}}{d(v^S)^2}\bigg|_{v_+} ({v^S} - v_+) - \frac{d^2{\theta^S}}{d(v^S)^2}\bigg|_{v^*} ({v^S} - v^*) \right) + O(\delta_S^2) \\
		&= \frac{1}{2} \frac{d^2{\theta^S}}{d(v^S)^2}\bigg|_{v^*} (v_+ - v^*) + \frac{1}{2} \left( \frac{d^2{\theta^S}}{d(v^S)^2}\bigg|_{v_+} - \frac{d^2{\theta^S}}{d(v^S)^2}\bigg|_{v^*} \right) ({v^S} - v_+) + O(\delta_S^2) \\
		&= \frac{1}{2} \frac{d^2{\theta^S}}{d(v^S)^2}\bigg|_{v^*} (v_+ - v^*) + O(\delta_S^2).
	\end{align*}
	Therefore, we have
	\begin{align} \label{equi_form_p}
		\frac{\hat{p} - \hat{p}_+}{{v^S} - v_+} - \frac{\hat{p} - \hat{p}^*}{{v^S} - v^*} 
		= \frac{R}{2v^*} \frac{d^2{\theta^S}}{d(v^S)^2}\bigg|_{v^*} (v_+ - v^*) + \left( \gamma p(v^*, {\theta^*}) + \frac{1}{v^*} \right) \frac{1}{(v^*)^2} (v_+ - v^*) + O(\delta_S^2). 
	\end{align}
	It remains to estimate $\dfrac{d^2{\theta^S}}{d(v^S)^2}\bigg|_{v^*}$.
	\\ 

	$\bullet$ \noindent{\bf Step 2 ( Derivation of $\dfrac{d^2{\theta^S}}{d(v^S)^2}\bigg|_{v^*}$ ):}
	Because $({v^S}, u^S, {\theta^S}, \phi^S)$ satisfies \reff{eq:viscous_shock}-\reff{far-field}, i.e., \reff{eq:v_xi}-\reff{eq:phi_xi_xi} and \reff{u_relation}, we derive
	\begin{align} \label{par_theta_v}
		\begin{array}{l}
			\dfrac{d{\theta^S}}{d{v^S}} = \dfrac{\theta^S_\xi}{v^S_\xi} \\
			= \dfrac{\mu \sigma^2}{\kappa} \dfrac{\frac{R({\theta^S} - {\theta^*})}{\gamma - 1} - \frac{\sigma^2}{2}({v^S} - v^*)^2 + p(v^*, {\theta^*})({v^S} - v^*) + \frac{\nu^2}{2}\frac{(\phi^S_\xi)^2}{{v^S}} - {v^S}(e^{\phi^S} - e^{\phi^*})+(\phi^S-\phi^*)}{\sigma^2({v^S} - v^*) + p({v^S}, {\theta^S}) - p(v^*, {\theta^*}) - \frac{\nu^2}{2}(\frac{\phi^S_\xi}{{v^S}})^2 + e^{\phi^S} - e^{\phi^*}}.
		\end{array}
	\end{align}
	Define
	\begin{align*}
		\begin{array}{l} 
			l^*:=\dfrac{d{\theta^S}}{d{v^S}}\bigg|_{v^*}= \underset{{v^S}\rightarrow v^*}{\lim} \dfrac{d{\theta^S}}{d{v^S}} = \underset{{v^S}\rightarrow v^*}{\lim} \dfrac{{\theta^S}-{\theta^*}}{{v^S}-v^*}, \\[3mm] 
			q^*:=\dfrac{d\phi^S}{d{v^S}}\bigg|_{v^*}= \underset{{v^S}\rightarrow v^*}{\lim} \dfrac{d\phi^S}{d{v^S}} = \underset{{v^S}\rightarrow v^*}{\lim} \dfrac{\phi^S-\phi^*}{{v^S}-v^*}.
		\end{array}
	\end{align*}
	
	In fact, by \reff{theta_t-}, 
	\begin{align} \label{theta_v}
		l^* := \underset{{v^S}\rightarrow v^*}{\lim} \frac{{\theta^S} - {\theta^*}}{{v^S} - v^*} = - \dfrac{(\gamma-1)p(v^*, {\theta^*})}{R} + O(\delta_S).
	\end{align}
	The Poisson equation \reff{eq:shock_d} gives 
	$$ e^{\phi^S} = \frac{1}{{v^S}} + \frac{\nu^2}{{v^S}} \left(\frac{\phi^S_\xi}{{v^S}}\right)_\xi. $$
	Then shifting the far-field value to the origin by the Taylor expansion, together with the second and fifth estimates of \reff{exist_all}, and $v^S_\xi(\xi=0) = C \delta_S^2$ in Lemma \ref{VS-pro}, we obtain 
	\begin{align} \label{phi_v}
		\begin{array}{ll}
			q^* 
			&= \dfrac{d\phi^S}{d{v^S}} \bigg|_{{v^S}(\xi=0)} + O(\delta_S) = \left( \dfrac{1}{e^{\phi^S}} \dfrac{d e^{\phi^S}}{d{v^S}} \right) \bigg|_{{v^S}(\xi=0)} + O(\delta_S) \\
			&= \left[ \dfrac{1}{e^{\phi^S}} \left( -\dfrac{1}{(v^S)^2} - \dfrac{\nu^2}{(v^S)^2} \left(\dfrac{\phi^S_\xi}{{v^S}}\right)_\xi + \dfrac{\nu^2}{{v^S}} \left(\dfrac{\phi^S_\xi}{{v^S}}\right)_{\xi\xi} \dfrac{1}{v^S_\xi} \right) \right] \Bigg|_{{v^S}(\xi=0)} + O(\delta_S) \\
		 	&= -\dfrac{1}{v^*} + O(\delta_S),
		\end{array}
	\end{align}
	and 
	\begin{align} \label{phi_v2}
		\begin{array}{ll}
			\dfrac{d^2 \phi^S}{d(v^S)^2} \bigg|_{v^*} 
			&= \dfrac{d^2 \phi^S}{d(v^S)^2} \bigg|_{{v^S}(\xi=0)} + O(\delta_S) = \left[ -\left(\dfrac{d\phi^S}{d{v^S}}\right) ^2 + \dfrac{1}{e^{\phi^S}} \dfrac{d^2 e^{\phi^S}}{d(v^S)^2} \right] \Bigg|_{{v^S}(\xi=0)} + O(\delta_S). \\[3mm]
			&= \left[ -\left(\dfrac{d\phi^S}{d{v^S}}\right)^2 + \dfrac{1}{e^{\phi^S}} \left( \dfrac{2}{(v^S)^3} + \dfrac{2\nu^2}{(v^S)^3} \left(\dfrac{\phi^S_\xi}{{v^S}}\right)_\xi - \dfrac{2\nu^2}{(v^S)^2} \left(\dfrac{\phi^S_\xi}{{v^S}}\right)_\xi \dfrac{1}{v^S_\xi} \right.\right. \\[3mm] 
			&\qquad\quad \left.\left. + \dfrac{\nu^2}{{v^S}} \left(\left(\dfrac{\phi^S_\xi}{{v^S}}\right)_{\xi\xi} \dfrac{1}{v^S_\xi} \right)_\xi \dfrac{1}{v^S_\xi} \right) \right] \Bigg|_{{v^S}(\xi=0)} + O(\delta_S) \\[5mm]
			&= \left[ -\left(\dfrac{d\phi^S}{d{v^S}}\right)^2 + \dfrac{2}{e^{\phi^S} (v^S)^3} + \dfrac{\nu^2}{{v^S}} \dfrac{\phi^S_{\xi\xi\xi\xi}}{{v^S} (v^S_\xi)^2} \right] \Bigg|_{{v^S}(\xi=0)} + O(\delta_S) \\[3mm]
			&= \dfrac{1}{(v^*)^2} + O(\delta_S).
		\end{array}
	\end{align}
	Furthermore, 
	\begin{align} \label{p_v}
		\begin{array}{ll}
			\underset{{v^S}\rightarrow v^*}{\lim} \dfrac{p({v^S}, {\theta^S}) - p(v^*, {\theta^*})}{{v^S} - v^*} 
			&=\underset{{v^S}\rightarrow v^*}{\lim} \dfrac{\frac{R{\theta^S}}{{v^S}} - \frac{R{\theta^*}}{v^*}}{{v^S} - v^*} 
			=\underset{{v^S}\rightarrow v^*}{\lim} \dfrac{\frac{R}{{v^S}}({\theta^S} - {\theta^*}) + R{\theta^*}(\frac{1}{{v^S}} - \frac{1}{v^*})}{{v^S} - v^*} \\[3mm]
			&=\underset{{v^S}\rightarrow v^*}{\lim} \left( \dfrac{R}{{v^S}} \dfrac{{\theta^S} - {\theta^*}}{{v^S} - v^*} - \dfrac{R{\theta^*}}{{v^S}v^*} \right) 
			= \dfrac{R l^*-p(v^*, {\theta^*})}{v^*},
		\end{array}
	\end{align}
	and 
	\begin{align} \label{e_phi_v}
		\begin{array}{l}
			\underset{{v^S}\rightarrow v^*}{\lim} \dfrac{(\phi^S_\xi)^2}{{v^S}} \dfrac{1}{{v^S} - v^*} = \underset{{v^S}\rightarrow v^*}{\lim} \left(\dfrac{(\phi^S_\xi)^2}{{v^S}}\right)_\xi \cdot \dfrac{1}{v^S_\xi} = 0, 
			\quad \underset{{v^S}\rightarrow v^*}{\lim} \left(\dfrac{\phi^S_\xi}{{v^S}}\right)^2 \dfrac{1}{{v^S} - v^*} = 0, \\[5mm]
			\underset{{v^S}\rightarrow v^*}{\lim} \dfrac{e^{\phi^S} - e^{\phi^*}}{{v^S} - v^*} = \underset{{v^S}\rightarrow v^*}{\lim} e^{\phi^S} \cdot \dfrac{d\phi^S}{d{v^S}} = \dfrac{q^*}{v^*} = -\dfrac{1}{(v^*)^2} + O(\delta_S).
		\end{array}
	\end{align}
	Thus, it follows from \reff{par_theta_v} that 
	\begin{align*}
		l^* 
		&= \underset{{v^S}\rightarrow v^*}{\lim} \frac{d{\theta^S}}{d{v^S}} \\
		&= \frac{\mu \sigma^2}{\kappa} \frac{\frac{R}{\gamma - 1} \underset{{v^S}\rightarrow v^*}{\lim} \frac{{\theta^S} - {\theta^*}}{{v^S} - v^*}  + p(v^*, {\theta^*}) - \underset{{v^S}\rightarrow v^*}{\lim} \frac{{v^S}(e^{\phi^S} - e^{\phi^*})}{{v^S} - v^*} + \underset{{v^S}\rightarrow v^*}{\lim} \frac{\phi^S - \phi^*}{{v^S} - v^*}}
		{\sigma^2 + \underset{{v^S}\rightarrow v^*}{\lim} \frac{p({v^S}, {\theta^S}) - p(v^*, {\theta^*})}{{v^S} - v^*} + \underset{{v^S}\rightarrow v^*}{\lim} \frac{e^{\phi^S} - e^{\phi^*}}{{v^S} - v^*}} \\
		&= \frac{\mu \sigma^2}{\kappa} \frac{\frac{R}{\gamma - 1} l^* + p(v^*, {\theta^*})}{\sigma^2 + \frac{R l^* - p(v^*, {\theta^*})}{v^*} + \frac{q^*}{v^*}},
	\end{align*}
	i.e.,
	\begin{align} \label{l-_satisfied}
		(l^*)^2 - \left( \frac{p(v^*, {\theta^*}) - \sigma^2 v^* - q^*}{R} + \frac{\mu \sigma^2 v^*}{\kappa(\gamma-1)} \right) l^* - \frac{\mu \sigma^2 v^*}{\kappa R} p(v^*, {\theta^*}) = 0.
	\end{align}
	
	With the above preparations, we will calculate the value of $\frac{d^2{\theta^S}}{d(v^S)^2}\big|_{v^*}$. From \reff{par_theta_v}, we have 
	\begin{align*}
		\frac{d^2{\theta^S}}{d(v^S)^2} 
		&= \frac{d}{d{v^S}} \left( \frac{d{\theta^S}}{d{v^S}} \right) \\ 
		&= \frac{\mu \sigma^2}{\kappa} \frac{\frac{R}{\gamma - 1} \frac{d{\theta^S}}{d{v^S}} - \sigma^2({v^S} - v^*) + p(v^*, {\theta^*}) + \frac{\nu^2}{2} \frac{d}{d{v^S}}(\frac{(\phi^S_\xi)^2}{{v^S}}) - \frac{d}{d{v^S}}({v^S}(e^{\phi^S} - e^{\phi^*})) + \frac{d\phi^S}{d{v^S}}}{\sigma^2({v^S} - v^*) + p({v^S}, {\theta^S}) - p(v^*, {\theta^*}) - \frac{\nu^2}{2}(\frac{\phi^S_\xi}{{v^S}})^2 + e^{\phi^S} - e^{\phi^*} } \\
		&\quad - \dfrac{d{\theta^S}}{d{v^S}} \frac{ \sigma^2 + \frac{d}{d{v^S}} p({v^S}, {\theta^S}) - \frac{\nu^2}{2} \frac{d}{d{v^S}} (\frac{\phi^S_\xi}{{v^S}})^2 + \frac{d}{d{v^S}} (e^{\phi^S} - e^{\phi^*}) }{ \sigma^2({v^S} - v^*) + p({v^S}, {\theta^S}) - p(v^*, {\theta^*}) - \frac{\nu^2}{2}(\frac{\phi^S_\xi}{{v^S}})^2 + e^{\phi^S} - e^{\phi^*} } \\
		&= - \frac{R}{{v^S} \left[ \sigma^2 ({v^S} - v^*) + p({v^S}, {\theta^S}) - p(v^*, {\theta^*}) - \frac{\nu^2}{2} (\frac{\phi^S_\xi}{{v^S}})^2 + e^{\phi^S} - e^{\phi^*} \right]} \\
		&\times \Bigg\{ \left( \frac{d{\theta^S}}{d{v^S}} \right)^2 - \left( \frac{p({v^S}, {\theta^S}) - \sigma^2 {v^S}}{R} + \frac{\nu^2{v^S}}{2R} \frac{d}{d{v^S}} (\frac{\phi^S_\xi}{{v^S}})^2 - \frac{{v^S}}{R} \frac{d}{d{v^S}} ( e^{\phi^S} - e^{\phi^*}) + \frac{\mu \sigma^2 {v^S}}{\kappa (\gamma-1)} \right) \frac{d{\theta^S}}{d{v^S}} \\
		&\qquad - \frac{\mu \sigma^2 {v^S}}{\kappa R} \left( - \sigma^2 ({v^S} - v^*) + p(v^*, {\theta^*}) + \frac{\nu^2}{2} \frac{d}{d{v^S}} (\frac{(\phi^S_\xi)^2}{{v^S}})- \frac{d}{d{v^S}} ({v^S} (e^{\phi^S} - e^{\phi^*})) + \frac{d\phi^S}{d{v^S}} \right) \\
		&\qquad - \left[ (l^*)^2 - \left( \frac{p(v^*, {\theta^*}) - \sigma^2 v^* - q^*}{R} + \frac{\mu \sigma^2 v^*}{\kappa(\gamma-1)} \right) l^* - \frac{\mu \sigma^2 v^*}{\kappa R} p(v^*, {\theta^*}) \right] \Bigg\}.
	\end{align*}
	Since
	\begin{align*}
		\left( \frac{d{\theta^S}}{d{v^S}} \right)^2 - (l^*)^2 = \left( \frac{d{\theta^S}}{d{v^S}} - l^* \right) \left( \frac{d{\theta^S}}{d{v^S}} + l^* \right),
	\end{align*}
	and 
	\begin{align*}
		&\left( \frac{p({v^S}, {\theta^S}) - \sigma^2 {v^S}}{R} + \frac{\nu^2{v^S}}{2R} \frac{d}{d{v^S}} (\frac{\phi^S_\xi}{{v^S}})^2 - \frac{{v^S}}{R} \frac{d}{d{v^S}} ( e^{\phi^S} - e^{\phi^*}) + \frac{\mu \sigma^2 {v^S}}{\kappa (\gamma-1)} \right) \frac{d{\theta^S}}{d{v^S}} \\
		&\qquad - \left( \frac{p(v^*, {\theta^*}) - \sigma^2 v^* - q^*}{R} + \frac{\mu \sigma^2 v^*}{\kappa (\gamma-1)} \right) l^* \\
		&= \left( \frac{ p({v^S}, {\theta^S}) - p(v^*, {\theta^*}) - \sigma^2 ({v^S} - v^*)}{R} + \frac{\nu^2{v^S}}{2R} \frac{d}{d{v^S}} (\frac{\phi^S_\xi}{{v^S}})^2 - \frac{{v^S} e^{\phi^S} - 1}{R} \frac{d\phi^S}{d{v^S}} - \frac{1}{R} ( \frac{d\phi^S}{d{v^S}} - q^* ) \right. \\
		&\qquad \left.  + \frac{\mu \sigma^2 ({v^S} - v^*)}{\kappa (\gamma-1)} \right) \frac{d{\theta^S}}{d{v^S}}  + \left( \frac{p(v^*, {\theta^*}) - \sigma^2 v^* - q^*}{R} + \frac{\mu \sigma^2 v^*}{\kappa (\gamma-1)} \right) \left( \frac{d{\theta^S}}{d{v^S}} - l^* \right), \\[4mm] 
	\end{align*}
	and 
	\begin{align*}
		&\frac{\mu \sigma^2 {v^S}}{\kappa R} \left( - \sigma^2 ({v^S} - v^*) + p(v^*, {\theta^*}) + \frac{\nu^2}{2} \frac{d}{d{v^S}} (\frac{(\phi^S_\xi)^2}{{v^S}})- \frac{d}{d{v^S}} ({v^S} (e^{\phi^S} - e^{\phi^*})) + \frac{d\phi^S}{d{v^S}} \right) - \frac{\mu \sigma^2 v^*}{\kappa R} p(v^*, {\theta^*}) \\
		&=\frac{\mu \sigma^2}{\kappa R} \left( (p(v^*, {\theta^*}) - \sigma^2 {v^S}) ({v^S} - v^*) + \frac{\nu^2 {v^S}}{2} \frac{d}{d{v^S}} (\frac{(\phi^S_\xi)^2}{{v^S}})- {v^S} (e^{\phi^S} - e^{\phi^*}) + {v^S} (1 - {v^S} e^{\phi^S}) \frac{d\phi^S}{d{v^S}} \right).
	\end{align*}
	Thus,
	\begin{align*}
		\frac{d^2{\theta^S}}{d(v^S)^2} 
		&= - \frac{R}{{v^S} [\sigma^2 ({v^S} - v^*) + p({v^S}, {\theta^S}) - p(v^*, {\theta^*}) - \frac{\nu^2}{2} (\frac{\phi^S_\xi}{{v^S}})^2 + e^{\phi^S} - e^{\phi^*}]} \\
		&\quad \times \Bigg\{ \left( \frac{d{\theta^S}}{d{v^S}} - l^* \right) \left( \frac{d{\theta^S}}{d{v^S}} + l^* \right) - \left( \frac{p(v, {\theta^*}) - \sigma^2 v^* - q^*}{R} + \frac{\mu \sigma^2 v^*}{\kappa(\gamma-1)} \right) \left( \frac{d{\theta^S}}{d{v^S}} - l^* \right) \\
		&\quad\qquad - \left( \frac{p({v^S}, {\theta^S}) - p(v^*, {\theta^*}) - \sigma^2 ({v^S} - v^*)}{R} + \frac{\nu^2{v^S}}{2R}  \frac{d}{d{v^S}} (\frac{\phi^S_\xi}{{v^S}})^2 - \frac{{v^S} e^{\phi^S} - 1}{R} \frac{d\phi^S}{d{v^S}} \right. \\
		&\quad\qquad\qquad \left. - \frac{1}{R} ( \frac{d\phi^S}{d{v^S}} - q^* ) + \frac{\mu \sigma^2 ({v^S} - v^*)}{\kappa (\gamma-1)} \right) \frac{d{\theta^S}}{d{v^S}} \\
		&\quad\qquad - \frac{\mu \sigma^2}{\kappa R} \left( (p(v^*, {\theta^*}) - \sigma^2 {v^S}) ({v^S} - v^*) + \frac{\nu^2 {v^S}}{2} \frac{d}{d{v^S}} (\frac{(\phi^S_\xi)^2}{{v^S}}) \right. \\
		&\qquad\qquad\qquad\qquad \left. - {v^S} (e^{\phi^S} - e^{\phi^*}) + {v^S} (1 - {v^S} e^{\phi^S}) \frac{d\phi^S}{d{v^S}} \right) \Bigg\}.
	\end{align*}
	Multiplying by ${v^S} \left[ \sigma^2 ({v^S} - v^*) + p({v^S}, {\theta^S}) - p(v^*, {\theta^*}) - \frac{\nu^2}{2} (\frac{\phi^S_\xi}{{v^S}})^2 + e^{\phi^S} - e^{\phi^*} \right]$ and
	dividing by ${v^S}-v^*$, then we get
	\begin{align*}
		&{v^S} \left[ \sigma^2 + \frac{p({v^S}, {\theta^S}) - p(v^*, {\theta^*})}{{v^S} - v^*} - \frac{\nu^2}{2} \frac{1}{{v^S} - v^*} (\frac{\phi^S_\xi}{{v^S}})^2 + \frac{e^{\phi^S} - e^{\phi^*}}{{v^S} - v^*} \right] \frac{d^2{\theta^S}}{d(v^S)^2} \\
		&= - R \Bigg\{ \frac{\frac{d{\theta^S}}{d{v^S}} - l^*}{{v^S}-v^*} \left( \frac{d{\theta^S}}{d{v^S}} + l^* \right) - \left( \frac{p(v^*, {\theta^*}) - \sigma^2 v^* - q^*}{R} + \frac{\mu \sigma^2 v^*}{\kappa(\gamma-1)} \right) \frac{\frac{d{\theta^S}}{d{v^S}} - l^*}{{v^S}-v^*} \\
		&\quad - \left( \frac{\frac{p({v^S}, {\theta^S}) - p(v^*, {\theta^*})}{{v^S} - v^*} - \sigma^2}{R} + \frac{\nu^2{v^S}}{2R}  \frac{1}{{v^S} - v^*} \frac{d}{d{v^S}} (\frac{\phi^S_\xi}{{v^S}})^2 \right.\\ 
		&\qquad\qquad \left. + \frac{1}{R} \frac{1 - {v^S} e^{\phi^S}}{{v^S} - v^*} \frac{d\phi^S}{d{v^S}} - \frac{1}{R} \frac{\frac{d\phi^S}{d{v^S}} - q^*}{{v^S}-v^*} + \frac{\mu \sigma^2}{\kappa(\gamma-1)} \right) \frac{d{\theta^S}}{d{v^S}} \\
		&\quad - \frac{\mu \sigma^2}{\kappa R} \left( (p(v^*, {\theta^*}) - \sigma^2 {v^S}) + \frac{\nu^2 {v^S}}{2} \frac{1}{{v^S} - v^*} \frac{d}{d{v^S}} (\frac{(\phi^S_\xi)^2}{{v^S}}) - {v^S} \frac{e^{\phi^S} - e^{\phi^*}}{{v^S} - v^*} + {v^S}\frac{ 1 - {v^S} e^{\phi^S}}{{v^S} - v^*} \frac{d\phi^S}{d{v^S}} \right) \Bigg\}.
	\end{align*}
	
	Due to \reff{p_v}-\reff{e_phi_v}, and $v^S_\xi(\xi=0) = C \delta_S^2$ in Lemma \ref{VS-pro}, we shift the far-field value to the origin (omitted for brevity) to obtain
	\begin{align*}
		\underset{{v^S}\rightarrow v^*}{\lim}\frac{1}{{v^S} - v^*} \frac{d}{d{v^S}} (\frac{\phi^S_\xi}{{v^S}})^2 
		&= \underset{{v^S}\rightarrow v^*}{\lim}\frac{d^2}{d(v^S)^2} (\frac{\phi^S_\xi}{{v^S}})^2 
		= \underset{{v^S}\rightarrow v^*}{\lim}\frac{d}{d{v^S}} \left[ \frac{2 \frac{\phi^S_\xi}{{v^S}} (\frac{\phi^S_\xi}{{v^S}})_\xi}{v^S_\xi} \right] \\
		&= 2 \underset{{v^S}\rightarrow v^*}{\lim}\left[ \frac{\frac{\phi^S_\xi}{{v^S}} (\frac{\phi^S_\xi}{{v^S}})_\xi}{v^S_\xi} \right]_\xi \frac{1}{v^S_\xi} = O(\delta_S), \\
	\end{align*}
	and 
	\begin{align*}
		&\underset{{v^S}\rightarrow v^*}{\lim}\frac{1 - {v^S} e^{\phi^S}}{{v^S} - v^*} = - \frac{d}{d{v^S}} ({v^S} e^{\phi^S}) \bigg|_{v^*} = -(e^{\phi^S} + {v^S} e^{\phi^S} \frac{d\phi^S}{d{v^S}}) \bigg|_{v^*} = O(\delta_S), \\[2mm]
		&\underset{{v^S}\rightarrow v^*}{\lim}\frac{\frac{d\phi^S}{d{v^S}} - q^*}{{v^S} - v^*} = \frac{d^2\phi^S}{d(v^S)^2} \bigg|_{v^*} = \frac{1}{(v^*)^2} + O(\delta_S).
	\end{align*}
	Then taking the limit ${v^S} \to v^*$ yields
	\begin{align*}
		&v^* \left[ \sigma^2 + \frac{R l^* - p(v^*, {\theta^*})}{v^*} - \frac{1}{(v^*)^2} + O(\delta_S) \right] \frac{d^2{\theta^S}}{d(v^S)^2}\bigg|_{v^*} \\
		&= - R \Bigg[ \left(2l^* - \frac{p(v^*, {\theta^*}) - \sigma^2 v^* + \frac{1}{v^*}}{R} - \frac{\mu \sigma^2 v^*}{\kappa(\gamma-1)} + O(\delta_S) \right) \frac{d^2{\theta^S}}{d(v^S)^2}\bigg|_{v^*} \\
		&\qquad\qquad - \left(\frac{1}{R}\big(\frac{R l^* - p(v^*, {\theta^*})}{v^*} - \sigma^2 - \frac{1}{(v^*)^2}\big) + \frac{\mu \sigma^2}{\kappa(\gamma-1)} + O(\delta_S)\right) \frac{d{\theta^S}}{d{v^S}}\bigg|_{v^*} \\
		&\qquad\qquad - \frac{\mu \sigma^2}{kR} \left(p(v^*, {\theta^*}) - \sigma^2 v^* + \frac{1}{v^*} + O(\delta_S) \right) \Bigg].
	\end{align*}
	Therefore, it follows from \reff{theta_v}-\reff{phi_v2} and \reff{sig_sig-} that 
	\begin{align} \label{theta_v_2}
		\frac{d^2{\theta^S}}{d(v^S)^2}\bigg|_{v^*} 
		= \frac{\dfrac{\gamma-1}{R v^*} \Big(\mu \gamma R - 2k(\gamma-1)\Big) p(v^*, {\theta^*}) \left( \gamma p(v^*, {\theta^*}) + \frac{1}{v^*} \right)}{ k(\gamma-1)^2 p(v^*, {\theta^*}) + R\mu \left( \gamma p(v^*, {\theta^*}) + \frac{1}{v^*} \right)} + O(\delta_S)
	\end{align}

	Applying \reff{theta_v_2} into \reff{equi_form_p}, we finally get
	\begin{align*}
		\frac{\hat{p} - \hat{p}_+}{{v^S} - v_+} - \frac{\hat{p} - \hat{p}^*}{{v^S} - v^*} 
		&= \frac{\left( \frac{\gamma(\gamma+1)}{2} p(v^*, {\theta^*}) + \frac{1}{v^*} \right) R\mu \left(\gamma p(v^*, {\theta^*}) + \frac{1}{v^*} \right)}{\kappa (\gamma-1)^2 p(v^*, {\theta^*}) + R\mu \left( \gamma p(v^*, {\theta^*}) + \frac{1}{v^*} \right)} \cdot \frac{1}{(v^*)^2} (v_+ - v^*) + O(\delta_S^2) \\
		&= \sigma^* \alpha^* \beta^* \delta_S + O(\delta_S^2),
	\end{align*}
	where
	\begin{align*}
		&\delta_S = v_+-v^*>0, \\ 
		&\sigma^* = \sqrt{\frac{\gamma R {\theta^*} + 1}{(v^*)^2}}, \\
		&\alpha^* = \left( \frac{\gamma(\gamma+1)p(v^*, {\theta^*})}{2} + \frac{1}{v^*}\right) \frac{1}{(v^*)^2 \sigma^*} \\,
		&\beta^* = \frac{R\mu \left(\gamma p(v^*, {\theta^*}) + \frac{1}{v^*} \right)}{\kappa(\gamma-1)^2 p(v^*, {\theta^*}) + R\mu \left( \gamma p(v^*, {\theta^*}) + \frac{1}{v^*} \right)}.
	\end{align*}
	Therefore, we complete the proof.
\end{proof}
\hspace*{\fill} 

\subsection*{Appendix C. Basic elliptic estimates on $\tilde{\phi}$}

\begin{proof}[\bf Proof of Lemma \ref{Basic estimates}]
	\vspace*{1.5em}
	The proof is given by the elemental energy method. Multiplying \reff{for basic} by $\tilde{\phi}$ and integrating over $\mathbb{R}$ yields
	\begin{align*}
		\int_{\mathbb{R}} \left( \bar{v} e^{\bar{\phi}} \tilde{\phi}^2 + \nu^2 \dfrac{\tilde{\phi}_\xi^2}{v} \right) d\xi = \int_{\mathbb{R}} \left[ - \nu^2 \tilde{\phi}_\xi \left( \dfrac{1}{v} - \dfrac{1}{\bar{v}} \right) \tilde{\phi}_\xi + \tilde{\phi} \left( - e^{\bar{\phi}} \bar{v} + e^{\bar{\phi}} N_p - Q_3 \right) \right] d\xi.
	\end{align*}
	From the assumption of Proposition \ref{A prior estimate} and \reff{Ass-Inf}, it follows $|N_p| \le C (|\tilde{\phi}| + |\tilde{\phi}| |\tilde{v}|)$. 
	By Lemma \ref{WI-est}, \reff{Q3R_Est} and \reff{Q3C_Est}, we have $\|Q_3\|_{L^2(\mathbb{R})} \le C\delta_0$. Thus, 
	\begin{align*}
		\begin{array}{l}
			\displaystyle \int_{\mathbb{R}} \left( \bar{v} e^{\bar{\phi}} \tilde{\phi}^2 + \nu^2 \dfrac{\tilde{\phi}_\xi^2}{v} \right) d\xi \\[4mm]
			\displaystyle \le C \int_{\mathbb{R}} \left( |\tilde{\phi}_\xi| \cdot |\tilde{v}| \cdot |\tilde{\phi}_\xi| + |\tilde{\phi}| \cdot (|\tilde{v}| + |\tilde{\phi}|^2 + |\tilde{\phi}| \cdot |\tilde{v}| + |Q_3|) \right) d\xi \\[4mm]
			\displaystyle \le \dfrac{1}{2} \int_{\mathbb{R}} \left( \bar{v} e^{\bar{\phi}} \tilde{\phi}^2 + \nu^2 \dfrac{\tilde{\phi}_\xi^2}{v} \right) d\xi + C \| \tilde{v} \|_{L^2(\mathbb{R})}^2 + C \delta_0.
		\end{array}
	\end{align*}
	Thus, 
	\begin{align*}
		\| \tilde{\phi} \|_{H^1(\mathbb{R})} \le C \| \tilde{v} \|_{L^2(\mathbb{R})} + C \delta_0.
	\end{align*}

	Similarly,
	\begin{align*}
		\| \tilde{\phi} \|_{H^k(\mathbb{R})} \le \| \tilde{v} \|_{H^{k-1}(\mathbb{R})} + C \delta_0, \quad k = 2, 3.
	\end{align*}
	Therefore, we complete the proof of this lemma.

\end{proof}

\bigskip

\noindent

\end{document}